# Optimization Based Solutions for Control and State Estimation in Non-holonomic Mobile Robots: Stability, Distributed Control, and Relative Localization

by

©Mohamed Walid Mehrez Elrafei Mohamed Said

A thesis submitted to the School of Graduate Studies
in partial fulfillment of the requirements for the degree of

**Doctor of Philosophy**

**Faculty of Engineering and Applied Science**
**Memorial University of Newfoundland**

December 2016 (Preprint)

St. John's, Newfoundland

# Abstract

Interest in designing, manufacturing, and using autonomous robots has been rapidly growing during the most recent decade. The main motivation for this interest is the wide range of potential applications these autonomous systems can serve in. The applications include, but are not limited to, area coverage, patrolling missions, perimeter surveillance, search and rescue missions, and situational awareness. In this thesis, the area of control and state estimation in non-holonomic mobile robots is tackled. Herein, optimization based solutions for control and state estimation are designed, analyzed, and implemented to such systems. One of the main motivations for considering such solutions is their ability of handling constrained and nonlinear systems such as non-holonomic mobile robots. Moreover, the recent developments in dynamic optimization algorithms as well as in computer processing facilitated the real-time implementation of such optimization based methods in embedded computer systems.

Two control problems of a single non-holonomic mobile robot are considered first; these control problems are point stabilization (regulation) and path-following. Here, a model predictive control (MPC) scheme is used to fulfill these control tasks. More precisely, a special class of MPC is considered in which terminal constraints and costs are avoided. Such constraints and costs are traditionally used in the literature to guarantee the asymptotic stability of the closed loop system. In contrast, we use a recently developed stability criterion in which the closed loop asymptotic stability can be guaranteed by appropriately choosing the prediction horizon length of the MPC controller. This method is based on finite time controllability as well as bounds on the MPC value function.

Afterwards, a regulation control of a multi-robot system (MRS) is considered. In this control problem, the objective is to stabilize a group of mobile robots to form a pattern. We achieve this task using a distributed model predictive control (DMPC) scheme based on a novel communication approach between the subsystems. This newly introduced method is based on the quantization of the robots' operating region. Therefore, the proposed communication technique allows for exchanging data in the form of integers instead of floating-point numbers. Additionally, we introduce a differential communication scheme to achieve a further reduction in the communication load.

Finally, a moving horizon estimation (MHE) design for the relative state estimation (relative localization) in an MRS is developed in this thesis. In this framework, robots with less payload/computational capacity, in a given MRS, are localized and tracked using robots fitted with high-accuracy sensory/computational means. More precisely, relative measurements between these two classes of robots are used to localize the less (computationally) powerful robotic members. As a complementary part of this study, the MHE localization scheme is combined with a centralized MPC controller to provide an algorithm capable of localizing and controlling an MRS based only on relative sensory measurements. The validity and the practicality of this algorithm are assessed by real-time laboratory experiments.

The conducted study fills important gaps in the application area of autonomous navigation especially those associated with optimization based solutions. Both theoretical as well as practical contributions have been introduced in this research work. Moreover, this thesis constructs a foundation for using MPC without stabilizing constraints or costs in the area of non-holonomic mobile robots.

iii

*To my caring parents and my lovely sisters ...*

# Acknowledgements

After thanking Almighty "ALLAH" for his blessing and guidance to complete this work, I would like to offer my sincere thanks to my supervisors Dr. George Mann and Dr. Ray Gosine for the amazing opportunity they provided me with. Their flexibility and the freedom they gave me in pursuing my research work allowed me to achieve a considerable number of contributions. Their help also allowed me to open collaboration channels with a number of universities in Europe.

I would like to acknowledge the financial support provided by my supervisors, Natural Sciences and Engineering Research Council of Canada (NSERC), the Research and Development Corporation (RDC), C-CORE J.I. Clark Chair, and Memorial University of Newfoundland. I would like to thank the staff of the Department of Mechanical Engineering. I thank also all of the ISLab colleagues. Beside it was a pleasure working next to them, their company made me always feel as if I am home with my family.

I would like also to thank all the members of the scientific network (DFG Grant WO 2056/1). In this network, I gladly worked with remarkable scholars in the field of optimal control. I specially thank Dr. Karl Worthmann (TU Ilmenau), who is the principal investigator of this network, for introducing me to such an amazing group of scientists.

Last but not least, I would like to thank my parents Mehrez and Magda; and my sisters Wedad and Walaa. Without their support, I wouldn't have had the stamina to finish this work.

# Table of Contents

# List of Tables

# List of Figures

# List of Abbreviations

| | |
|---|---|
| DMPC | Distributed model predictive control |
| DOF | Degrees of freedom |
| EKF | Extended Kalman filter |
| GPS | Global Positioning System |
| IRRM | Inter robot relative measurement |
| LSQ | Least squares |
| MHE | Moving horizon estimation |
| MIMO | Multiple-input multiple-output |
| MPC | Model predictive control |
| MPFC | Model predictive path-following control |
| MRS | Multi-robot system |
| OCP | Optimal control problem |
| PF | Particle filter |
| RL | Relative localization |
| RMSE | Root mean square error |
| RTI | Real-time iteration |
| SISO | Single-input single-output |
| SLAM | Simultaneous localization and mapping |
| TCP/IP | Transfer control protocol/internet protocol |
| UAV | Unmanned aerial vehicle |
| UGV | Unmanned ground vehicle |
| UKF | Unscented Kalman filter |

# Mathematical Notation

| | |
|---|---|
| $\mathbb{N}$ | Natural numbers. |
| $\mathbb{R}$ | Real numbers. |
| $\mathbb{Z}$ | Integer numbers. |
| $\mathbb{N}_0$ | Non-negative natural numbers, i.e. $\mathbb{N}_0 := \mathbb{N} \cup \{0\}$. |
| $\mathbb{R}_{\geq 0}$ | Non-negative real numbers, where $\mathbb{R}_{>0} := \{\mathbb{R}_{\geq 0} \setminus 0\}$. |
| $t$ | time variable. |
| $\delta$ | Sampling time (period). |
| $k$ | Open-loop time index $k \in \mathbb{N}_0$ unless otherwise defined. |
| $n$ | Closed-loop time index $n \in \mathbb{N}_0$ unless otherwise defined. |
| $\mathbf{x}$ | A state vector $\mathbf{x} \in \mathbb{R}^{n_x}$ with number of elements $n_x \in \mathbb{N}$. |
| $\mathbf{u}$ | A control input vector $\mathbf{u} \in \mathbb{R}^{n_u}$ with number of elements $n_u \in \mathbb{N}$. |
| $\mathbf{y}$ | An output vector $\mathbf{y} \in \mathbb{R}^{n_y}$ with number of elements $n_y \in \mathbb{N}$. |
| $\mathbf{f}$ | An analytic vector mapping. |
| $\mathbf{h}$ | An output function. |
| $\mathbf{x}_0$ | initial condition of the system dynamics. |
| $\mathbf{x_u}(\cdot; \mathbf{x}_0)$ | A discrete-time state trajectory starting from the initial condition $\mathbf{x}_0$ driven by the control sequence $\mathbf{u}$. |
| $\mathbf{x}(\cdot; \mathbf{x}_0, \mathbf{u})$ | A continuous-time state trajectory starting from the initial condition $\mathbf{x}_0$ driven by the control function $\mathbf{u}$. |
| $c$ | Cell width of a discretized map. |
| $(\cdot)^r$ | A superscript used to denote references. |
| $(\cdot)^\star$ | A superscript used to denote optimal solutions of OCP's. |
| $N_E$ | Estimation horizon length. |
| $N$ or $N_C$ | Prediction horizon length. |
| $T$ | Prediction horizon length in continuous time settings, i.e. $T = \delta N$. |
| $n_{mc}$ | Number of Monte-Carlo simulations. |
| $\mathcal{K}$ | Set of class $\mathcal{K}$ functions. A continuous function $\eta : \mathbb{R}_{\geq 0} \to \mathbb{R}_{\geq 0}$ is said to be of class $\mathcal{K}$ if it is zero at zero and strictly monotonically increasing. |
| $\mathcal{K}_\infty$ | Set of class $\mathcal{K}_\infty$ functions, i.e. a class $\mathcal{K}$ function that is additionally unbounded. |
| $\mathcal{KL}$ | Set of class $\mathcal{KL}$ functions. A function $\beta : \mathbb{R}_{\geq 0} \times \mathbb{N}_0 \to \mathbb{R}_{\geq 0}$ is said to be of class $\mathcal{KL}$ if $\beta(\cdot, n) \in \mathcal{K}_\infty$ for all $n \in \mathbb{N}_0$ and $\beta(r, \cdot)$ is strictly monotonically decaying to zero for each $r > 0$. |
| $\mathcal{PC}([0,T], \mathbb{R}^m)$ | Class of piecewise continuous functions for $T \in \mathbb{R}_{>0} \cup \{\infty\}$ and $m \in \mathbb{N}$. |
| $\|\mathbf{x}\|$ | 2-norm (Euclidean norm) of a vector $\mathbf{x} \in \mathbb{R}^n$, $n \in \mathbb{N}$. |
| $\|\mathbf{x}\|_\mathbf{M}^2$ | A short notation for $\mathbf{x}^\top \mathbf{M} \mathbf{x}$. |
| $\mathrm{diag}(q_1, \ldots, q_n)$ | Diagonal matrix with entries $q_1, \ldots, q_n$. |
| $\lceil \cdot \rceil$ | The Ceil operator, which is defined as $\lceil x \rceil := \min\{n \in \mathbb{Z} : n \geq x\}$. |
| $\lfloor \cdot \rfloor$ | The Floor operator, which is defined as $\lfloor x \rfloor := \max\{m \in \mathbb{Z} : m \leq x\}$. |

# Chapter 1

# Introduction

In this chapter, the targeted area of this study; the related control/estimation problems; and the proposed solution methods are introduced. Afterwards, the thesis statement as well as the expected outcomes/contributions of this work are presented. Finally, a brief outlook and organization of the thesis is presented.

## 1.1 Motivation

Autonomy is the broad umbrella under which most of the industrial/academic research interests meet. As a matter of fact, the recent technological advances in industry showed that autonomous navigation is, certainly, the only way to explore the outer space part of our universe. A very good example of such a mission can be seen in the recently launched and landed Mars rover autonomous vehicle *Curiosity*, see, e.g. [1]. Moreover, a numerous number of applications exploit autonomous navigation in order to reduce/eliminate human endangering and increase efficiency of the executed mission. Examples of these applications are area coverage and patrolling missions [2], aerial/perimeter surveillance [3–5], search and rescue missions [6], searching operations [7], situational awareness [8], and grid searching [9].

## 1.2 Autonomous Systems Considered in the Thesis

Autonomous robots can be classified into several groups, e.g. legged and wheeled, holonomic and non-holonomic; aerial, ground and underwater; homogeneous multi-robot systems and heterogeneous; and large scale robots, macro and micro, see, e.g. [10], for details. The scope of this thesis includes both non-holonomic individual ground robots as well as heterogeneous multi-robot systems (MRS's). On the one hand, the interest in non-holonomic robots arises due to practical and theoretical concerns. The non-holonomic nature of this class of mobile robots makes their number of control inputs limited. This feature imposes a controllability issue that is directly linked to the point-stabilization task, i.e. the task cannot be accomplished with a pure feedback control. This is because Brockett's condition [11] implies that the linearized non-holonomic model is not stabilizable, see, e.g. [12, 13], for details. On the other hand, multi-robot systems (MRS's) are interesting due to their service in a wide range of applications, e.g. education, military, field, demining, and surveillance, see the review articles [14–16]. In such applications, some given tasks may not be feasible by depending on a single working robot. Furthermore, for a give autonomous task, using multiple robots may require a shorter execution time with a possibly higher performance when compared with using a single working robot. Finally, MRS's have flexibility in mission execution and tolerance to possible robots' faults, i.e. a failure of a robotic member in a given MRS team should not necessarily lead to the failure of the overall mission. Therefore, MRS's have been realized, in the literature, for both homogeneous systems, e.g. ground mobile robots [17], underwater vehicles [18] and unmanned aerial vehicles [19]; as well as teams of heterogeneous robots [20–23].

In fact, heterogeneous MRS's are receiving even more attention in the robotics community. In these systems, both aerial and ground robotic members work cooperatively. This framework allows a heterogeneous MRS to perform complex tasks/missions that cannot be achieved by a homogeneous MRS. For example, in an area searching application, imagine a ground robot (carrying an aerial vehicle) that navigates in a cluttered environment while scanning it by cameras, laser scanners, and/or sonar sensors. At a particular point of the mission, the terrain may become impassable; thus, the aerial vehicle can be deployed and take aerial pictures of the area inaccessible to the ground robot. Moreover, as ground robots have in general longer motion time, they can assist the limited flight time of the aerial vehicles by deploying them only when necessary [24].

## 1.3 Control/Estimation Problems Considered

### 1.3.1 Control problems of a single non-holonomic robot

For a given robotic application, there is a number of basic tasks, which includes, but not limited to mapping, localization, planning, and control [25]. Usually, the first two tasks, i.e. mapping and localization are connected such that a robotic member can localize itself once a map of the navigated environment is possessed. A typical mapping technique is known as simultaneous localization and mapping (SLAM) [26]. After robotic members in a given application are *accurately* localized, the mission planning task starts. In this task, a trajectory/path for possibly more than one robotic member are designed. This task requires the satisfaction of the environmental constraints of the operating region, e.g. (static/dynamic) obstacles, and map margins. Finally, a controller is used in order to achieve the requirements of the planning task, i.e. in the control task, appropriate control actions are calculated to achieve the overall mission in a given application. For a single non-holonomic robot, the motion control problems can be classified into three problems as [27]

1. **point stabilization** (regulation), where the objective is to steer a non-holonomic mobile robot to a desired position and orientation, i.e. a desired posture,

2. **trajectory tracking**, where the vehicle is required to track a time-parameterized reference, i.e. a time varying reference, and

3. **path following**, where the vehicle is required to approach and follow a desired path-parameterized reference, without a fixed timing law.

The path following control gives an additional degree of freedom in the controller design than the trajectory tracking control. More precisely, the choice of the timing law is not fixed, see, e.g. [28], for details. Therefore, in this thesis, we focus on the point stabilization and the path following control problems of a single non-holonomic robot.

### 1.3.2 Distributed control of a multi-robot system

Next, we consider the stabilization of a group of multiple robots to a permissible equilibria. In this framework, each robot in the given group is to be stabilized to a reference posture while avoiding collision with other robots. To this end, centralized as well as distributed control approaches can be applied. While centralized control schemes lead normally to superior performance, they become real-time infeasible when considering robotic groups with a large number of robots (scalability problem). Distributed control techniques solves the scalability problem by distributing the control tasks among the subsystems of the robotic team. However, this approach requires data communication among the subsystems essentially for their locations, see, e.g. [27]. In this thesis, we give a particular attention to relaxing the communication load in distributed control of systems with multiple robots.

### 1.3.3 Relative-localization in a multi-robot system

In a multi-robot system (MRS), the localization task can be achieved by equipping all the involved robots with the state-of-the-art sensory means. However, this increases the overall cost of such systems. Moreover, some robotic members, e.g. flying/hovering robots may not have enough payload/computational capacity to operate such means. Therefore, relative localization (RL) has been developed as a practical solution for effective and accurate execution of multi-robot collaborative tasks. The objective of relative localization is to detect and locate robots with limited sensory capabilities (observed/child robots) with respect to another robots with accurate localization means (observing/leader robots). This is achieved by using the relative observation between the two robots categories, see [20–23, 29] for more details.

## 1.4 Proposed Solution Methods

In this thesis, the control task is considered first for the single non-holonomic robot control problems, i.e. point stabilization (regulation), and path following. In this context, an optimization based controller, i.e. nonlinear model predictive control (MPC) *without stabilizing constraints or costs* is applied, for the first time, to the mentioned control problems while the closed loop stability requirements are rigorously proven. In MPC, we first measure the state of the system, which serves as a basis for solving a finite horizon optimal control problem (OCP) that characterize the control objective. This results in a sequence of future control values. Then, the first element of the computed sequence is applied before the process is repeated, see, e.g. [30]. In the settings of MPC without stabilizing constraints or costs, stability can be guaranteed by *appropriately* choosing the prediction horizon length. MPC without stabilizing constraints or costs with guaranteed closed loop asymptotic stability has recently emerged in the control Engineering community and it is of a particular interest due to many reasons particularly the following: first, this MPC scheme is the easiest one to implement and it is widely used in the industry. Second, the knowledge of a local Lyapunov function–a design requirement in other MPC stabilizing designs–is not a requirement in this MPC scheme, see, e.g. [31, 32].

Next, we implement a distributed model predictive control (DMPC) scheme to stabilize a multi-robot system of ground robots to an equilibria. In DMPC, the control task is distributed among the subsystems, i.e. the robots solve their own (local) OCP's and communicate with the other robots to avoid collisions. The robots share information in order to allow other robots to formulate proper coupling constraints and, thus, avoid collision with each other. To formulate the coupling constraints, previous studies are based on the communication of predicted trajectories, see, e.g. [33]. In contrast, we first partition the operating region into a grid and derive an estimate on the minimum width of a grid cell. Then, the predicted trajectories are projected onto the grid resulting in an occupancy grid, which serves as quantization of the communication data. Based on a data exchange of these projections, each robot formulates suitable coupling constraints. Utilizing the occupancy grid reduces bandwidth limits and congestion issues since a more compact data representation (integers instead of floating point values) is employed. Moreover, based on an introduced differential communication method, we show that the communication load can be reduced significantly.

Finally, the relative localization task in MRS is solved using an optimization based state estimator, i.e. nonlinear moving horizon estimation (MHE) scheme. MHE considers the evolution of a constrained and possibly nonlinear model on a fixed time horizon, and minimizes the deviation of the model from a past window of measurements [34]. The main motivation of using MHE is the limitations associated with the traditionally used relative localization estimators, which lie mainly in their relatively long estimation settling time, i.e. the time required such that the estimation error reaches acceptable error margins, see, e.g. [20–23].

Indeed, considering optimization based solutions, e.g. MPC and MHE, is generically interesting due to their ability of handling constrained and (possibly) nonlinear systems. Moreover, although optimization based solutions were initially criticized in the literature by their real-time applicability limitations, the recent developments in optimization algorithms as well as in digital computing fulfilled the real-time applicability requirements. Real-time applicability is also demonstrated in the thesis.

## 1.5 Thesis Statement

Two main research areas in mobile robots, i.e. control and state estimation, are tackled in this study, where optimization based solutions are used. As highlighted in the literature, the control problems presented in Section 1.3, for a single non-holonomic robot, are resolved by MPC controllers utilizing either stabilizing constraints and/or terminal costs in order to guarantee the closed loop stability of the system. In contrast, in this study, we explore the MPC designs without stabilizing constraints or costs for these control problems, and investigate the required stability conditions. Afterwards, we consider the distributed implementation of MPC to regulation of systems with multiple robots. To this end, we introduce a novel approach of reducing the communication load via quantization techniques and differential communication. Finally, a moving horizon estimation (MHE) observer is designed, investigated, and tested for the relative localization problem involved in multi-robot systems.

## 1.6 Objectives and Expected Contributions

The proposed research here has both theoretical as well as practical objectives. These objectives are summarized in the following, where the expected contributions are highlighted.

**Objective 1:** A novel design of an MPC control algorithm, without stabilizing constraints or costs, for the point stabilization of a non-holonomic mobile robot.

- A novel verification of the controllability assumptions for point stabilization.

- Determination of the necessary requirements that guarantee the closed loop system stability, e.g. MPC running costs, and prediction horizon.

- Investigation of the controller components' effects on its performance regarding the stability.

**Objective 2:** A novel design of an MPC control algorithm, without stabilizing constraints or costs, for the path following control of a non-holonomic mobile robot.

- A novel verification of the controllability assumptions for path following.

- Determination of the necessary requirements that guarantee the closed loop system stability, e.g. MPC running costs, and prediction horizon.

- Investigation of the controller components' effects on its performance regarding the stability.

**Objective 3:** A novel design of a DMPC based on occupancy grid communication.

- Introduction of a grid generation technique and quantization of the communicated data.

- Introduction of a differential communication technique.

- Novel formulation of coupling constraints based on occupancy grid.

- Numerical validation and demonstration of the introduced technique.

**Objective 4:** A novel implementation of MHE observers in relative localization in multi-robot systems.

- Problem formulation of the relative localization estimation problem in the MHE frame work.

- Investigation of the MHE performance by the traditionally used estimators in order to highlight its advantages.

- Development of an optimization based algorithm in which relative localization as well as relative trajectory tracking are achieved by MHE and MPC, respectively.

- Numerical as well as experimental validation of the developed algorithm.

As can be noticed, the thesis combines both theoretical as well as practical aspects related to optimization based solutions for control/state-estimation of non-holonomic mobile robots. Moreover, both individual mobile robots and multi-robotic systems are considered.

## 1.7 Organization of the Thesis

Four different control/estimation problems in the area of non-holonomic robots are considered in this thesis. Therefore, a general background is given in chapter 2. Moreover, the core chapters of the thesis, i.e. Chapters 3, 4, 5, 6 are written in the **manuscript format** in which the background, the motivation, and the conclusion to each problem considered are discussed. The outline of the thesis is as follows.

**Chapter 1:** In this chapter, the research area is motivated first. Then, the considered systems and their associated control/estimation problems as well as the proposed solution methods are presented.

**Chapter 2:** In this chapter, a brief background of the considered control/estimation problems is provided. Next, the main drawbacks of the previously utilized solution methods are discussed. Finally, the main challenges of the used solution methods are outlined.

**Chapter 3:** In this chapter, the *discrete time* formulation of MPC control without stabilizing constraints and costs are presented. Then, the fundamental stability theorem and its necessary assumptions are shown. Next, these assumptions are verified for the **point stabilization** problem of a non-holonomic mobile robot. Afterwards, these results are verified by means of numerical experiments. Finally, the concluding remarks are drawn.

**The content of this chapter appeared in the publication:**

Karl Worthmann, Mohamed W. Mehrez, Mario Zanon, George K.I. Mann, Raymond G. Gosine, and Moritz Diehl, "Model Predictive Control of Nonholonomic Mobile Robots Without Stabilizing Constraints and Costs", in IEEE Transactions on Control Systems Technology, 2016.

**Chapter 4:** In this chapter, the *continuous time* formulation of MPC control without stabilizing constraints and costs are presented. Then, the fundamental stability theorem and its necessary assumptions are shown. Next, these assumptions are verified for the **path following** problem of a non-holonomic mobile robot. Afterwards, these results are verified by means of numerical experiments. Finally, the concluding remarks are drawn.

**The content of this chapter appeared in the publication:**

Mohamed W. Mehrez, Karl Worthmann, George K.I. Mann, Raymond G. Gosine, and Timm Faulwasser, "Predictive Path Following of Mobile Robots without Terminal Stabilizing Constraints", in Proceedings of the IFAC 2017 World Congress, Toulouse, France, 2017, accepted for publication.

**Chapter 5:** In this chapter, a gird generation is introduced while a DMPC scheme is presented. Afterwards, a differential communication scheme is presented. Then, a suitable representation of the collision avoidance constraints are derived. Next, the proposed method is investigated by means of numerical simula-tions. Finally, conclusions are drawn.

**The content of this chapter appeared partially in the publication:**

Mohamed W. Mehrez, Tobias Sprodowski, Karl Worthmann, George K.I. Mann, Raymond G. Gosine, Juliana K. Sagawa, and Jürgen Pannek, "Occupancy Grid based Distributed MPC of Mobile Robots", Submitted.

**Chapter 6:** In this chapter, the relative localization problem in the MHE framework is formulated with demonstration of MHE capability of achieving the localization task with high accuracy. Then, the development of an optimization based algorithm is highlighted in which relative localization as well as relative trajectory tracking are achieved by means of MHE and MPC with numerical/experimental demonstration. Finally, conclusions are drawn.

**The content of this chapter appeared in the publication:**

Mohamed W. Mehrez, George. K.I. Mann, and Raymond G. Gosine, "An Optimization Based Approach for Relative Localization and Relative Tracking Control in Multi-Robot Systems", in Journal of Intelligent and Robotic Systems, 2016.

**Chapter 7:** In this chapter, the thesis findings are discussed and analyzed and the related conclusions are drawn. Finally, possible future research directions are presented.

### 1.7.1 Note to readers

The discrete time formulation of MPC and MHE is used throughout the thesis except for Chapter 4, where the continuous time formulation of MPC is followed. This is because of the functional derivatives required in the path following control problem, i.e. the control problem considered in Chapter 4.

As mentioned earlier, the core chapters of the thesis, i.e. Chapters 3, 4, 5, 6 are written in the manuscript format. Nonetheless, the notations of these chapters have been modified from their original format such that the notation is consistent throughout the whole thesis. Additionally, more details/figures have been added in order to increase the clarity of the presentation.

# References

[1] A. D. Steltzner, A. M. S. Martin, T. P. Rivellini, A. Chen, and D. Kipp, "Mars science laboratory entry, descent, and landing system development challenges," *Journal of Spacecraft and Rockets*, vol. 51, no. 4, pp. 994–1003, 2014.

[2] J. Acevedo, B. Arrue, I. Maza, and A. Ollero, "Distributed approach for coverage and patrolling missions with a team of heterogeneous aerial robots under communication constraints," *International Journal of Advanced Robotic Systems*, vol. 10, no. 28, pp. 1–13, 2013.

[3] R. Beard, T. McLain, D. Nelson, D. Kingston, and D. Johanson, "Decentralized cooperative aerial surveillance using fixed-wing miniature uavs," *Proceedings of the IEEE*, vol. 94, no. 7, pp. 1306–1324, July 2006.

[4] D. Kingston, R. Beard, and R. Holt, "Decentralized perimeter surveillance using a team of uavs," *Robotics, IEEE Transactions on*, vol. 24, no. 6, pp. 1394–1404, Dec 2008.

[5] J. Acevedo, B. Arrue, I. Maza, and A. Ollero, "Cooperative large area surveillance with a team of aerial mobile robots for long endurance missions," *Journal of Intelligent and Robotic Systems*, vol. 70, no. 1-4, pp. 329–345, 2013.

[6] M. Bernard, K. Kondak, I. Maza, and A. Ollero, "Autonomous transportation and deployment with aerial robots for search and rescue missions," *Journal of Field Robotics*, vol. 28, no. 6, pp. 914–931, 2011.

[7] I. Maza and A. Ollero, "Multiple uav cooperative searching operation using polygon area decomposition and efficient coverage algorithms," in *Distributed Autonomous Robotic Systems 6*, R. Alami, R. Chatila, and H. Asama, Eds. Birkhäuser Basel, 2007, pp. 221–230.

[8] M. A. Hsieh, A. Cowley, J. F. Keller, L. Chaimowicz, B. Grocholsky, V. Kumar, C. J. Taylor, Y. Endo, R. C. Arkin, B. Jung, D. F. Wolf, G. S. Sukhatme, and D. C. MacKenzie, "Adaptive teams of autonomous aerial and ground robots for situational awareness," *Journal of Field Robotics*, vol. 24, no. 11-12, pp. 991–1014, 2007.

[9] W. Dunbar and R. Murray, "Model predictive control of coordinated multi-vehicle formations," in *Decision and Control, 2002, Proceedings of the 41st IEEE Conference on*, vol. 4, 2002, pp. 4631–4636 vol.4.

[10] R. Siegwart, I. R. Nourbakhsh, and D. Scaramuzza, *Introduction to autonomous mobile robots*. MIT press, 2011.

[11] R. W. Brockett, "Asymptotic stability and feedback stabilization," in *Differential Geometric Control Theory*, R. W. Brockett, R. S. Millman, and H. J. Sussmann, Eds. Birkhäuser, Boston, MA, 1983, pp. 181–191.

[12] J. P. Laumond, S. Sekhavat, F. Lamiraux, J. paul Laumond (editor, J. P. Laumond, S. Sekhavat, and F. Lamiraux, "Guidelines in nonholonomic motion planning for mobile robots," in *Robot Motion Plannning and Control*. Springer-Verlag, 1998, pp. 1–53.

[13] M. Michalek and K. Kozowski, "Vector-field-orientation feedback control method for a differentially driven vehicle," *IEEE Transactions on Control Systems Technology*, vol. 18, no. 1, pp. 45–65, 2010.

[14] M. Lewis and K.-H. Tan, "High precision formation control of mobile robots using virtual structures," *Autonomous Robots*, vol. 4, no. 4, pp. 387–403, 1997.

[15] T. Arai, E. Pagello, and L. Parker, "Guest editorial advances in multirobot systems," *Robotics and Automation, IEEE Transactions on*, vol. 18, no. 5, pp. 655–661, Oct 2002.

[16] D. Portugal and R. Rocha, "A survey on multi-robot patrolling algorithms," in *Technological Innovation for Sustainability*, ser. IFIP Advances in Information and Communication Technology, L. Camarinha-Matos, Ed. Springer Berlin Heidelberg, 2011, vol. 349, pp. 139–146.

[17] S. Carpin and L. E. Parker, "Cooperative motion coordination amidst dynamic obstacles," in *in Distributed Autonomous Robotic Systems*, 2002, pp. 145–154.

[18] D.-U. Kong and J. An, "Modular behavior controller for underwater robot teams: A biologically inspired concept for advanced tasks," in *Intelligent Robotics and Applications*, ser. Lecture Notes in Computer Science, S. Jeschke, H. Liu, and D. Schilberg, Eds. Springer Berlin Heidelberg, 2011, vol. 7102, pp. 536–547.

[19] C. Yokoyama, M. Takimoto, and Y. Kambayashi, "Cooperative control of multi-robots using mobile agents in a three-dimensional environment," in *Systems, Man, and Cybernetics (SMC), 2013 IEEE International Conference on*, Oct 2013, pp. 1115–1120.

[20] O. De Silva, G. Mann, and R. Gosine, "Development of a relative localization scheme for ground-aerial multi-robot systems," in *Intelligent Robots and Systems (IROS), 2012 IEEE/RSJ International Conference on*, 2012, pp. 870–875.

[21] T. Wanasinghe, G. I. Mann, and R. Gosine, "Relative localization approach for combined aerial and ground robotic system," *Journal of Intelligent and Robotic Systems*, vol. 77, no. 1, pp. 113–133, 2015.

[22] T. Wanasinghe, G. Mann, and R. Gosine, "Distributed leader-assistive localization method for a heterogeneous multirobotic system," *Automation Science and Engineering, IEEE Transactions on*, vol. 12, no. 3, pp. 795–809, July 2015.

[23] M. W. Mehrez, G. K. Mann, and R. G. Gosine, "Nonlinear moving horizon state estimation for multi-robot relative localization," in *Electrical and Computer Engineering (CCECE), 2014 IEEE 27th Canadian Conference on*, May 2014, pp. 1–5.

[24] O. D. Silva, "Relative localization of gps denied ground aerial multi-robot systems," Ph.D. dissertation, Memorial University of Newfoundland, 2015.

[25] P. Morin and C. Samson, "Motion control of wheeled mobile robots," in *Springer Handbook of Robotics*, B. Siciliano and O. Khatib, Eds.   Springer Berlin Heidelberg, 2008, pp. 799–826.

[26] H. Durrant-Whyte and T. Bailey, "Simultaneous localization and mapping: part i," *Robotics Automation Magazine, IEEE*, vol. 13, no. 2, pp. 99–110, June 2006.

[27] K. Kanjanawanishkul, "Coordinated path following control and formation control of mobile robots," Ph.D. dissertation, University of Tübingen, 2010.

[28] T. Faulwasser and R. Findeisen, "Nonlinear model predictive path-following control," in *Nonlinear Model Predictive Control*, ser. Lecture Notes in Control and Information

Sciences, L. Magni, D. Raimondo, and F. Allgöwer, Eds. Springer Berlin Heidelberg, 2009, vol. 384, pp. 335–343.

[29] F. Rivard, J. Bisson, F. Michaud, and D. Letourneau, "Ultrasonic relative positioning for multi-robot systems," in *Robotics and Automation, 2008. ICRA 2008. IEEE International Conference on*, 2008, pp. 323–328.

[30] J. B. Rawlings and D. Q. Mayne, *Model Predictive Control: Theory and Design*. Nob Hill Publishing, 2009.

[31] N. Altmüller, L. Grüne, and K. Worthmann, "Performance of nmpc schemes without stabilizing terminal constraints," in *Recent Advances in Optimization and its Applications in Engineering*, M. Diehl, F. Glineur, E. Jarlebring, and W. Michiels, Eds. Springer Berlin Heidelberg, 2010, pp. 289–298.

[32] K. Worthmann, M. Reble, L. Grüne, and F. Allgöwer, "The role of sampling for stability and performance in unconstrained nonlinear model predictive control," *SIAM Journal on Control and Optimization*, vol. 52, no. 1, pp. 581–605, 2014.

[33] L. Grüne and K. Worthmann, "A distributed nmpc scheme without stabilizing terminal constraints," in *Distributed Decision Making and Control*, ser. Lecture Notes in Control and Information Sciences, R. Johansson and A. Rantzer, Eds. Springer London, 2012, vol. 417, pp. 261–287.

[34] C. Rao and J. Rawlings, "Nonlinear moving horizon state estimation," in *Nonlinear Model Predictive Control*, ser. Progress in Systems and Control Theory, F. Allgöwer and A. Zheng, Eds. Birkhäuser Basel, 2000, vol. 26, pp. 45–69.

# Chapter 2

# Background and Literature Review

In this chapter, we first present a background on the systems considered in the thesis as well as the proposed solution methods. Then, we highlight the main issues related to the control/state-estimation problems considered and how they are tackled by the proposed solutions.

## 2.1 Non-holonomic Mobile Robots

Non-holonomic mobile robots constitute a large class of mobile robots used in autonomous navigation. In terms of holonomy, holonomic mobile robots also exist, see Figure 2.1 for examples of research holonomic and non-holonomic robots. The term non-holonomic is directly related to the kinematic constraints this class of robots possess, i.e. these robots are constrained to move laterally. Under this class of mobile robots, a number of kinematic configurations can be observed, e.g. unicycle (differential drive robots), bicycle, and car-like robots, see Figure 2.2 for details on these configurations. The non-holonomic nature of this class of mobile robot imposes limitations on acceptable system velocities [1]. Nevertheless, non-holonomy becomes useful as it limits the number of control inputs, while maintaining the full controllability of the system in the state space [2]. This advantage,

*Nonholonomic Research Platforms*

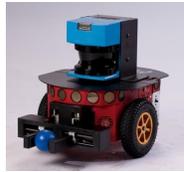
**Pioneer 3DX**

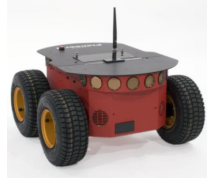
**Pioneer 3AT**

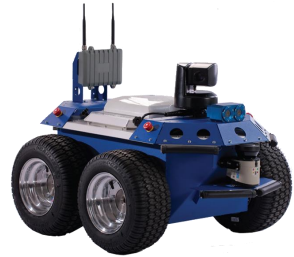
**Seekur Jr**

*Holonomic Research Platforms*

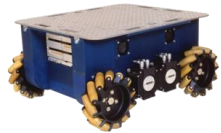
**URANUS**

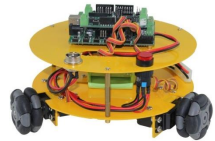
**3WD Compact Omni-Directional**

Fig. 2.1: Examples of non-holonomic and holonomic research platforms.

however, imposes a difficulty that is connected with the point-stabilization task, i.e. the task cannot be accomplished with a pure feedback control algorithm. This is due to Brockett's condition [3], which implies that the linearized non-holonomic model is not stabilizable, see, e.g. [1, 2], for more details. In the following subsection, we present briefly the mathematical formulation of the non-holonomic mobile robots control problems.

### 2.1.1 Control problems of a non-holonomic robot

Three fundamental control problems are realized for non-holonomic mobile robots, i.e. point stabilization, trajectory tracking, and path following [5]. In order to differentiate between these control problems, we make use of the following general discrete dynamical model.

$$\mathbf{x}(k+1) = \mathbf{f}(\mathbf{x}(k), \mathbf{u}(k)), \qquad \mathbf{x}(0) = \mathbf{x}_0 \qquad (2.1)$$

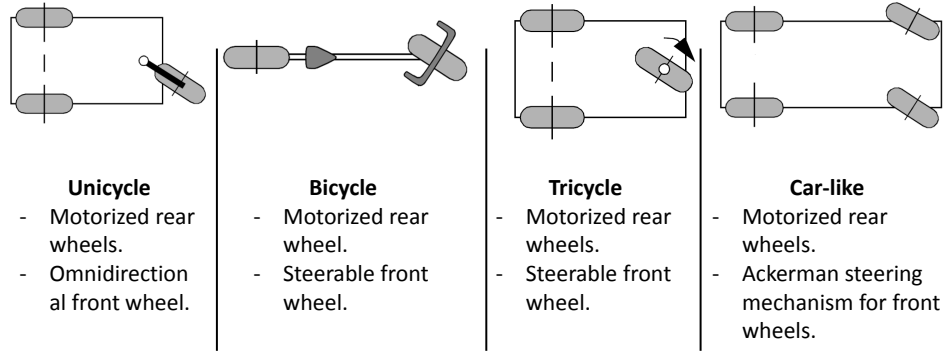

Fig. 2.2: Different robotic configurations with non-holonomic constraints based on [4].

where $k \in \mathbb{N}_0$ is the current sampling instant, $\mathbf{x} \in \mathbb{R}^{n_x}$ is the system state vector, $\mathbf{u} \in \mathbb{R}^{n_u}$ is the system control input, and $\mathbf{f}(\cdot) : \mathbb{R}^{n_x} \times \mathbb{R}^{n_u} \to \mathbb{R}^{n_x}$ is the nonlinear dynamical mapping. In addition, state constraints $X \subset \mathbb{R}^{n_x}$ as well as control constraints $U \subset \mathbb{R}^{n_u}$ can be defined for the considered system. In point stabilization control problem, a feedback control $\boldsymbol{\mu}(\mathbf{x}(k)) : X \to U$ is designed such that the solution of (2.1) starting from the initial condition $\mathbf{x}_0 := \mathbf{x}(0) \in X$ stays close to a desired set point, $\mathbf{x}^r \in X$, and converges, i.e.

$$\lim_{k \to \infty} \|\mathbf{x}(k) - \mathbf{x}^r\| = 0. \tag{2.2}$$

As can be noticed, the set point $\mathbf{x}^r$ is a constant reference in the point stabilization case. When this reference is time dependent (varying), i.e. $\mathbf{x}^r(k) : k \to X$, the considered control problem is realized as a trajectory tracking problem. Here, the feedback control $\boldsymbol{\mu}(\mathbf{x}(k)) : X \to U$ task is to steer the solution of (2.1) to track the time varying reference

such that

$$\lim_{k \to \infty} \|\mathbf{x}(k) - \mathbf{x}^r(k)\| = 0. \tag{2.3}$$

Indeed, the point stabilization task can be considered as a special case of the trajectory tracking problem. As can be observed, in trajectory tracking, the reference defines a time-state mapping for the system state $\mathbf{x}$, i.e. the reference carries the information of when to be where in the state space. On the other hand, in path following, both the system and the time-state mapping itself are affected, i.e. the mapping is not known a priori and it is a by-product of the controller. This feature makes the controller more flexible in following paths that a trajectory tracking controller cannot follow, see, e.g. [6], for a numerical example. In path following, a geometric reference curve $\mathcal{P} \in X$ parameterized in $\lambda \in \mathbb{R}$ has to be followed, i.e.

$$\mathcal{P} = \left\{ p(\lambda) \in \mathbb{R}^{n_x} : \lambda \in [\bar{\lambda}, 0] \to p(\lambda) \right\},$$

where the map $p : [\bar{\lambda}, 0] \to \mathbb{R}^{n_x}$ specifies the path $\mathcal{P}$ over the interval $[\bar{\lambda}, 0]$. The path-following control problem involves the following tasks: first, ensure the convergence of the system state $\mathbf{x}(k)$ to the path $\mathcal{P}$. Second, determine the rate of path evolution, i.e. to determine the mapping $\lambda : \mathbb{R} \to [\bar{\lambda}, 0]$. This mapping is not given in advance, but rather is used as an extra degree of freedom in the controller design, see, e.g. [7].

### 2.1.2 Non-predictive controllers used

Since it is relatively simple, a considerable number of trajectory tracking control laws have been pursued, e.g. backstepping [8], dynamic feedback linearization [9], and sliding mode control [10]; however, these techniques impose constraints on the reference speeds within which tracking can be achieved, e.g. reference speeds have to be persistently excited, i.e. their values cannot be set to zero. On the other hand, when the point stabilization (regulation) problem is considered, a linearized non-holonomic model loses stabilizability; thus, according to Brockett's theorem [3], a smooth time-invariant feedback control law does not exist. Various studies solving the regulation problem are reported in [11]. These studies include piecewise-continuous feedback control [12], smooth time-varying control [13], Lyapunov control [14], and dynamic feedback linearization control [9]. However, these methods do not limit the controller design to the feasible sets of the mobile robot's states, which is an important controller feature especially when considering domains with obstacles.

Control techniques achieving simultaneous regulation and tracking include differential kinematic control [15], and feedback linearization control [9]; however, these methods do not provide a single controller architecture capable of simultaneously achieving tracking and regulation without the need to switch between the two control modes. This problem has been resolved using backstepping control [16], and vector field orientation feedback control [2]. Still, the later controllers do not consider actuator saturation limits in a straightforward manner, i.e. a post processing step is required to scale the calculated control signals to their saturation values, see [2] for details. Moreover, the mentioned controllers are not designed, in general, in a straightforward manner; thus, the associated process of choosing suitable tuning parameters achieving an acceptable performance is difficult, see, e.g. [17], for details. Different approaches have been also established in the literature to tackle path-following problems, e.g. back stepping [18] and feedback linearization [19]. In the mentioned approaches, the consideration of state and input constraints is in general difficult. For extensive reviews on control of non-holonomic robots and methods of implementations, see, e.g., [20, 21].

As can be observed, most of the controllers presented above ***do not*** combine, in one design, many important features, e.g. consideration of the natural saturation limits on control inputs and states; a single control architecture capable of simultaneously achieving regulation and tracking control without the need to switch between the two controllers; and a controller with well understood tuning parameters designed in a straightforward manner, see, e.g. [22], for more details.

## 2.2 Model Predictive Control (MPC)

Model predictive control (MPC), also known as receding horizon control, is considered to be one of the most attractive control strategies due to its applicability to constrained nonlinear multiple inputs multiple outputs (MIMO) systems. In MPC, a future control sequence minimizing an objective function is computed over a finite prediction horizon. Then, the first control in this sequence is applied to the system. This process is repeated at every sampling instant [23, 24].

### 2.2.1 MPC mathematical formulation

Here, we present a brief mathematical formulation of MPC. As discussed earlier, in MPC, the control input applied to the system is obtained by solving a finite horizon open-loop optimal control problem (OCP), every decision instant. We consider now the discrete dynamics (2.1) at a time instant $n$. Then, for the prediction horizon $N_C \in \mathbb{N}$ and the open-loop control sequence

$$\mathbf{u} = (\mathbf{u}(0), \mathbf{u}(1), \ldots, \mathbf{u}(N-1)) \in U^{N_C},$$

the online optimal control problem (OCP) of MPC can be outlined as the following:

$$\min_{\mathbf{u}\in\mathbb{R}^{n_u\times N_C}} \quad J_{N_C}(\mathbf{x}_0, \mathbf{u}) \tag{2.4}$$

$$\text{subject to} \quad \mathbf{x}(0) = \mathbf{x}_0,$$

$$\mathbf{x}(k+1) = \mathbf{f}(\mathbf{x}(k), \mathbf{u}(k)) \quad \forall k \in \{0, 1, \ldots, N_C - 1\},$$

$$\mathbf{x}(k) \in X \quad \forall k \in \{1, 2, \ldots, N_C\},$$

$$\mathbf{u}(k) \in U \quad \forall k \in \{0, 1, \ldots, N_C - 1\},$$

where the objective function $J_{N_C}(\mathbf{x}_0, \mathbf{u}) : X \times U^N \to \mathbb{R}_{\geq 0}$ is *generally* given by:

$$J_{N_C}(\mathbf{x}_0, \mathbf{u}) = \sum_{k=0}^{N_C - 1} \ell(\mathbf{x}(k), \mathbf{u}(k)) + F(\mathbf{x}(N_C)). \tag{2.5}$$

The first term of the objective (cost) function (2.5) is referred to as the *running costs* and is computed by penalizing the deviation of the prediction of the system state $\mathbf{x}(\cdot)$ and its reference state $\mathbf{x}^r(\cdot)$, as well as penalizing the deviation of the control $\mathbf{u}(\cdot)$ from its reference $\mathbf{u}^r(\cdot)$. The term $F(\mathbf{x}(N_C))$ is referred to as the *terminal cost* and it corresponds the deviation of the last entry of the predicted trajectory from the reference. The running cost $\ell(\cdot) : X \times U \to \mathbb{R}_{\geq 0}$ is generally given by:

$$\ell(\mathbf{x}(\cdot), \mathbf{u}(\cdot)) = \|\mathbf{x}(\cdot) - \mathbf{x}^r(\cdot)\|_{\mathbf{Q}}^2 + \|\mathbf{u}(\cdot) - \mathbf{u}^r(\cdot)\|_{\mathbf{R}}^2, \tag{2.6}$$

where $\mathbf{Q}$ and $\mathbf{R}$ are positive definite symmetric weighting-matrices of the appropriate dimensions. It has to be mentioned here that, in the case of a static reference (regulating MPC), the control reference $\mathbf{u}^r(\cdot) = 0$. Moreover, in the tracking case, the control reference $\mathbf{u}^r(\cdot)$ is the nominal control to steer the system state along the reference trajectory. As can be noticed from (2.6), the deviation of the state $\mathbf{x}$ to the reference $\mathbf{x}^r$ as well as the deviation of the control $\mathbf{u}$ to the reference $\mathbf{u}^r$ are penalized along the prediction trajectory. While penalizing the state deviation from its reference is intuitive, penalizing the deviation of the control has computational advantages, i.e. penalizing the control variable may render the optimal control problem solution easier. Moreover, when the control variable is penalized, the control values with expensive energy can be avoided [25].

Finally, the terminal cost term $F(\cdot) : X \to \mathbb{R}_{\geq 0}$ is given by:

$$F(\mathbf{x}(N_C)) = \|\mathbf{x}(N_C)) - \mathbf{x}^r(N_C)\|_{\mathbf{P}}^2 , \qquad (2.7)$$

where $\mathbf{P}$ is a positive definite weighting-matrix penalizing the deviation of the last entry of the state prediction, i.e. $\mathbf{x}(N_C)$, from its reference $\mathbf{x}^r(N_C)$. As will be discussed later in this chapter, the terminal cost (2.7) is used primarily for ensuring MPC closed-loop stability [25].

The minimizing control sequence resulting from solving OCP (3.6) is denoted by

$$\mathbf{u}^\star := (\mathbf{u}^\star(0), \mathbf{u}^\star(1), \ldots, \mathbf{u}^\star(N_C - 1)) \in U^{N_C},$$

where $\mathbf{u}^\star(0)$ is the control action to be applied on the system. Algorithm 2.1 summarizes the MPC scheme for nonlinear systems, e.g. (2.1). See also Figure 2.3 for an example of two MPC iterations applied to a simple single input single output (SISO) system.

Based on the introduced cost function (2.5), a corresponding (optimal) value function $V_{N_C} : X \to \mathbb{R}_{\geq 0} \cup \infty$ is defined for a given prediction horizon $N_C \in \mathbb{N}$ as the following:

$$V_{N_C}(\mathbf{x}_0) := \min_{\mathbf{u} \in \mathbb{R}^{2 \times N_C}} J_{N_C}(\mathbf{x}_0, \mathbf{u}) = J_{N_C}(\mathbf{x}_0, \mathbf{u}^\star). \qquad (2.8)$$

In addition to the state and control constraints presented in OCP (2.4), which account mainly for the physical limits of state and control, the OCP can be also subject to terminal equality or inequality constraints. The terminal equality constraint can be written as

$$\mathbf{x}(N_C) - \mathbf{x}^r(N_C) = 0,$$

where this constraint requires that the last entry of the predicted trajectory to be equal to its reference. Moreover, the terminal inequality constraint can be written as

$$\mathbf{x}(N_C) \in \Omega(\mathbf{x}^r(N_C)),$$

where $\Omega(\cdot) \subset X$. This constraint requires that the last entry of the predicted trajectory to be within a region around the reference $\mathbf{x}^r(N_C)$. Both terminal equality and inequality constraints are used in the literature to guarantee the stability of MPC closed-loop; they are referred to as terminal stabilizing constraints. This will be discussed in more details later in this chapter.

---
**Algorithm 2.1** MPC scheme steps
---
1: **for** every sampling instant $n = 0, 1, 2, \ldots$ **do**
2:     Measure the current state $\hat{\mathbf{x}}(n) := \mathbf{x}(n) \in X$ of the system (2.1).
3:     Set $\mathbf{x}_0 = \hat{\mathbf{x}}(n)$
4:     Find the minimizing control sequence $\mathbf{u}^\star = (\mathbf{u}^\star(0), \cdots, \mathbf{u}^\star(N_C - 1)) \in U^{N_C}$,
       which satisfies $J_{N_C}(\hat{\mathbf{x}}, \mathbf{u}^\star) = V_{N_C}(\hat{\mathbf{x}})$.
5:     Define the MPC-feedback control law $\boldsymbol{\mu}_{N_C} : X \to U$ at $\hat{\mathbf{x}}$ by $\boldsymbol{\mu}_{N_C}(\hat{\mathbf{x}}) := \mathbf{u}^\star(0)$.
6:     Apply $\boldsymbol{\mu}_{N_C}(\hat{\mathbf{x}})$ to the system (2.1).
7: **end for**
---

### 2.2.2 MPC related applications

Beside the natural consideration of physical limits, MPC has well understood tuning parameters, i.e. prediction horizon length and optimization problem weighting parameters, see, e.g. [23, 26] for details. Therefore, the use of MPC has been reported in a considerable number of applications. Both linear and nonlinear MPC have been used in a wide range of application, e.g. process control [27, 28]; heating systems [29]; wind turbines [30]; power electronics and drives [31, 32]; unmanned aerial vehicles (UAV's) including both quadrotors and helicopters [33–38], and fixed wing aircrafts [39, 40]; marine surface vessels [41, 42]; missiles and projectiles guidance [43, 44]; parallel robots [45]; and holonomic mobile robots [46].

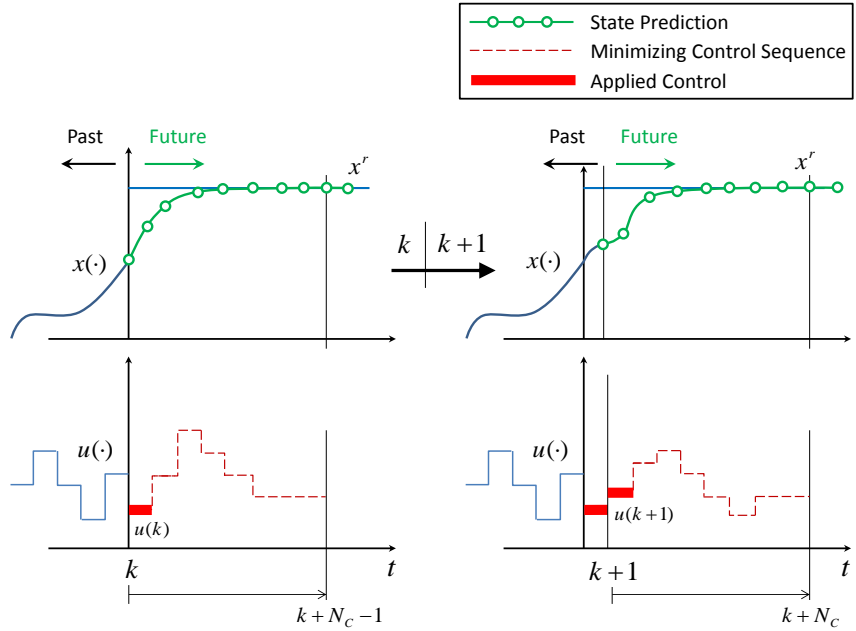

Fig. 2.3: Illustration of two MPC iterations for a simple SISO system. $k$ denotes the time step and $N_C$ denotes the prediction horizon length.

For non-holonomic robot control problems, studies utilizing linear MPC in order to achieve the trajectory tracking objective, which adopt a linearized non-holonomic motion model, are presented in [47, 48]. Nonlinear MPC, which uses the *nonlinear* motion model of robots, has been used for tracking problems [22, 49]; regulation problems [11, 17, 50, 51]; and both [52–54]. The path following control problem has been also solved using linear MPC in [55, 56]. Additionally, a primary numerical investigation of path following control to non-holonomic robots was done in [6].

### 2.2.3 MPC and stability

In this thesis, we focus on a major concern in MPC for non-holonomic mobile robots, which is **stability**. Indeed, MPC has also other issues, e.g. ***feasibility*** and ***real-time applicability***. Feasibility refers to the problem of ensuring the existence of a minimizing control sequence (function) in the involved MPC optimization problem. In general, the constrained optimization problem in MPC involves both state constraints as well as input (control) constraints, e.g. $X$ and $U$ presented in Section 2.1.1. Typically, state constraints components are treated as soft constraints, i.e. while solving the MPC optimization problem these constraints are relaxed numerically by adding a slack variable, e.g. $e << 0$. However, the physical input constraints are maintained as hard because of their direct relation to the considered system actuator limits. Relaxing the state constraints solves the feasibility problem at least for stable systems. As a matter of fact, tightening the state constraints is not, generally, practical because of the inherent presence of noise, disturbances, and numerical errors, see [57] for more details. Moreover, for a given nonlinear system, e.g. (2.1), if $(\mathbf{0}, \mathbf{0}) \in X \times U$, the recursive feasibility of OCP (2.4) holds trivially, i.e. $\mathbf{u}^\star(\cdot) \neq \emptyset$, see [25] for details.

As MPC requires a repeated solution of an optimization problem, it was realized in the early literature as a time consuming control strategy. This constructed a main obstacle for MPC to be applied to physical systems especially those with fast dynamics, e.g. mobile robots. Nonetheless, due to the advancements in computing machines and development of efficient numerical algorithms, a considerable number of dynamic optimization and MPC implementation packages can be found in the literature, e.g. DONLP2 [58], APOPT [59], and NLOPT [60]. Beside the standard MATLAB optimization toolbox, we employed in this thesis recently developed MPC implementation packages, e.g. CasADi [61]. CasADi is a symbolic framework for algorithmic (also known as automatic) differentiation and numeric optimization. This tool provides a low-level framework for quick, yet highly efficient implementation of algorithms for nonlinear numerical optimization. Additionally, we also adopted (ACADO-Toolkit) [62], which implements real-time MPC routines. ACADO-toolkit has been proven to be open source, user friendly, extensible, self-contained, and computationally versatile, when compared with the previously developed optimization packages [63].

Since only finite horizon problems are solved in each MPC step, **_stability_** of the MPC closed loop can not be trivially guaranteed [64]. Several studies have been reported in the literature studying how stability can be ensured, e.g. by imposing terminal equality constraints [65]. The terminal point constraint can be relaxed to a terminal set constraint and a terminal cost added to the objective function [66]. Using stabilizing terminal constraints introduces computational complexity to the online optimization problem. Moreover, when systems with large feasible sets are considered, these methods, typically, require long prediction horizons and hence become more computationally demanding, see [25] for details. Nonetheless, further studies proved stability characteristics of MPC without stabilizing constraints, but with terminal costs, see, e.g. [67]. In this framework, the key parameter that guarantees stability is the prediction horizon length.

Although MPC has been implemented for non-holonomic robots control problems in a considerable number of studies, as observed in Section 2.2, closed loop stability has been considered by only a few. A practical investigation of MPC with terminal equality constraints have been demonstrated in [53]. Additionally, stabilizing MPC using terminal inequality constraints and costs has been pursued in [11] for regulation problems. This has been extended to tracking problems in [22]. Both results have been experimentally demonstrated in [54]. In [52], stability has been studied for regulation and tracking by a contractive MPC design; in this context, stability has been ensured by adding a contractive constraint to the first state of the prediction horizon. Moreover, for the path following problem, stability has been resolved in [6, 7, 68]. Yet, for non-holonomic robots, most of the considered stabilizing MPC designs use stabilizing constraints and/or costs to ensure stability.

### 2.2.4 MPC without stabilizing constraints or costs

As observed in Section 2.2.3, all the studies cited there adopt either stabilizing constraints and/or terminal costs to ensure the stability and performance estimates of the MPC controller. Although such techniques are, in theory, capable of improving the performance of MPC designs, they are normally avoided in practice. Reasons for this are that stabilizing constraints limit the operating region of the controller and may pose problems in numerically solving the optimization problem involved in each MPC step. Moreover, design of terminal costs may be complicated in particular for systems having domains with dynamic obstacles (a typical feature of non-holonomic systems when formation control problems are considered). Moreover, although terminal costs may in principle be used without terminal constraints, they typically provide only a local approximation to the true cost-to-go and thus require terminal constraints in order to ensure that the optimized trajectories end up in a region, where the terminal cost attains meaningful values. Finally, and most importantly, stability-like-behaviour and good performance are often observed without any terminal conditions [69]. Therefore, a new generation of stabilizing MPC controllers, without stabilizing constraints or costs, has been emerged [24].

In the framework of MPC without stabilizing constraints or costs, the prediction horizon length is the main factor in guaranteeing the stability. Therefore, estimates of the prediction horizon length such that stability is guaranteed were investigated in [70–72] based on controllability assumptions, i.e. a finite time controllability of a given system has to be verified. These MPC schems without additional stability enforcing constraints are attractive due to their easiness of implementation. As a matter of fact, this MPC scheme is widely used in the industry. Moreover, the knowledge of a local Lyapunov function–a design requirement in other MPC stabilizing designs–is not a requirement in this MPC scheme. Finally, this MPC design facilitates the derivation of bounds on the infinite horizon performance of the MPC closed loop, see [73, 74] for more details. MPC without stabilizing constraints or costs has been investigated for a number of applications, e.g. non-holonomic integrators [74], two dimensional tank reactors [69], semilinear parabolic partial differential equation (PDE's) [73], and synchronous generators [75].

As an early result of this thesis, MPC without stabilizing constraints or costs has been investigated for the first time for the point stabilization control problem in [50], where the stability results have been presented in the discrete time settings. These results have been extended to the continuous time settings in [51]. In this thesis, we also show an extension of the stability results to the path following control problem.

### 2.2.5 Distributed model predictive control (DMPC)

Although we will present more details on systems with multiple robots later in this chapter, we present a distributed model predictive control (DMPC) here for completeness of the current section. As a motivation for a multi-robot control task, consider a group of robots which first has to set up a particular formation and, then, maintain this pattern while working on the actual control task. As will be seen later in this chapter, this control problem is known as formation regulation.

DMPC provides an appropriate framework for the formation regulation task. In DMPC, the subsystems (robots), in the formation, are not combined into a large overall system (centralized control case), but rather are regulated as independent systems with common (coupled) state-constraints. These constraints account for the inter-robot collision avoidance. In this context, the individual subsystems communicate, every sampling instant, their future predictions to other subsystems. A distributed predictive control provides a flexible and an efficient control algorithm for formation control when compared to the centralized predictive approaches on one hand (due to the distribution of the optimal control problem solution and *scalability*, i.e. the controller applicability is independent of the number of robots in the formation group), and with other distributed non-predictive controllers on the other hand (due to the consideration of the future prediction when calculating the current control actions). See, e.g. [76], for different numerical examples of DMPC.

The DMPC tasks are to: first, provide sufficient conditions under which the formation subsystems are stabilized to their references. Second, determine the required conditions such that the executed optimization problems are feasible, i.e. when optimal solutions exist. Several MPC solutions are presented in the literature for the formation control of mobile robots, see, e.g. [77, 78]

To formulate the coupling constraints in DMPC, previous studies were based on the full communication of predicted trajectories, see [78]. However, relying on sharing all of the prediction data may lead to bandwidth limits and congestion issues in communication. Therefore, in this thesis, we present a method by which the communication load in DMPC is reduced. This method is based on partitioning an operating region into a grid. Then, the predicted trajectories are projected onto the grid resulting in an occupancy grid, which serves as quantization of the communication data. Finally, these data are communicated via a differential communication scheme, which reduces the communication load significantly.

## 2.3 Multi Robot Systems (MRS's)

One of the systems that is considered in this thesis is a heterogeneous MRS, which consists of both ground mobile robots as well as aerial vehicles. In a GPS denied environment, operating such systems is very challenging due to the lack of the localization data, i.e.

the information of the robots' locations. This particular configuration of multi-robot systems attracted a considerable attention in the scientific community, mainly, because of two reasons: first, the high control capabilities and aerial surveillance of modern aerial vehicles, i.e. this class of robots has demonstrated versatile navigation capabilities even in cluttered indoor environment and multi-floor buildings. However, these aerial vehicles have limited payload, computational, and flight time capacities; thus, equipping such systems with high accuracy localization means is, in general, impractical. Second, ground mobile robots are remarked in the literature to have high payload capacity. This allows them to operate for longer periods of time and to be equipped with powerful sensing and actuation devices essential for interacting with the environment. Moreover, ground mobile robots are characterized by a more stable behaviour and relatively slower dynamics when compared to their aerial counterparts. Therefore, these unique complementary characteristics of aerial and ground agents makes a heterogeneous MRS team structures advantageous for missions, see [79] for more details.

In summary, precise localization of the robotic members in a given MRS is a necessity in order to guarantee the success of their missions. As can be inferred from the previous discussions, high accuracy localization can be achieved by equipping all the involved robots with the state-of-the-art sensory means. However, this increases the overall cost of such systems. Moreover, as already observed, some robotic members, e.g. flying/hovering robots may not have enough payload/computational capacity to operate such means. Therefore, relative localization (RL) has been developed as a practical solution for effective and accurate execution of multi-robot collaborative tasks. The objective of relative localization is to detect and locate robots with limited sensory capabilities (referred to as **observed robots**) with respect to other robots with accurate localization means (referred to as **observing robots**). This is achieved by using the relative observation between the two robots categories as will be shown in details in the following subsections, see [80–84]

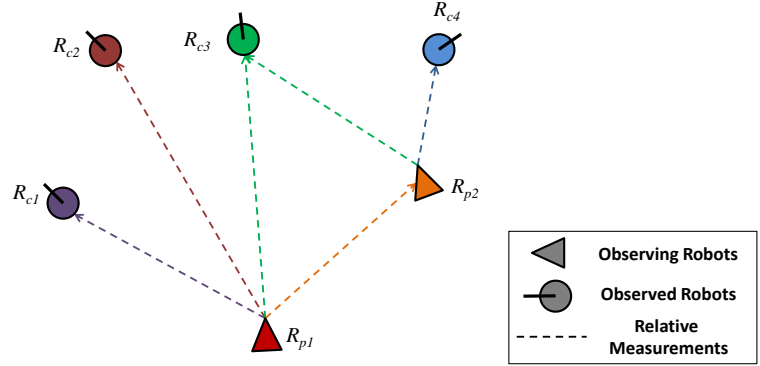

Fig. 2.4: Configuration space of a multi-robot system at a given time instant.

for more details.

### 2.3.1 Relative localization (RL)

In multi-robots systmes, inter-robot relative measurements (IRRM) becomes available in which we consider the range and bearing form of such IRRM. Relative range refers to the Euclidean distance between two robots in a given MRS, while relative bearing refers to the angles that the range vector makes with respect to an observing robot, i.e. azimuth and elevation angles, see Figure 2.4 for an example of an MRS configuration space. IRRM has other different forms, e.g. relative range only, relative bearing only, and full relative pose observation (relative range, relative bearing and relative orientation), see [83] for more details.

Using the IRRM, the subsystems (robots) of a given MRS can perform "Cooperative Localization" [85]. The term cooperative localization generally refers to the sharing of information for localization purposes regardless the choice of navigation frame. A more specific scenario of cooperative localization (known as relative localization (RL)) is illustrated in Figure 2.4. In this framework, less capable (observed) robots rely on the IRRM in order to be detected and localized with respect to the more capable (observing) robots frames or references. The term Relative Localization, which has different counterpart notations, e.g. Mutual Localization, and Network Localization, will be used throughout the text of the thesis, see [79] for details.

In all attempts at solving the relative localization problem, extended Kalman filter (EKF) is the most commonly adopted nonlinear estimator [81]. However, inappropriate initialization of EKF can generally lead to instability of the estimator as well as longer estimation settling time, i.e. time required to reach acceptable estimation error levels, and this results in misleading relative pose estimates. Indeed, this erroneous localization causes undesirable behaviour and possibly failure in collaborative missions. In order to avoid these issues, the majority of the past work assumed a ***known*** initial relative pose between two arbitrary robots, see, e.g. [86]. In addition to EKF, particle filter (PF) [87], unscented Kalman filter (UKF) [88], and pseudolinear Kalman filter [82, 89] have been also used to achieve relative localization.

### 2.3.2 Moving horizon estimation (MHE) and RL

The previously highlighted estimators, except the particle filter, use the Gaussian probability distribution to approximate a given state noise; the PF instead approximates the distribution via Monte Carlo sampling [90]. Moreover, these filters employ Markov property in order to reduce the computational complexity, i.e. the current state estimate is based only on the most recent measurement and the previous state estimate. In contrast to the previous methods, it is proposed in this thesis to employ a moving horizon estimation (MHE) scheme to solve the RL problem. MHE considers the evolution of a constrained and possibly nonlinear model on a fixed time horizon, and minimizes the deviation of the model from a past window of measurements [91]. Therefore, MHE relaxes the Markov assumption highlighted above. Although MHE does not generally rely on any specific error distribution, tuning the estimator becomes easier when a probabilistic insight is considered, see [92, 93] for details, and [94] for a relative localization example using

MHE in autonomous underwater vehicles. MHE has not been realized in the literature when applied to the relative localization problem in MRS. Therefore, an early contribution, resulted from this thesis, showed its applicability to MRS relative localization with real-time constraints satisfaction, see [84].

### 2.3.3 MHE mathematical formulation

Here, we show the generic formulation of a moving horizon estimation (MHE) scheme. To this end, we define the disturbed form of system (2.1) and we include an output, i.e $\mathbf{y} \in \mathbb{R}^{n_y}$, equation as the following

$$\mathbf{x}(k+1) = \mathbf{f}(\mathbf{x}(k), \mathbf{u}(k)) + \boldsymbol{\nu}_{\mathbf{x}}, \tag{2.9}$$
$$\mathbf{y}(k) = \mathbf{h}(\mathbf{x}(k)) + \boldsymbol{\nu}_{\mathbf{y}},$$

where in the above model we defined the zero mean Gaussian noises for state $\boldsymbol{\nu}_{\mathbf{x}}$ and output $\boldsymbol{\nu}_{\mathbf{y}}$, respectively. Then, the MHE state estimator is formulated as a least squares (LSQ) cost function $J_{N_E} : X \times U^{N_E} \to \mathbb{R}_{\geq 0}$ shown in (2.10). Hence, using (2.10), the MHE repeatedly solves the constrained nonlinear dynamic optimization problem (2.11) over a fixed estimation horizon of length $N_E \in \mathbb{N}$ [93]. In (2.11), for $q \in \mathbb{N}_0$, the control sequence $\mathbf{u}$ is defined by

$$\mathbf{u} = (\mathbf{u}(q - N_E), \mathbf{u}(q - N_E + 1), \ldots, \mathbf{u}(q - 1)) \in U^{N_E}.$$

$$J_{N_E}(\mathbf{x}(q - N_E), \mathbf{u}) = \|\mathbf{x}(q - N_E) - \mathbf{x}^{\text{est}}_{q - N_E}\|^2_{\mathbf{A}} \\ + \sum_{k=q-N_E}^{q} \|\mathbf{h}(\mathbf{x}(k)) - \tilde{\mathbf{y}}(k)\|^2_{\mathbf{B}} + \sum_{k=q-N_E}^{q-1} \|\mathbf{u}(k) - \tilde{\mathbf{u}}(k)\|^2_{\mathbf{C}} \tag{2.10}$$

$$\min_{\mathbf{x}(q-N_E), \mathbf{u} \in \mathbb{R}^{n_u \times N_E}} J_{N_E}(\mathbf{x}(q-N_E), \mathbf{u}) \qquad (2.11)$$

subject to:

$$\mathbf{x}(k+1) = \mathbf{f}(\mathbf{x}(k), \mathbf{u}(k)) \qquad \forall k \in \{q-N_E, \ldots, q-1\},$$

$$\mathbf{x}(k) \in X \qquad \forall k \in \{q-N_E, \ldots, q\},$$

$$\mathbf{u}(k) \in U \qquad \forall k \in \{q-N_E, \ldots, q-1\},$$

In (2.10), $\tilde{\mathbf{y}}$ and $\tilde{\mathbf{u}}$ denote the actually measured system outputs and inputs, respectively. The first term in (2.10) is known as the arrival cost and it penalizes the deviation of the first state in the moving horizon window and its priori estimate $\mathbf{x}^{\text{est}}_{q-N_E}$ by the diagonal positive definite matrix $\mathbf{A} \in \mathbb{R}^{n_x \times n_x}$. Normally, the estimate $\mathbf{x}^{\text{est}}_{q-N_E}$ is adopted from the previous MHE estimation step. Moreover, the weighting matrix $\mathbf{A}$ is chosen as a smoothed EKF-update based on sensitivity information gathered while solving the most recent MHE step. Therefore, $\mathbf{x}^{\text{est}}_{q-N_E}$ and $\mathbf{A}$ initializations are required, see [93, 95] for more details.

The second term in (2.10) penalizes the change in the system predicted outputs $\mathbf{h}(\mathbf{x})$ from the actually measured outputs $\tilde{\mathbf{y}}$ by the diagonal positive-definite matrix $\mathbf{B} \in \mathbb{R}^{n_y \times n_y}$. Similarly, the change in the applied control inputs $\mathbf{u}$ from the measured inputs $\tilde{\mathbf{u}}$ is penalized using the diagonal positive-definite matrix $\mathbf{C} \in \mathbb{R}^{n_u \times n_u}$. The latter term is included in the cost function (2.10) to account for actuator noise and/or inaccuracy, see [92] for details. $\mathbf{B}$ and $\mathbf{C}$ are chosen to match the applied motion and measurement noise covariances. Furthermore, all estimated quantities in (2.11) are subject bounds, which signify the system physical limitations. The solution of the optimization problem (2.11) leads mainly to the relative-state estimate sequence $\hat{\mathbf{x}}(k), k = (q - N_E, \cdots, q)$, where $\hat{\mathbf{x}}(0) := \hat{\mathbf{x}}(q)$ denotes the current estimate of the relative state vector of a given observed robot; $\hat{\mathbf{x}}(0)$ can be later used as a measurement feedback. Moreover, estimate sequence of the actually applied control input, i.e. $\hat{\mathbf{u}}$ is a byproduct of the MHE estimator. Algorithm 2.2 summarizes the MHE scheme for nonlinear systems. See also Figure 2.5 for an illustration of two successive steps of state estimation using MHE for a simple single input single output (SISO) system.

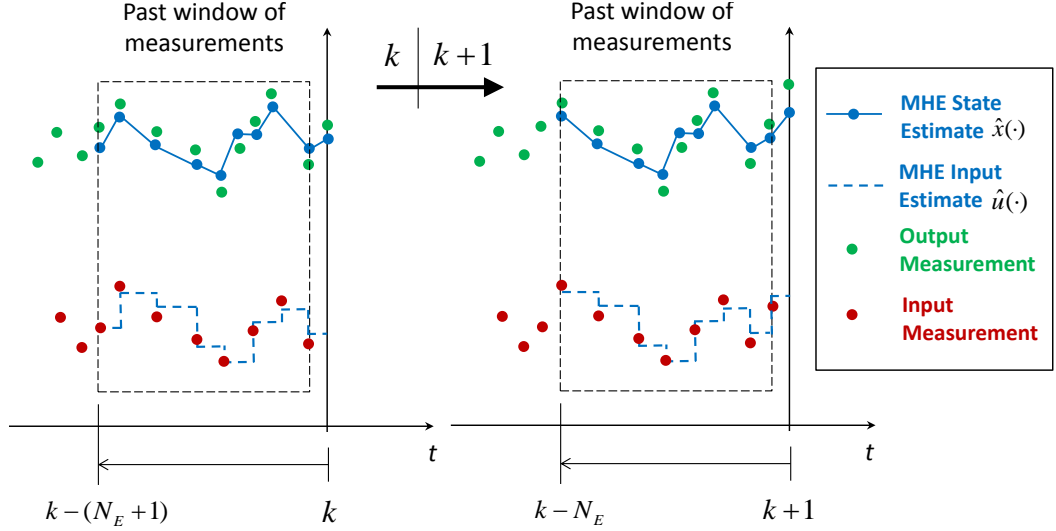

Fig. 2.5: Illustration of two MHE iterations for a simple SISO system in which the measured output is the system state. $k$ denotes the time step and $N_E$ denotes the estimation horizon length.

---
**Algorithm 2.2** MHE scheme steps
---
1: **for** every sampling instant $n = 0, 1, 2, \ldots$ **do**
2:     Get the past $N_E$ input measurements, i.e. $\tilde{\mathbf{u}}$.
3:     Get the past $N_E + 1$ relative measurements, i.e. $\tilde{\mathbf{y}}$.
4:     Solve the optimization problem (2.11) over estimation horizon $N_E$, and find the best current state estimate, i.e. $\hat{\mathbf{x}}(0)|_n$.
5: **end for**
---

### 2.3.4 Control problems in MRS

As a complementary part of this section, we present the control problems involved in multi-robot systems. The control objectives within a given MRS formation structure are classified under virtual structure, behaviour based, and leader follower approaches. In the virtual structure approach, the individual robots in the MRS formation are treated as points on a rigid body where their individual motion are determined based on the formation overall motion. In behavior-based control, several desired behaviors are assigned to each robot in the MRS; the behaviours include fromation keeping (stabilization), goal seeking, and obstacle avoidance. Finally, in the approach of the leader follower structure, a follower robot is assigned to follow the pose of a leader robot with an offset, see the overview articles [96–98] for more details.

The solution to the virtual structure approach is presented in [99], where authors designed an iterative algorithm that fits the virtual structure to the positions of the individual robots in the formation, then displaces the virtual structure in a prescribed direction and then updates the robot positions. Solutions to the behaviour based approach include motor scheme control [100], null-space-based behaviour control [101], and social potential fields control [102]. The leader-follower control approach is achieved using feedback linearization [103], backstepping control [104], sliding mode control [105], and model predictive control [106].

In this thesis, based on relative localization achieved using MHE, it is proposed to use a centralized nonlinear model predictive control (MPC) to perform a relative trajectory tracking task where one or more robots in a given MRS are commanded to follow time varying trajectories with respect to another robot reference frame. In addition, inter-robot possible collisions must be avoided using the adopted controller. This control problem can be categorized under the leader follower structure, with the exception that the relative references of the following robots are time varying. In fact, this study was carried out with both simulations and experimental results and was published in [107].

# References

[1] J. P. Laumond, S. Sekhavat, F. Lamiraux, J. paul Laumond (editor, J. P. Laumond, S. Sekhavat, and F. Lamiraux, "Guidelines in nonholonomic motion planning for mobile robots," in *Robot Motion Plannning and Control*. Springer-Verlag, 1998, pp. 1–53.

[2] M. Michalek and K. Kozowski, "Vector-field-orientation feedback control method for a differentially driven vehicle," *IEEE Transactions on Control Systems Technology*, vol. 18, no. 1, pp. 45–65, 2010.

[3] R. W. Brockett, "Asymptotic stability and feedback stabilization," in *Differential Geometric Control Theory*, R. W. Brockett, R. S. Millman, and H. J. Sussmann, Eds. Birkhäuser, Boston, MA, 1983, pp. 181–191.

[4] R. Siegwart, I. R. Nourbakhsh, and D. Scaramuzza, *Introduction to autonomous mobile robots*. MIT press, 2011.

[5] T. Faulwasser, "Optimization-based solutions to constrained trajectory-tracking and path-following problems," Ph.D. dissertation, Otto-von-Guericke-Universität Magdeburg, 2012.

[6] T. Faulwasser, B. Kern, and R. Findeisen, "Model predictive path-following for constrained nonlinear systems," in *Proceedings of the 48th IEEE Conference on Decision and Control (CDC)*, 2009, pp. 8642–8647.

[7] T. Faulwasser, *Optimization-based Solutions to Constrained Trajectory-tracking and Path-following Problems*. Shaker, Aachen, Germany, 2013.

[8] Z.-P. Jiang and H. Nijmeijer, "Tracking control of mobile robots: A case study in backstepping," *Automatica*, vol. 33, no. 7, pp. 1393 – 1399, 1997.

[9] G. Oriolo, A. De Luca, and M. Vendittelli, "Wmr control via dynamic feedback linearization: design, implementation, and experimental validation," *Control Systems Technology, IEEE Transactions on*, vol. 10, no. 6, pp. 835–852, 2002.

[10] D. Chwa, "Sliding-mode tracking control of nonholonomic wheeled mobile robots in polar coordinates," *Control Systems Technology, IEEE Transactions on*, vol. 12, no. 4, pp. 637–644, July 2004.

[11] D. Gu and H. Hu, "A Stabilizing Receding Horizon Regulator for Nonholonomic Mobile Robots," *IEEE Transactions on Robotics*, vol. 21, no. 5, pp. 1022–1028, 2005.

[12] C. de Wit and O. Sordalen, "Exponential stabilization of mobile robots with nonholonomic constraints," *Automatic Control, IEEE Transactions on*, vol. 37, no. 11, pp. 1791–1797, Nov 1992.

[13] R. M'Closkey and R. Murray, "Exponential stabilization of driftless nonlinear control systems using homogeneous feedback," *Automatic Control, IEEE Transactions on*, vol. 42, no. 5, pp. 614–628, May 1997.

[14] G. Indiveri, "Kinematic time-invariant control of a 2d nonholonomic vehicle," in *Decision and Control, 1999. Proceedings of the 38th IEEE Conference on*, vol. 3, 1999, pp. 2112–2117 vol.3.

[15] W. Dixon, D. Dawson, F. Zhang, and E. Zergeroglu, "Global exponential tracking control of a mobile robot system via a pe condition," in *Proceedings of the 38th IEEE Conference on Decision and Control*, vol. 5, 1999, pp. 4822–4827.

[16] T.-C. Lee, K.-T. Song, C.-H. Lee, and C.-C. Teng, "Tracking Control of Unicycle-Modeled Mobile Robots Using a Saturation Feedback Controller," *IEEE Transactions on Control Systems Technology*, vol. 9, no. 2, pp. 305–318, 2001.

[17] F. Kühne, W. F. Lages, and J. M. Gomes da Silva Jr., "Point Stabilization of Mobile Robots with Nonlinear Model Predictive Control," in *Proceedings of the IEEE International Conference on Mechatronics and Automation*, vol. 3, 2005, pp. 1163–1168 Vol. 3.

[18] A. P. Aguiar, J. P. Hespanha, and P. V. Kokotović, "Performance limitations in reference tracking and path following for nonlinear systems," *Automatica*, vol. 44, no. 3, pp. 598 – 610, 2008.

[19] C. Nielsen, C. Fulford, and M. Maggiore, "Path following using transverse feedback linearization: Application to a maglev positioning system," *Automatica*, vol. 46, no. 3, pp. 585–590, Mar. 2010.

[20] W. Dixon, *Nonlinear Control of Wheeled Mobile Robots*, ser. Lecture Notes in Control and Information Sciences.   Springer, 2001.

[21] G. Cook, *Mobile Robots: Navigation, Control and Remote Sensing.*   Wiley, 2011.

[22] D. Gu and H. Hu, "Receding horizon tracking control of wheeled mobile robots," *IEEE Transactions on Control Systems Technology*, vol. 14, no. 4, pp. 743–749, 2006.

[23] J. B. Rawlings and D. Q. Mayne, *Model Predictive Control: Theory and Design.* Nob Hill Publishing, 2009.

[24] L. Grüne and J. Pannek, *Nonlinear Model Predictive Control: Theory and Algorithms*, ser. Communications and Control Engineering.   Springer London Dordrecht Heidelberg New York, 2011.

[25] ——, *Nonlinear Model Predictive Control: Theory and Algorithms*, ser. Communications and Control Engineering, A. Isidori, J. H. van Schuppen, E. D. Sontag,

M. Thoma, and M.Krstic, Eds. Springer London Dordrecht Heidelberg New York, 2011.

[26] M. Alamir, *Stabilization of Nonlinear Systems Using Receding-horizon Control Schemes: A Parametrized Approach for Fast Systems*, ser. Lecture notes in control and information sciences. Springer, 2006.

[27] J. Bartee, P. Noll, C. Axelrud, C. Schweiger, and B. Sayyar-Rodsari, "Industrial application of nonlinear model predictive control technology for fuel ethanol fermentation process," in *American Control Conference, 2009. ACC '09.*, June 2009, pp. 2290–2294.

[28] P. Mhaskar, "Robust model predictive control design for fault-tolerant control of process systems," *Industrial & Engineering Chemistry Research*, vol. 45, no. 25, pp. 8565–8574, 2006.

[29] J. Siroky, F. Oldewurtel, J. Cigler, and S. Prívara, "Experimental analysis of model predictive control for an energy efficient building heating system," *Applied Energy*, vol. 88, no. 9, pp. 3079 – 3087, 2011.

[30] W. H. Lio, J. Rossiter, and B. Jones, "A review on applications of model predictive control to wind turbines," in *Control (CONTROL), 2014 UKACC International Conference on*, July 2014, pp. 673–678.

[31] P. Cortes, M. Kazmierkowski, R. Kennel, D. Quevedo, and J. Rodriguez, "Predictive control in power electronics and drives," *Industrial Electronics, IEEE Transactions on*, vol. 55, no. 12, pp. 4312–4324, Dec 2008.

[32] S. Kouro, P. Cortes, R. Vargas, U. Ammann, and J. Rodriguez, "Model predictive control – a simple and powerful method to control power converters," *Industrial Electronics, IEEE Transactions on*, vol. 56, no. 6, pp. 1826–1838, June 2009.

[33] K. Alexis, G. Nikolakopoulos, and A. Tzes, "Switching model predictive attitude control for a quadrotor helicopter subject to atmospheric disturbances," *Control Engineering Practice*, vol. 19, no. 10, pp. 1195 – 1207, 2011.

[34] K. Dalamagkidis, K. Valavanis, and L. Piegl, "Nonlinear model predictive control with neural network optimization for autonomous autorotation of small unmanned helicopters," *Control Systems Technology, IEEE Transactions on*, vol. 19, no. 4, pp. 818–831, July 2011.

[35] K. Alexis, G. Nikolakopoulos, and A. Tzes, "Model predictive quadrotor control: attitude, altitude and position experimental studies," *Control Theory Applications, IET*, vol. 6, no. 12, pp. 1812–1827, Aug 2012.

[36] P. Bouffard, A. Aswani, and C. Tomlin, "Learning-based model predictive control on a quadrotor: Onboard implementation and experimental results," in *Robotics and Automation (ICRA), 2012 IEEE International Conference on*, May 2012, pp. 279–284.

[37] M. Abdolhosseini, Y. Zhang, and C. Rabbath, "An efficient model predictive control scheme for an unmanned quadrotor helicopter," *Journal of Intelligent & Robotic Systems*, vol. 70, no. 1-4, pp. 27–38, 2013.

[38] K. Alexis, G. Nikolakopoulos, and A. Tzes, "On trajectory tracking model predictive control of an unmanned quadrotor helicopter subject to aerodynamic disturbances," *Asian Journal of Control*, vol. 16, no. 1, pp. 209–224, 2014.

[39] J. Eklund, J. Sprinkle, and S. Sastry, "Implementing and testing a nonlinear model predictive tracking controller for aerial pursuit/evasion games on a fixed wing aircraft," in *American Control Conference, 2005. Proceedings of the 2005*, June 2005, pp. 1509–1514 vol. 3.

[40] Y. Kang and J. Hedrick, "Linear tracking for a fixed-wing uav using nonlinear model predictive control," *Control Systems Technology, IEEE Transactions on*, vol. 17, no. 5, pp. 1202–1210, Sept 2009.

[41] F. Fahimi, "Non-linear model predictive formation control for groups of autonomous surface vessels," *International Journal of Control*, vol. 80, no. 8, pp. 1248–1259, 2007.

[42] S.-R. Oh and J. Sun, "Path following of underactuated marine surface vessels using line-of-sight based model predictive control," *Ocean Engineering*, vol. 37, no. 2–3, pp. 289 – 295, 2010.

[43] M. Gross and M. Costello, "Impact point model predictive control of a spin-stabilized projectile with instability protection," *Journal of Aerospace Engineering*, 2013.

[44] H. B. Oza and R. Padhi, "Impact-angle-constrained suboptimal model predictive static programming guidance of air-to-ground missiles," *Journal of Guidance, Control, and Dynamics*, vol. 35, no. 1, pp. 153–164, 2012.

[45] A. Vivas and P. Poignet, "Predictive functional control of a parallel robot," *Control Engineering Practice*, vol. 13, no. 7, pp. 863 – 874, 2005, control Applications of Optimisation.

[46] K. Kanjanawanishkul and A. Zell, "Path following for an omnidirectional mobile robot based on model predictive control," in *Robotics and Automation, 2009. ICRA '09. IEEE International Conference on*, May 2009, pp. 3341–3346.

[47] P. Falcone, M. Tufo, F. Borrelli, J. Asgari, and H. Tsengz, "A linear time varying model predictive control approach to the integrated vehicle dynamics control prob-

lem in autonomous systems," in *Decision and Control, 2007 46th IEEE Conference on*, Dec 2007, pp. 2980–2985.

[48] G. Klancar and I. Skrjanc, "Tracking-error model-based predictive control for mobile robots in real time," *Robotics and Autonomous Systems*, vol. 55, no. 6, pp. 460 – 469, 2007.

[49] H. Lim, Y. Kang, C. Kim, J. Kim, and B.-J. You, "Nonlinear model predictive controller design with obstacle avoidance for a mobile robot," in *Mechtronic and Embedded Systems and Applications, 2008. MESA 2008. IEEE/ASME International Conference on*, Oct 2008, pp. 494–499.

[50] K. Worthmann, M. W. Mehrez, M. Zanon, G. K. Mann, R. G. Gosine, and M. Diehl, "Model predictive control of nonholonomic mobile robots without stabilizing constraints and costs," *Control Systems Technology, IEEE Transactions on*, vol. PP, no. 99, pp. 1–13, 2015.

[51] ——, "Regulation of differential drive robots using continuous time mpc without stabilizing constraints or costs," in *Proceedings of the 5th IFAC Conference on Nonlinear Model Predictive Control (NPMC'15), Sevilla, Spain*, 2015, pp. 129–135.

[52] F. Xie and R. Fierro, "First-state contractive model predictive control of nonholonomic mobile robots," in *Proceedings of the American Control Conference*, 2008, pp. 3494–3499.

[53] M. W. Mehrez, G. K. I. Mann, and R. G. Gosine, "Stabilizing nmpc of wheeled mobile robots using open-source real-time software," in *Advanced Robotics (ICAR), 2013 16th International Conference on*, Nov 2013, pp. 1–6.

[54] ——, "Comparison of stabilizing nmpc designs for wheeled mobile robots: an experimental study," in *Proceedings of Moratuwa Engineering Research Conference (MERCon)*, 2015.

[55] G. Raffo, G. Gomes, J. Normey-Rico, C. Kelber, and L. Becker, "A predictive controller for autonomous vehicle path tracking," *Intelligent Transportation Systems, IEEE Transactions on*, vol. 10, no. 1, pp. 92–102, March 2009.

[56] J. Backman, T. Oksanen, and A. Visala, "Navigation system for agricultural machines: Nonlinear model predictive path tracking," *Computers and Electronics in Agriculture*, vol. 82, no. 0, pp. 32 – 43, 2012.

[57] K. Kanjanawanishkul, "Coordinated path following control and formation control of mobile robots," Ph.D. dissertation, University of Tübingen, 2010.

[58] P. Spellucci, "An sqp method for general nonlinear programs using only equality constrained subproblems," *Mathematical Programming*, vol. 82, no. 3, pp. 413–448, 1998.

[59] L. Advanced Process Solutions, "Advanced process optimizer (apopt)," http://apopt.com/index.php.

[60] S. G. Johnson, "The nlopt nonlinear-optimization package," http://ab-initio.mit.edu/nlopt.

[61] J. Andersson, "A General-Purpose Software Framework for Dynamic Optimization," PhD thesis, Arenberg Doctoral School, KU Leuven, Department of Electrical Engineering (ESAT/SCD) and Optimization in Engineering Center, Kasteelpark Arenberg 10, 3001-Heverlee, Belgium, October 2013.

[62] B. Houska, H. Ferreau, and M. Diehl, "ACADO Toolkit – An Open Source Framework for Automatic Control and Dynamic Optimization," *Optimal Control Applications and Methods*, vol. 32, no. 3, pp. 298–312, 2011.

[63] B. Houska, H. J. Ferreau, and M. Diehl, "An auto-generated real-time iteration algorithm for nonlinear {MPC} in the microsecond range," *Automatica*, vol. 47, no. 10, pp. 2279 – 2285, 2011.

[64] T. Raff, S. Huber, Z. K. Nagy, and F. Allgöwer, "Nonlinear model predictive control of a four tank system: An experimental stability study," in *Proceedings of the IEEE Conference on Control Applications*, Munich, Germany, 2006, pp. 237–242.

[65] S. Keerthi and E. Gilbert, "Optimal infinite-horizon feedback laws for a general class of constrained discrete-time systems: Stability and moving-horizon approximations," *Journal of Optimization Theory and Applications*, vol. 57, no. 2, pp. 265–293, 1988.

[66] H. Michalska and D. Mayne, "Robust receding horizon control of constrained nonlinear systems," *Automatic Control, IEEE Transactions on*, vol. 38, no. 11, pp. 1623–1633, Nov 1993.

[67] G. Grimm, M. Messina, S. Tuna, and A. Teel, "Model predictive control: for want of a local control Lyapunov function, all is not lost," *IEEE Transactions on Automatic Control*, vol. 50, no. 5, pp. 546–558, 2005.

[68] T. Faulwasser and R. Findeisen, "Nonlinear model predictive control for constrained output path following," *Cond. accepted for publication in IEEE Transactions on Automatic Control*, 2015.

[69] L. Grüne and M. Stieler, "Asymptotic stability and transient optimality of economic {MPC} without terminal conditions," *Journal of Process Control*, vol. 24, no. 8, pp. 1187 – 1196, 2014.

[70] S. Tuna, M. Messina, and A. Teel, "Shorter horizons for model predictive control," in *Proceedings of the American Control Conference*, 2006, pp. 863–868.

[71] L. Grüne, "Analysis and design of unconstrained nonlinear MPC schemes for finite and infinite dimensional systems," *SIAM J. Control Optim.*, vol. 48, no. 2, pp. 1206–1228, 2009.

[72] L. Grüne, J. Pannek, M. Seehafer, and K. Worthmann, "Analysis of unconstrained nonlinear MPC schemes with varying control horizon," *SIAM J. Control Optim.*, vol. Vol. 48 (8), pp. 4938–4962, 2010.

[73] N. Altmüller, L. Grüne, and K. Worthmann, "Performance of nmpc schemes without stabilizing terminal constraints," in *Recent Advances in Optimization and its Applications in Engineering*, M. Diehl, F. Glineur, E. Jarlebring, and W. Michiels, Eds. Springer Berlin Heidelberg, 2010, pp. 289–298.

[74] K. Worthmann, M. Reble, L. Grüne, and F. Allgöwer, "The role of sampling for stability and performance in unconstrained nonlinear model predictive control," *SIAM Journal on Control and Optimization*, vol. 52, no. 1, pp. 581–605, 2014.

[75] P. Braun, J. Pannek, and K. Worthmann, "Predictive Control Algorithms: Stability Despite Shortened Optimization Horizons," in *Proceedings on the 15th IFAC Workshop on Control Applications of Optimization CAO 2012, Rimini, Italy*, 2012, pp. 284–289.

[76] M. A. Müller, M. Reble, and F. Allgöwer, "Cooperative control of dynamically decoupled systems via distributed model predictive control," *International Journal of Robust and Nonlinear Control*, vol. 22, no. 12, pp. 1376–1397, 2012.

[77] S.-M. Lee, H. Kim, and H. Myung, "Cooperative coevolution-based model predictive control for multi-robot formation," in *Proceedings of IEEE International Conference on Robotics and Automation (ICRA)*, 2013, pp. 1890–1895.

[78] L. Grüne and K. Worthmann, "A distributed nmpc scheme without stabilizing terminal constraints," in *Distributed Decision Making and Control*, ser. Lecture Notes in Control and Information Sciences, R. Johansson and A. Rantzer, Eds. Springer London, 2012, vol. 417, pp. 261–287.

[79] O. D. Silva, "Relative localization of gps denied ground aerial multi-robot systems," Ph.D. dissertation, Memorial University of Newfoundland, 2015.

[80] F. Rivard, J. Bisson, F. Michaud, and D. Letourneau, "Ultrasonic relative positioning for multi-robot systems," in *Robotics and Automation, 2008. ICRA 2008. IEEE International Conference on*, 2008, pp. 323–328.

[81] O. De Silva, G. Mann, and R. Gosine, "Development of a relative localization scheme for ground-aerial multi-robot systems," in *Intelligent Robots and Systems (IROS), 2012 IEEE/RSJ International Conference on*, 2012, pp. 870–875.

[82] T. Wanasinghe, G. I. Mann, and R. Gosine, "Relative localization approach for combined aerial and ground robotic system," *Journal of Intelligent and Robotic Systems*, vol. 77, no. 1, pp. 113–133, 2015.

[83] T. Wanasinghe, G. Mann, and R. Gosine, "Distributed leader-assistive localization method for a heterogeneous multirobotic system," *Automation Science and Engineering, IEEE Transactions on*, vol. 12, no. 3, pp. 795–809, July 2015.

[84] M. W. Mehrez, G. K. Mann, and R. G. Gosine, "Nonlinear moving horizon state estimation for multi-robot relative localization," in *Electrical and Computer Engineering (CCECE), 2014 IEEE 27th Canadian Conference on*, May 2014, pp. 1–5.

[85] A. Mourikis and S. Roumeliotis, "Performance analysis of multirobot cooperative localization," *Robotics, IEEE Transactions on*, vol. 22, no. 4, pp. 666–681, Aug 2006.

[86] J. Fenwick, P. Newman, and J. Leonard, "Cooperative concurrent mapping and localization," in *Robotics and Automation, 2002. Proceedings. ICRA '02. IEEE International Conference on*, vol. 2, 2002, pp. 1810–1817.

[87] L. Carlone, M. Kaouk Ng, J. Du, B. Bona, and M. Indri, "Simultaneous localization and mapping using rao-blackwellized particle filters in multi robot systems," *Journal of Intelligent and Robotic Systems*, vol. 63, no. 2, pp. 283–307, 2011.

[88] S. Xingxi, W. Tiesheng, H. Bo, and Z. Chunxia, "Cooperative multi-robot localization based on distributed ukf," in *Computer Science and Information Technology (ICCSIT), 2010 3rd IEEE International Conference on*, vol. 6, 2010, pp. 590–593.

[89] C. Pathiranage, K. Watanabe, B. Jayasekara, and K. Izumi, "Simultaneous localization and mapping: A pseudolinear kalman filter (plkf) approach," in *Information and Automation for Sustainability, 2008. ICIAFS 2008. 4th International Conference on*, 2008, pp. 61–66.

[90] J. B. Rawlings and B. R. Bakshi, "Particle filtering and moving horizon estimation," *Computers and Chemical Engineering*, vol. 30, no. 10–12, pp. 1529 – 1541, 2006.

[91] C. Rao and J. Rawlings, "Nonlinear moving horizon state estimation," in *Nonlinear Model Predictive Control*, ser. Progress in Systems and Control Theory, F. Allgöwer and A. Zheng, Eds. Birkhäuser Basel, 2000, vol. 26, pp. 45–69.

[92] M. Zanon, J. Frasch, and M. Diehl, "Nonlinear moving horizon estimation for combined state and friction coefficient estimation in autonomous driving," in *Control Conference (ECC), 2013 European*, 2013, pp. 4130–4135.

[93] H. Ferreau, T. Kraus, M. Vukov, W. Saeys, and M. Diehl, "High-speed moving horizon estimation based on automatic code generation," in *Decision and Control (CDC), 2012 IEEE 51st Annual Conference on*, 2012, pp. 687–692.

[94] S. Wang, L. Chen, D. Gu, and H. Hu, "An optimization based moving horizon estimation with application to localization of autonomous underwater vehicles," *Robotics and Autonomous Systems*, vol. 62, no. 10, pp. 1581 – 1596, 2014.

[95] E. Kayacan, E. Kayacan, H. Ramon, and W. Saeys, "Learning in centralized nonlinear model predictive control: Application to an autonomous tractor-trailer system," *Control Systems Technology, IEEE Transactions on*, vol. 23, no. 1, pp. 197–205, Jan 2015.

[96] M. Saffarian and F. Fahimi, "Non-iterative nonlinear model predictive approach applied to the control of helicopters group formation," *Robotics and Autonomous Systems*, vol. 57, no. 6-7, pp. 749 – 757, 2009.

[97] M. W. Mehrez, G. K. Mann, and R. G. Gosine, "Formation stabilization of nonholonomic robots using nonlinear model predictive control," in *Electrical and Computer Engineering (CCECE), 2014 IEEE 27th Canadian Conference on*, May 2014, pp. 1–6.

[98] W. Dunbar and R. Murray, "Model predictive control of coordinated multi-vehicle formations," in *Decision and Control, 2002, Proceedings of the 41st IEEE Conference on*, vol. 4, 2002, pp. 4631–4636 vol.4.

[99] M. Lewis and K.-H. Tan, "High precision formation control of mobile robots using virtual structures," *Autonomous Robots*, vol. 4, no. 4, pp. 387–403, 1997.

[100] T. Balch and R. Arkin, "Behavior-based formation control for multirobot teams," *Robotics and Automation, IEEE Transactions on*, vol. 14, no. 6, pp. 926–939, 1998.

[101] G. Antonelli, F. Arrichiello, and S. Chiaverini, "Flocking for multi-robot systems via the null-space-based behavioral control," in *Intelligent Robots and Systems, 2008. IROS 2008. IEEE/RSJ International Conference on*, 2008, pp. 1409–1414.

[102] T. Balch and M. Hybinette, "Social potentials for scalable multi-robot formations," in *Robotics and Automation, 2000. Proceedings. ICRA '00. IEEE International Conference on*, vol. 1, 2000, pp. 73–80 vol.1.

[103] G. Mariottini, F. Morbidi, D. Prattichizzo, G. Pappas, and K. Daniilidis, "Leader-follower formations: Uncalibrated vision-based localization and control," in *Robotics and Automation, 2007 IEEE International Conference on*, 2007, pp. 2403–2408.

[104] X. Li, J. Xiao, and Z. Cai, "Backstepping based multiple mobile robots formation control," in *Intelligent Robots and Systems, 2005. (IROS 2005). 2005 IEEE/RSJ International Conference on*, 2005, pp. 887–892.

[105] J. Sanchez and R. Fierro, "Sliding mode control for robot formations," in *Intelligent Control. 2003 IEEE International Symposium on*, 2003, pp. 438–443.

[106] D. Gu and H. Hu, "A model predictive controller for robots to follow a virtual leader," *Robotica*, vol. 27, pp. 905–913, 10 2009.

[107] M. W. Mehrez, G. K. I. Mann, and R. G. Gosine, "An optimization based approach for relative localization and relative tracking control in multi-robot systems," *Journal of Intelligent and Robotic Systems*, pp. 1–24, 2015.

# Chapter 3

# Model Predictive Control of Non-holonomic Mobile Robots without Stabilizing Constraints and Costs

## 3.1 Abstract

The problem of steering a non-holonomic mobile robot to a desired position and orientation is considered. In this chapter, a model predictive control (MPC) scheme based on tailored non-quadratic stage cost is proposed to fulfil this control task. We rigorously prove asymptotic stability while neither stabilizing constraints nor costs are used. To this end, we first design suitable maneuvers to construct bounds on the value function. Second, these bounds are exploited to determine a prediction horizon length such that asymptotic stability of the MPC closed-loop is guaranteed. Finally, numerical simulations are conducted to explain the necessity of having non-quadratic running costs.

## 3.2 Introduction

Unmanned ground vehicles (UGVs) have attracted considerable interest in the recent decades due to their wide range of applicability, see [1] or [2] for a thorough review. Non-holonomic differential-drive models, such as unicycle models, are commonly used to describe the kinematics of UGVs. Typically, the control objective is to drive the robot between two static poses, which can be identified as set-point stabilization (regulation) [3]. For this problem, Brockett's condition [4] implies that neither the linearized model is stabilizable nor a smooth time-invariant feedback control law exists – a typical characteristic of non-holonomic systems, see also [5]. Nonetheless, various solution strategies like piecewise-continuous feedback control or smooth time-varying control have been reported, see the overview paper [6]. Further control approaches based on differential kinematic control [7], backstepping [8], and vector field orientation feedback [9] have also been proposed. However, these control strategies ignore natural input saturation limits and, thus, require a post processing step in order to scale the calculated control signals to their physical bounds, see [9] for details. In addition, determining suitable tuning parameters in order to achieve an acceptable performance remains a challenging task [10]. In contrast, several successful case studies using model predictive control (MPC) were conducted, see, e.g. [3, 6, 11–13].

MPC is considered to be one of the most attractive control strategies due to its applicability to constrained nonlinear multiple input multiple output (MIMO) systems. In MPC, a sequence of control inputs minimizing an objective function is computed over a finite prediction horizon; then, the first element of this (optimal) control sequence is applied to the plant. This process is repeated every sampling instant, see, e.g. [14], for further details. Since only finite horizon problems are solved in each MPC step, closed-loop stability may not hold [15]. Nonetheless, stability can be ensured, e.g. by imposing terminal constraints [16, 17], or by using bounds on the value function in order to determine a stabilizing prediction horizon length, see, e.g. [18–20].

For regulation of non-holonomic robots, stabilizing MPC using terminal region constraints and costs has been pursued in [6] while a contraction constraint on the first state in the prediction horizon was used in [3]. Moreover, in [11] a non-quadratic terminal cost was constructed on a terminal region for car-like non-holonomic robots. Here, the desired set point was located at the boundary of the closed terminal region, see [21] for a robust version. MPC without stabilizing constraints but with terminal costs has been first studied for non-holonomic systems in [22]. For the regulation of differential drive robots, MPC without stabilizing constraints is particularly attractive since computing (possibly time varying) terminal regions for large feasible sets can be an extremely challenging task [23]. This is especially true if the results shall be generalized to multi robot systems or domains with obstacles.

In this chapter, a stability analysis of MPC schemes without stabilizing constraints or costs, for regulation of non-holonomic mobile robots, is performed. Herein, a methodology is proposed, which allows to determine a prediction horizon length such that asymptotic stability of the MPC closed-loop is guaranteed. To this end, a proof of concept for verifying the controllability assumption introduced in [19] is presented. Herein, the running costs are tailored to the design specification of controlling both the position and the orientation. Then, the less conservative technique of [20, 24–26] is applied in order to rigorously prove asymptotic stability.

While the construction of particular open-loop maneuvers used to derive the growth condition of [19] heavily relies on the kinematic unicycle model, the pursued approach is outlined such that it can be used as a framework for verifying the above mentioned controllability assumption and, thus, being able to conclude asymptotic stability of the MPC closed-loop also for other systems. In particular, the insight provided by our analysis yields guidelines for the design of MPC controllers also for more accurate models of differential drive robots – a topic for future research. An extension of our discrete time results to the continuous time domain based on the presented results can be found in [27].

Finally, we numerically demonstrate that the canonical choice of quadratic running costs is not suited for regulation of non-holonomic mobile robots without (stabilizing) terminal constraints and/or costs. Moreover, the effectiveness of our approach is shown by means of numerical simulations.

This chapter is organized as follows: Section 3.3 outlines the regulation problem of non-holonomic mobile robots as well as the MPC algorithm. The stability results presented in [20, 26] are revisited in Section 3.4. In Section 3.5, bounds on the value function are derived by constructing appropriate feasible open-loop trajectories. Based on these bounds, a suitable prediction horizon length can be determined such that the MPC closed-loop is asymptotically stable. Our findings are illustrated by numerical simulations in Section 3.6. Finally, conclusions are drawn in Section 3.7.

## 3.3 Problem Setup

In this section, a differential drive mobile robot is described by an ordinary differential equation. Then, a corresponding discrete time model is presented and a model predictive control scheme is proposed in order to asymptotically stabilize the robot.

### 3.3.1 Non-holonomic mobile robot

The kinematic model of the mobile robot is given by

$$\begin{pmatrix} \dot{x}(t) \\ \dot{y}(t) \\ \dot{\theta}(t) \end{pmatrix} = \dot{\mathbf{x}}(t) = \mathbf{f}(\mathbf{x}(t), \mathbf{u}(t)) = \begin{pmatrix} v(t)\cos(\theta(t)) \\ v(t)\sin(\theta(t)) \\ \omega(t) \end{pmatrix} \quad (3.1)$$

with an analytic vector field $\mathbf{f} : \mathbb{R}^3 \times \mathbb{R}^2 \to \mathbb{R}^3$. The first two (spatial) components of the state $\mathbf{x} = (x, y, \theta)^\top$ (m,m,rad)$^\top$ represent the position in the plane while the angle $\theta$ corresponds to the orientation of the robot. The control input is $\mathbf{u} = (v, \omega)^\top$ (m/s,rad/s)$^\top$, where $v$ and $\omega$ are the linear and the angular speeds of the robot, respectively. Assuming piecewise constant control inputs on each interval $[k\delta, (k+1)\delta)$, $k \in \mathbb{N}_0$, with sampling period $\delta$ (seconds) and using direct integration, the (exact) discrete time dynamics $\mathbf{f}_\delta : \mathbb{R}^3 \times \mathbb{R}^2 \to \mathbb{R}^3$ are given by [28]

$$\mathbf{x}(k+1) = \mathbf{f}_\delta(\mathbf{x}(k), \mathbf{u}(k)) = \begin{pmatrix} x(k) \\ y(k) \\ \theta(k) \end{pmatrix} + \begin{pmatrix} \frac{v(k)}{\omega(k)} (\sin(\theta(k) + \delta\omega(k)) - \sin(\theta(k))) \\ \frac{v(k)}{\omega(k)} (\cos(\theta(k)) - \cos(\theta + \delta\omega(k))) \\ \delta\omega(k) \end{pmatrix} \quad (3.2)$$

for $\omega \neq 0$. When the robot moves in a straight line (angular speed $\omega = 0$), the right hand side of (3.2) becomes

$$\mathbf{x}(k) + \lim_{\omega \to 0} \begin{pmatrix} \frac{v(k)}{\omega(k)} (\sin(\theta(k) + \delta\omega(k)) - \sin(\theta(k))) \\ \frac{v(k)}{\omega(k)} (\cos(\theta(k)) - \cos(\theta(k) + \delta\omega(k))) \\ \delta\omega(k) \end{pmatrix} = \mathbf{x}(k) + \delta \cdot v(k) \begin{pmatrix} \cos(\theta(k)) \\ \sin(\theta(k)) \\ 0 \end{pmatrix}.$$

The movement is restricted to a rectangle, which is modelled by the box constraints

$$\begin{pmatrix} -\bar{x} \\ -\bar{y} \end{pmatrix} \leq \begin{pmatrix} x(k) \\ y(k) \end{pmatrix} \leq \begin{pmatrix} \bar{x} \\ \bar{y} \end{pmatrix} \quad \forall\, k \in \mathbb{N}_0. \quad (3.3)$$

The control inputs are limited by

$$\begin{pmatrix} -\bar{v} \\ -\bar{\omega} \end{pmatrix} \leq \begin{pmatrix} v(k) \\ \omega(k) \end{pmatrix} \leq \begin{pmatrix} \bar{v} \\ \bar{\omega} \end{pmatrix} \quad \forall\, k \in \mathbb{N}_0 \quad (3.4)$$

with $\bar{x}, \bar{y}, \bar{v}, \bar{\omega} > 0$. Then, admissibility of a sequence of input signals can be defined as follows.

**Definition 3.1.** *Let $X := [-\bar{x}, \bar{x}] \times [-\bar{y}, \bar{y}] \times \mathbb{R} \subset \mathbb{R}^3$ and $U := [-\bar{v}, \bar{v}] \times [-\bar{\omega}, \bar{\omega}] \subset \mathbb{R}^2$ be given. Then, for a given state $\mathbf{x}_0 := \mathbf{x}(0) \in X$, a sequence of control values $\mathbf{u} = (\mathbf{u}(0), \mathbf{u}(1), \ldots, \mathbf{u}(N-1)) \in U^N$ of length $N \in \mathbb{N}$ is called admissible, denoted by $\mathbf{u} \in \mathcal{U}^N(\mathbf{x}_0)$, if the state trajectory*

$$\mathbf{x_u}(\cdot; \mathbf{x}_0) = (\mathbf{x_u}(0; \mathbf{x}_0), \mathbf{x_u}(1; \mathbf{x}_0), \ldots, \mathbf{x_u}(N; \mathbf{x}_0))$$

*iteratively generated by system dynamics (3.2) and $\mathbf{x_u}(0; \mathbf{x}_0) = \mathbf{x}_0$ satisfies $\mathbf{x_u}(k; \mathbf{x}_0) \in X$ for all $k \in \{0, 1, \ldots, N\}$. An infinite sequence of control values $\mathbf{u} = (\mathbf{u}(k))_{k \in \mathbb{N}_0} \subset U$ is said to be admissible for $\mathbf{x}_0 \in X$, denoted by $\mathbf{u} \in \mathcal{U}^\infty(\mathbf{x}_0)$, if the truncation to its first $N$ elements is contained in $\mathcal{U}^N(\mathbf{x}_0)$ for all $N \in \mathbb{N}$.*

### 3.3.2 Model predictive control (MPC)

The goal is to steer the mobile robot to a desired (feasible) state $\mathbf{x}^r \in X$, which is without loss of generality chosen to be the origin, i.e. $\mathbf{x}^r = 0_{\mathbb{R}^3}$.[1] Indeed, $\mathbf{x}^r$ is a (controlled) equilibrium since $\mathbf{f}_\delta(\mathbf{x}^r, 0) = \mathbf{x}^r$. More precisely, our goal is to find a static state feedback law $\boldsymbol{\mu} : X \to U$ such that, for each $\mathbf{x}_0 \in X$, the resulting closed-loop system $\mathbf{x}_{\boldsymbol{\mu}}(\cdot; \mathbf{x}_0)$ generated by

$$\mathbf{x}_{\boldsymbol{\mu}}(n+1; \mathbf{x}_0) = \mathbf{f}_\delta(\mathbf{x}_{\boldsymbol{\mu}}(n; \mathbf{x}_0), \boldsymbol{\mu}(\mathbf{x}_{\boldsymbol{\mu}}(n; \mathbf{x}_0)))$$

and $\mathbf{x}_{\boldsymbol{\mu}}(0; \mathbf{x}_0) = \mathbf{x}_0$, satisfies the constraints $\mathbf{x}_{\boldsymbol{\mu}}(n; \mathbf{x}_0) \in X$ and $\boldsymbol{\mu}(\mathbf{x}_{\boldsymbol{\mu}}(n; \mathbf{x}_0)) \in U$ for all closed-loop time indices $n \in \mathbb{N}_0$ and is asymptotically stable, i.e. there exists a $\mathcal{KL}$- function $\beta : \mathbb{R}_{\geq 0} \times \mathbb{N}_0 \to \mathbb{R}_{\geq 0}$ such that, for each $\mathbf{x}_0 \in X$, the closed-loop trajectory obeys the inequality

$$\|\mathbf{x}_{\boldsymbol{\mu}}(n;\mathbf{x}_0)\| \leq \beta(\|\mathbf{x}_0\|, n) \qquad \forall\, n \in \mathbb{N}_0,$$

see Figure 3.1 for an illustration sketch of a $\mathcal{KL}$-function.

---
[1] $\mathbf{x}^r$ is supposed to be in the interior of the state constraint set $X$.

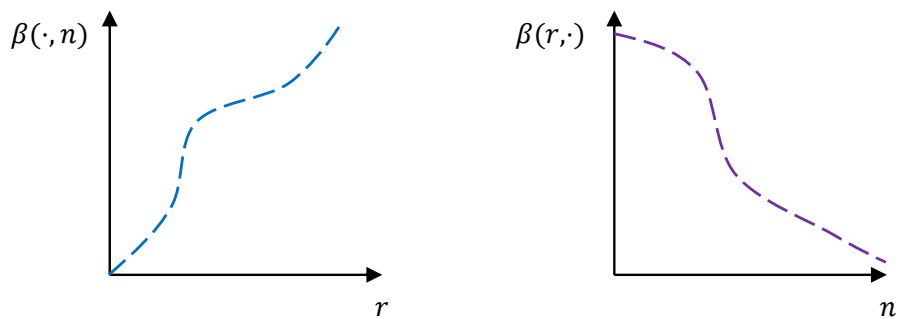

Fig. 3.1: An illustration of a typical class $\mathcal{KL}$-function $(\beta(r,n) : \mathbb{R}_{\geq 0} \times \mathbb{N}_0 \to \mathbb{R}_{\geq 0})$ based on [23].

As briefly discussed in the introduction, several control techniques have been developed for this purpose. In this chapter we use MPC, which makes use of the system dynamics in order to design a control strategy minimizing a cost function. This cost function sums up given stage costs along predicted (feasible) trajectories. We propose to deploy the running (stage) costs $\ell : X \times U \to \mathbb{R}_{\geq 0}$ defined as

$$\ell(\mathbf{x}, \mathbf{u}) = q_1 x^4 + q_2 y^2 + q_3 \theta^4 + r_1 v^4 + r_2 \omega^4 \tag{3.5}$$

with $q_1, q_2, q_3, r_1, r_2 \in \mathbb{R}_{>0}$. In (3.5), small deviations in the $y$-direction are penalized more than deviations with respect to $x$ or $\theta$. The motivation behind this particular choice becomes clear in Section 3.5.2, where the different order of $y$ is exploited in order to verify Assumption 3.2 and, thus, to ensure asymptotic stability. Moreover, in Section 3.6.2, we explain why quadratic running costs $\ell(\mathbf{x}, \mathbf{u}) = \mathbf{x}^\top \mathbf{Q} \mathbf{x} + \mathbf{u}^\top \mathbf{R} \mathbf{u}$, $\mathbf{Q} \in \mathbb{R}^{3\times 3}$ and $\mathbf{R} \in \mathbb{R}^{2\times 2}$, are not suited for our example by conducting numerical simulations.

Based on the introduced running costs, a cost function $J_N : X \times U^N \to \mathbb{R}_{\geq 0}$ is defined as

$$J_N(\mathbf{x}_0, \mathbf{u}) := \sum_{k=0}^{N-1} \ell(\mathbf{x}_\mathbf{u}(k), \mathbf{u}(k)).$$

Then, at every time step $n$, the following finite horizon Optimal Control Problem (OCP) is solved for the initial condition $\mathbf{x}_0 := \mathbf{x}_\mathbf{u}(n)$

$$\min_{\mathbf{u} \in \mathbb{R}^{2 \times N}} \quad J_N(\mathbf{x}_0, \mathbf{u}) \tag{3.6}$$

$$\text{subject to } \mathbf{x}_\mathbf{u}(0) = \mathbf{x}_0,$$

$$\mathbf{x}_\mathbf{u}(k+1) = \mathbf{f}_\delta(\mathbf{x}_\mathbf{u}(k), \mathbf{u}(k)) \quad \forall k \in \{0, 1, \ldots, N-1\},$$

$$\mathbf{x}_\mathbf{u}(k) \in X \quad \forall k \in \{1, 2, \ldots, N\},$$

$$\mathbf{u}(k) \in U \quad \forall k \in \{0, 1, \ldots, N-1\},$$

for $N \in \mathbb{N} \cup \{\infty\}$. The corresponding (optimal) value function $V_N : X \to \mathbb{R}_{\geq 0} \cup \{\infty\}$ of OCP (3.6) is defined as

$$V_N(\mathbf{x}_0) := \min_{\mathbf{u} \in \mathcal{U}^N(\mathbf{x}_0)} J_N(\mathbf{x}_0, \mathbf{u}),$$

where $V_N(\mathbf{x}_0) = \infty$ if $\mathcal{U}^N(\mathbf{x}_0) = \emptyset$ holds. The resulting minimizing control sequence of OCP (3.6) is denoted by

$$\mathbf{u}^\star := (\mathbf{u}^\star(0), \mathbf{u}^\star(1), \ldots, \mathbf{u}^\star(N-1)) \in \mathcal{U}^N(\mathbf{x}_0),$$

where $\mathbf{u}^\star(0)$ is the control action to be applied on the system. The overall steps of the used MPC controller are summarized in Algorithm 3.1, which is an MPC scheme without stabilizing constraints or costs. For a detailed discussion on MPC we refer to [14, 23].

---
**Algorithm 3.1** MPC
---
**Initialization**: set prediction horizon $N$ and time index $n := 0$.
1: Measure the current state $\hat{\mathbf{x}} := \mathbf{x}(n)$.
2: Compute $\mathbf{u}^\star = (\mathbf{u}^\star(0), \mathbf{u}^\star(1), \ldots, \mathbf{u}^\star(N-1)) \in \mathcal{U}^N(\hat{\mathbf{x}})$ satisfying $J_N(\hat{\mathbf{x}}, \mathbf{u}^\star) = V_N(\hat{\mathbf{x}})$.
3: Define the MPC feedback law $\boldsymbol{\mu}_N : X \to U$ at $\hat{\mathbf{x}}$ by $\boldsymbol{\mu}_N(\hat{\mathbf{x}}) := \mathbf{u}^\star(0)$ and implement $\mathbf{u}(n) := \boldsymbol{\mu}_N(\hat{\mathbf{x}})$ at the plant. Then, increment the time index $n$ and goto step 1.

---

Since $0_{\mathbb{R}^2} \in U$ holds, $\mathcal{U}^N(\mathbf{f}_\delta(\hat{\mathbf{x}}, \boldsymbol{\mu}_N(\hat{\mathbf{x}}))) \neq \emptyset$ holds, i.e. recursive feasibility of the MPC closed-loop is ensured. Existence of an admissible sequence of control values minimizing $J_N(\hat{\mathbf{x}}, \cdot)$ can be infered from compactness of the nonempty domain and continuity of the cost function by applying the Weierstrass theorem, see [29] for details. However, since neither stabilizing constraints nor terminal costs are incorporated in our MPC formulation, asymptotic stability is far from being trivial and does, in general, not hold, see, e.g. [15]. In the following, we will show how to ensure asymptotic stability by appropriately choosing the MPC prediction horizon $N$.

## 3.4 Stability of MPC Without Stabilizing Constraints or Costs

In this section, known results from [20, 25] are recalled. Later, these results are exploited in order to rigorously prove asymptotic stability of the exact discrete time model of the mobile robot governed by (3.2). The following assumption, introduced in [19], is a key ingredient in order to show asymptotic stability of the MPC closed-loop.

**Assumption 3.2.** *Let a monotonically increasing and bounded sequence $(\gamma_i)_{i \in \mathbb{N}_0}$ be given*

and suppose that, for each $\mathbf{x}_0 \in X$, the estimate

$$V_i(\mathbf{x}_0) \leq \gamma_i \cdot \ell^\star(\mathbf{x}_0) \quad \forall\, i \in \mathbb{N}, \tag{3.7}$$

where

$$\ell^\star(\mathbf{x}_0) := \min_{\mathbf{u} \in \mathcal{U}^1(\mathbf{x}_0)} \ell(\mathbf{x}_0, \mathbf{u})$$

holds. Furthermore, let there exist two $\mathcal{K}_\infty$-functions $\underline{\eta}, \bar{\eta} : \mathbb{R}_{\geq 0} \to \mathbb{R}_{\geq 0}$ satisfying

$$\underline{\eta}(\|\mathbf{x} - \mathbf{x}^r\|) \leq \ell^\star(\mathbf{x}) \leq \bar{\eta}(\|\mathbf{x} - \mathbf{x}^r\|) \qquad \forall\, \mathbf{x} \in X. \tag{3.8}$$

Based on Assumption 3.2 and the fact that recursive feasibility trivially holds for our example, as observed in the preceding section, asymptotic stability of the MPC closed-loop can be established, see [20, Theorems 4.2 and 5.3] and [26].

**Theorem 3.3.** *Let Assumption 3.2 hold and let the performance index $\alpha_N$ be given by the formula*

$$\alpha_N := 1 - \frac{(\gamma_N - 1)\prod_{k=2}^N (\gamma_k - 1)}{\prod_{k=2}^N \gamma_k - \prod_{k=2}^N (\gamma_k - 1)}. \tag{3.9}$$

*Then, if $\alpha_N > 0$, the relaxed Lyapunov inequality*

$$V_N(\mathbf{f}_\delta(\mathbf{x}, \boldsymbol{\mu}_N(\mathbf{x}))) \leq V_N(\mathbf{x}) - \alpha_N \ell(\mathbf{x}, \boldsymbol{\mu}_N(\mathbf{x})) \tag{3.10}$$

*holds for all $\mathbf{x} \in X$ and the MPC closed-loop with prediction horizon $N$ is asymptotically stable.*

Note that the relaxed Lyapunov inequality (3.10) requires that the value function of MPC to be monotonically decreasing.

While Condition (3.8) holds trivially for the chosen running costs, the derivation of the growth bounds $\gamma_i$, $i \in \mathbb{N}_0$, of Condition (3.7), is, in general, difficult. One option to derive $\gamma_i$, is the following proposition.

**Proposition 3.4.** *Let a sequence $(c_n)_{n \in \mathbb{N}_0} \subseteq \mathbb{R}_{\geq 0}$, be given and assume that $\sum_{n=0}^{\infty} c_n < \infty$ holds. In addition, suppose that for each $\mathbf{x}_0 \in X$ an admissible sequence of control values $\mathbf{u}_{\mathbf{x}_0} = (\mathbf{u}_{\mathbf{x}_0}(n))_{n \in \mathbb{N}_0} \in \mathcal{U}^{\infty}(\mathbf{x}_0)$ exists such that the inequality*

$$\ell(\mathbf{x}_{\mathbf{u}_{\mathbf{x}_0}}(n; \mathbf{x}_0), \mathbf{u}_{\mathbf{x}_0}(n)) \leq c_n \cdot \ell^{\star}(\mathbf{x}_0) \qquad \forall\, n \in \mathbb{N}_0 \tag{3.11}$$

*holds. Then, the growth bounds $\gamma_i$, $i \in \mathbb{N}_0$, of Condition (3.7) are given by $\gamma_i = \sum_{n=0}^{i-1} c_n$, $i \in \mathbb{N}_0$.*

**Proof**: Let $\mathbf{x}_0 \in X$ and $\mathbf{u}_{\mathbf{x}_0} \in \mathcal{X}^{\infty}(\mathbf{x}_0)$ be given such that Inequality (3.11) holds. Then, the definition of the value function $V_i$ yields

$$V_i(\mathbf{x}_0) \leq \sum_{n=0}^{i-1} \ell(\mathbf{x}_{\mathbf{u}_{\mathbf{x}_0}}(n; \mathbf{x}_0), \mathbf{u}_{\mathbf{x}_0}(n)) \leq \sum_{n=0}^{i-1} c_n \ell^{\star}(\mathbf{x}_0) = \gamma_i \ell^{\star}(\mathbf{x}_0).$$

While monotonicity of the sequence $(\gamma_i)_{i \in \mathbb{N}}$ results from $c_n \geq 0$, $n \in \mathbb{N}_0$, boundedness follows from the assumed summability of the sequence $(c_n)_{n \in \mathbb{N}_0}$. $\square$

In order to illustrate these results, a simple example taken from [30] is presented for which Condition (3.11) is deduced.

**Example 3.5.** *The system dynamics are given by the simple model $\mathbf{x}(k+1) = \mathbf{x} + \mathbf{u}$ with state and control constraints $X = [-1, 1]^2$ and $U = [-\bar{u}, \bar{u}]^2$ for some $\bar{u} > 0$, respectively. The desired equilibrium $\mathbf{x}^r$ is supposed to be contained in $X$. The running costs are $\ell(\mathbf{x}, \mathbf{u}) = \|\mathbf{x} - \mathbf{x}^r\|^2 + \lambda \|\mathbf{u}\|^2$ with weighting factor $\lambda \geq 0$.*

*Let $c := \max_{\mathbf{x} \in X} \|\mathbf{x} - \mathbf{x}^r\|$, i.e. the maximal distance of a feasible point from the desired state $\mathbf{x}^r$. We define inductively a control $\mathbf{u}_{\mathbf{x}_0} \in \mathcal{U}^N(\mathbf{x}_0)$ for some design parameter*

$\rho \in (0, 1)$

$$\mathbf{u}(k) = \kappa(\mathbf{x}^r - \mathbf{x}_{\mathbf{u}_{\mathbf{x}_0}}(k; \mathbf{x}_0)) \quad \text{with } \kappa = \min\{\bar{u}/c, \rho\}.$$

*The choice of $\kappa$ implies $\mathbf{u}(k) \in U$ for $\mathbf{x}_{\mathbf{u}_{\mathbf{x}_0}}(k; \mathbf{x}_0) \in X$. Since*

$$\mathbf{x}_{\mathbf{u}_{\mathbf{x}_0}}(k+1; \mathbf{x}_0) = \mathbf{x}_{\mathbf{u}_{\mathbf{x}_0}}(k; \mathbf{x}_0) + \kappa(\mathbf{x}^r - \mathbf{x}_{\mathbf{u}_{\mathbf{x}_0}}(k; \mathbf{x}_0))$$

*holds, we obtain*

$$\|\mathbf{x}_{\mathbf{u}_{\mathbf{x}_0}}(k+1; \mathbf{x}_0) - \mathbf{x}^r\| = (1-\kappa)\|\mathbf{x}_{\mathbf{u}_{\mathbf{x}_0}}(k; \mathbf{x}_0) - \mathbf{x}^r\|$$

*and due to convexity of $X$ and $\kappa \in (0, 1)$, feasibility of the state trajectory $(\mathbf{x}_{\mathbf{u}_{\mathbf{x}_0}}(k; \mathbf{x}_0))_{k \in \mathbb{N}_0}$ is ensured. Then, Condition (3.11) can be deduced by*

$$\begin{aligned}
\ell(\mathbf{x}_{\mathbf{u}_{\mathbf{x}_0}}(k), \mathbf{u}_{\mathbf{x}_0}(k)) &= \|\mathbf{x}_{\mathbf{u}_{\mathbf{x}_0}}(k) - \mathbf{x}^r\|^2 + \lambda \|\mathbf{u}_{\mathbf{x}_0}(k)\|^2 \\
&= (1 + \lambda \kappa^2)\|\mathbf{x}_{\mathbf{u}_{\mathbf{x}_0}}(k) - \mathbf{x}^r\|^2 \\
&= (1 + \lambda \kappa^2)(1 - \kappa)^{2k} \underbrace{\|\mathbf{x}_{\mathbf{u}_{\mathbf{x}_0}}(0) - \mathbf{x}^r\|^2}_{= \ell^\star(\mathbf{x}_0)}
\end{aligned}$$

*with $\mathbf{x}_{\mathbf{u}_{\mathbf{x}_0}}(k) = \mathbf{x}_{\mathbf{u}_{\mathbf{x}_0}}(k; \mathbf{x}_0)$, i.e. Condition (3.11) with $c_n = C\sigma^n$ where the parameters $C = 1 + \lambda \kappa^2$ and $\sigma = (1 - \kappa)^2$ are used. Hence, an exponential decay is shown which implies the summability of the sequence $(c_n)_{n \in \mathbb{N}_0}$.*

Based on the sequence $(c_n)_{n \in \mathbb{N}_0}$ computed in Example 3.5, and for $N = 2$, Formula (3.9) yields $\alpha_2 = 1 - (C + \sigma C - 1)^2$. Then, supposing $\lambda \in (0, 1)$ and using the

estimate

$$C + \sigma C = (1 + \lambda \kappa^2)(1 - \kappa^2) \leq (1 + \kappa)(1 - \kappa)^2 = (1 - \kappa^2)(1 - \kappa) < 1$$

$\alpha_2 > 0$ is implied. Hence, Theorem 3.3 can be used to conclude asymptotic stability for prediction horizon $N = 2$. For the parameters $\lambda = 0.1$ and $\rho = 0.5$, the performance index $\alpha_2$ is approximately 0.9209.

**Remark 3.6.** *A direct verification of Assumption 3.2 yields, in general, less conservative bounds on the required prediction horizon in order to ensure that $\alpha_N \in (0, 1]$ is satisfied. However, Proposition 3.4 is instructive for the construction in the subsequent section.*

## 3.5 Stability Analysis of the Unicycle Mobile Robot

In this section, a bounded sequence $(\gamma_i)_{i \in \mathbb{N}_{\geq 2}}$ is constructed such that Assumption 3.2 holds. For this purpose, first an open set $\mathcal{N}_1 = \mathcal{N}_1(s)$ of initial conditions depending on a parameter $s \in [0, \infty)$ is defined by

$$\mathcal{N}_1 = \left\{ \mathbf{x} = \begin{pmatrix} x \\ y \\ \theta \end{pmatrix} \in \mathbb{R}^3 : \mathbf{x} \in X \text{ and } \ell^\star \left( \begin{pmatrix} x \\ y \\ 0 \end{pmatrix} \right) := q_1 x^4 + q_2 y^2 < s \right\}. \quad (3.12)$$

Based on this definition, the feasible set $X$ is split up into $\mathcal{N}_1$ and $\mathcal{N}_2 := X \setminus \mathcal{N}_1$ such that $X = \mathcal{N}_1 \cup \mathcal{N}_2$ holds. Then, bounded sequences $(\gamma_i^{\mathcal{N}_j})_{i \in \mathbb{N}_{\geq 2}}$, $j \in \{1, 2\}$, are derived such that

$$V_i(\mathbf{x}_0) \leq \gamma_i^{\mathcal{N}_j} \cdot \ell^\star(\mathbf{x}_0) \qquad \forall \ \mathbf{x}_0 \in \mathcal{N}_j \quad (3.13)$$

holds for all $i \in \mathbb{N}$. In conclusion, taking into account that the input sequence $(\mathbf{u}(k))_{k \in \mathbb{N}_0}$, $\mathbf{u}(k) = 0_{\mathbb{R}^2}$, is admissible on the infinite horizon and implies Inequality (3.7) with $\gamma_i = i$,

Inequality (3.7) holds for all $\mathbf{x}_0 \in X$ with

$$\gamma_i := \min\{i, \max\{\gamma_i^{\mathcal{N}_1}, \gamma_i^{\mathcal{N}_2}\}\}, \qquad i \in \mathbb{N}_{\geq 2}. \tag{3.14}$$

The motivation behind partitioning the set $X$ is that we design two different maneuvers in order to deduce bounded sequences $(\gamma_i^{\mathcal{N}_j})_{i \in \mathbb{N}_{\geq 2}}$, $j \in \{1, 2\}$. While in principle one strategy could be sufficient, one of the proposed maneuvers works for initial states close to the origin (inside the set $\mathcal{N}_1$) while the other becomes more advantageous outside $\mathcal{N}_1$. In this vein, the vehicle is just turned towards the origin $0 \in \mathbb{R}^2$ and, then, drives in that direction before the angle is set to zero if $\mathbf{x}_0 \in \mathcal{N}_2$. However, this move does not allow to derive a bounded $\gamma_i$-sequence for initial positions $\mathbf{x}_0 = (0, y_0, 0)^\top$ whose distance $\ell^\star(\mathbf{x}_0)$ tends to zero. But, boundedness is essential in order to deduce asymptotic stability of the MPC closed-loop via Theorem 3.3.

Before we present the (technical) details in following Sections 3.5.1 and 3.5.2, we explain briefly the strategy used to construct $(\gamma_i^{\mathcal{N}_j})_{i \in \mathbb{N}_{\geq 2}}$, $j \in \{1, 2\}$. First, for initial values $\mathbf{x}_0 = (x_0, y_0, 0)^\top \in \mathcal{N}_j$, a family of particular control sequences $\mathbf{u}_{\mathbf{x}_0} := (\mathbf{u}(k; \mathbf{x}_0))_{k \in \mathbb{N}_0} \in \mathcal{U}^\infty(\mathbf{x}_0)$ is proposed such that the robot is steered to the origin in a finite number of steps. These input sequences $\mathbf{u}_{\mathbf{x}_0}$ yield (suboptimal) running costs $\ell(\mathbf{x}_{\mathbf{u}_{\mathbf{x}_0}}(k; \mathbf{x}_0), \mathbf{u}(k; \mathbf{x}_0))$ such that, by definition of optimality, the following quotients can be estimated uniformly with respect to $\mathbf{x}_0 = (x_0, y_0, 0)^\top \in \mathcal{N}_j$ by

$$\frac{\ell(\mathbf{x}_{\mathbf{u}_{\mathbf{x}_0}}(k; \mathbf{x}_0), \mathbf{u}(k; \mathbf{x}_0))}{\ell^\star(\mathbf{x}_0)} \leq c_k \quad \forall\, k \in \mathbb{N}_0 \tag{3.15}$$

with coefficients $c_k = c_k^{\mathcal{N}_j}$, $k \in \mathbb{N}_0$, i.e. a coefficient sequence $(c_k)_{k \in \mathbb{N}_0}$ such that Inequality (3.11) holds. Moreover, as highlighted in (3.7), for $\mathbf{x}_0 = (x_0, y_0, 0)^\top$, $\ell^\star(\mathbf{x}_0)$ is given

by

$$\ell^\star(\mathbf{x}_0) = q_1 x_0^4 + q_2 y_0^2.$$

Since also the number of steps needed in order to steer the considered initial states $\mathbf{x}_0$ to the origin exhibits a uniform upper bound, there exists $\bar{k}$ such that $c_k = 0$ holds for all $k \geq \bar{k}$. Then, the coefficients $c_1, c_2, \ldots, c_{\bar{k}-1}$ are rearranged in a descending order denoted by $(\bar{c}_k)_{k \in \mathbb{N}_0}$ with $\bar{c}_0 = c_0$, which still implies Condition (3.7) with $\gamma_i := \sum_{n=0}^{i-1} \bar{c}_n$. Finally, these $\gamma_i$-sequences are used in order to ensure Condition (3.7) for *all* initial states contained in $\mathcal{N}_j$, i.e. those initial conditions with $\theta_0 \neq 0$. Due to symmetries (the robot can go back and forth), it is sufficient to consider initial positions with $(x_0, y_0)^\top \geq 0_{\mathbb{R}^2}$.

### 3.5.1  Trajectory generation for $\mathbf{x}_0 \in \mathcal{N}_2$

In this section, we first consider initial conditions inside $\mathcal{N}_2$ with $\theta_0 = 0$. Subsequently, we will prove that the derived bounds also hold for the case $\theta_0 \neq 0$.

**Initial Conditions $\mathbf{x}_0 \in \mathcal{N}_2$ with $\theta_0 = 0$:** For initial conditions $\mathbf{x}_0 = (x_0, y_0, 0)^\top$ in the set $\mathcal{N}_2$, the following simple maneuver can be employed:

a) choose an angle $\bar{\theta} \in [-\pi, \pi)$ such that the vehicle points towards (or in the opposite direction to) the origin $(0,0)^\top \in \mathbb{R}^2$,

b) drive directly towards the origin,

c) turn the vehicle to the desired angle $\theta^r = 0$.

See Figure 3.2 for an illustration sketch of the proposed maneuver. The number of steps needed in order to carry out this maneuver depends on the constraints and the sampling time $\delta$, which is supposed to satisfy $i \cdot \delta = 1$ for some integer $i \in \mathbb{N}$. We define the minimal

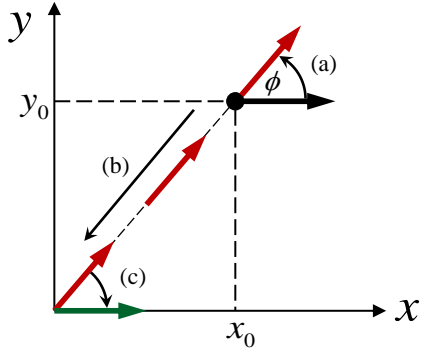

Fig. 3.2: Illustration of the maneuver for initial value $\mathbf{x}_0 = (x_0, y_0, 0) \in \mathcal{N}_2$.

number of steps required to turn the vehicle by 90 degrees as

$$k_\delta^\star := \left\lceil \frac{\pi/2}{\min\{\bar{\omega}, \pi/2\} \cdot \delta} \right\rceil$$

assuming reasonable bounds control constraints. We define also the minimal number of steps required to drive the vehicle from the farthest corner of the box defined by the constraints (3.3) to the origin as

$$l_\delta^\star := \left\lceil \frac{\sqrt{\bar{y}^2+, \bar{x}^2}}{\bar{v} \cdot \delta} \right\rceil,$$

respectively. Additionally, the inequality

$$r_2 \leq \frac{q_3 \cdot \delta}{2} \tag{3.16}$$

is assumed to hold in order to avoid technical difficulties resulting from not reflecting the sampling time $\delta$ in the running costs.

Initial values $\mathbf{x}_0 = (x_0, y_0, 0)^\top \geq 0$ are considered first. Let the angle $\arctan(y_0/x_0) \in [0, \pi/2)$ be denoted by $\phi$. The vehicle stays at the initial position without moving for $k_\delta^\star$ steps, i.e $(v_i, w_i)^\top = (0,0)^\top$, $i \in \{0, 1, \ldots, k_\delta^\star - 1\}$, which yields Inequality (3.15) with $c_i^{\mathcal{N}_2} = 1$, $i = 0, 1, \ldots, k_\delta^\star - 1$. This artificially added phase is introduced here in order to facilitate the treatment of initial positions with $\theta_0 \neq 0$.

Next, the vehicle turns $k_\delta^\star$ steps such that $\theta_{\mathbf{u}}(2k_\delta^\star; z_0) = \phi$ holds by applying the input $\mathbf{u}(k_\delta^\star + i) = (0, \phi \cdot (k_\delta^\star \delta)^{-1})^\top \in U$ for all $i \in \{0, 1, \ldots, k_\delta^\star - 1\}$. This control action yields the running costs

$$\ell(\mathbf{x}_{\mathbf{u}}(k_\delta^\star + i; \mathbf{x}_0), \mathbf{u}(k_\delta^\star + i)) = q_1 x_0^4 + q_2 y_0^2 + q_3 \left(\frac{i\phi}{k_\delta^\star}\right)^4 + r_2 \left(\frac{\phi}{k_\delta^\star \delta}\right)^4. \tag{3.17}$$

Since $\phi \in [0, \pi/2)$, $\ell^\star(\mathbf{x}_0) \geq s$, and Assumption (3.16) hold, Inequality (3.15) is ensured with the coefficients

$$c_{k_\delta^\star + i}^{\mathcal{N}_2} := 1 + \frac{q_3 \pi^4}{16 k_\delta^{\star 4} \cdot s} \left(i^4 + \frac{1}{2\delta^3}\right), \tag{3.18}$$

$i = 0, 1, \ldots, k_\delta^\star - 1$. Then, the vehicle drives towards the origin in $l_\delta^\star$ steps with constant backward speed $\mathbf{u}(2k_\delta^\star + i) = (-\|(x_0, y_0)^\top\| \cdot (l_\delta^\star \delta)^{-1}, 0)^\top \in X$, $i \in \{0, 1, \ldots, l_\delta^\star - 1\}$. This leads to the running costs

$$\ell(\mathbf{x}_{\mathbf{u}}(2k_\delta^\star + i; \mathbf{x}_0), \mathbf{u}(2k_\delta^\star + i)) = \left(\frac{l_\delta^\star - i}{l_\delta^\star}\right)^2 \left[q_1 \left(\frac{l_\delta^\star - i}{l_\delta^\star}\right)^2 x_0^4 + q_2 y_0^2\right] + q_3 \phi^4 + r_1 \left(\frac{\|(x_0, y_0)\|}{l_\delta^\star \delta}\right)^4.$$

Since $\phi \leq \pi/2$ and the control effort is smaller than $\bar{v}$, the respective coefficients for Inequality (3.15) can be chosen as

$$c_{2k_\delta^\star + i}^{\mathcal{N}_2} := \left(\frac{l_\delta^\star - i}{l_\delta^\star}\right)^2 + \frac{q_3 (\pi/2)^4 + r_1 \bar{v}^4}{s} \tag{3.19}$$

for $i \in \{0, 1, \ldots, l_\delta^\star - 1\}$. Finally, the vehicle turns $k_\delta^\star$ steps in order to reach $\theta_{\mathbf{u}}(3k_\delta^\star + l_\delta^\star; z_0) = 0$ using the input $\mathbf{u}(2k_\delta^\star + l_\delta^\star + i) = (0, -\phi \cdot (k_\delta^\star \delta)^{-1})^\top$, $i \in \{0, 1, \ldots, k_\delta^\star - 1\}$. Thus, we have the running costs

$$\ell(\mathbf{x_u}(2k_\delta^\star + l_\delta^\star + i; \mathbf{x}_0), \mathbf{u}(2k_\delta^\star + l_\delta^\star + i)) = \left[q_3\left(\frac{k_\delta^\star - i}{k_\delta^\star}\right)^4 + r_2\left(\frac{1}{k_\delta^\star \delta}\right)^4\right]\phi^4. \quad (3.20)$$

Then, invoking (3.16) ensures Inequality (3.15) with

$$c_{2k_\delta^\star + l_\delta^\star + i}^{\mathcal{N}_2} := \frac{q_3 \pi^4}{16 k_\delta^{\star 4} \cdot s}\left[(k_\delta^\star - i)^4 + \frac{1}{2\delta^3}\right] \quad (3.21)$$

for $i \in \{0, 1, \ldots, k_\delta^\star - 1\}$. The calculated coefficients $c_i^{\mathcal{N}_2}$, $i = 1, 2, \ldots, 3k_\delta^\star + l_\delta^\star - 1$, are ordered descendingly resulting in a new sequence $(\bar{c}_i^{\mathcal{N}_2})_{i=1}^{3k_\delta^\star + l_\delta^\star - 1}$, satisfying $\bar{c}_i^{\mathcal{N}_2} \leq \bar{c}_{i-1}^{\mathcal{N}_2}$ for $i \in \{2, 3, \ldots, 3k_\delta^\star + l_\delta^\star - 1\}$. Then, setting $\bar{c}_0^{\mathcal{N}_2} = c_0^{\mathcal{N}_2}$ and $\bar{c}_i^{\mathcal{N}_2} = 0$ for all $i \geq 3k_\delta^\star + l_\delta^\star$ yields $(\bar{c}_i^{\mathcal{N}_2})_{i=0}^\infty$. Hence, the accumulated bounds $(\gamma_i^{\mathcal{N}_2})_{i\in\mathbb{N}_{\geq 2}}$ of Condition (3.7) for the first maneuver are given by

$$\gamma_i^{\mathcal{N}_2} := \sum_{n=0}^{i-1} \bar{c}_i^{\mathcal{N}_2}, \qquad i \in \mathbb{N}_{\geq 2}. \quad (3.22)$$

**Initial Conditions $\mathbf{x}_0 \in \mathcal{N}_2$ with $\theta_0 \neq 0$**: In this subsection, we show that Condition (3.7) holds for arbitrary initial conditions $\mathbf{x}_0 \in \mathcal{N}_2$, i.e. $\theta_0 \in [-\pi, 0) \cup (0, \pi)$, using the bounds defined in (3.22). To this end, we distinguish four intervals in dependence of the initial angular deviation $\theta_0$, see Figure 3.3. While the basic ingredients are similar to the described maneuver for $\theta_0 = 0$, the order of the involved *motions* differs as summarized in Figure 3.4 in order to facilitate the accountability of the upcoming presentation.

**Case 1**: let $\theta_0$ be contained in the interval $(0, \phi)$. The robot stays at the initial position without moving for $k_\delta^\star$ steps; thus, Inequality (3.15) holds with the coefficients $c_i = 1$, $i = 0, 1, \ldots, k_\delta^\star - 1$. Then, the control input $\mathbf{u}(k_\delta^\star) = (0, \omega_{k_\delta^\star})^\top$, $\omega_{k_\delta^\star} \in (0, \phi \cdot (k_\delta^\star \delta)^{-1}]$ is

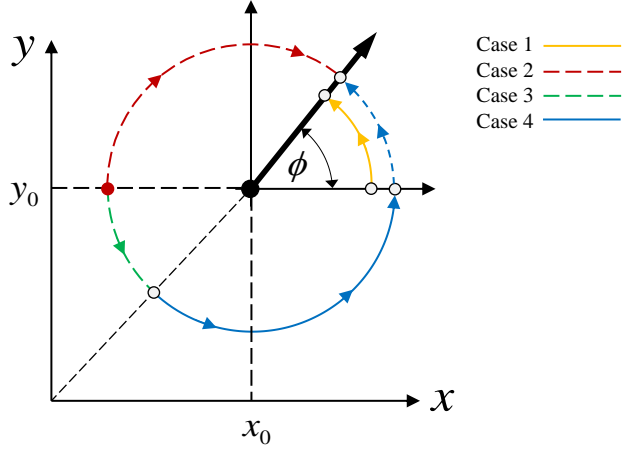

Fig. 3.3: Classification of the four different cases for $\mathbf{x}_0 \in \mathcal{N}_2$ with $\theta_0 \neq 0$.

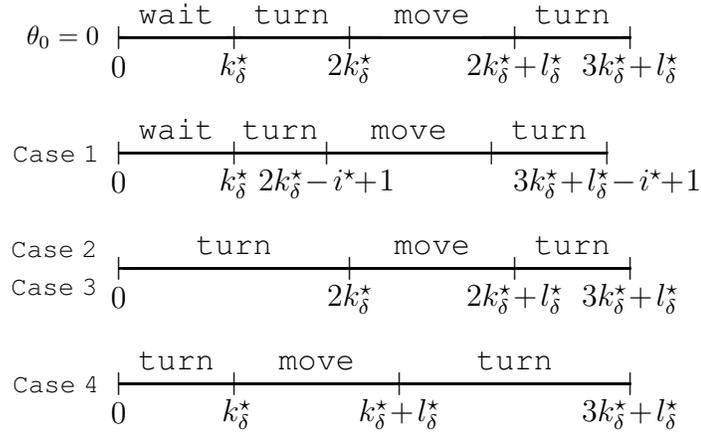

Fig. 3.4: The maneuver for initial conditions $\mathbf{x}_0 \in \mathcal{N}_2$ consists of waiting, turning, and moving the differential drive robot. However, the order of these *motions* depends on the initial angular deviation, see Figure 3.3.

adjusted such that

$$\exists!\ i^\star \in \{1, \ldots, k_\delta^\star - 1\} : \theta_0 + \delta\omega_{k_\delta^\star} = \theta_{k_\delta^\star + i^\star}$$

where $\theta_{k_\delta^\star + i^\star}$ is one of the achieved angles during the maneuver for $\theta_0 = 0$. Then, the robot turns $k_\delta^\star - i^\star$ steps such that $\theta_\mathbf{u}(2k_\delta^\star - i^\star + 1) = \phi$ is achieved using the input $\mathbf{u}(k_\delta^\star + i) = (0, \phi \cdot (k_\delta^\star \delta)^{-1})^\top$, $i \in \{1, 2, \ldots, k_\delta^\star - i^\star\}$. Hence, Inequality (3.15) is valid with the coefficient $c_{k_\delta^\star + i + (i^\star - 1)}^{\mathcal{N}_2}$, $i \in \{0, 1, \ldots, k_\delta^\star - i^\star\}$. The remaining parts of the maneuver are performed as for $\theta_0 = 0$. Since we have $\ell^\star(\mathbf{x}_0) > s$, but precisely the same running costs, the growth bounds given by (3.22) can be used to ensure condition (3.7) for the considered case.

**Case 2**: let $\theta_0 \in (\phi, \pi]$ hold. The first part of the maneuver is performed by turning the robot $2k_\delta^\star$ steps, such that $\theta_\mathbf{u}(2k_\delta^\star; z_0) = \phi$ is achieved using the input $\mathbf{u}(i) = (0, -\Delta\theta \cdot (k_\delta^\star \delta)^{-1})^\top$, $i = 0, 1, \ldots, 2k_\delta^\star - 1$, $\Delta\theta = (\theta_0 - \phi)/2$. Hence, the corresponding running costs are

$$\ell(\mathbf{x}_\mathbf{u}(i; \mathbf{x}_0), \mathbf{u}(i)) = q_1 x_0^4 + q_2 y_0^2 + q_3 \left[\theta_0 - i\left(\frac{\Delta\theta}{k_\delta^\star}\right)\right]^4 + r_2 \left(\frac{\Delta\theta}{k_\delta^\star \delta}\right)^4. \tag{3.23}$$

for $i \in \{0, 1, \ldots, 2k_\delta^\star - 1\}$. Then, using $\theta_0 \geq h := i\left(\frac{\Delta\theta}{k_\delta^\star}\right)\theta > 0$, the following gives an overestimate of the third term in (3.23)

$$\begin{aligned}\left[\theta_0 - i\left(\frac{\Delta\theta}{k_\delta^\star}\right)\right]^4 &= \theta_0^4 - h\theta_0^3 - 3h\theta_0(\theta_0 - h)^2 - h^3(\theta_0 - h) \\ &\leq \theta_0^4 - i\left(\frac{\Delta\theta}{k_\delta^\star}\right)\theta_0^3\end{aligned} \tag{3.24}$$

and, thus, invoking Assumption (3.16), i.e. $r_2 \leq \frac{q_3 \cdot \delta}{2}$, leads to

$$\ell(\mathbf{x_u}(i; \mathbf{x}_0), \mathbf{u}(i)) \leq \ell^\star(\mathbf{x}_0) - q_3 \left(\frac{\Delta\theta}{k_\delta^\star}\right) \left[i\theta_0^3 - \frac{1}{2}\left(\frac{\Delta\theta}{k_\delta^\star \delta}\right)^3\right] \tag{3.25}$$

for $i = 0, 1, \ldots, 2k_\delta^\star - 1$. In conclusion, the right hand side of this inequality is always less than or equal to $\ell^\star(\mathbf{x}_0)$ for $i > 0$. Hence, Inequality (3.15) holds with

$$c_{k_\delta^\star}^{\mathcal{N}_2}, c_0^{\mathcal{N}_2}, c_1^{\mathcal{N}_2}, \ldots, c_{k_\delta^\star - 1}^{\mathcal{N}_2}, c_{k_\delta^\star + 1}^{\mathcal{N}_2}, \ldots, c_{2k_\delta^\star - 1}^{\mathcal{N}_2}$$

for $i = 0, 1, 2, \ldots, 2k_\delta^\star - 1$. In particular, the construction of the sequence $(\bar{c}_i^{\mathcal{N}_2})_{i=0}^\infty$ yields $c_{k_\delta^\star}^{\mathcal{N}_2} + c_0^{\mathcal{N}_2} \leq \bar{c}_0^{\mathcal{N}_2} + \bar{c}_1^{\mathcal{N}_2} = \gamma_2^{\mathcal{N}_2}$. Finally, the remaining parts of the maneuver can be dealt with analogously to case 1 showing that Condition (3.7) is ensured with the accumulated bounds defined by (3.22).

**Case 3**: let $\theta_0 \in (-\pi, -\pi + \phi)$ hold. First, the robot is turned $k_\delta^\star$ steps such that $\theta_{k_\delta^\star} = -\pi + \phi$ holds using the input $\mathbf{u}(i) = (0, \Delta\theta \cdot (k_\delta^\star \delta)^{-1})^\top$, $i = 0, 1, \ldots, k_\delta^\star - 1$, with $\Delta\theta = |\theta_0| - \pi + \phi$. The respective running costs are given by (3.23), which also satisfy Inequality (3.25) with $\theta_0$ replaced by $|\theta_0|$. Like in Case 2, the inequality $\ell(\mathbf{x_u}(i; \mathbf{x}_0), \mathbf{u}(i)) \leq \ell^\star(\mathbf{x}_0)$ holds for $i = 1, 2, \ldots, k_\delta^\star - 1$. During the second part of the maneuver the robot is driven to the origin in $l_\delta^\star$ steps; thus, Inequality (3.15) holds with the coefficients defined by (3.19) — indeed, $q_3(\pi/2)^4$ could have been dropped. Next, the robot is turned $k_\delta^\star$ steps until $\theta_\mathbf{u}(2k_\delta^\star + l_\delta^\star; \mathbf{x}_0) = -\pi/2$ holds using $\mathbf{u}(k_\delta^\star + l_\delta^\star + i) = (0, \Delta\theta \cdot (k_\delta^\star \delta)^{-1})^\top$ with $\Delta\theta = \pi/2 - \phi$ for $i \in \{0, 1, \ldots, k_\delta^\star - 1\}$. Hence, for $i = 0, 1, \ldots, k_\delta^\star - 1$, we have the running costs

$$\ell(\mathbf{x_u}(k_\delta^\star + l_\delta^\star + i; \mathbf{x}_0), \mathbf{u}(k_\delta^\star + l_\delta^\star + i)) = q_3 \left(\phi - \pi + \frac{i\Delta\theta}{k_\delta^\star}\right)^4 + r_2 \left(\frac{\Delta\theta}{k_\delta^\star \delta}\right)^4.$$

Similar to (3.24), we have the estimate

$$\left((\phi - \pi) + \frac{i\Delta\theta}{k_\delta^\star}\right)^4 \leq (\phi - \pi)^4 - \frac{i\Delta\theta}{k_\delta^\star}(\pi - \phi)^3.$$

Therefore, using Assumption (3.16), $\pi - \phi \geq \Delta\theta \cdot (k_\delta^\star \delta)^{-1}$, and $|\theta_0| \geq |\phi - \pi|$, one obtains the inequality

$$\ell(\mathbf{x_u}(k_\delta^\star + l_\delta^\star + i; z_0), \mathbf{u}(k_\delta^\star + l_\delta^\star + i)) \leq q_3 \theta_0^4 \leq \ell^\star(\mathbf{x}_0)$$

for $i = 1, 2, \ldots, k_\delta^\star - 1$. Then, the robot turns another $k_\delta^\star$ steps such that $\theta_\mathbf{u}(3k_\delta^\star + l_\delta^\star; z_0) = 0$ holds using the input $\mathbf{u}(2k_\delta^\star + l_\delta^\star + i) = (0, \pi \cdot (2k_\delta^\star \delta)^{-1})^\top$, $i \in \{0, 1, \ldots, k_\delta^\star - 1\}$. The resulting running costs for this part of the maneuver are given by (3.20) with $\phi = \pi/2$ and, thus, also satisfy Inequality (3.15) with coefficients $c_{2k_\delta^\star + l_\delta^\star + i}$, $i = 0, 1, \ldots, k_\delta^\star - 1$ defined by (3.21), respectively. We show that Case 3 is less costly than the reference case $\theta_0 = 0$ by the following calculations, in which Assumption (3.16), i.e. $r_2 \leq \frac{q_3 \cdot \delta}{2}$, is used:

$$\ell(\mathbf{x_u}(0; \mathbf{x}_0), \mathbf{u}(0)) + \ell(\mathbf{x_u}(k_\delta^\star + l_\delta^\star; \mathbf{x}_0), \mathbf{u}(k_\delta^\star + l_\delta^\star))$$
$$= \ell^\star(\mathbf{x}_0) + q_3(\pi - \phi)^4 + \frac{r_2 \overbrace{\left[(|\theta_0| - \pi + \phi)^4 + (\frac{\pi}{2} - \phi)^4\right]}^{\leq (|\theta_0| - \frac{\pi}{2})^4 \leq (\frac{\pi}{2})^4}}{(k_\delta^\star \delta)^4}$$
$$\leq \ell^\star(\mathbf{x}_0) + q_3 \theta_0^4 + \frac{q_3 \pi^4 \delta}{32(k_\delta^\star \delta)^4}$$
$$\leq \left(2 + \frac{q_3 \pi^4 \delta}{32(k_\delta^\star \delta)^4 \cdot s}\right) \ell^\star(\mathbf{x}_0) = \left(c_0^{\mathcal{N}_2} + c_{k_\delta^\star}^{\mathcal{N}_2}\right) \cdot \ell^\star(\mathbf{x}_0)$$

In conclusion, the accumulated bounds given by (3.22) can be used to ensure Condition (3.7) for the case considered here.

**Case 4:** let $\theta_0 \in (-\pi + \phi, 0)$ hold. First, for $i = 0, 1, \ldots, k_\delta^\star - 1$, the robot uses the control inputs $\mathbf{u}(i) = (0, \Delta\theta \cdot (k_\delta^\star \delta)^{-1})^\top$ with $\Delta\theta$ defined as $\max\{0, \phi - \pi/2 - \theta_0\}$ in order to achieve that $\phi - \theta_\mathbf{u}(k_\delta^\star; z_0) \leq \pi/2$ holds. Then, the robot employs $\mathbf{u}(k_\delta^\star + i) = $

$(0, (\phi - \theta_{\mathbf{u}}(k_\delta^\star; \mathbf{u}_0)) \cdot (k_\delta^\star \delta)^{-1})^\top$ for all $i \in \{0, 1, \ldots, k_\delta^\star - 1\}$, which yields $\theta_{\mathbf{u}}(2k_\delta^\star; \mathbf{x}_0) = \phi$.

Proceeding analogously to Case 2 leads to Estimate (3.25) for all $i\{0, 1, \ldots, k_\delta^\star - 1\}$ while the running costs $\ell(\mathbf{x_u}(k_\delta^\star + i; \mathbf{x}_0), \mathbf{u}(k_\delta^\star + i))$ for the next $k_\delta^\star$ steps are given by

$$q_1 x_0^4 + q_2 y_0^2 + q_3 \left(\theta_{k_\delta^\star} + \frac{i(\phi - \theta_{k_\delta^\star})}{k_\delta^\star}\right)^4 + r_2 \left(\frac{\phi - \theta_{k_\delta^\star}}{k_\delta^\star \delta}\right)^4$$

with $\theta_{k_\delta^\star} = \theta_{\mathbf{u}}(k_\delta^\star; \mathbf{x}_0)$ for all $i \in \{0, 1, \ldots, k_\delta^\star - 1\}$. Invoking Assumption (3.16), i.e. $r_2 \leq \frac{q_3 \cdot \delta}{2}$, yields the bound

$$2\ell^\star(\mathbf{x}_0) + q_3 \left[\underbrace{\theta_{k_\delta^\star}^4 - \theta_0^4}_{\leq -\Delta\theta \cdot |\theta_0|^3} + \frac{\Delta\theta}{2k_\delta^\star}\underbrace{\left(\frac{\Delta\theta}{k_\delta^\star \delta}\right)^3}_{\leq |\theta_0|^3} + \frac{\phi - \theta_{k_\delta^\star}}{2k_\delta^\star}\left(\frac{\phi - \theta_{k_\delta^\star}}{k_\delta^\star \delta}\right)^3\right]$$

for the running costs $\ell(\mathbf{x_u}(0; \mathbf{x}_0), \mathbf{u}(0)) + \ell(\mathbf{x_u}(k_\delta^\star; \mathbf{x}_0), \mathbf{u}(k_\delta^\star))$ and, thus, allows to derive the estimate

$$\ell(\mathbf{x_u}(0; \mathbf{x}_0), \mathbf{u}(0)) + \ell(\mathbf{x_u}(k_\delta^\star; \mathbf{x}_0), \mathbf{u}(k_\delta^\star)) \tag{3.26}$$
$$\leq \left(2 + \frac{q_3 \pi^4}{32 k_\delta^{\star 4} \delta^3 \cdot s}\right)\ell^\star(\mathbf{x}_0) = \left(c_0^{\mathcal{N}_2} + c_{k_\delta^\star}^{\mathcal{N}_2}\right) \cdot \ell^\star(\mathbf{x}_0).$$

The running costs $\ell(\mathbf{x_u}(i; \mathbf{x}_0), \mathbf{u}(i))$ can be estimated by $c_i^{\mathcal{N}_2} \ell^\star(\mathbf{x}_0)$ for all

$$i \in \{1, 2, \ldots, k_\delta^\star - 1\} \cup \{k_\delta^\star + 1, \ldots, 2k_\delta^\star - 1\},$$

see Case 2 and the derivation of the coefficients (3.18) for details while taking $\theta_{\mathbf{u}}(k_\delta^\star; \mathbf{x}_0) + i(\phi - \theta_{k_\delta^\star})/k_\delta^\star \leq i\phi/k_\delta^\star$ into account. Since the remaining parts of the maneuver are performed precisely as in Case 1, combining this with Inequality (3.26) shows that the accumulated bounds given by (3.22) can be used to ensure Condition (3.7) for the considered case.

### 3.5.2 Trajectory generation for $x_0 \in \mathcal{N}_1$

We consider initial conditions inside $\mathcal{N}_1$ with $\theta_0 = 0$ and construct a suitable coefficient sequence $(c_n^{\mathcal{N}_1})_{n \in \mathbb{N}_0}$ satisfying Inequality (3.15). Here, the particular choice of the stage costs $\ell$ is heavily exploited in order to successfully steer the robot from a position $(0, y_0, 0)^\top$ with $y_0 \neq 0$ to the origin while simultaneously deriving finitely many bounds $c_n^{\mathcal{N}_1}$, $n \in \mathbb{N}_0$. These bounds are on the one hand uniform in $y_0$, i.e. Inequality (3.15) holds independently of $y_0$ and, thus, also for $y_0 \to 0$. On the other hand, the number of coefficients $c_n^{\mathcal{N}_1}$, which are strictly greater than zero, is uniformly bounded. Combining these two properties ensures that the sequence remains summable – an important ingredient to make Proposition 3.4 applicable in order to ensure Assumption 3.2. Subsequently, we reorder this sequence in order to get $(\bar{c}_n^{\mathcal{N}_1})_{n \in \mathbb{N}_0}$ and prove that the resulting bounds $\gamma_i^{\mathcal{N}_1} := \sum_{n=0}^{i-1} \bar{c}_n^{\mathcal{N}_1}$ also yield Inequality (3.13) for the case $\theta_0 \neq 0$.

**Initial Conditions $x_0 \in \mathcal{N}_1$ with $\theta_0 = 0$:** The following maneuver is used in order to derive bounds $\gamma_i^{\mathcal{N}_1}$, $i \in \mathbb{N}_{\geq 2}$, satisfying Inequality (3.7) for initial condition whose angular deviation is equal to zero:

a) drive towards the $y$-axis until $(0, y_0, 0)^\top$ is reached.

b) drive forward while slightly steering in order to reduce the $y$-component to $y_0/2$; a position $(\tilde{x}, y_0/2, 0)^\top$ for some $\tilde{x} > 0$ is reached.

c) carry out a symmetric maneuver while driving backward so that the origin $0_{\mathbb{R}^3}$ is reached.

See Figure 3.5 for a visualization for parts b) and c) of the maneuver. The number of steps needed in order to perform this maneuver depends on the constraints and the sampling time $\delta$, which is supposed to satisfy $i \cdot \delta = 1$ for some integer $i \in \mathbb{N}$ as in Section 3.5.1.

To this end, we define

$$k_\delta^\star := \left\lceil \frac{\pi}{\min\{\bar{\omega}, \pi\} \cdot \delta} \right\rceil, l_\delta^\star := \left\lceil \frac{\sqrt[4]{s/q_1}}{\min\{\bar{v}, \sqrt[4]{s/q_1}\} \cdot \delta} \right\rceil,$$

where the vehicle can turn by 180 degrees in $k_\delta^\star$ steps and drive to the $y$-axis in $l_\delta^\star$ steps, respectively. In addition to Inequality (3.16), the Condition

$$r_1 \leq \frac{q_1 \cdot \delta}{2} \tag{3.27}$$

is assumed in order to keep the presentation technically simple.

Initial conditions $\mathbf{x}_0 = (x_0, y_0, 0)^\top \geq 0$ are considered first. Firstly, the vehicle does not move for $k_\delta^\star$ steps. Hence, Inequality (3.15) holds with $c_i^{\mathcal{N}_1} = 1$ for $i \in \{0, 1, \ldots, k_\delta^\star - 1\}$. Then, the robot drives towards the $y$-axis in $l_\delta^\star$ steps using $\mathbf{u}(k_\delta^\star + i) = (-x_0 \cdot (l_\delta^\star \delta)^{-1}, 0)^\top \in U$, $i = 0, 1, \ldots, l_\delta^\star - 1$, which allows to estimate $\ell(\mathbf{x_u}(k_\delta^\star + i; \mathbf{x}_0), \mathbf{u}(k_\delta^\star + i))$ first by

$$q_1 x_0^4 (1 - i/l_\delta^\star)^4 + q_2 y_0^2 + q_1 x_0^4 / (2 l_\delta^\star (l_\delta^\star \delta)^3)$$

using (3.27) and, then, by

$$\ell^\star(\mathbf{x}_0) - q_1 \left(\frac{x_0^4}{l_\delta^\star}\right) \left[i - \frac{1}{2(l_\delta^\star \delta)^3}\right].$$

Hence, Inequality (3.15) holds with $c_0^{\mathcal{N}_1} = 1 + (2l_\delta^\star (l_\delta^\star \delta)^3)^{-1}$ and $c_{k_\delta^\star + i}^{\mathcal{N}_1} = 1$ for all $i \in \{1, 2, \ldots, l_\delta^\star - 1\}$.

The next part of the maneuver is performed in four seconds with constant control effort $\|\mathbf{u}(\cdot)\|$ such that the angle is decreased to $-\arctan(\sqrt{y_0})$ during the first second and then put back to zero while the $y$-position of the robot decreases to $y_0/2$. Afterwards, these two moves are carried out backwards in order to reach the origin, see Figure 3.5 on

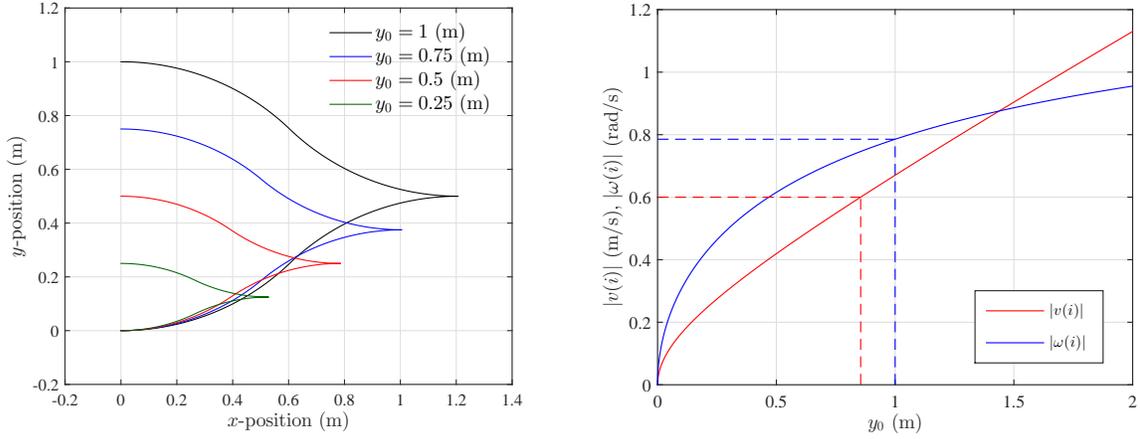

Fig. 3.5: Trajectories of parts b) and c) of the maneuver starting from different initial conditions on the $y$-axis for $\mathbf{x}_0 \in \mathcal{N}_1$ (left). The respective controls are displayed on the right.

the left. To this end, the controls

$$\omega(i) = -\omega(\delta^{-1} + i) = -\omega(2\delta^{-1} + i) = \omega(3\delta^{-1} + i),$$
$$v(i) = v(\delta^{-1} + i) = -v(2\delta^{-1} + i) = -v(3\delta^{-1} + i),$$

$i \in \{k_\delta^\star + l_\delta^\star, k_\delta^\star + l_\delta^\star + 1, \ldots, k_\delta^\star + l_\delta^\star + \delta^{-1} - 1\}$, with

$$\omega(i) = -\arctan(\sqrt{y_0}) \quad \text{and} \quad v(i) = -\frac{y_0 \arctan(\sqrt{y_0})}{\frac{4}{\sqrt{y_0+1}} - 4} \tag{3.28}$$

are employed. Note that this strategy ensures not to move when starting at the origin. The resulting $y$-positions are given by $y(k_\delta^\star + l_\delta^\star + n\delta^{-1}) = (1 - n/4)y_0$ for $n \in \{0, 1, 2, 3\}$ while the $x$-positions are, for $i = k_\delta^\star + l_\delta^\star$, given by $x(i) = 0$, $x(i + \delta^{-1}) = x(i + 3\delta^{-1}) = \sin(\omega(i)) \cdot v(i)/\omega(i)$, and $x(i + 2\delta^{-1}) = 2x(i + \delta^{-1})$. The maneuver has to be suitably adapted if either control constraints enforce $v(\cdot)$ or $\omega(\cdot)$ to be smaller or $x(k_\delta^\star + l_\delta^\star + 2\delta^{-1})$ violates the state constraints. However, since this maneuver is constructed for small $y_0$, constraints can be neglected.

Next, we evaluate the running costs and determine coefficients $c^{\mathcal{N}_1}_{k^\star_\delta+l^\star_\delta+i}$ such that Inequality (3.15) holds for all $i \in \{0, 1, \ldots, 4\delta^{-1} - 1\}$. To this end, the estimates

$$\arctan^2 \sqrt{y_0} \leq y_0 \quad \text{and} \quad v(i)^4 \leq (2 + 3y_0 + y_0^2)^2 y_0^2 / 64$$

for $v(i)$ from (3.28) are employed where, for the derivation of the latter, the two auxiliary inequalities

$$(2 + y_0)^2 \cdot (y_0 + 1)^2 \leq (2 + 3y_0 + y_0^2)^2, \text{ and}$$

$$y_0^2 \leq 2(2 + y_0) \cdot (\sqrt{y_0 + 1} - 1)^2$$

were exploited. Hence, using $(2 + 3y_0 + y_0^2) \leq (y_0 + 1.5)^2$ and $q_2 y_0^2 \leq \ell^\star(\mathbf{x}_0)$, we obtain that the running costs $\ell(\mathbf{x_u}(k^\star_\delta + l^\star_\delta); \mathbf{x}_0), \mathbf{u}(k^\star_\delta + l^\star_\delta))$ are bounded by

$$\left(1 + \frac{(\sqrt{s/q_2} + 1.5)^4 r_1}{64 q_2} + r_2/q_2\right) \ell^\star(\mathbf{x}_0) =: c^{\mathcal{N}_1}_{k^\star_\delta+l^\star_\delta} \ell^\star(\mathbf{x}_0)$$

and, thus, Inequality (3.15) holds. Then, since $\sin^2(\omega(i)) \leq \omega(i)^2$ holds for $\omega(i)$ from (3.28), $\ell^\star(\mathbf{x_u}(k^\star_\delta + l^\star_\delta + \delta^{-1})) \leq q_1 v(i)^4 + (9/16) q_2 y_0^2 + q_3 y_0^2$ this yields Inequality (3.15) with $c^{\mathcal{N}_1}_{k^\star_\delta+l^\star_\delta+\delta^{-1}}$ given by

$$9/16 + \left(q_3 + r_2 + (q_1 + r_1)(\sqrt{s/q_2} + 1.5)^4/64\right) q_2^{-1}.$$

Analogously, the coefficients $c^{\mathcal{N}_1}_{k^\star_\delta+l^\star_\delta+2\delta^{-1}}$ and $c^{\mathcal{N}_1}_{k^\star_\delta+l^\star_\delta+3\delta^{-1}}$ defined by

$$1/4 + \left(r_2 + (16 q_1 + r_1)(\sqrt{s/q_2} + 1.5)^4/64\right) q_2^{-1} \quad \text{and}$$

$$1/16 + \left(q_3 + r_2 + (q_1 + r_1)(\sqrt{s/q_2} + 1.5)^4/64\right) q_2^{-1}$$

are derived. Forسampling time $\delta < 1$, further coefficients have to be determined. To this end, the running costs $\ell(\cdot, \cdot)$ at time

$$k_\delta^\star + l_\delta^\star + n\delta^{-1} + i, \qquad (n, i) \in \{0, 1, 2, 3\} \times \{1, 2, \ldots, \delta^{-1} - 1\},$$

are overestimated by plugging in the state

$$\begin{pmatrix} x_{\mathbf{u}}(k_\delta^\star + l_\delta^\star + (2.125 - 0.5(n - 1.5)^2) \cdot \delta^{-1}; \mathbf{x}_0) \\ (1 - 0.25n)y_0 \\ \theta_{\mathbf{u}}(k_\delta^\star + l_\delta^\star + \delta^{-1}; \mathbf{x}_0) \end{pmatrix}$$

instead of $\mathbf{x_u}(k_\delta^\star + l_\delta^\star + n\delta^{-1} + i)$ while leaving the control as it is. This yields, for $i \in \{1, 2, \ldots, \delta^{-1} - 1\}$, Inequality (3.15) with the coefficients $c_{k_\delta^\star + l_\delta^\star + n\delta^{-1} + i}^{\mathcal{N}_1}$ defined by

$$c_{k_\delta^\star + l_\delta^\star + n\delta^{-1}}^{\mathcal{N}_1} + \begin{cases} \left(q_1(\sqrt{s/q_2} + 1.5)^4/64 + q_3\right) q_2^{-1}, & n = 0 \\ \left(15q_1(\sqrt{s/q_2} + 1.5)^4/64\right) q_2^{-1}, & n = 1 \\ q_3 q_2^{-1}, & n = 2 \\ 0, & n = 3 \end{cases}. \qquad (3.29)$$

The coefficients $c_i^{\mathcal{N}_1}$, $i = 1, 2, \ldots, k_\delta^\star + l_\delta^\star + 4\delta^{-1} - 1$, are ordered descendingly in order to construct a new sequence $(\bar{c}_i^{\mathcal{N}_1})_{i=1}^{k_\delta^\star + l_\delta^\star + 4\delta^{-1} - 1}$ such that the property $\bar{c}_{i-1}^{\mathcal{N}_1} \geq \bar{c}_i^{\mathcal{N}_1}$ holds for all $i \in \{2, 3, \ldots, k_\delta^\star + l_\delta^\star + 4\delta^{-1} - 1\}$. Then, setting $\bar{c}_0^{\mathcal{N}_1} = c_0^{\mathcal{N}_1}$ and $\bar{c}_i^{\mathcal{N}_1} = 0$ for all $i \geq k_\delta^\star + l_\delta^\star + 4\delta^{-1}$ yields $(\bar{c}_i^{\mathcal{N}_1})_{i=0}^\infty$. In conclusion, the accumulated bounds $(\gamma_i^{\mathcal{N}_1})_{i \in \mathbb{N}_{\geq 2}}$ of Condition (3.7) for the second maneuver are given by

$$\gamma_i^{\mathcal{N}_1} := \sum_{n=0}^{i-1} \bar{c}_i^{\mathcal{N}_1}, \qquad i \in \mathbb{N}_{\geq 2}. \qquad (3.30)$$

**Initial Conditions $\mathbf{x}_0 \in \mathcal{N}_1$ with $\theta_0 \neq 0$:** Next, we show that Condition (3.7) with $\gamma_i^{\mathcal{N}_1}$, $i \in \mathbb{N}_{\geq 2}$, holds also for $\mathbf{x}_0$ with $\theta_0 \in [-\pi, 0) \cup (0, \pi)$ and, thus, for all initial conditions $\mathbf{x}_0 \in \mathcal{N}_1$. Firstly, the robot turns $k_\delta^\star$ steps using $\mathbf{u}(i) = (0, -\theta_0/k_\delta^\star \delta)^\top$, $i = 0, 1, \ldots, k_\delta^\star - 1$, such that $\theta_\mathbf{u}(k_\delta^\star; \mathbf{x}_0) = 0$ is attained. This yields the running costs

$$\ell(\mathbf{x}_\mathbf{u}(i; \mathbf{x}_0), \mathbf{u}(i)) = q_1 x_0^4 + q_2 y_0^2 + q_3 \theta_0^4 \left[1 - \left(\frac{i}{k_\delta^\star}\right)\right]^4 + r_2 \left(\frac{\theta_0}{k_\delta^\star \delta}\right)^4$$

Using $1 - (i/k_\delta^\star) \in [0, 1]$ and Assumption (3.16) leads to

$$\ell(\mathbf{x}_\mathbf{u}(i; \mathbf{x}_0), \mathbf{u}(i)) \leq \ell^\star(\mathbf{x}_0) - q_3 \left(\frac{\theta_0^4}{k_\delta^\star}\right) \left[i - \frac{1}{2(k_\delta^\star \delta)^3}\right]$$

for $i \in \{0, 1, \ldots, k_\delta^\star - 1\}$. Hence, the right hand side of this inequality is always less or equal $\ell^\star(\mathbf{x}_0)$ for $i > 0$. The remaining parts of the maneuver are performed as before. We show that this case is less costly than its counterpart $\theta_0 = 0$ by the following calculations, in which the abbreviation $\Xi := r_1 v(k_\delta^\star + l_\delta^\star)^4 + r_2 \omega(k_\delta^\star + l_\delta^\star)^4$ is used:

$$\sum_{i \in \{0, k_\delta^\star + l_\delta^\star\}} \ell(\mathbf{x}_\mathbf{u}(i; \mathbf{x}_0), \mathbf{u}(i)) = \ell^\star(\mathbf{x}_0) + q_2 y_0^2 + r_2 \left(\frac{\theta_0}{k_\delta^\star \delta}\right)^4 + \Xi$$

$$\overset{(3.16)}{\leq} 2 \cdot \ell^\star(\mathbf{x}_0) + \Xi \leq \left(c_0^{\mathcal{N}_1} + c_{k_\delta^\star + l_\delta^\star}^{\mathcal{N}_1}\right) \cdot \ell^\star(\mathbf{x}_0).$$

In conclusion, the accumulated bounds given by (3.30) can be used to ensure Condition (3.7) for initial conditions $\mathbf{x}_0$ with $\theta_0 \neq 0$.

## 3.6 Numerical Results

In the preceding section, bounds

$$\gamma_i = \min\{i, \max\{\gamma_i^{\mathcal{N}_1}, \gamma_i^{\mathcal{N}_2}\}\}, \text{ for } i \in \mathbb{N}_{\geq 2}$$

satisfying Assumption 3.2 were deduced, see (3.14). Here, the bounds $\gamma_i^{\mathcal{N}_2} = \sum_{n=0}^{i-1} \bar{c}_n^{\mathcal{N}_2}$ were constructed according to the procedure presented in the paragraph before (3.22) based on the coefficients $\bar{c}_n^{\mathcal{N}_2}$ displayed in (3.18), (3.19), and (3.21). Similarly, $\gamma_i^{\mathcal{N}_1} = \sum_{n=0}^{i-1} \bar{c}_n^{\mathcal{N}_1}$ are derived using (3.29). In the following, a prediction horizon $N$ is determined such that the resulting MPC closed-loop is asymptotically stable – based on these bounds $\gamma_i$, $i \in \mathbb{N}_{\geq 2}$. To this end, the minimal stabilizing horizon $\hat{N}$ is defined as

$$\min \left\{ N \in \mathbb{N}_{\geq 2} : \alpha_N = 1 - \frac{(\gamma_N - 1) \prod_{k=2}^{N}(\gamma_k - 1)}{\prod_{k=2}^{N} \gamma_k - \prod_{k=2}^{N}(\gamma_k - 1)} > 0 \right\},$$

a quantity, which depends on the sampling rate $\delta$ and the weighting coefficients of the running cost $\ell(\cdot,\cdot)$. Then, a comparison with quadratic running costs is presented in Subsection 3.6.2 before, in Section 3.6.3, numerical simulations are conducted in order to show that MPC without stabilizing constraints or costs steers differential drive robots to a desired equilibrium.

### 3.6.1 Computation of the minimal stabilizing horizon

Let the sets $U = [-0.6, 0.6] \times [-\pi/4, \pi/4]$ and $X := [-2, 2]^2 \times \mathbb{R}$ be given. Moreover, the weighting parameters $q_1 = 1$, $q_3 = 0.1$, $r_1 = q_1 \delta/2$, and $r_2 = q_3 \delta/2$ of the running costs $\ell(\cdot,\cdot)$ are defined depending on the sampling time $\delta$. Then, for a given sampling time $\delta$ and weighting coefficient $q_2$, the minimal stabilizing horizon $\hat{N}$ can be computed by Algorithm 3.2.

The only optimization is carried out in Step 3 of Algorithm 3.2. This can be done by a line search. Here, the optimization variable $s$ can be restricted to a compact interval depending on the size of the state and control constraints, and the weighting parameter $q_2$.

The results of Algorithm 3.2 for sampling time $\delta = 1$ and weighting coefficient $q_2 \in$

**Algorithm 3.2** Calculating the minimal stabilizing horizon $\hat{N}$

**Given:** Control bounds $\bar{v}$, $\bar{\omega}$, box constraints $\bar{x}$, $\bar{y}$, weighting coefficients $q_1$, $q_2$, $q_3$, $r_1$, $r_2$, and sampling time $\delta$.

**Initialization:** Set $N = 1$ and $\alpha = 0$.

1: **while** $\alpha = 0$ **do**
2:     Increment $N$.
3:     Minimize $\gamma_N^\star := \max\{\gamma_N^{\mathcal{N}_2}, \gamma_N^{\mathcal{N}_1}\}$ subject to $s \in \mathbb{R}_{\geq 0}$, (3.22), and (3.30).
4:     Define $\gamma_N := \min\{N, \gamma_N^\star\}$.
5:     Set $\alpha = 1 - \frac{(\gamma_N - 1)\prod_{k=2}^N (\gamma_k - 1)}{\prod_{k=2}^N \gamma_k - \prod_{k=2}^N (\gamma_k - 1)}$.
6: **end while**

**Output** Minimal stabilizing horizon length $\hat{N} = N$ and the performance index $\alpha_{\hat{N}} = \alpha$.

Table 3.1: Minimal stabilizing horizon $\hat{N}$ in dependence on the sampling time $\delta$ and the weighting parameter $q_2$ for $q_1 = 1$, $q_3 = 0.1$, $r_1 = q_1\delta/2$, and $r_2 = q_3\delta/2$.

| Sampling time $\delta$ | $\hat{N}(\hat{N} \cdot \delta(\text{seconds}))$ | | | |
| (seconds) | $q_2 = 2$ | $q_2 = 5$ | $q_2 = 10$ | $q_2 = 100$ |
|---|---|---|---|---|
| 1.00 | 12(12) | 10(10) | 8(8) | 8(8) |
| 0.50 | 25(12.5) | 19(9.5) | 16(8) | 15(7.5) |
| 0.25 | 48(12) | 37(9.25) | 32(8) | 29(7.25) |
| 0.10 | 122(12.2) | 93(9.3) | 79(7.9) | 70(7) |

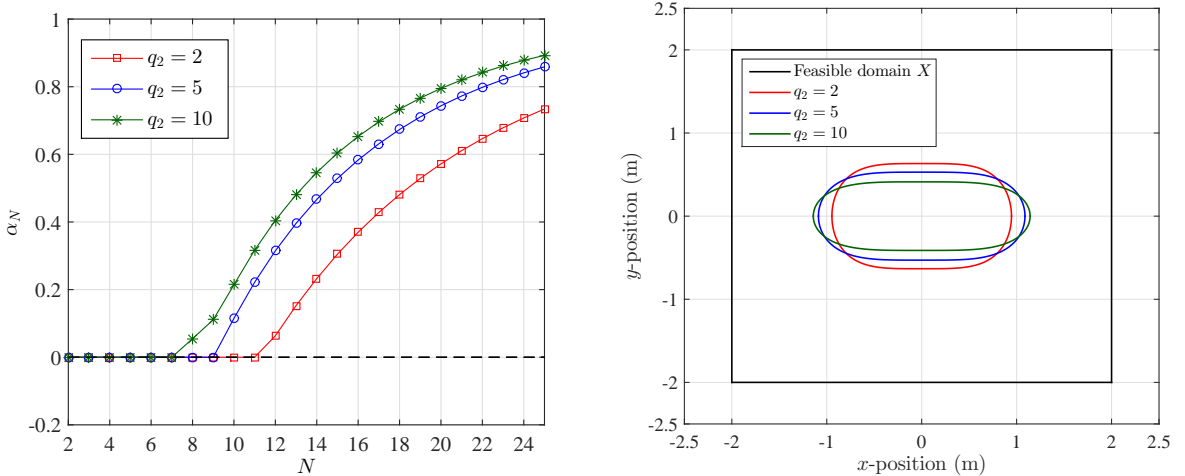

Fig. 3.6: Left: Dependence of the performance bound $\alpha_N$ on the prediction horizon $N$ and different values of $q_2$. Right: Dependence of the set $\mathcal{N}_1$ radius, i.e. $s$, on different values of $q_2$. For both sub figures, we have sampling time $\delta = 1$ and $q_1 = 1$, $q_3 = 0.1$, $r_1 = q_1\delta/2$, $r_2 = q_3\delta/2$.

$\{2, 5, 10\}$ are presented in Figure 3.6 (left). Using a larger weighting coefficient $q_2$ results in smaller minimal stabilizing horizons $\hat{N}$. Moreover, it can be seen that the suboptimality index $\alpha$ converges to one for prediction horizon $N$ tending to infinity. Furthermore, it is observed that the radius $s$ of the set $\mathcal{N}_1$ increases for larger $q_2$, i.e. $s = 0.8$ ($q_2 = 2$), $s = 1.4$ ($q_2 = 5$), and $s = 1.7$ ($q_2 = 10$), see Figure 3.6 (right) for visualization of the set $\mathcal{N}_1$. In contrast to that, the influence of the sampling time $\delta$ is negligible, see Table 3.1.

### 3.6.2 Comparison with quadratic running costs

Here, the proposed MPC scheme without stabilizing constraints or costs, i.e. Algorithm 3.1, is applied in order to stabilize a unicycle mobile robot to the origin. The constraints and weighting coefficients of the running costs $\ell$ in this subsection and the subsequent one are the same as in the preceding subsection with $q_2 = 5$ and sampling time $\delta = 0.25$. In this case, the theoretically calculated minimal stabilizing horizon is given by $\hat{N} = 37$, see Table 3.1. All simulations have been run using the *Matlab* routine fmincon to solve the optimal control problem in each MPC step. However, for a real time implementation we recommend ACADO toolkit [31]. The MPC performance is investigated through two sets of numerical simulations – on the one hand under the proposed running costs (3.5); on the other hand using the standard quadratic running costs with weighting matrices $\mathbf{Q} = \mathrm{diag}(q_1, q_2, q_3)$ and $\mathbf{R} = \mathrm{diag}(r_1, r_2)$.

First, the initial state of the robot is chosen to be $\mathbf{x}_0 = (0, 0.1, 0)^\top$, i.e. located on the $y$-axis, close to the origin, and with an orientation angle of zero. Both controllers steer the robot close to the origin, but only the MPC controller with the proposed running costs fulfils the control objective of steering the robot to the origin, see Figure 3.7. This conclusion can be also inferred from the scaled value function $V_N(\mathbf{x}_{\boldsymbol{\mu}_N}(n; \mathbf{x}_0)) \cdot \ell^\star(\mathbf{x}_0)^{-1}$, $n \in \mathbb{N}_0$, evaluated along the closed-loop trajectories as depicted in Figure 3.8 (left). Since the value function does not decrease anymore after a few ($n \approx 12$) time steps, MPC with quadratic running costs fails to ensure asymptotic stability for the chosen prediction horizon $N = 37$.

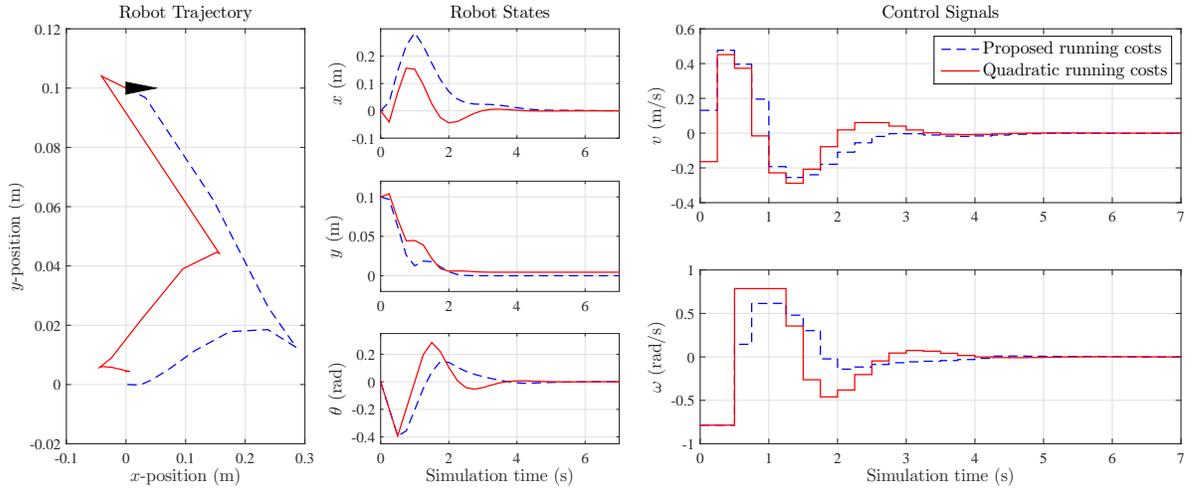

Fig. 3.7: MPC Closed loop state trajectory and employed controls for sampling time $\delta = 0.25$ and prediction horizon $N = 37$ under the proposed and quadratic running costs with weighting matrices $\mathbf{Q} = \text{diag}(q_1, q_2, q_3)$ and $\mathbf{R} = \text{diag}(r_1, r_2)$.

Moreover, since uniform boundedness of $\sup_{\mathbf{x}_0 \in X} V_N(\mathbf{x}_0) \cdot \ell^\star(\mathbf{x}_0)^{-1}$ with respect to the prediction horizon $N$ is a necessary condition for asymptotic stability of the MPC closed-loop, we further investigate this quantity. To this end, three initial conditions $\mathbf{x}_0 = (0, y_0, 0)^\top, y_0 \in \{0.1, 0.01, 0.001\}$, are considered, see Figure 3.8 (right). Under the proposed stage costs, the quantity $V_N(\mathbf{x}_0) \cdot \ell^\star(\mathbf{x}_0)^{-1}$ is bounded for all chosen initial conditions. In contrary to this, for quadratic running costs, the quantity $V_N(\mathbf{x}_0) \cdot \ell^\star(\mathbf{x}_0)^{-1}$ grows unboundedly for decreasing $y_0$-component, e.g. $y_0 = 10^{-i}, i \in \mathbb{N}$, in our numerical simulations. Indeed, this observation was also made for different weighting coefficients and prediction horizons. Even in the setting with stabilizing terminal constraints and costs [11, 21], non-quadratic terminal costs were deployed. We conjecture that Assumption 3.2 cannot be satisfied for quadratic running costs. In conclusion, using a non-quadratic running cost $\ell(\cdot, \cdot)$ like (3.5) seems to be necessary in order to ensure asymptotic stability of the MPC closed-loop without stabilizing constraints or costs.

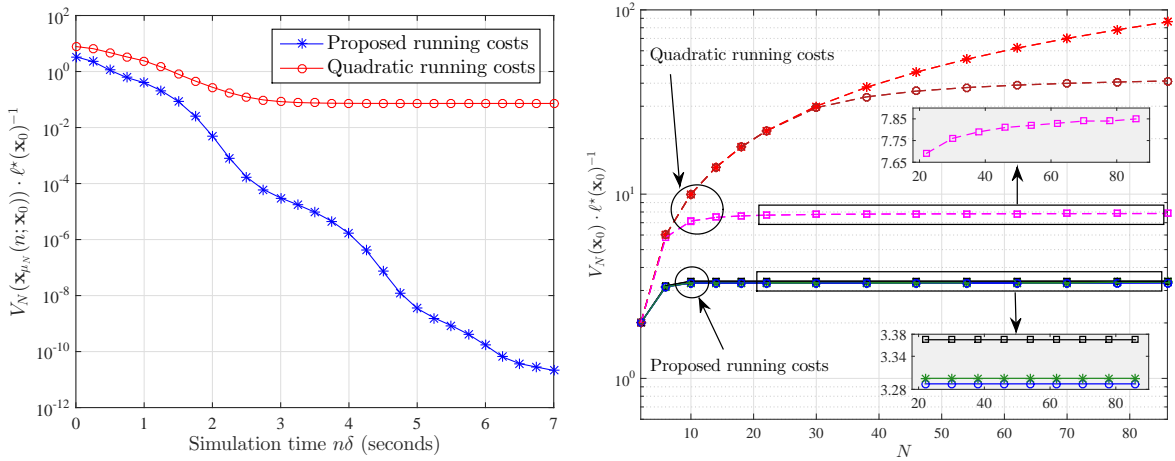

Fig. 3.8: Left: Evolution of $V_N(\mathbf{x}_{\boldsymbol{\mu}_N}(n;\mathbf{x}_0)) \cdot \ell^\star(\mathbf{x}_0)^{-1}$, $n \in \{0,1,\ldots,28\}$, for $\mathbf{x}_0 = (0,0.1,0)^\top$, $\delta = 0.25$, and $N = 37$. Right: Evaluation of $V_N(\mathbf{x}_0) \cdot \ell^\star(\mathbf{x}_0)^{-1}$ for $N = 2,3,\ldots,86$ for the proposed and quadratic running costs with initial conditions: $\mathbf{x}_0 = (0,0.1,0)^\top$ ($\square$), ($\square$); $\mathbf{x}_0 = (0,0.01,0)^\top$ ($\circ$), ($\circ$); and $\mathbf{x}_0 = (0,0.001,0)^\top$ ($\ast$), ($\ast$). $\delta = 0.25$.

### 3.6.3 Numerical investigation of the required horizon length

In this subsection, the minimal stabilizing horizon $\hat{N}$ is numerically examined for the MPC controller. To this end, the evolution of the value function $V_N(\mathbf{x}_{\boldsymbol{\mu}_N}(n;\mathbf{x}_0))$, $n \in \mathbb{N}_0$, along the MPC closed-loop, using the proposed running costs (3.5) for initial conditions $\mathbf{x}_0 = (0, 10^{-i}, 0)^\top$, $i \in \{0,1,2,3,4,5\}$, is considered, see Figure 3.9 (left). If the value function decays strictly, the relaxed Lyapunov inequality (3.10) holds — a sufficient stability condition, see [32]. Hence, we compute the minimal prediction horizon such that this stability condition is satisfied until a numerical tolerance is reached, i.e. $V_N(\mathbf{x}_{\boldsymbol{\mu}_N}(n;\mathbf{x}_0)) \leq 3 \cdot 10^{-11}$ as shown in Figure 3.9 (right).

So far, we concentrated on very particular initial conditions. Now, the ability of the proposed MPC controller to stabilize a unicycle mobile robot to an equilibrium point is demonstrated. To this end, eight initial positions evenly distributed along a large circle of 1.9 (m) radius, as well as five initial positions distributed along a small circle of 0.1 (m) radius, are selected. The initial orientation angle $\theta_0$ is randomly chosen from the set $\{i \cdot \pi/4 | i \in \{0, 1, 2, 3, 4, 5, 6, 7\}\}$. The prediction horizon $N$ is chosen such that the value function $V_N(\mathbf{x}_{\boldsymbol{\mu}_N}(n;\mathbf{x}_0))$, $n \in \mathbb{N}_0$, evaluated along the closed-loop reaches a neighbourhood of the origin corresponding to a reference magnitude of $10^{-9}$ for initial conditions on the large circle depicted and $10^{-11}$ for initial conditions on the small circle, which is illustrated in Figure 3.10. It is observed that stabilizing horizons of $N = 7$ and $N = 15$ are required for the initial conditions located on the large and small circles, respectively.

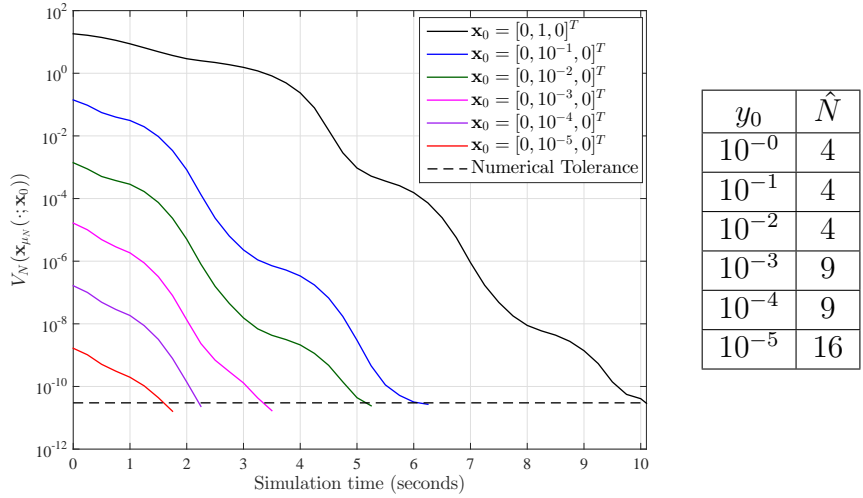

Fig. 3.9: Evolution of $V_N(\mathbf{x}_{\boldsymbol{\mu}_N}(\cdot;\mathbf{x}_0))$ along the closed-loop trajectories (left) and numerically computed stabilizing prediction horizons $\hat{N}$ (right) for sampling time $\delta = 0.25$ and different initial conditions.

Our numerical simulations show that the required prediction horizon $N$ rapidly grows if the initial condition is located (very) close to the origin. Otherwise, much shorter horizons $N$ are sufficient to steer the robot (very close) to the desired equilibrium. Independently of this observation, the numerically calculated stabilizing prediction horizon is shorter than its theoretically derived bound $\hat{N} = 37$. However, the calculated stabilizing horizon $\hat{N}$ holds for all initial states $\mathbf{x}_0$ in the feasible domain $X$. Moreover, both the estimates and the maneuvers used in order to derive $\gamma_N^{\mathcal{N}_2}$ and $\gamma_N^{\mathcal{N}_1}$ given by (3.22) and (3.30), respectively, are not optimal as highlighted in Section 3.5. Hence, the derived estimate

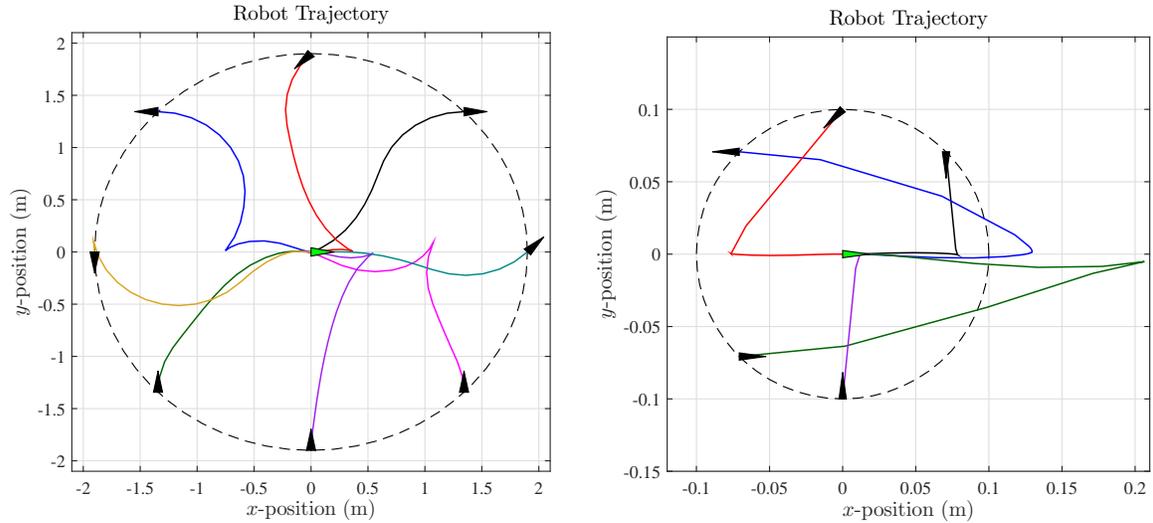

Fig. 3.10: MPC closed-loop trajectories emanating from initial conditions $(x_0, y_0)^\top$ on the circle of radius 1.9 (left; $N = 7$) and 0.1 (right; $N = 15$) using the proposed running costs. The initial state and orientation is indicated by the filled (black) triangles ($\delta = 0.25$).

of $\hat{N}$ can be further improved.

## 3.7 Conclusions

In this chapter, a stabilizing MPC controller is developed for the regulation problem of unicycle non-holonomic mobile robots. Unlike the common stabilizing schemes presented in the literature, where terminal constraints and/or costs are adopted, asymptotic stability of the developed controller is guaranteed by the combination of suitably chosen running costs and prediction horizon. Herein, the design of the running costs reflects the task to control both the position and the orientation of the robot and, thus, penalizes the direction orthogonal to the desired orientation more than other directions. Then, open loop trajectories are constructed in order to derive bounds on the value function and to determine the length of the prediction horizon such that asymptotic stability of the MPC closed-loop can be rigorously proven. The presented proof of concept can serve as a blueprint for deducing stability properties of similar applications. Finally, numerical simulations are conducted in order to examine the proposed controller and assess its performance in comparison with a controller based on quadratic running costs.

# References

[1] R. Siegwart, I. R. Nourbakhsh, and D. Scaramuzza, *Introduction to autonomous mobile robots*. MIT press, 2011.

[2] R. Velazquez and A. Lay-Ekuakille, "A review of models and structures for wheeled mobile robots: Four case studies," in *Proceedings of the 15th International Conference on Advanced Robotics (ICAR)*, 2011, pp. 524–529.

[3] F. Xie and R. Fierro, "First-state contractive model predictive control of nonholonomic mobile robots," in *Proceedings of the American Control Conference*, 2008, pp. 3494–3499.

[4] R. W. Brockett, "Asymptotic stability and feedback stabilization," in *Differential Geometric Control Theory*, R. W. Brockett, R. S. Millman, and H. J. Sussmann, Eds. Birkhäuser, Boston, MA, 1983, pp. 181–191.

[5] A. Astolfi, "Discontinuous control of nonholonomic systems," *Systems & Control Letters*, vol. 27, no. 1, pp. 37 – 45, 1996.

[6] D. Gu and H. Hu, "A Stabilizing Receding Horizon Regulator for Nonholonomic Mobile Robots," *IEEE Transactions on Robotics*, vol. 21, no. 5, pp. 1022–1028, 2005.

[7] W. Dixon, D. Dawson, F. Zhang, and E. Zergeroglu, "Global exponential tracking control of a mobile robot system via a pe condition," in *Proceedings of the 38th IEEE Conference on Decision and Control*, vol. 5, 1999, pp. 4822–4827.

[8] T.-C. Lee, K.-T. Song, C.-H. Lee, and C.-C. Teng, "Tracking Control of Unicycle-Modeled Mobile Robots Using a Saturation Feedback Controller," *IEEE Transactions on Control Systems Technology*, vol. 9, no. 2, pp. 305–318, 2001.

[9] M. Michalek and K. Kozowski, "Vector-field-orientation feedback control method for a differentially driven vehicle," *IEEE Transactions on Control Systems Technology*, vol. 18, no. 1, pp. 45–65, 2010.

[10] F. Kühne, W. F. Lages, and J. M. Gomes da Silva Jr., "Point Stabilization of Mobile Robots with Nonlinear Model Predictive Control," in *Proceedings of the IEEE International Conference on Mechatronics and Automation*, vol. 3, 2005, pp. 1163–1168 Vol. 3.

[11] F. A. Fontes, "A general framework to design stabilizing nonlinear model predictive controllers," *Systems & Control Letters*, vol. 42, no. 2, pp. 127–143, 2001.

[12] D. Gu and H. Hu, "Receding horizon tracking control of wheeled mobile robots," *IEEE Transactions on Control Systems Technology*, vol. 14, no. 4, pp. 743–749, 2006.

[13] Y. Zhu and U. Ozguner, "Robustness analysis on constrained model predictive control for nonholonomic vehicle regulation," in *Proceedings of the American Control Conference*, 2009, pp. 3896–3901.

[14] J. B. Rawlings and D. Q. Mayne, *Model Predictive Control: Theory and Design*. Nob Hill Publishing, 2009.

[15] T. Raff, S. Huber, Z. K. Nagy, and F. Allgöwer, "Nonlinear model predictive control of a four tank system: An experimental stability study," in *Proceedings of the IEEE Conference on Control Applications*, 2006, pp. 237–242.

[16] S. Keerthi and E. Gilbert, "Optimal infinite-horizon feedback laws for a general class of constrained discrete-time systems: Stability and moving-horizon approximations," *Journal of Optimization Theory and Applications*, vol. 57, no. 2, pp. 265–293, 1988.

[17] H. Chen and F. Allgöwer, "A quasi-infinite horizon nonlinear model predictive control scheme with guaranteed stability," *Automatica*, vol. 34, no. 10, pp. 1205–1218, 1998.

[18] J. A. Primbs and V. Nevistić, "Feasibility and stability of constrained finite receding horizon control," *Automatica*, vol. 36, no. 7, pp. 965–971, 2000.

[19] S. Tuna, M. Messina, and A. Teel, "Shorter horizons for model predictive control," in *Proceedings of the American Control Conference*, 2006, pp. 863–868.

[20] L. Grüne, J. Pannek, M. Seehafer, and K. Worthmann, "Analysis of unconstrained nonlinear MPC schemes with varying control horizon," *SIAM J. Control Optim.*, vol. Vol. 48 (8), pp. 4938–4962, 2010.

[21] F. Fontes and L. Magni, "Min-max model predictive control of nonlinear systems using discontinuous feedbacks," *IEEE Transactions on Automatic Control*, vol. 48, no. 10, pp. 1750–1755, 2003.

[22] G. Grimm, M. Messina, S. Tuna, and A. Teel, "Model predictive control: for want of a local control lyapunov function, all is not lost," *IEEE Transactions on Automatic Control*, vol. 50, no. 5, pp. 546–558, 2005.

[23] L. Grüne and J. Pannek, *Nonlinear Model Predictive Control: Theory and Algorithms*, ser. Communications and Control Engineering, A. Isidori, J. H. van Schuppen, E. D. Sontag, M. Thoma, and M.Krstic, Eds.  Springer London Dordrecht Heidelberg New York, 2011.

[24] L. Grüne, "Analysis and design of unconstrained nonlinear MPC schemes for finite and infinite dimensional systems," *SIAM J. Control Optim.*, vol. 48, no. 2, pp. 1206–1228, 2009.

[25] K. Worthmann, "Stability Analysis of unconstrained Receding Horizon Control," Ph.D. dissertation, University of Bayreuth, 2011.

[26] ——, "Estimates on the Prediction Horizon Length in Model Predictive Control," in *Proceedings of the 20th International Symposium on Mathematical Theory of Networks and Systems, CD–ROM, MTNS2012_0112_paper.pdf*, 2012.

[27] K. Worthmann, M. W. Mehrez, M. Zanon, R. G. Gosine, G. K. I. Mann, and M. Diehl, "Regulation of Differential Drive Robots using Continuous Time MPC without Stabilizing Constraints or Costs," in *Proceedings of the 5th IFAC Conference on Nonlinear Model Predictive Control (NPMC'15), Sevilla, Spain*, 2015, pp. 129–135.

[28] S. Thrun, W. Burgard, and D. Fox, "Probabilistic robotics," 2005.

[29] J. Brinkhuis and V. Tikhomirov, *Optimization: Insights and Applications*. Princeton University Press, 2011.

[30] L. Grüne and K. Worthmann, *Distributed Decision Making and Control*, ser. Lecture Notes in Control and Information Sciences. Springer Verlag, 2012, no. 417, ch. A distributed NMPC scheme without stabilizing terminal constraints.

[31] B. Houska, H. Ferreau, and M. Diehl, "ACADO Toolkit – An Open Source Framework for Automatic Control and Dynamic Optimization," *Optimal Control Applications and Methods*, vol. 32, no. 3, pp. 298–312, 2011.

[32] L. Grüne and A. Rantzer, "On the infinite horizon performance of receding horizon controllers," *IEEE Transactions on Automatic Control*, vol. 53, no. 9, pp. 2100–2111, 2008.

# Chapter 4

# Predictive Path Following of Mobile Robots without Terminal Stabilizing Constraints

## 4.1 Abstract

This chapter considers model predictive path-following control for differentially driven (non-holonomic) mobile robots and state-space paths. In contrast to previous works, we analyze stability of model predictive path-following control without stabilizing terminal constraints or terminal costs. To this end, we verify cost controllability assumptions and compute bounds on the stabilizing horizon length. Finally, we draw upon simulations to verify our stability results.

## 4.2 Introduction

Recently, non-holonomic mobile robots have attracted considerable interest as they are increasingly used in industry, for discovery and observation purposes, and in autonomous services. Often, the differential drive model, i.e. the unicycle, is used to describe the kinematics of non-holonomic robots. In applications, different control tasks, such as set-point stabilization/regulation, trajectory tracking, and path following, arise. Set-point stabilization refers to the control task of stabilizing a given setpoint. In case the reference is time-varying, the control task is referred to as trajectory tracking. Control tasks in which a geometric reference is to be followed, while the speed to move along the reference is not given apriori, are commonly referred to as path-following problems [1, 2]. We refer to [3] for detailed overview on path-following control methods.

Different approaches have been also established in the literature to tackle path-following problems, e.g. back stepping [4] and feedback linearization [5]. The mentioned approaches have in common that the consideration of state and input constraints is in general difficult. Model predictive control (MPC) is of a particular interest in robot control as it can handle constrained nonlinear systems. In MPC, a finite-horizon optimal control problem (OCP) is solved in a receding horizon fashion, and the first part of the optimal control is applied to the plant. Several successful MPC approaches to path-following problems have been presented in the literature. An early numerical investigation for non-holonomic robots is presented in [6]. Generalizations and extensions are discussed in follow-up papers such as [2, 3, 7–9]. In all these path-following MPC schemes, stability and path convergence are enforced using additional stabilizing terminal constraints and/or terminal costs.

The present chapter discusses the stability of model predictive path-following control (MPFC) as proposed in [2]. Similar to [6], we investigate paths defined in the state space. However, we try to explicitly avoid the use of stabilizing terminal constraints and terminal penalties. Instead, we work in the framework of controllability assumptions proposed in [10], which allows to guarantee stability by appropriately choosing the prediction horizon length. Recently, these techniques have been extended to a continuous time setting, see [11, 12]. The pursued design also allows the estimation of bounds on the infinite horizon performance of the closed loop.

Here, we first reformulate the path-following problem as the set-point stabilization of an augmented system, which combines the robot model and a timing law. The augmented state is subject to a specific constraint such that close to the end of the path the robot has to be exactly on the geometric reference. We show that this particular structure of the state constraint simplifies the verification of the controllability assumption. Specifically, we verify it by designing an open-loop control maneuver such that a robot is steered from any arbitrary initial position to the end point of the reference path. This control maneuver leads to bounds on the value function, which allow to derive a stabilizing horizon length. In essence, we extend the set-point stabilization analysis of [13, 14], which embodies Chapter 3 results, to path-following problems.

The chapter is organized as follows: in Section 4.3, a brief description of the considered model predictive path-following control scheme is given. In the subsequent Section 4.4, stability results from [11] are recalled. Then, in Section 4.5, a growth bound on the value function is derived based on feasible open-loop trajectories, which is used in Section 4.6 to determine a stabilizing prediction horizon length. Finally, we draw upon numerical simulations to assess the closed-loop performance.

## 4.3 Problem Formulation

We recall the model of a differentially driven robot as well as the path-following problem and the MPFC formulation from [6].

The continuous-time kinematic model of a differentially driven robot at time $t \in \mathbb{R}_{\geq 0}$

is given by

$$\begin{pmatrix} \dot{x}(t) \\ \dot{y}(t) \\ \dot{\theta}(t) \end{pmatrix} = \dot{\mathbf{x}}(t) = \mathbf{f}(\mathbf{x}(t), \mathbf{u}(t)) = \begin{pmatrix} v(t)\cos(\theta(t)) \\ v(t)\sin(\theta(t)) \\ \omega(t) \end{pmatrix} \quad (4.1)$$

with vector field $\mathbf{f} : \mathbb{R}^3 \times \mathbb{R}^2 \to \mathbb{R}^3$ and initial condition $\mathbf{x}(0) = \mathbf{x}_0$. The state vector $\mathbf{x} = (x, y, \theta)^\top$ (m,m,rad)$^\top$ contains the robot's posture variables, i.e. the spatial components $x$, $y$, and the orientation $\theta$. The control input $\mathbf{u} = (v, \omega)^\top$ consists of the linear and the angular speeds of the robot $v$ (m/s) and $\omega$ (rad/s), respectively. The state and input constraints are given by

$$X_\varepsilon = [-\bar{x}, -\varepsilon] \times [-\bar{y}, \bar{y}] \times \mathbb{R} \subset \mathbb{R}^3,$$
$$U = [-\bar{v}, \bar{v}] \times [-\bar{\omega}, \bar{\omega}] \subset \mathbb{R}^2$$

with $\bar{x}, \bar{y}, \bar{v}, \bar{w} > 0$. Later, the design parameter $\varepsilon \in [0, \bar{x})$ is employed to construct an extended set $Z_\varepsilon$ (introduced later) taking the particular structure of the path-following problem into account.

### 4.3.1 Path-following problem

The state-space path-following problem is to steer system (4.1) along a geometric reference curve $\mathcal{P} \subset X_0$. To this end, a parametrization $p : [\bar{\lambda}, 0] \to \mathbb{R}^3$, $\bar{\lambda} < 0$, with $\bar{\lambda} = -\bar{x}$ and satisfying $p(0) = 0$, specifies $\mathcal{P}$. Here, the scalar variable $\lambda \in \mathbb{R}$ is called the path parameter. In other words, the reference is given by

$$\mathcal{P} = \left\{ p(\lambda) \in \mathbb{R}^3 : \lambda \in [\bar{\lambda}, 0] \to p(\lambda) \right\}.$$

In path-following problems, *when to be where* on $\mathcal{P}$ is not a strict requirement. Nonetheless, the path parameter $\lambda$ is time dependent with unspecified time evolution $t \to \lambda(t)$. Therefore, the control function $\mathbf{u} \in \mathcal{PC}([0,\infty), \mathbb{R}^2)$ and the timing of $\lambda : [0,\infty) \to [\bar{\lambda}, 0]$ are chosen such that the path $\mathcal{P}$ is followed as closely as possible while maintaining feasibility with respect to the state and control constraints. The state-space path-following problem is summarized as follows, see [3, Chap. 5.1]:

**Problem 4.1.** *(State-space path following)*

1. *Convergence to path: The robot state* $\mathbf{x}$ *converges to the path* $\mathcal{P}$ *such that*

$$\lim_{t \to \infty} \|\mathbf{x}(t) - p(\lambda(t))\| = 0.$$

2. *Convergence on path: The robot moves along* $\mathcal{P}$ *in the direction of increasing* $\lambda$ *such that* $\dot{\lambda}(t) \geq 0$ *holds and* $\lim_{t \to \infty} \lambda(t) = 0$.

3. *Constraint satisfaction: The state and control constraints* $X$ *and* $U$ *are satisfied for all* $t \geq 0$.

Similar to [6], we treat the path parameter $\lambda$ as a virtual state, whose time evolution $t \to \lambda(t)$ is governed by an additional (virtual) control input $g \in \mathbb{R}$. Therefore, the dynamics of $\lambda$, which is an extra degree of freedom in the controller design, is described by a differential equation referred to as *timing law*. Here, we define the timing law as a single integrator

$$\dot{\lambda}(t) = g(t), \qquad \lambda(0) = \lambda_0, \tag{4.2}$$

$\lambda_0 \in [\bar{\lambda}, 0]$. The virtual input of the timing law is assumed to be piecewise continuous and bounded, i.e. for all $t \geq 0$, $g(t) \in G := [0, \bar{g}] \subset \mathbb{R}$. Using systems (4.1) and (4.2), the

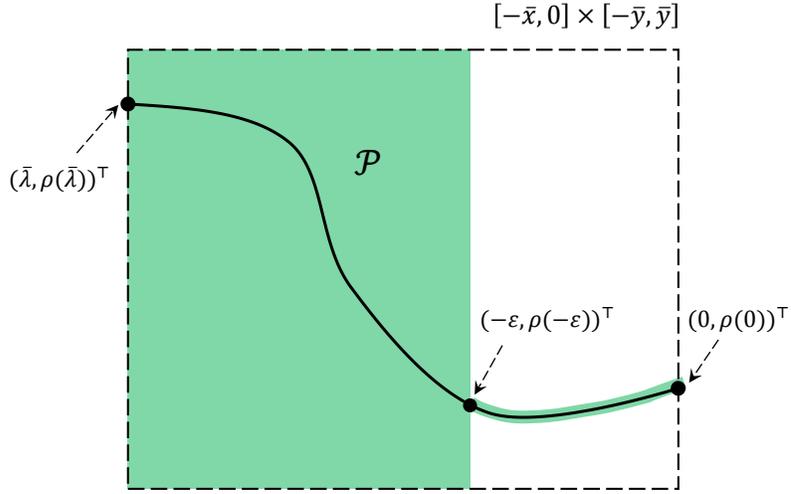

Fig. 4.1: Visualization of the spatial components (highlighted in green) of the set $Z_\varepsilon$ given by (4.4).

path-following problem is analyzed via the following augmented system

$$\dot{\mathbf{z}} = \begin{pmatrix} \dot{\mathbf{x}} \\ \dot{\lambda} \end{pmatrix} = \begin{pmatrix} \mathbf{f}(\mathbf{x}, \mathbf{u}) \\ g \end{pmatrix} = \mathbf{f}_g(\mathbf{z}, \mathbf{w}). \tag{4.3}$$

The augmented state vector $\mathbf{z} := (\mathbf{x}^\top, \lambda)^\top \in \mathbb{R}^4$ embodies the state of the robot $\mathbf{x}$ as well as the virtual state $\lambda$. For a constant $\varepsilon > 0$, the constraint set of the augmented state variable $\mathbf{z}$ is defined as

$$Z_\varepsilon = X_\varepsilon \times [\bar{\lambda}, -\varepsilon] \cup \left\{ (p(\lambda)^\top, \lambda)^\top \mid \lambda \in (-\varepsilon, 0] \right\} \subset \mathbb{R}^4. \tag{4.4}$$

In other words, the robot is forced to be on the reference path $\mathcal{P}$ if the path parameter $\lambda$ satisfies $\lambda > -\varepsilon$. This particular structure of $Z_\varepsilon$ serves as a kind of *stabilizing constraint* mainly used to simplify the later derivations, see Figure 4.1 for a visualization of the spatial components of the set $Z_\varepsilon$. We remark that while $Z_\varepsilon$ imposes a constraint on the robot motion close to the end of the path $\mathcal{P}$, it will *not* be enforced at the end of each prediction horizon, i.e. it is not a stabilizing terminal constraint in the classical sense.

Additionally, we define the vector of augmented control actions $\mathbf{w} := (\mathbf{u}^\top, g)^\top \in \mathbb{R}^3$, which contains the robot control input as well as the virtual control. The input constraint of $\mathbf{w}$ is

$$W = U \times G.$$

The path-following problem is reduced now to the point stabilization of the augmented system (4.3) [6]. This allows us to directly use the stability results of [11] and the techniques presented in [14] as will be shown in the following section.

For system (4.3), we use $\mathbf{z}(\cdot; \mathbf{z}_0, \mathbf{w})$ to denote a trajectory originating at $\mathbf{z}_0$ and driven by the input $\mathbf{w} \in \mathcal{PC}([0,T), \mathbb{R}^3)$, $T \in \mathbb{R}_{>0} \cup \{\infty\}$. Additionally, the control function $\mathbf{w}$ is called admissible on the interval $[0,T)$ if

$$\mathbf{w}(t) \in W, t \in [0,T) \quad \text{and} \quad \mathbf{z}(t; \mathbf{z}_0, \mathbf{w}) \in Z_\varepsilon, t \in [0,T]$$

hold. We denote the set of all admissible control functions $w$ for initial value $\mathbf{z}_0$ on $[0,T)$ by $\mathcal{W}_T(\mathbf{z}_0)$.

Here, we consider paths parameterized by

$$p(\lambda) = \left(\lambda, \rho(\lambda), \arctan\left(\frac{\partial \rho}{\partial \lambda}\right)\right)^\top, \tag{4.5}$$

whereby $\rho$ is at least twice continuously differentiable. A basic requirement for the path-following problem to be feasible, is the consistency of the path with the state constraints, i.e. we assume that $\mathcal{P} \subset X_0$. System (4.1) is *differentially flat*, and $(x,y)^\top$ is one of its flat outputs.[1] Hence, an input $\mathbf{u}_{ref}$ ensuring that the system follows the path (4.5) exactly for a given timing $\lambda(t)$ can be obtained using ideas from [17]. First, by investigating the first two equations of (4.1), $v_{ref}$ is obtained as

$$v_{ref}(\lambda, \dot\lambda) = \dot\lambda \cdot \sqrt{1 + \left(\frac{\partial \rho(\lambda)}{\partial \lambda}\right)^2}. \tag{4.6a}$$

---

[1] We refer to [15, 16] for more details on differential flatness.

Similarly, using the last equation in model (4.1), $\omega_{ref}$ is computed as

$$\begin{aligned}\omega_{ref}(\lambda, \dot\lambda) &= \frac{d}{dt}\left(\arctan\left(\frac{\partial \rho(\lambda)}{\partial \lambda}\right)\right) \\ &= \dot\lambda \cdot \left(1 + \left(\frac{\partial \rho(\lambda)}{\partial \lambda}\right)^2\right)^{-1}\left(\frac{\partial^2 \rho(\lambda)}{\partial \lambda^2}\right).\end{aligned} \tag{4.6b}$$

Note that due to (4.2), $\mathbf{u}_{ref}$ can also be regarded as a function of $\lambda$ and the virtual control $g$. The largest value of $g$ for which $\mathbf{u}_{ref}(\lambda, g)$ is admissible is given by

$$\hat{g} := \max\{g \in [0, \bar{g}] \,|\, \mathbf{u}_{ref}(\lambda, g) \in U \quad \forall \lambda \in [\bar\lambda, 0]\}. \tag{4.7}$$

It is readily seen that, for any twice continuously differentiable $\rho(\lambda)$ in (4.5), the timing law $\dot\lambda = g$, the input constraint $U$ and the structure of (4.6) imply $\hat{g} > 0$.

### 4.3.2 Model predictive path following control (MPFC)

Here, we recall the MPFC scheme for state-space paths as proposed in [6]. We refer to [2, 3] for extension to paths defined in output spaces. In MPFC Problem 4.1 is solved via a continuous-time sampled-data MPC.

We consider continuous running (stage) costs $\ell : Z_0 \times W \to \mathbb{R}_{\geq 0}$ satisfying

$$\ell(\mathbf{z}^r, 0) = 0 \quad \text{and} \quad \min_{\mathbf{w} \in W} \ell(\mathbf{z}, \mathbf{w}) > 0 \quad \forall\, \mathbf{z} \in \mathbb{R}^4 \setminus \mathbf{z}^r,$$

where $\mathbf{z}^r$ is the path end point we would like to stabilize dynamics (4.3) at, i.e. $\mathbf{z}^r$ is given by

$$\mathbf{z}^r := \left(0, \rho(0), \arctan\left(\frac{\partial \rho(\lambda)}{\partial \lambda}\right)|_{\lambda=0}, 0\right)^\top. \tag{4.8}$$

Similar to [6], $\ell$ is chosen as

$$\ell(\mathbf{z}, \mathbf{w}) = \left\|\begin{matrix} \mathbf{x} - p(\lambda) \\ \lambda \end{matrix}\right\|_{\mathbf{Q}}^2 + \left\|\begin{matrix} \mathbf{u} - \mathbf{u}_{ref}(\lambda, g) \\ g \end{matrix}\right\|_{\mathbf{R}}^2. \tag{4.9}$$

$\mathbf{Q} = \operatorname{diag}(q_1, q_2, q_3, \hat{q}), \mathbf{R} = \operatorname{diag}(r_1, r_2, \hat{r})$ are positive definite weighting matrices and $\mathbf{u}_{ref}(\lambda, g)$ is from (4.6). The objective functional to be minimized in the MPFC scheme reads

$$J_T(\mathbf{z}_k, \mathbf{w}) := \int_{t_k}^{t_k+T} \ell(\mathbf{z}(\tau; \mathbf{z}_k, \mathbf{w}), \mathbf{w}(\tau))\, \mathrm{d}\tau$$

with the prediction horizon $T \in \mathbb{R}_{>0}$. Hence, the MPFC scheme is based on repeatedly solving the following optimal control problem (OCP)[2]

$$V_T(\mathbf{z}_k) = \min_{\mathbf{w} \in \mathcal{PC}([t_k, t_k+T), \mathbb{R}^3)} J_T(\mathbf{z}_k, \mathbf{w}) \tag{4.10}$$

$$\text{subject to } \mathbf{z}(t_k) = \mathbf{z}_k,$$

$$\dot{\mathbf{z}}(\tau) = \mathbf{f}_g(\mathbf{z}(\tau), \mathbf{w}(\tau)) \quad \forall \tau \in [t_k, t_k + T]$$

$$\mathbf{z}(\tau) \in Z_\varepsilon \qquad \forall \tau \in [t_k, t_k + T]$$

$$\mathbf{w}(\tau) \in W \qquad \forall \tau \in [t_k, t_k + T)$$

For a given time $t_k$ and initial value $\mathbf{z}(t_k) = \mathbf{z}_k \in Z_\varepsilon$, the minimum of this OCP, i.e.

---

[2]To avoid cumbersome technicalities, we assume that the minimum exists and is attained.

$V_T(\mathbf{z}_k)$, does actually not depend on $t_k$. Moreover, $V_T : Z_0 \to \mathbb{R}_{\geq 0} \cup \{\infty\}$ is the optimal value function.

The solution of (4.10) results in the optimal control function $\mathbf{w}^\star \in \mathcal{W}_T(\mathbf{z}_0)$. Then, for a sampling period $\delta \in (0, T)$, the MPFC feedback applied to the robot is given by

$$\mathbf{u}(t) = (v^\star(t, \mathbf{z}(t_k)), \omega^\star(t, \mathbf{z}(t_k)))^\top, \quad t \in [t_k, t_k + \delta). \tag{4.11}$$

Note that the initial condition $\mathbf{z}_k$ in the optimization is composed of the robot state $\mathbf{x}(t_k)$ and the path parameter $\lambda(t_k)$. Similar to [6], the initial condition is $\lambda(t_k) = \lambda(t_k; \lambda_{k-1}, g^\star)$, i.e. the corresponding value of the last predicted trajectory. Note that $Z_\varepsilon$ from (4.4) requires the robot to follow the final part of the path exactly. However, observe that OCP (4.10) does not involve any terminal constraint.

## 4.4 Stability and performance bounds

We first recall an MPC stability result presented in [11] in a slighly reformulated version, see [14]. Then, we show that, if this stability result is verified for the MPFC scheme, for the considered initial conditions, Problem 4.1 is solved.

We consider the nominal case (no plant-model mismatch). For a specific choice of $\delta$ and $T$, let the subscript $(\cdot)_{T,\delta}$ denote the MPC closed-loop variables (states and inputs) of the augmented system (4.3), where $\mathbf{w}^\star$ is applied in the sense of (4.11).

**Theorem 4.2.** *Assume existence of a monotonically increasing and bounded function $B : \mathbb{R}_{\geq 0} \to \mathbb{R}_{\geq 0}$ satisfying, for all $\mathbf{z}_0 \in Z_\varepsilon$,*

$$V_t(\mathbf{z}_0) \leq B(t) \cdot \ell^\star(\mathbf{z}_0) \quad \forall\, t \geq 0, \tag{4.12}$$

*where $\ell^\star(\mathbf{z}_0) := \min_{\mathbf{w} \in W} \ell(\mathbf{z}_0, \mathbf{w})$. Let the sampling period $\delta > 0$ and the prediction horizon*

$T > \delta$ be chosen such that the condition $\alpha_{T,\delta} > 0$ holds for

$$\alpha_{T,\delta} = 1 - \frac{e^{-\int_\delta^T B(t)^{-1}\,dt} \cdot e^{-\int_{T-\delta}^T B(t)^{-1}\,dt}}{\left[1 - e^{-\int_\delta^T B(t)^{-1}\,dt}\right]\left[1 - e^{-\int_{T-\delta}^T B(t)^{-1}\,dt}\right]}.$$

*Then, for all $\mathbf{z} \in Z_\varepsilon$, the relaxed Lyapunov inequality*

$$V_T(\mathbf{z}_{T,\delta}(\delta;\mathbf{z})) \leq V_T(\mathbf{z}) - \alpha_{T,\delta} \int_0^\delta \ell(\mathbf{z}_{T,\delta}(t;\mathbf{z}), \mathbf{w}_{T,\delta}(t,\mathbf{z}))\,dt$$

*as well as the performance estimate*

$$V_\infty^{\mathbf{w}_{T,\delta}}(\mathbf{z}) \leq \alpha_{T,\delta}^{-1} \cdot V_\infty(\mathbf{z}) \tag{4.13}$$

*are satisfied, whereby $V_\infty^{\mathbf{w}_{T,\delta}(\mathbf{z})}$ is the MPC closed-loop costs*

$$V_\infty^{\mathbf{w}_{T,\delta}}(\mathbf{z}) := \int_0^\infty \ell(\mathbf{z}_{T,\delta}(t;\mathbf{z}), \mathbf{w}_{T,\delta}^{MPC}(t;\mathbf{z}))\,dt,$$

*and $\mathbf{w}_{T,\delta}^{MPC}(t;\mathbf{z})$ is the closed-loop control. If, in addition, there exists $\mathcal{K}_\infty$-functions $\underline{\eta}, \bar{\eta} : \mathbb{R}_{\geq 0} \to \mathbb{R}_{\geq 0}$ satisfying*

$$\underline{\eta}(\|\mathbf{z} - \mathbf{z}^r\|) \leq \ell^\star(\mathbf{z}) \leq \bar{\eta}(\|\mathbf{z} - \mathbf{z}^r\|) \qquad \forall\, \mathbf{z} \in Z_\varepsilon, \tag{4.14}$$

*the MPC closed loop is asymptotically stable.* □

While condition (4.14) holds trivially for the chosen running costs (4.9), it is crucial to verify the growth condition (4.12).

The next result provides the connection between Theorem 4.2 and Problem 4.1.

**Proposition 4.3.** *Let $\delta, T$ in MPFC scheme based on OCP (4.10) be chosen such that the conditions of Theorem 4.2 hold. Then, the MPFC feedback (4.11) solves Problem 4.1.*

**Proof.** Recursive feasibility of the optimization follows from the observation that the state constraint set $Z_\varepsilon$ can be rendered forward invariant[3] by means of inputs satisfying the input constraint $W$. Furthermore, the chosen running costs (4.9) implies that, as $\mathbf{z}$ goes to $\mathbf{z}^r$, the robot state converges to the path and $\lambda$ converges to 0. Hence, stabilizing $\mathbf{z}$ at $\mathbf{z}^r$ in an admissible way implies solving Problem 4.1. ∎

In the next section, we show a method to construct a function $B : \mathbb{R}_{\geq 0} \to \mathbb{R}_{\geq 0}$ such that all assumptions of Theorem 4.2 hold. In our construction, we exploit the structure of the state constraint $Z_\varepsilon$ from (4.4), i.e. we rely on the restriction of the robot motion close to the final equilibrium point. We remark that the simplification induced by the structure of $Z_\varepsilon$ is inline with observations in [19], wherein it is shown that the local analysis of growth bounds can be difficult for initial conditions in the neighborhood of the reference.

## 4.5 Cost Controllability Analysis of MPFC

Here, we present a method to compute a stabilizing horizon length for model predictive path following control (MPFC). This is achieved by deriving the growth bound (4.12) of Theorem 4.2 for the MPFC scheme of Section 4.3.2. The following theorem summarizes the main result of this section.

**Theorem 4.4.** *Let $\varepsilon > 0$ be given and let the relation $r_2 \leq \frac{q_3}{2}$ hold[4]. Then, for all $\mathbf{z}_0 \in Z_\varepsilon$, the function $B \in \mathcal{PC}([0, \infty), \mathbb{R}_{\geq 0})$, $t \to \min\{t, \int_0^t c(s) \, ds\}$ with the function $c(t)$*

---

[3] We say that a set $\tilde{Z} \subseteq Z_\varepsilon$ is *forward invariant* for system (4.3) if $\mathbf{z}(\delta; \mathbf{z}, \mathbf{w}) \in \tilde{Z}$ holds for all $\mathbf{z} \in \tilde{Z}$, see, e.g. [18].

[4] The relation $r_2 \leq \frac{q_3}{2}$ is only imposed to keep the presentation technically simple.

*defined by*

$$c(t) = \begin{cases} 1 + \frac{q_3 \pi^2 (t^2 + 6 t_r t + 0.5)}{4 t_r^2 \cdot \hat{q} \varepsilon^2} & t \in [0, t_r) \\ 1 + \frac{q_3 (3\pi/2)^2 + r_1 \bar{u}_1^2}{\hat{q} \varepsilon^2} & t \in [t_r, t_r + t_l) \\ 1 + \frac{q_3 \pi^2 ((4 t_r + t_l - t)^2 + 0.5)}{4 t_r^2 \cdot \hat{q} \varepsilon^2} & t \in [t_r + t_l, 2 t_r + t_l) \\ \frac{(4 t_r + t_l + t_g - t)^2 + \hat{r}\hat{q}^{-1}}{t_g^2} & t \in [2 t_r + t_l, 2 t_r + t_l + t_g) \\ 0 & \textit{otherwise} \end{cases} \quad (4.15)$$

*whereby*

$$t_r = \frac{\pi/2}{\bar{\omega}}, \quad t_l = \frac{\sqrt{4 \bar{y}^2 + (\bar{\lambda} + \varepsilon)^2}}{\bar{v}}, \quad t_g = \frac{|\bar{\lambda}|}{\hat{g}}$$

*with $\hat{g}$ given by (4.7), satisfies the growth bound (4.12) of Theorem 4.2 for $\ell$ from (4.9).*

**Proof.** Firstly, we construct a bounded function $\tilde{B}$, which already satisfies Inequality (4.12). Moreover, since the constant input function $\mathbf{w} \equiv 0_{\mathbb{R}^3}$ is admissible on the infinite time horizon, Inequality (4.12) is also satisfied for the identity $\hat{B}(t) = t$. In summary, the function $B$ defined in Theorem 4.4 yields Inequality (4.12) for all $\mathbf{z}_0 \in Z_\varepsilon$.

We first consider initial conditions

$$\mathbf{z}_0 = (x_0, y_0, \theta_0, \lambda_0)^\top \in X_\varepsilon \cup [\bar{\lambda}, -\varepsilon] \subsetneq Z_\varepsilon. \quad (4.16)$$

We derive an admissible (open-loop) control function $\mathbf{w}_{\mathbf{z}_0} \in \mathcal{W}_\infty(\mathbf{z}_0)$ steering the robot, in finite time, to the end point of the considered path, i.e. to $p(0) = 0$. This function $\mathbf{w}_{\mathbf{z}_0}$ yields (suboptimal) running costs $\ell(\mathbf{z}(t; \mathbf{z}_0, \mathbf{w}_{\mathbf{z}_0}), \mathbf{w}_{\mathbf{z}_0}(t))$. Then, we uniformly estimate the quotient

$$\frac{\ell(\mathbf{z}(t; \mathbf{z}_0, \mathbf{w}_{\mathbf{z}_0}), \mathbf{w}_{\mathbf{z}_0}(t))}{\ell^\star(\mathbf{z}_0)} \leq c(t) \qquad \forall\, \mathbf{z} \in X_\varepsilon \cup [\bar{\lambda}, -\varepsilon] \quad (4.17)$$

for all $t \geq 0$ with $c \in \mathcal{PC}(\mathbb{R}_{>0}, \mathbb{R}_{\geq 0})$. Then, since the maneuver is carried out in finite time, there exists a $\bar{t} > 0$ such that $c(t) = 0$ for all $t \geq \bar{t}$, and a bounded growth function $\tilde{B}$ given by

$$\tilde{B}(t) := \int_0^t c(\tau) \, \mathrm{d}\tau.$$

As the robot can move forward and backward, it suffices to consider the case $y_0 < \rho(\lambda_0)$. Moreover, note that $\ell^\star(\mathbf{z}_0)$ is given by

$$\ell^\star(\mathbf{z}_0) = q_1(x_0 - \lambda_0)^2 + q_2(y_0 - \rho(\lambda_0))^2 + q_3(\theta_0 - \hat{\theta})^2 + \hat{q}\lambda_0^2,$$

with angle $\hat{\theta} = \arctan\left(\frac{\partial \rho}{\partial \lambda}(\lambda_0)\right) \in (-\pi/2, \pi/2)$, i.e. $\hat{\theta}$ represents the desired orientation along the path $\mathcal{P}$ at $\lambda = \lambda_0$. We emphasize that using the virtual input $g = 0$ implies that $\mathbf{u}_{ref}(\lambda, g) = (0, 0)^\top$. It's important to keep this in mind while reading the following construction of the function (4.15) satisfying (4.17). All initial conditions (4.16) satisfy $\lambda_0 \leq -\varepsilon$. Hence, $\ell^\star(\mathbf{z}_0)$ is uniformly bounded from below by $\hat{q}\varepsilon^2$.

The blueprint for the construction of the control function $\mathbf{w}_{\mathbf{z}_0}$ can be summarized as follows, see also Figure 4.2:

I. Define the angle

$$\phi_{\mathbf{z}_0} = \operatorname{atan2}\left(\frac{\rho(\lambda_0) - y_0}{\lambda_0 - x_0}\right) \in [-\pi, \pi)$$

and turn the robot such that the condition $\theta(t) = \phi_{\mathbf{z}_0}$ or $\theta(t) = \phi_{\mathbf{z}_0} - \pi$ is satisfied. At the end of this step, the robot points towards (or in the opposite direction to) the path point $(\lambda_0, \rho(\lambda_0))^\top$.

II. Drive directly (forward or backwards) until the robot reaches the path at $(\lambda_0, \rho(\lambda_0))^\top$.

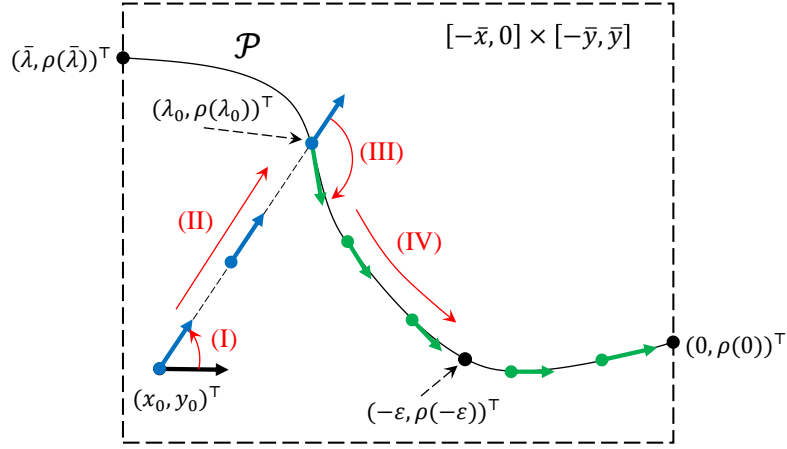

Fig. 4.2: Sketch of the proposed maneuver for initial conditions characterized by (4.16).

III. Turn the robot until its orientation becomes tangent to the path at the selected point, i.e. until the angle $\hat{\theta}$ is reached.

IV. Drive the robot along the path until the end of the path, i.e. until $\lambda = 0$.

The time needed to complete this maneuver depends on the constraints $X_0$ and $W$: the minimal time required to turn the vehicle by 90 degrees is $t_r$, the minimal time required to drive the vehicle between the two corners of the box $[\bar{\lambda}, -\varepsilon] \times [-\bar{y}, \bar{y}]$ is $t_l$, and the minimum time to drive the path variable $\lambda$, with a constant control $g$, from its limit $\bar{\lambda}$ to $0$ is $t_g$. The proposed maneuver actions I–IV are as follows.

<u>Part I.</u> First, the vehicle turns until time $t_r$ such that $\theta(t_r) = \phi_{\mathbf{z}_0}$ or $\theta(t_r) = \phi_{\mathbf{z}_0} - \pi$. This is achieved by applying the constant input $\mathbf{w}_{\mathbf{z}_0}(t) \in W$, $t \in [0, t_r)$, given by

$$\mathbf{w}_{\mathbf{z}_0} \equiv \begin{cases} \pm(0, (\phi_{\mathbf{z}_0} - \theta_0)t_1^{-1}, 0)^\top & \text{if } |\phi_{\mathbf{z}_0} - \theta_0| \leq \frac{\pi}{2}, \\ \pm(0, (\phi_{\mathbf{z}_0} - \pi - \theta_0)t_1^{-1}, 0)^\top & \text{otherwise,} \end{cases}$$

where the control sign is adjusted such that the robot rotation is achieved in the proper direction, i.e. clock wise or counter clock wise. This yields, for $\ell(\mathbf{z}(t; \mathbf{z}_0, \mathbf{w}_{\mathbf{z}_0}), \mathbf{w}_{\mathbf{z}_0}(t))$ and

$t \in [0, t_r)$, the expression[5]

$$\ell^\star(\mathbf{z}_0) + q_3 \left[ \left( \frac{(t_r - t)\theta_0 + t\phi_{\mathbf{z}_0}}{t_r} - \hat{\theta} \right)^2 - (\theta_0 - \hat{\theta})^2 \right] + r_2 \left( \frac{\phi_{\mathbf{z}_0} - \theta_0}{t_r} \right)^2.$$

When expanding the square bracket, the terms $(t/t_r)^2(\phi_{\mathbf{z}_0} - \theta_0)^2$ and $(2t/t_r)(\phi_{\mathbf{z}_0} - \theta_0)(\theta_0 - \hat{\theta})$ appear. Since $|\phi_{\mathbf{z}_0} - \theta_0| \leq \pi/2$, the term $(\phi_{\mathbf{z}_0} - \theta_0)^2$ is upper bounded by $(\pi/2)^2$. Moreover, it is straightforward to show that the product $(\phi_{\mathbf{z}_0} - \theta_0)(\theta_0 - \hat{\theta})$ is bound from above by $(\pi/2 \cdot 3\pi/2)$. Using these upper bounds along with our estimate for $\ell^\star(\mathbf{z}_0)$ and the assumption on $r_2$, Inequality (4.17) holds for $t \in [0, t_r)$ with $c$ defined by (4.15).

<u>Part II.</u> Next, the vehicle drives to the point $(\lambda_0, \rho(\lambda_0))^\top$ until time $t_r + t_l$ with a constant control $\mathbf{w}_{\mathbf{z}_0}(t)$, $t \in [t_r, t_r + t_l)$, defined as

$$\mathbf{w}_{\mathbf{z}_0} \equiv \begin{cases} \left( Dt_l^{-1}, 0, 0 \right)^\top & \text{if } |\phi_{\mathbf{z}_0} - \theta_0| \leq \frac{\pi}{2} \\ -\left( Dt_l^{-1}, 0, 0 \right)^\top & \text{otherwise} \end{cases}$$

with $D = \sqrt{(x_0 - \lambda_0)^2 + (y_0 - \rho(\lambda_0))^2}$. This yields, for $t \in [t_r, t_r + t_l)$,

$$\ell(\mathbf{z}(t; \mathbf{z}_0, \mathbf{w}_{\mathbf{z}_0}), \mathbf{w}_{\mathbf{z}_0}(t)) =$$
$$\left( \frac{t_r + t_l - t}{t_l} \right)^2 \left[ q_1(x_0 - \lambda_0)^2 + q_2(y_0 - \rho(\lambda_0))^2 \right] + q_3(\phi_{\mathbf{z}_0} - \hat{\theta})^2 + \hat{q}\lambda_0^2 + \frac{r_1}{t_l^2} D^2.$$

Here, we have $|\phi_{\mathbf{z}_0} - \hat{\theta}| \leq 3\pi/2$. Moreover, the magnitude of the control effort is upper bounded by $\bar{v}$. Thus, using the assumption on $r_2$ and $\ell^\star(\mathbf{z}_0) \geq \hat{q}\varepsilon^2$, ensures Inequality (4.17) with $c(t)$ defined by (4.15).

<u>Part III.</u> The vehicle turns until its orientation becomes $\hat{\theta}$ at time $4t_r + t_l$, i.e. a turn

---

[5] Without loss of generality, the term $(\phi_{\mathbf{z}_0} - \theta_0)$ can be replaced by $(\phi_{\mathbf{z}_0} - \pi - \theta_0)$ in the considered $\ell(\cdot, \cdot)$ because the norm of either term is upper bounded by $\pi/2$.

of a maximum 270 degrees. This is achieved using the (constant) control

$$\mathbf{w}_{\mathbf{z}_0} = \left(0, -\frac{(\phi_{\mathbf{z}_0} - \hat{\theta})}{3t_r}, 0\right)^\top, \quad t \in [t_r + t_l, 4t_r + t_l).$$

Thus, we have

$$\ell(\mathbf{z}(t; \mathbf{z}_0, \mathbf{w}_{\mathbf{z}_0}), \mathbf{w}_{\mathbf{z}_0}(t)) = q_3 \left[(\phi_{z_0} - \hat{x}_3)\left(\frac{4t_r + t_l - t}{3t_r}\right)\right]^2 + \hat{q}\theta_0^2 + r_2 \left[\frac{-(\phi_{z_0} - \hat{x}_3)}{3t_r}\right]^2.$$

Similar to the considerations above, we use $|\phi_{\mathbf{z}_0} - \hat{\theta}| \leq 3\pi/2$. Thus, using the assumption on $r_2$ and $\ell^\star(\mathbf{z}_0) \geq \hat{q}\varepsilon^2$, ensures (4.17) with $c(t)$ defined by (4.15).

<u>Part IV.</u> Finally, the robot drives along the path until it reaches its end point, i.e. $\mathbf{z}^r$, at time $4t_r + t_l + t_g$. This is achieved using the control

$$\mathbf{w}_{\mathbf{z}_0}(t) = \left(v_{ref}, \omega_{ref}, \frac{-\lambda_0}{t_g}\right)^\top, \quad t \in [4t_r + t_l, 4t_r + t_l + t_g),$$

where $\mathbf{u}_{ref} = (v_{ref}, \omega_{ref})^\top$, has the argument $(\mathbf{z}_4(t; \mathbf{z}_0, \mathbf{w}_{\mathbf{z}_0}), -\lambda_0 t_g^{-1})$ at time $t$. This results in running costs

$$\ell(\mathbf{z}(t; \mathbf{z}_0, \mathbf{w}_{\mathbf{z}_0}), \mathbf{w}_{\mathbf{z}_0}(t)) = \hat{q}\left(\frac{\lambda_0(2t_r + t_l + t_v - t)}{t_g}\right)^2 + \hat{r}\left(\frac{\lambda_0}{t_g}\right)^2.$$

Then, using the lower bound on $\ell^\star(\mathbf{z}_0)$, Inequality (4.17) is ensured with $c(t)$ given by (4.15).

Moreover, for $t \geq 2t_r + t_l + t_g$, the function $c(t)$ is defined as $c(t) = 0$. Note that the function $c$ is independent of the particular initial condition $\mathbf{z}_0$. Additionally, the proposed maneuver ensures the satisfaction of the state constraints, i.e. $\mathbf{z}(t; \mathbf{z}_0, \mathbf{w}_{\mathbf{z}_0}) \in Z_\varepsilon$ and the selection of the times $t_l$, $t_r$, and $t_g$ ensures satisfaction of the control constraints. Clearly, $c$ is a piecewise continuous function with at most four points of discontinuity. Moreover, since $c(t) = 0$ for all $t \geq 2t_r + t_l + t_g$ it is bounded and integrable on $t \in [0, \infty)$. Therefore, the growth bound $B$ is well defined.

Additionally, for initial conditions $\mathbf{z}_0 \in Z_0$ with $\lambda_0 > -\varepsilon$, the proposed maneuver corresponds to waiting ($\mathbf{w}_{\mathbf{z}_0} \equiv 0$ on $[0, 4t_r + t_l)$) before the mobile robots travels along the final segment of the path. Again, the constructed function $c$ satisfies Inequality (4.17). In conclusion, the deduced function $B$ satisfies Condition (4.12) for all $\mathbf{z}_0 \in Z_0$. ∎

The idea of using an integrable function, which is zero after a finite time interval goes back to [20] and is presented in a continuous time setting in [11]. We emphasize that, for given $\delta > 0$, the existence of a prediction horizon $T > 0$ such that the stability condition $\alpha_{T,\delta} > 0$ holds, is shown in [10].

## 4.6 Numerical Results

Next, using the bound $B$ derived in Section 4.5, we first determine a prediction horizon length $\hat{T}$ such that the MPFC closed loop is asymptotically stable. $\hat{T}$ is then employed for closed-loop simulations.

The derived growth function $B$ of the previous section is employed to determine $\hat{T}$ which is defined as

$$\hat{T} = \min \{T > 0 \,|\, \exists N \in \mathbb{N} : T = N\delta \text{ and } \alpha_{T,\delta} > 0\}.$$

The horizon $\hat{T}$ depends on the constraints $Z_0$ and $W$, the sampling period $\delta$, the weights of the running costs (4.9), the parameter $\varepsilon$, and the shape of the path to be followed. The path to be followed is given by

$$\rho(\lambda) = 0.6 \sin(0.25 \cdot \lambda).$$

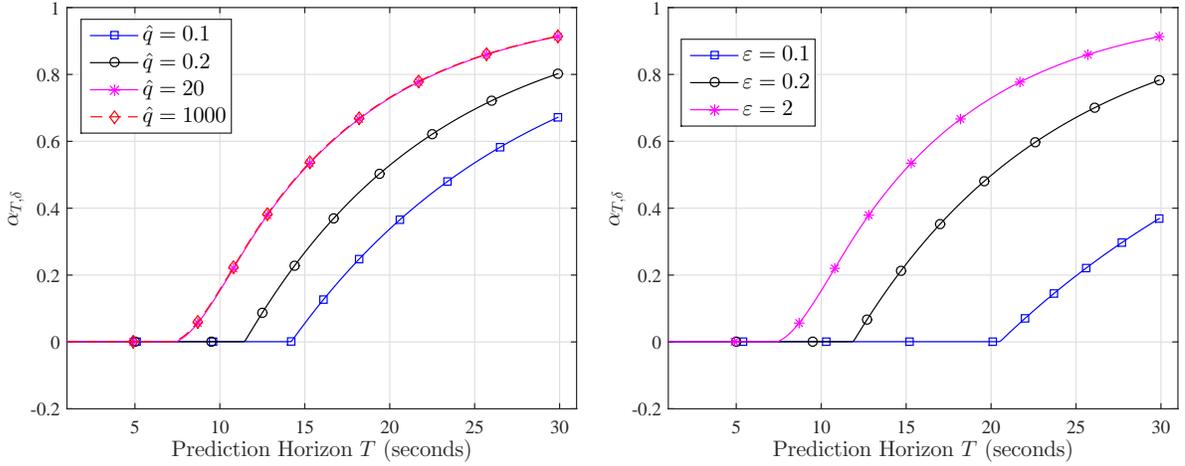

Fig. 4.3: Left: Effect of changing $\hat{q}$ on $\alpha_{T,\delta}$ for $\varepsilon = 2$. Right: Effect of changing $\varepsilon$ on $\alpha_{T,\delta}$ for $\hat{q} = 20$.

The sampling rate $\delta$ is set as $\delta = 0.1$ (seconds). The constraints are

$$X_0 = [-20, 0] \times [-1, 1] \times \mathbb{R}, \quad [\bar{\lambda}, 0] = [-20, 0],$$
$$W = [-4, 4] \times [-\pi/2, \pi/2] \times [0, 4].$$

Additionally, we choose the parameters of $\ell$ as $q_1 = q_2 = 10^4$, $q_3 = 0.01$, $r_1 = 0.1$, $r_2 = q_3/2$, and $\hat{r} = 0.1$ while $\hat{q}$ is investigated over the set $\{0.1, 0.2, 20\}$. For $\varepsilon = 2$, Figure 4.3 (left) shows the effect of changing $\hat{q}$ on $\alpha_{T,\delta}$. As one can see, increasing $\hat{q}$ reduces the stabilizing prediction horizon length. However, setting $\hat{q}$ above 20 does not make a noticeable improvement to the stabilizing prediction horizon length. For $\hat{q} = 20$, Figure 4.3 (right) shows the effect of changing $\varepsilon$. As shown, increasing $\varepsilon$ reduces the stabilizing prediction horizon length.

Next, we use $\hat{T} = 7.5$ (seconds) for MPFC closed-loop simulations. Here, we investigate the MPFC scheme performance by considering 6 initial conditions of the robot.

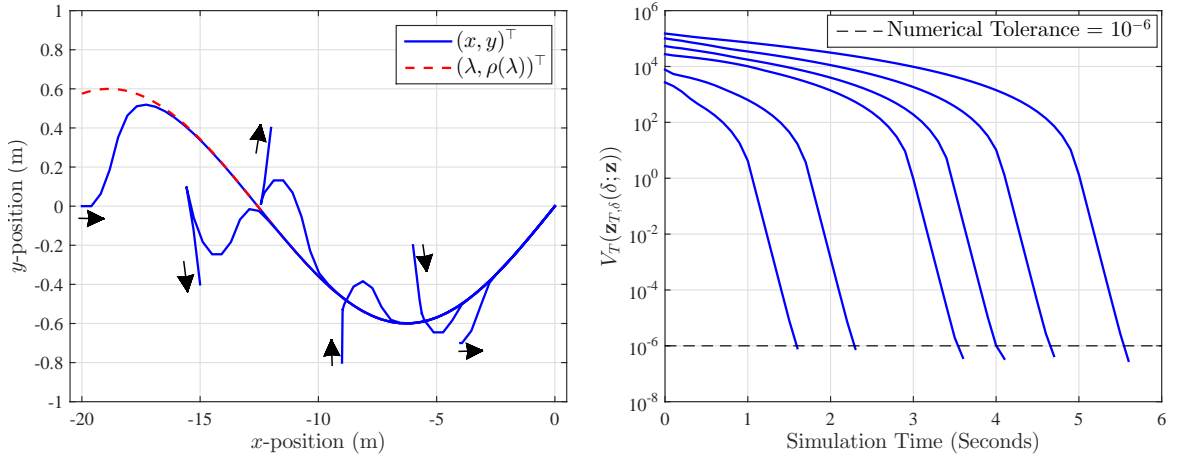

Fig. 4.4: Left: Robot's closed-loop trajectories for the considered initial conditions. Right: Value function $V_T$ along closed-loop trajectories. In all simulations, we set $\hat{q} = 20$, $\varepsilon = 2$, and $T = 7.5$.

Following [6], $\lambda_0$ is set as

$$\lambda_0 = \underset{\lambda \in [\bar{\lambda}, -\varepsilon]}{\mathrm{argmin}} \quad \|\mathbf{x}_0 - p(\lambda)\|,$$

i.e. $\lambda_0$ is chosen such that the distance $\|\mathbf{x}_0 - p(\lambda_0)\|$ is minimized. All simulations have been run utilizing the interior-point optimization method provided by the IPOPT package, see [21], coupled with MATLAB via the CasADi toolbox, see [22]. Closed-loop trajectories were considered until the condition $V_T(\mathbf{z}_{T,\delta}(\delta; \mathbf{z})) \leq 10^{-6}$ is met. Figure 4.4 (left) shows the closed-loop trajectories exhibited by the robot for the considered initial conditions. Moreover, for 2 initial conditions out of 6, Figure 4.5 shows the closed loop time evolution of all the elements of the state $\mathbf{z}$ and all the elements of the control vector $\mathbf{w}$.

In all cases, MPFC steered successfully the robot to the end point on the considered path, while satisfying all the conditions presented in Problem 4.1. Moreover, as can be noticed from Figure 4.4 (right), the evolution of the value function $V_T$ exhibited a strictly monotonically decreasing behavior for the all considered initial conditions. This demonstrates that the relaxed Lyapunov inequality presented in Theorem 4.2 is satisfied

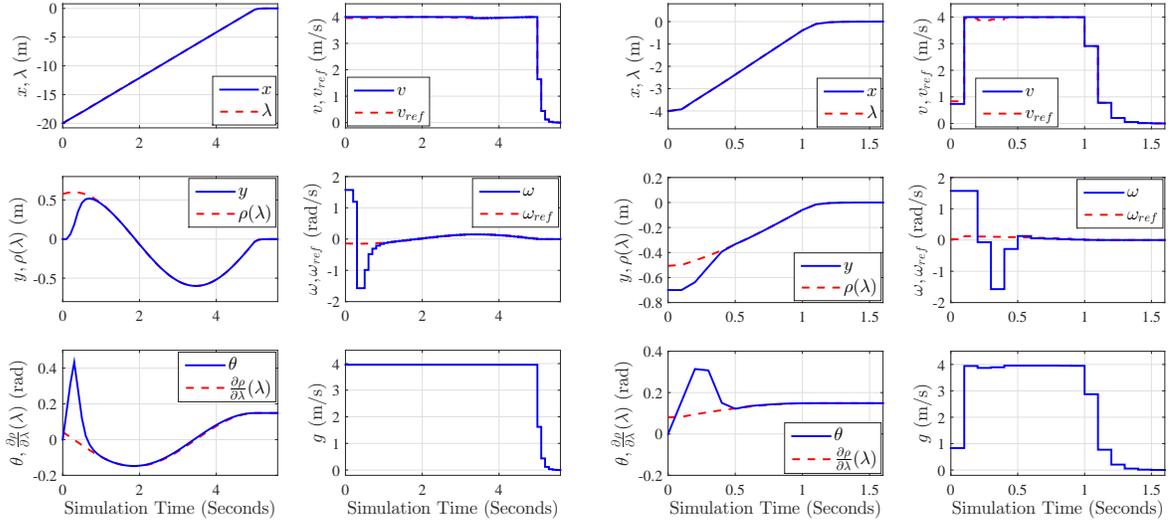

Fig. 4.5: Closed loop time evolution of the elements of the state $\mathbf{z}$ and control $\mathbf{w}$ for the initial conditions $\mathbf{x}_0 = (-20, 0, 0)^\top$ (Left) and $\mathbf{x}_0 = (-4, -0.7, 0)^\top$ (Right). In all simulations, we set $\hat{q} = 20$, $\varepsilon = 2$, and $T = 7.5$.

and, thus, the closed loop stability is verified under the chosen prediction horizon length.

## 4.7 Conclusions

In this chapter, asymptotic stability of model predictive path-following control (MPFC) for non-holonomic mobile robots and state-space paths has been rigorously proven without stabilizing terminal constraints or costs. We have demonstrated how the main assumptions, i.e. the growth bound on the value function, of the stability result can be derived. We also studied the effect of the tuning parameters of MPFC on the stabilizing horizon length.

# References

[1] A. P. Aguiar, J. P. Hespanha, and P. V. Kokotović, "Performance limitations in reference tracking and path following for nonlinear systems," *Automatica*, vol. 44, no. 3, pp. 598 – 610, 2008.

[2] T. Faulwasser and R. Findeisen, "Nonlinear model predictive control for constrained output path following," *IEEE Trans. Automat. Control*, vol. 61, no. 4, pp. 1026–1039, April 2016.

[3] T. Faulwasser, *Optimization-based Solutions to Constrained Trajectory-tracking and Path-following Problems.*  Shaker, Aachen, Germany, 2013.

[4] A. P. Aguiar, J. P. Hespanha, and P. V. Kokotović, "Performance limitations in reference tracking and path following for nonlinear systems," *Automatica*, vol. 44, no. 3, pp. 598 – 610, 2008.

[5] C. Nielsen, C. Fulford, and M. Maggiore, "Path following using transverse feedback linearization: Application to a maglev positioning system," *Automatica*, vol. 46, no. 3, pp. 585–590, Mar. 2010.

[6] T. Faulwasser and R. Findeisen, "Nonlinear model predictive path-following control," in *Nonlinear Model Predictive Control*, ser. Lecture Notes in Control and Inform. Sciences, L. Magni, D. Raimondo, and F. Allgöwer, Eds.  Springer Berlin Heidelberg, 2009, vol. 384, pp. 335–343.

[7] A. Alessandretti, A. P. Aguiar, and C. Jones, "Trajectory-tracking and Path-following Controllers for Constrained Underactuated Vehicles using Model Predictive Control," in *Proc. European Control Conf.*, Zürich, Switzerland, 2013, pp. 1371–1376.

[8] D. Lam, C. Manzie, and M. C. Good, "Multi-axis model predictive contouring control," *Int. J. Control*, vol. 86, no. 8, pp. 1410–1424, 2013.

[9] S. Yu, X. Li, H. Chen, and F. Allgöwer, "Nonlinear model predictive control for path following problems," *Int. J. Robust and Nonlinear Control*, vol. 25, no. 8, pp. 1168–1182, 2015.

[10] L. Grüne, J. Pannek, M. Seehafer, and K. Worthmann, "Analysis of unconstrained nonlinear MPC schemes with varying control horizon," *SIAM J. Control and Optimization*, vol. 48, no. 8, pp. 4938–4962, 2010.

[11] M. Reble and F. Allgöwer, "Unconstrained model predictive control and suboptimality estimates for nonlinear continuous-time systems," *Automatica*, vol. 48, no. 8, pp. 1812–1817, 2012.

[12] K. Worthmann, M. Reble, L. Grüne, and F. Allgöwer, "The role of sampling for stability and performance in unconstrained nonlinear model predictive control," *SIAM Journal on Control and Optimization*, vol. 52, no. 1, pp. 581–605, 2014.

[13] K. Worthmann, M. W. Mehrez, M. Zanon, G. K. I. Mann, R. G. Gosine, and M. Diehl, "Model predictive control of nonholonomic mobile robots without stabilizing constraints and costs," *IEEE Trans. Control Syst. Technol.*, vol. 24, no. 4, pp. 1394–1406, July 2016.

[14] ——, "Regulation of Differential Drive Robots using Continuous Time MPC without Stabilizing Constraints or Costs," *IFAC-PapersOnLine: 5th IFAC Conf. Nonlinear Model Predictive Control, Seville, Spain*, vol. 48, no. 23, pp. 129–135, 2015.

[15] J. Lévine, *Analysis and control of nonlinear systems: a flatness-based approach*, ser. Math. Eng.   Springer, Berlin, 2009.

[16] P. Martin, R. M. Murray, and P. Rouchon, "Flat systems, equivalence and trajectory generation," 2003.

[17] T. Faulwasser, V. Hagenmeyer, and R. Findeisen, "Optimal exact path-following for constrained differentially flat systems," in *Proc. 18th IFAC World Congress, Milano, Italy*, 2011, pp. 9875–9880.

[18] L. Grüne and J. Pannek, *Nonlinear Model Predictive Control: Theory and Algorithms*, ser. Communications and Control Engineering.  Springer London Dordrecht Heidelberg New York, 2011.

[19] A. Boccia, L. Grüne, and K. Worthmann, "Stability and feasibility of state constrained MPC without stabilizing terminal constraints," *Syst. Control Lett.*, vol. 72(8), pp. 14–21, 2014.

[20] L. Grüne, "Analysis and design of unconstrained nonlinear MPC schemes for finite and infinite dimensional systems," *SIAM J. Control Optimization*, vol. 48, pp. 1206–1228, 2009.

[21] A. Wächter and T. L. Biegler, "On the implementation of an interior-point filter line-search algorithm for large-scale nonlinear programming," *Math. Programming*, vol. 106, no. 1, pp. 25–57, 2006.

[22] J. Andersson, "A General-Purpose Software Framework for Dynamic Optimization," PhD thesis, Arenberg Doctoral School, KU Leuven, Dept. Elect. Eng., Belgium, October 2013.

# Chapter 5

# Occupancy Grid based Distributed Model Predictive Control of Mobile Robots

## 5.1 Abstract

In this chapter, we introduce a novel approach for reducing the communication load in distributed model predictive control (DMPC) for non-holonomic mobile robots. The key idea is to project the state prediction into a grid resulting in an occupancy grid, which can be communicated as tuples of integer values instead of floating point values. Additionally, we employ a differential communication technique to further reduce the communication load. The proposed approach has the advantage of utilizing continuous optimization methods while the communication is conducted in a quantized setting. We explore and evaluate this method numerically when applied to a group of non-holonomic robots.

## 5.2 Introduction

Formation control of mobile robots has attracted significant interest in the last decade due to its potential application areas, e.g. patrolling missions [1], search and rescue operations [2], and situational awareness [3]. The control objective is to coordinate a group of robots to first form, and then maintain a prescribed formation. This can be done either in a cooperative or non-cooperative manner, see, e.g. [4] or [5]. In the latter case, each robot has its own target with respect to dynamic (e.g. collision avoidance) or static constraints (e.g. restrictions of the operating region).

Most of formation control approaches presented in the literature can be categorized under virtual structure, behavior-based, and leader follower [6]. In the virtual structure framework, a robotic team is considered as a rigid body and the individual robots are treated as points of this body. Thus, the overall motion of the formation defines the motion of the individual robots [7, 8]. In behavior-based control, several kinds of behaviors are distinguished, e.g. formation keeping (stabilization), goal seeking, and obstacle avoidance [9]. Solution methods used in behavior-based control include motor scheme control [10], null-space-based behavior control [11], potential fields functions [12], and reciprocal velocity obstacles approach [13]. The basic approach of the leader follower structure is that a follower is assigned to track another vehicle with an offset. Here, proposed solution methods include feedback linearization [14], backstepping control [15], and sliding mode control [16].

In this chapter, we utilize a model predictive control (MPC) scheme to stabilize a group of autonomous non-cooperative vehicles (formation stabilization) while avoiding collisions. In MPC, we first measure the state of the system, which serves as a basis for solving a finite horizon optimal control problem (OCP). This results in a sequence of future control values. Then, the first element of the computed sequence is applied before the process is repeated, see [17] for an overview on non-linear MPC. Within this setting, the cost function of the OCP allows us to encode a control objective and to evaluate the closed loop performance. Moreover, static and dynamic constraints can be directly taken into account [18]. MPC has been used in several studies considering formation control, see, e.g. [19–21]. Here, we consider a distributed implementation of MPC in which the system is split into subproblems regarding a single robot each. The robots solve their own (local) OCP's and communicate with each other to avoid collisions. Since the strong coupling among the robots may render the control task infeasible in decentralized control (no communication), see, e.g. [22], this data exchange is necessary. Here, the individual optimization tasks are executed in a fixed order similar to the method proposed by Richards and How [23].

To formulate the coupling constraints, previous studies were based on the communication of predicted trajectories, see [24]. In contrast, we first partition the operating region into a grid and derive an estimate on the minimum width of a grid cell. Then, the predicted trajectories are projected onto the grid resulting in an occupancy grid, which serves as quantization of the communication data. Based on a data exchange of these projections, each robot formulates suitable coupling constraints. Utilizing the occupancy grid reduces bandwidth limits and congestion issues since a more compact data representation (integers instead of floating point values) is employed. In summary, the optimization is performed in a spatial set to make use of sensitivity information while the communication is conducted in a quantized set to reduce the transmission load. To construct the constraints, we consider squircles as an approximation for a grid cell to use gradient-based algorithms for optimization. Furthermore, we numerically compute the minimal horizon length such that a desired closed-loop performance is achieved and analyse the influence of the grid width on the communication effort and total cost.

Since the predictions are now in a discrete set, a differential communication scheme can be utilized. The differential communication scheme ensures that only altered cells in a given occupancy grid prediction are communicated in comparison to predictions of previous time steps. As a result, this method reduces further the communication capacity requirements significantly.

The chapter is organized as follows: we first introduce the grid generation and present the distributed MPC scheme in Section 5.3. Then, we derive a suitable representation of the collision avoidance constraints. A differntial communication scheme is introduced in Section 5.5. Afterwards, we investigate the proposed method by means of numerical simulations for a group of non-holonomic robots in Section 5.6, before conclusions are drawn.

## 5.3 Problem Setting

Firstly, we present the dynamics of the considered robots and a respective grid representation of the operating region. Thereafter, we introduce the optimal control problem (OCP) of each robot before these components are integrated in a distributed model predictive control (DMPC) scheme.

The control objective is to stabilize a group of $P$ non-holonomic robots, $P \in \mathbb{N}$, at reference equilibria while avoiding inter-robot collisions. Each robot is considered as an agent of a multi-agent system with state $\mathbf{x}_i \in X \subset \mathbb{R}^3$, $i \in \{1, 2, \ldots, P\}$. For the $i^{th}$ robot, the Euler discrete time model is given by

$$\begin{pmatrix} x_i(k+1) \\ y_i(k+1) \\ \theta_i(k+1) \end{pmatrix} = \mathbf{x}_i(k+1) = \mathbf{f}_{\delta,i}(\mathbf{x}_i(k), \mathbf{u}_i(k)) = \begin{pmatrix} x_i(k) \\ y_i(k) \\ \theta_i(k) \end{pmatrix} + \delta \cdot \begin{pmatrix} v_i(k) \cos(\theta_i(k)) \\ v_i(k) \sin(\theta_i(k)) \\ \omega_i(k) \end{pmatrix} \tag{5.1}$$

with vector field $\mathbf{f}_{\delta,i} : \mathbb{R}^3 \times \mathbb{R}^2 \to \mathbb{R}^3$. Here, $k$ denotes the current time step, and $\delta$ denotes the sampling period. The state $\mathbf{x}_i \in X \subset \mathbb{R}^3$ contains position and orientation coordinates, i.e. $\mathbf{x}_i = (x_i, y_i, \theta_i)^\top$ (m,m,rad)$^\top$. The control $\mathbf{u}_i \in U \subset \mathbb{R}^2$ is given by $\mathbf{u}_i = (v_i, \omega_i)^\top$. $v_i$ (m/s) is the linear speed of the robot and $\omega_i$ (rad/s) is the angular speed. Each robot has an initial state $\mathbf{x}_i^0 := \mathbf{x}_i(0)$ and a reference $\mathbf{x}_i^r$. State and control constraints are incorporated in

$$X := [-\bar{x}, \bar{x}] \times [-\bar{y}, \bar{y}] \times \mathbb{R}, \qquad U := [-\bar{v}, \bar{v}] \times [-\bar{\omega}, \bar{\omega}]$$

for $\bar{x}, \bar{y}, \bar{v}, \bar{\omega} > 0$.

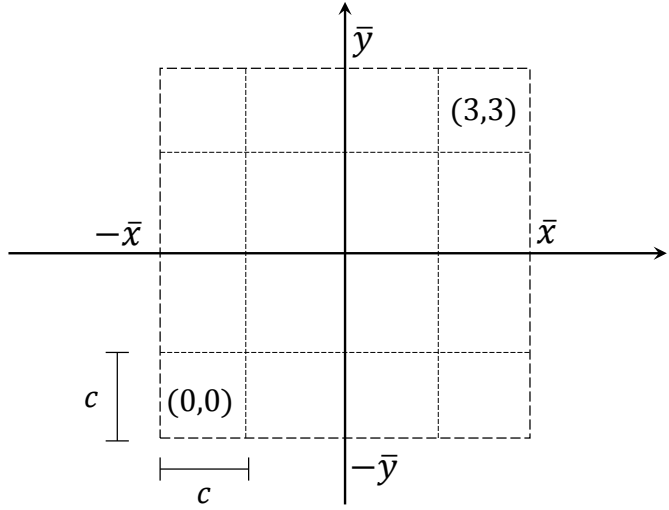

Figure 5.1: Example of a discretized set $\mathcal{G}$ with $a_{\max} = b_{\max} = 4$

We partition the set $[-\bar{x}, \bar{x}] \times [-\bar{y}, \bar{y}]$ into a grid of squared cells, see Figure 5.1. Each cell has a cell width $c$, where $c$ satisfies $2\bar{x} = a_{\max} c$ and $2\bar{y} = b_{\max} c$ with $a_{\max}, b_{\max} \in \mathbb{N}$. Therefore, each cell in the grid has a unique index $(a, b) \in \mathcal{G}$, where

$$\mathcal{G} := \{0, 1, \ldots, a_{\max} - 1\} \times \{0, 1, \ldots, b_{\max} - 1\} \subset \mathbb{N}_0^2.$$

Then, a state $\mathbf{x}_i$ (and its location $(x_i, y_i)$) can be mapped to the discrete set $\mathcal{G}$ using the quantization $q : X \to \mathcal{G}$ defined as

$$q(\mathbf{x}_i) = (a_i, b_i) := \left( \left\lfloor \frac{x_i + \bar{x}}{c} \right\rfloor, \left\lfloor \frac{y_i + \bar{y}}{c} \right\rfloor \right). \quad (5.2)$$

Note that a location $(x, y)$ can be also expressed by a unique single index instead of $(a, b)$ by utilizing a full enumeration.

We next present the required notation to set up the DMPC algorithm. To this end, we define our quantized trajectory, communication scheme and running costs to compose the optimal control problem (OCP). For a finite control sequence

$$\mathbf{u}_i = (\mathbf{u}_i(0), \mathbf{u}_i(1), \ldots, \mathbf{u}_i(N-1)) \in U^N$$

and an initial value $\mathbf{x}_i^0$, the predicted state trajectory of robot $i$ over prediction horizon $N \in \mathbb{N}$ is given by

$$\mathbf{x}_i^{\mathbf{u}}(\cdot; \mathbf{x}_i^0) := \left( \mathbf{x}_i^{\mathbf{u}}(0; \mathbf{x}_i^0), \mathbf{x}_i^{\mathbf{u}}(1; \mathbf{x}_i^0), \ldots, \mathbf{x}_i^{\mathbf{u}}(N; \mathbf{x}_i^0) \right).$$

Here, we removed the subscript $i$ from $\mathbf{u}_i$ in $\mathbf{x}_i^{\mathbf{u}}(\cdot; \mathbf{x}_i^0)$ to simplify our notation. The trajectory together with the introduced quantization $q$ allow us to define the occupancy grid $\mathcal{I}_i(n) \in (\mathbb{N}_0 \times \mathcal{G})^{N+1}$ at time $n$ as

$$\mathcal{I}_i(n) := \left( n+k, q\left( \mathbf{x}_i^{\mathbf{u}}\left(k; \mathbf{x}_i^0\right) \right) \right) =: \left( n+k, a_i^{\mathbf{u}}(k; \mathbf{x}_i^0), b_i^{\mathbf{u}}(k; \mathbf{x}_i^0) \right), \quad k \in \{0, 1, \ldots, N\}, \quad (5.3)$$

where $(a_i^{\mathbf{u}}(k; \mathbf{x}_i^0), b_i^{\mathbf{u}}(k; \mathbf{x}_i^0))$ denotes the quantized state. While each robot $i$ sends only one such package $\mathcal{I}_i(n)$, it collects all received occupancy grids in $I_i$, i.e.

$$I_i(n) := (\mathcal{I}_1(n), \ldots, \mathcal{I}_{i-1}(n), \mathcal{I}_{i+1}(n), \ldots, \mathcal{I}_P(n)).$$

Later, we use this data to construct the coupling constraints to ensure collision avoidance. The main motivation of communicating the occupancy grid instead of the predicted trajectories is to transmit integers instead of floating point values to reduce the communication load.

Here, we utilize a distributed model predictive control (DMPC) scheme, which implements a local controller for each robot $i$. To this end, for a given reference state $\mathbf{x}_i^r$, we choose the (local) running costs $\ell_i : X \times U \to \mathbb{R}_{\geq 0}$, which are supposed to be positive definite with respect to $\mathbf{x}_i^r$. Using the results of [25, 26], which embody the analysis presented in Chapter 3, the running costs $\ell_i$ are chosen as

$$\ell_i(\mathbf{x}_i, \mathbf{u}_i) := q_1(x_i - x_i^r)^4 + q_2(y_i - y_i^r)^2 + q_3(\theta_i - \theta_i^r)^4 + r_1 v_i^4 + r_2 \omega_i^4 \tag{5.4}$$

to ensure stabilization of non-holonomic robots without terminal constraints or costs. Then, at every time step $n$, each robot $i$ solves the following finite horizon optimal control problem for initial condition $\mathbf{x}_i^0 = \mathbf{x}_i(n)$ and given data $I_i(n)$:

$$\min_{\mathbf{u}_i \in \mathbb{R}^{2 \times N}} J_i^N\left(\mathbf{u}_i; \mathbf{x}_i^0, I_i(n)\right) := \sum_{k=0}^{N-1} \ell_i\left(\mathbf{x}_i^{\mathbf{u}}\left(k; \mathbf{x}_i^0\right), \mathbf{u}_i(k)\right) \tag{5.5}$$

$$\text{subject to } \mathbf{x}_i^{\mathbf{u}}(0) = \mathbf{x}_i^0,$$

$$\mathbf{x}_i^{\mathbf{u}}(k+1) = \mathbf{f}_{\delta,i}(\mathbf{x}_i(k), \mathbf{u}_i(k)) \quad \forall k \in \{0, 1, \ldots, N-1\},$$

$$G\left(\mathbf{x}_i^{\mathbf{u}}(k), I_i(n)\right) \leq 0 \quad \forall k \in \{1, 2, \ldots, N\},$$

$$\mathbf{x}_i^{\mathbf{u}}(k) \in X \quad \forall k \in \{1, 2, \ldots, N\},$$

$$\mathbf{u}_i(k) \in U \quad \forall k \in \{0, 1, \ldots, N-1\},$$

As a result, an optimal control sequence

$$\mathbf{u}_i^\star = \left(\mathbf{u}_i^\star(0), \ldots, \mathbf{u}_i^\star(N-1)\right) \tag{5.6}$$

is obtained, which minimizes the cost function (5.5).[1] The corresponding value function $V_i^N : X \times (\mathbb{N}_0 \times \mathcal{G})^{(P-1)(N+1)} \to \mathbb{R}_{\geq 0}$ is defined as

$$V_i^N \left( \mathbf{x}_i^0, I_i(n) \right) := J_i^N(\mathbf{u}_i^\star; \mathbf{x}_i^0, I_i(n)). \quad (5.7)$$

The detailed construction of the collision avoidance constraints $G$ will be explicated in the ensuing Section 5.4. Using this formulation of the OCP, in analogy to Richards and How [23, 27] the DMPC algorithm is as follows:

---
**Algorithm 5.1** DMPC-algorithm for the overall system
---
1: **Obtain and communicate** admissible control sequences $\mathbf{u}_i$ for initial values $\mathbf{x}_i(0)$ for all $i \in \{1, 2, \ldots, P\}$
2: **for** $n = 0, 1, \ldots$ **do**
3:     **for** $i$ from 1 to $P$ **do**
4:         **Receive** $I_i(n)$
5:         **Solve** OCP (5.5) and **Apply** $\mathbf{u}_i^*(0)$ to the $i^{th}$ robot.
6:         **Broadcast** $\mathcal{I}_i(n)$
7:     **end for**
8: **end for**
---

In the generic case, the calculation of an admissible set of controls for all robots in the initialization of Algorithm 5.1 is particularly demanding. In our setting, however, we avoid the usage of terminal constraints and follow an approach outlined in [24, 28], which allows us to employ $\mathbf{u}_i \equiv 0$ to ensure initial feasibility. We like to note that in contrast to [23, 27] initial feasibility is not sufficient for stability in our case, but instead stability conditions similar to [24] need to be verified, which is outside the scope of this thesis for the distributed control case.

---
[1] To avoid technical difficulties, existence of a minimizer is assumed.

## 5.4 Formulation of the Constraints

In this section, we derive an appropriate cell width $c$ based on the dynamics of the system and specify a safety distance between two robots. Thereafter, we expound a method to set up collision avoidance constraints $G$ for the OCP (5.5). This method relies on a squircle approximation.

To formulate the coupling constraints, we first define the backward mapping $f_c : \mathbb{N}_0 \times \mathcal{G} \to \mathbb{R}^2$ given by

$$f_c((n, a, b)) = \underbrace{(c \cdot (a + 0.5) - \bar{x}, c \cdot (b + 0.5) - \bar{y})^\top}_{=:(x^c, y^c)_{(a,b)}}, \tag{5.8}$$

which transforms a cell index $(a, b) \in \mathcal{G}$ to the corresponding cell center location, i.e.

$$(x^c, y^c)_{(a,b)} \in [-\bar{x}, \bar{x}] \times [-\bar{y}, \bar{y}].$$

### 5.4.1 Grid generation and safety margins

The minimum cell width $\underline{c}$ has to be chosen large enough such that a cell can not be skipped during one time step. Therefore, a minimum cell width is defined by

$$\underline{c} = \lceil \bar{v} \cdot \delta \rceil.$$

With this minimum cell width, we ensure that each robot cannot skip a cell by applying $\mathbf{u}_i \in U$. Moreover, to guarantee a minimum distance $d_{\min}$ between the robots and to avoid overlapping of predicted trajectories, see Figure 5.2 (left) for an illustration, we additionally inflate each (occupied) cell. The inflation margin is given by $\max\{d_{\min}, \underline{c}\}$. Furthermore, we ensure that constraint violations due to sampling cannot occur. Here, the worst case is an intermediate step with angle $\frac{\pi}{4}$ as illustrated in Figure 5.2 (right). Therefore, cells are further inflated to compensate for constraint violation due to sampling by $\underline{c}\cos(\pi/4)/2$. Hence, for a given cell width $c \geq \underline{c}$, the safety margins of an occupied cell are specified by a square with side length $\Psi$ given by

$$\Psi := c + 2 \cdot \max\{d_{\min}, \underline{c}\} + \underline{c} \cdot \cos(\pi/4) + \epsilon$$

where $\epsilon \ll 1$ is a numerical safety margin.

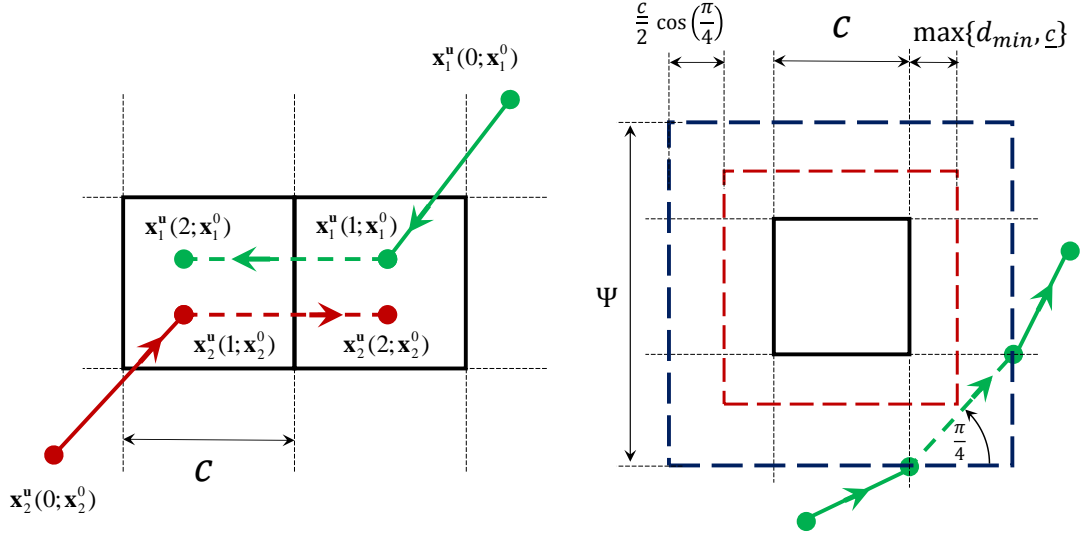

Figure 5.2: Left: overlapping of prediction trajectories. Right: cell inflation to avoid violation of a state constraint due to quantization.

### 5.4.2 Constraint construction

As robot $i$ obtains the data $I_i(n)$ from the other robots, the constraints $G(\mathbf{x}_i^\mathbf{u}(k;\mathbf{x}_i^0), I_i(n))$, $k \in \{1, 2, \ldots, N\}$, are assembled by combining $g_{q,k}^i$, $q \in \{1, 2, \ldots, P-1\}$, as

$$\left(g_{1,k}^i, \ldots, g_{i-1,k}^i, g_{i+1,k}^i, \ldots, g_{P,k}^i\right),$$

which resemble the collision avoidance. The functions

$$g_{q,k}^i = g_{q,k}^i \left( \mathbf{x}_i^{\mathbf{u}}\left(k; \mathbf{x}_i^0\right), \mathcal{I}_q(n)(k) \right)$$

are constructed utilizing squircles [29]. A squircle, which is a special case of a superellipsoid, is a geometric shape between a square and a circle, see Figure 5.3 for an illustration of apporximating a square by a squircle; as can be noticed, a squircle approximation is less conservative than the circular approximation. The polar equations of a squircle are given by

$$x(\phi) = 2^{-3/4}h \cdot |\cos(\phi)|^{\frac{1}{2}} \cdot \text{sgn}(\cos(\phi)),$$
$$y(\phi) = 2^{-3/4}h \cdot |\sin(\phi)|^{\frac{1}{2}} \cdot \text{sgn}(\sin(\phi))$$

with $\phi \in [0, 2\pi)$. $h$ is the width of the square, which a squircle is approximating from

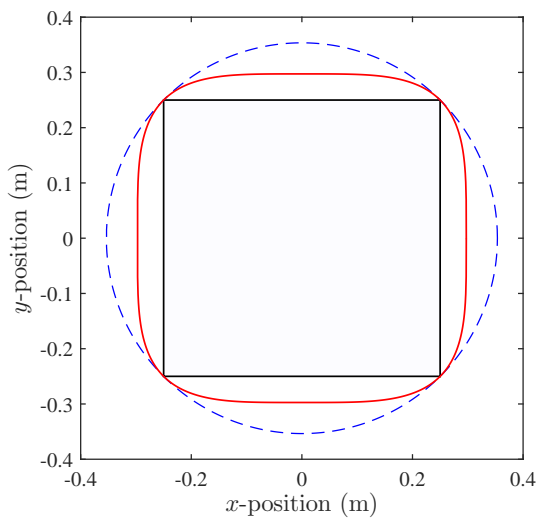

Figure 5.3: Approximating a cell constraints (black) using a circle (blue) and a squircle (red).

outside. Now, the coupling constraints $g_{q,k}^i$ of the $i^{th}$ robot can be formulated using the

Euclidean norm as[2]

$$\left\| \begin{pmatrix} x_i(k) \\ y_i(k) \end{pmatrix} - f_c\left(\mathcal{I}_q(n)(k)\right) \right\|_2 - \frac{\Psi\sqrt{|\cos(\beta)|+|\sin(\beta)|}}{2^{3/4}}$$

with $k \in \{1, 2, \ldots, N\}$ and $q \in \{1, 2, \ldots, P\} \setminus \{i\}$. The angle $\beta$ is calculated via

$$\beta = \arctan\left(\frac{y_i(k) - y_i^c(k)}{x_i(k) - x_i^c(k)}\right),$$

where $(x_i^c(k), y_i^c(k))$ is obtained via the mapping (5.8) given the predictions $\mathcal{I}_q(n)(k)$.

## 5.5 Differential Communication

In this section, we present the communication scheme, which utilizes a differential update approach. Each occupancy grid prediction $\mathcal{I}_i$, i.e. (5.3), is composed of a sequence of time stamps and their corresponding indices of the occupied cells prediction. $\mathcal{I}_i(n)$ returns the corresponding tuple $(n, a_n, b_n)$. The fundamental idea of differential communication is, that each robot $i$ keeps a memory of the most recent communicated tuple of its own prediction denoted as $\mathcal{I}_i^-$, and uses it in the preparation of the broadcasted tuple at the current time step to create a differential update denoted as $\mathcal{I}_i^c$. Hence, the differential update $\mathcal{I}_i^c$ considers only altered cells from a previous time step or new occupied cells obtained from the current occupancy grid prediction. Algorithm 5.2 shows how the communicated tuple $\mathcal{I}_i^c$ is prepared to be broadcasted. In Algorithm 5.2, after the complete occupancy grid $\mathcal{I}_i$ is computed in step 3, the differential update $\mathcal{I}_i^c$ is created by comparing $\mathcal{I}_i$ with $\mathcal{I}_i^-$ using a For-loop. Thereafter, the differential update $\mathcal{I}_i^c$ is broadcasted in step 9 before the first data entry in $\mathcal{I}_i^-$ is discarded in step 10 because it will not be used

---

[2]The absolute function $|\cdot|$ is approximated by the differentiable expression $|e| \approx \sqrt{e^2 + \epsilon}$, for $e \in \mathbb{R}$ and $\epsilon \ll 1$.

in the following time step. In the same step, $\mathcal{I}_i^c$ is additionally emptied.

---

**Algorithm 5.2** Preparation of the communicated tuples for a single robot $i$
---
1: **Set** $\mathcal{I}_i^c = \mathcal{I}_i^- := \emptyset$
2: **for** each time step $n = 0, 1, \ldots$ **do**
3:     **Compute** $\mathcal{I}_i := ((n, a_n, b_n), \ldots, (n+N, a_{n+N}, b_{n+N}))$
4:     **for** $m = n : n + N$ **do**
5:         **if** $\mathcal{I}_i(m) \notin \mathcal{I}_i^-$ **then**
6:             $\mathcal{I}_i^c := \mathcal{I}_i^c \cup \mathcal{I}_i(m)$
7:         **end if**
8:     **end for**
9:     **Broadcast** $\mathcal{I}_i^c$
10:     **Set** $\mathcal{I}_i^- := \mathcal{I}_i \setminus \mathcal{I}_i(n), \mathcal{I}_i^c := \emptyset$
11: **end for**

---

Since each robot communicates only the altered cells in its prediction, each other robot has to additionally keep a memory of the most recently received information from other robots. This information is later used to assemble the complete occupancy grid predictions used in formulating the coupling constraints of robot $i$. Algorithm 5.3 shows how the complete prediction of a robot $i$ is assembled by robot $j$ after receiving the differential update $\mathcal{I}_i^c$. In Algorithm 5.3, after receiving the differential update $\mathcal{I}_i^c$ in step 3, the occupancy gird $\mathcal{I}_i$ is assembled by comparing the data in $\mathcal{I}_i^c$ and $\mathcal{I}_i^-$ using a For-loop with two If-conditions; the first one updates the entries of $\mathcal{I}_i^-$, which have the same time-stamps in $\mathcal{I}_i^c$. In the second If-condition, $\mathcal{I}_i$ is assembled by combining $\mathcal{I}_i^-$ with the extra entries from $\mathcal{I}_i^c$. Thereafter, in step 12, the first data entry in $\mathcal{I}_i$ is discarded because it will not be used in the following time step. In addition, at this step, $\mathcal{I}_i^-$ is updated by $\mathcal{I}_i$.

**Algorithm 5.3** Assembly of the complete prediction of robot $i$ by another robot $j \in \{1, 2, \ldots, P\} \setminus i$

1: **Set** $\mathcal{I}_i = \mathcal{I}_i^- := \emptyset$
2: **for** each time step $n = 0, 1, \ldots$ **do**
3:     **Receive** $\mathcal{I}_i^c$
4:     **for** $m = n : n + N$ **do**
5:         **if** $\mathcal{I}_i^c(m)_{ts} \in \mathcal{I}_i^-$ **then**
6:             $\mathcal{I}_i^-(m) := \mathcal{I}_i^c(m)$
7:         **end if**
8:         **if** $\mathcal{I}_i^c(m) \notin \mathcal{I}_i^-$ **then**
9:             $\mathcal{I}_i := \mathcal{I}_i^- \cup \mathcal{I}_i^c(m)$
10:        **end if**
11:     **end for**
12:     **Set** $\mathcal{I}_i := \mathcal{I}_i \setminus \mathcal{I}_i(n)$, $\mathcal{I}_i^- := \mathcal{I}_i$
13: **end for**

$\mathcal{I}_i^c(m)_{ts}$ returns only the time stamp at index $m$.

These assembled tuples are used to generate the coupling constraints among the subsystem, which were explored in the Section 5.4.

As a result, Algorithm 5.1 can be reformulated based on the differential communication, i.e. in Algorithm 5.1, line 4 (the receive command) will be replaced with

**Receive** $\mathcal{I}_i^c, i \in \{1, 2, \ldots, P\} \setminus \{i\}$

**Assemble** $\mathcal{I}_i, i \in \{1, 2, \ldots, P\} \setminus \{i\}$ using Algorithm 5.3.

Moreover, (the broadcast command), i.e. line 6 will be replaced with

**Broadcast** $\mathcal{I}_i^c$ using Algorithm 5.2.

## 5.6 Numerical Experiments

In this section, we explore the performance of Algorithm 5.1 through numerical simulations of non-holonomic mobile robots. To this end, we consider a group of $P = 4$ mobile robots and set the state and the control constraints as

$$X := [-6, 6]^2 \times \mathbb{R}, \qquad U := [-1, 1] \times [-1, 1].$$

The sampling time is chosen as $\delta = 0.5$ (seconds), while the weighting parameters of the running costs (5.4) are chosen as

$$q_1 = 1, q_2 = 25, q_3 = 1, r_1 = 0.2, \text{ and } r_2 = 0.2.$$

The minimum distance to be ensured between the robots is chosen as $d_{\min} = 0.5$ (m). Based on the given constraints, a minimum cell width $\underline{c} = 0.5$ (m) is computed, see Section 5.4.1. We perform simulations with different cell widths $c \in \{0.5, 1, 1.5, 2\}$. The initial conditions of the robots as well as their references can be found in Table 5.1. Moreover, closed loop simulations are executed until all robots meet the condition

$$\|\mathbf{x}_i(n) - \mathbf{x}_i^r\| \leq 0.01. \tag{5.9}$$

The total number of closed loop iterations is denoted by $n_\#$.

We investigate the performance of Algorithm 5.1 in terms of communication effort as

Table 5.1: Initial Conditions and references of the robots

| System ($i$) | Initial Conditions ($\mathbf{x}_i^0$) | References ($\mathbf{x}_i^r$) |
|---|---|---|
| 1 | $(4.5, 4.5, -3\pi/4)^\top$ | $(-4.5, -4.5, \pi)^\top$ |
| 2 | $(-4.5, 4.5, -\pi/4)^\top$ | $(4.5, -4.5, 0)^\top$ |
| 3 | $(4.5, -4.5, 3\pi/4)^\top$ | $(-4.5, 4.5, \pi)^\top$ |
| 4 | $(-4.5, -4.5, \pi/4)^\top$ | $(4.5, 4.5, 0)^\top$ |

well as accumulated closed loop costs:

- The communication effort is given by

$$K := \sum_{n=0}^{n_\#} \sum_{p=1}^{P} \#\mathcal{I}_i^c|_n, \tag{5.10}$$

i.e. $K$ is the accumulated number of the broadcasted tuples among all robots over the closed loop simulation.

- The accumulated closed loop costs are given by

$$M := \sum_{n=0}^{n_\#} \sum_{i=1}^{P} \ell_i\left(\mathbf{x}_i^{MPC}(n), \mathbf{u}_i^{MPC}(n)\right) \tag{5.11}$$

where $\mathbf{u}_i^{MPC}(n)$ denotes the control signal applied in Step 5) of Algorithm 5.1 at time instant $n$ and $\mathbf{x}_i^{MPC}(\cdot)$ the corresponding closed loop trajectory.

Simulations were performed utilizing the (interior-point) optimization method provided by IPOPT [30] coupled with MATLAB via CasADi toolbox [31].

Figure 5.4 illustrates the required communication effort $K$ computed via (5.10) and the accumulated costs (5.11), respectively. The growth of these key figures is caused by robots, which must take detours to avoid collisions and reach their reference states, see Figure 5.5 for simulation snap shots[3]. As the cell sizes are growing, the length of these detours as well as the simulation time are increasing, which leads to the observed increases in communication load and costs. It can be also noticed from Figure 5.4 (left) that adopting the differential communication scheme, i.e Algorithm 5.2, reduces the communication load significantly. To this end, a communication reduction of at least 72% has been observed, see Table 5.2. Within our simulations, we also observed that the minimal prediction horizon $N$ required to meet Condition (5.9) grows with increasing cell size, see Table 5.3. The cause is the quantization together with the collision avoidance constraints as the robots need a longer prediction horizon to observe a decrease in the target distance while circumventing occupied cells. Hence, if chosen too coarse, the quantization may lead to a larger prediction horizon $N$ and raise the computational effort.

---
[3]Click here for a visualization video of the numerical simulations.

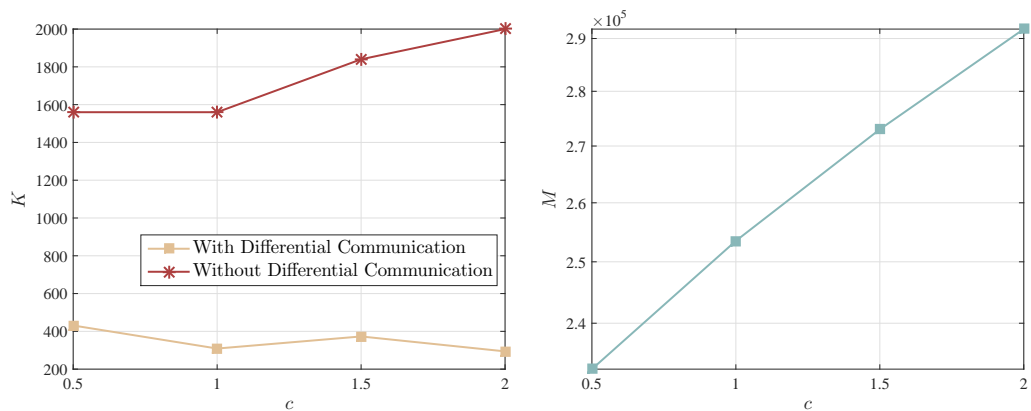

Figure 5.4: Left: Cell width $c$ vs. communication $K$. Right: closed loop costs $M$ (right). All simulations were run for $N = 9$.

Table 5.2: The number of closed loop iterations as well as the communication reduction percentage for $N = 9$ and for each cell size $c$.

| Cell Width $c$ | Number of Iterations for $N = 9$ | Communication Reduction Percentage (%) $N = 9$ |
|---|---|---|
| 0.5 | 40 | 72.4 |
| 1 | 40 | 80 |
| 1.5 | 47 | 80 |
| 2 | 51 | 85 |

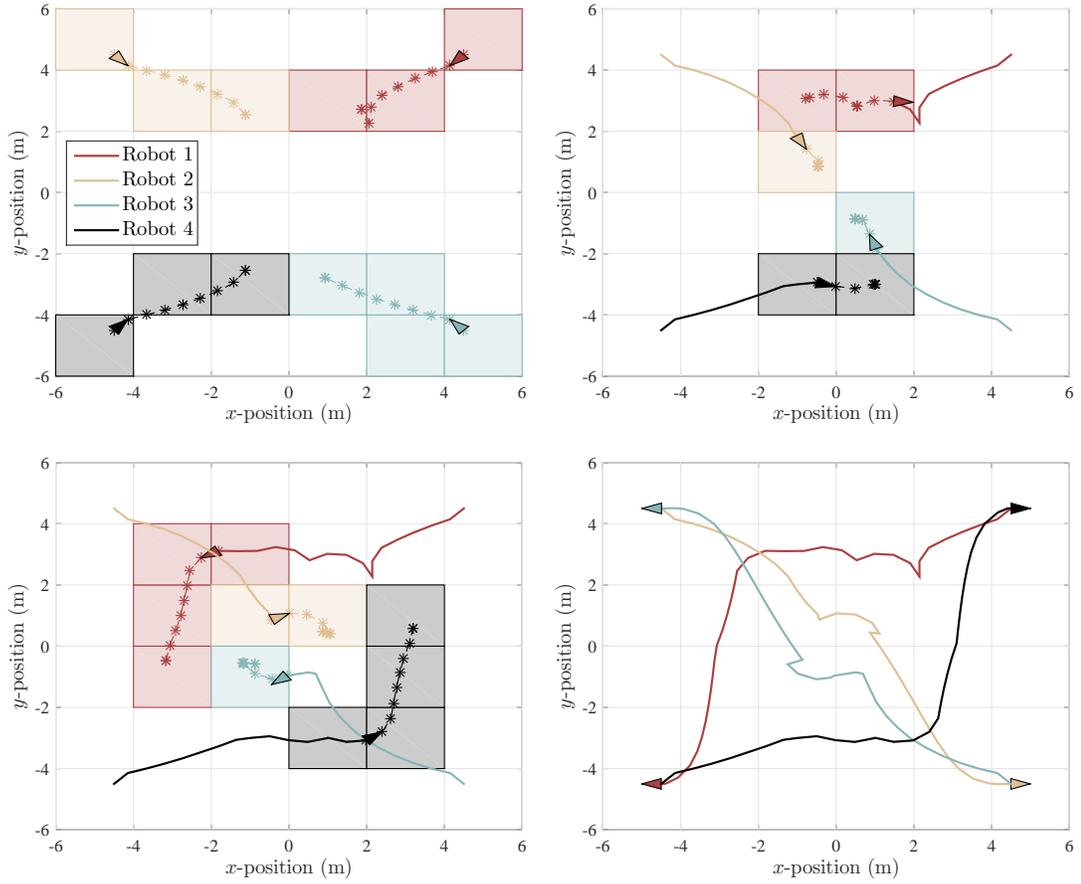

Figure 5.5: Snap shots of closed loop behavior for $c = 2$ and $N = 9$ taken at $n = 1$ (top left), $n = 10$ (top right), $n = 20$ (bottom left) and $n = 51$ (bottom right). Trajectories are presented as continuous lines and predictions as dashed lines with star markers.

Table 5.3: Impact of cell width $c$ on minimal $N$ and $n_{\#}$.

| $c$ | $N$ | $n_{\#}$ |
|---|---|---|
| 0.5 | 7 | 40 |
| 1 | 7 | 45 |
| 1.5 | 9 | 47 |
| 2 | 9 | 51 |

Figure 5.6 (left) shows the evolution of the sum

$$M_P(n) = \sum_{i=1}^{P} \ell_i \left( \mathbf{x}_i^{MPC}(n), \mathbf{u}_i^{MPC}(n) \right)$$

for each considered cell width, whereas Figure 5.6 (right) depicts the development of $\ell_i \left( \mathbf{x}_i^{MPC}(n), \mathbf{u}_i^{MPC}(n) \right)$ for each robot with $c = 2$. It can be observed that the overall instantaneous costs $M_P(n)$ are monotonically decreasing for any quantization $c$, yet, on a local level, increases in the costs $\ell_i$ may occur.

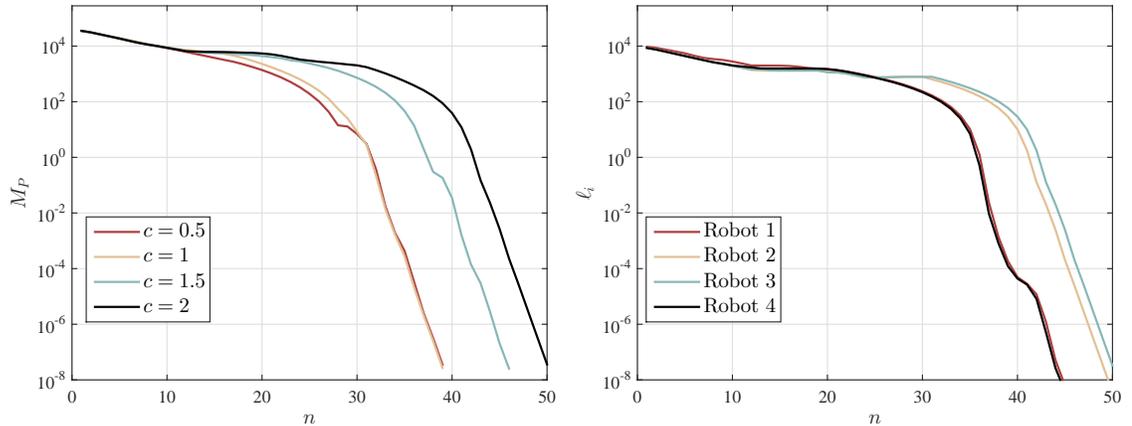

Figure 5.6: Left: Development of $M_P$ with $c = \{0.5, 1.0, 1.5, 2.0\}$ and $N = 9$. Right: Development of $\ell_i \left( \mathbf{x}_i^{MPC}(n), \mathbf{u}_i^{MPC}(n) \right)$ for each robot with $c = 2$ and $N = 9$.

### 5.6.1 Discussion of results

The advantage of the proposed approach in relation to previous models based on the communication of predicted trajectories [24] is given by its low communication requirements. The latter is induced by the quantization of broadcasted data, which is transmitted as integers instead of floating point values. Using the proposed differential communication technique allowed for reducing the communication load even further. Moreover, Figure 5.6 (left) illustrates the decreasing costs for the robots for the distributed control method proposed here. This can be interpreted as convergence of the overall system towards its target. Hence, we have numerically shown that the proposed approach is applicable.

From the individual costs, see Figure 5.6 (right), we can conclude that the decrease of the costs depends on two aspects: the ordering in which robots compute their optimal control law, and the convergence of the robots to their reference locations while obeying the collision avoidance. The first one is straightforward because the individual optimization tasks are executed by each robot in a fixed order (priority rule), see [28] for an in-depth analysis of the influence of priority rules in distributed MPC. Therefore, the solution space is smaller for the last robot when all previous robots have reserved their occupancy grid. Consequently, the fixed ordering may temporarily increase the costs of single robots. The influence of these two aspects on the closed loop simulations can be seen as follows. In Figure 5.5 (top right), robot 1 computed its optimal solution first, which reveals an initial backward motion. This allowed robot 2 to approach the center of the map before robot 1. Moreover, since robot 3 performs its optimization before robot 4, robot 3 approaches the map center before robot 4. However, since robot 1 is optimizing first, robot 3 had to wait for robot 1 (see Figure 5.5 (bottom left)). Only after robot 2 proceeded a number of steps while avoiding collisions with robot 3 (see Figure 5.5 (top right)), the way of robot 4 to its reference became completely clear and the robot reached its reference first (nearly with robot 1). After the way is now clear for robots 2 and 3, they reach their references almost simultaneously, see Figure 5.6 (right).

Generally, adopting the minimum cell width $\underline{c}$ in forming the grid, lead to the most efficient closed loop performance of the proposed approach, see Figure 5.6. Nonetheless, increasing the cell width can reduce the communication load, see Table 5.3. Therefore, the cell width $c$ can be adapted to favor either performance or communication where, in all cases of $c$, the data-type change of the transmitted information is guaranteed.

## 5.7 Conclusions

In this chapter, we proposed a distributed MPC scheme based on communication of occupancy grids instead of predicted trajectories for mobile robots. Each occupancy grid is composed of a sequence of tuples describing the corresponding occupied cells. Changing the data representation of the communicated stream into a quantized space leads to alleviated requirements on the communication capacity of the controller. It has been also shown that the communication load can be further relaxed by adopting the proposed differential communication scheme, where only altered cells from previous communication steps are broadcasted. Moreover, despite the quantization of the communicated information, sensitivity information can be exploited because the optimization is performed in a continuous setting.

In this chapter, we provided a possible formulation method of the coupling constraints utilized by the proposed controller. This method is based on utilizing a squircle shape as an approximation of a square. This method has the advantage of being less conservative than a circular approximation as well as being suitable for derivative-based optimization methods. The proposed controller with the constraints formulation method proposed led to satisfactory results in terms of communication as well as convergence.

## References

[1] J. Acevedo, B. Arrue, I. Maza, and A. Ollero, "Distributed approach for coverage and patrolling missions with a team of heterogeneous aerial robots under communication constraints," *International Journal of Advanced Robotic Systems*, vol. 10, no. 28, pp. 1–13, 2013.

[2] M. Bernard, K. Kondak, I. Maza, and A. Ollero, "Autonomous transportation and

deployment with aerial robots for search and rescue missions," *Journal of Field Robotics*, vol. 28, no. 6, pp. 914–931, 2011.

[3] M. A. Hsieh, A. Cowley, J. F. Keller, L. Chaimowicz, B. Grocholsky, V. Kumar, C. J. Taylor, Y. Endo, R. C. Arkin, B. Jung, D. F. Wolf, G. S. Sukhatme, and D. C. MacKenzie, "Adaptive teams of autonomous aerial and ground robots for situational awareness," *Journal of Field Robotics*, vol. 24, no. 11-12, pp. 991–1014, 2007.

[4] B. T. Stewart, S. J. Wright, and J. B. Rawlings, "Cooperative distributed model predictive control for nonlinear systems," *Journal of Process Control*, vol. 21, no. 5, pp. 698–704, 2011.

[5] M. Farina and R. Scattolini, "An output feedback distributed predictive control algorithm," in *IEEE Conference on Decision and Control and European Control Conference*, December 2011, pp. 8139–8144.

[6] M. Saffarian and F. Fahimi, "Non-iterative nonlinear model predictive approach applied to the control of helicopters group formation," *Robotics and Autonomous Systems*, vol. 57, no. 6–7, pp. 749–757, 2009.

[7] M. A. Lewis and K.-H. Tan, "High precision formation control of mobile robots using virtual structures," *Autonomous Robots*, vol. 4, no. 4, pp. 387–403, 1997.

[8] W. Ren and R. Beard, "Decentralized scheme for spacecraft formation flying via the virtual structure approach," *Journal of Guidance, Control, and Dynamics*, vol. 27, no. 1, pp. 73–82, 2004.

[9] M. Farina, A. Perizzato, and R. Scattolini, "Application of distributed predictive control to motion and coordination problems for unicycle autonomous robots," *Robotics and Autonomous Systems*, vol. 72, pp. 248–260, 2015.

[10] T. Balch and R. Arkin, "Behavior-based formation control for multirobot teams," *IEEE Transactions on Robotics and Automation*, vol. 14, no. 6, pp. 926–939, 1998.

[11] G. Antonelli, F. Arrichiello, and S. Chiaverini, "Flocking for multi-robot systems via the null-space-based behavioral control," in *IEEE/RSJ International Conference on Intelligent Robots and Systems*, 2008, pp. 1409–1414.

[12] T. Balch and M. Hybinette, "Social potentials for scalable multi-robot formations," in *IEEE International Conference on Robotics and Automation*, vol. 1, 2000, pp. 73–80 vol.1.

[13] J. Snape, S. J. Guy, J. van den Berg, and D. Manocha, *Smooth Coordination and Navigation for Multiple Differential-Drive Robots.* Berlin, Heidelberg: Springer Berlin Heidelberg, 2014, pp. 601–613.

[14] G. Mariottini, F. Morbidi, D. Prattichizzo, G. Pappas, and K. Daniilidis, "Leader-follower formations: Uncalibrated vision-based localization and control," in *IEEE International Conference on Robotics and Automation*, 2007, pp. 2403–2408.

[15] X. Li, J. Xiao, and Z. Cai, "Backstepping based multiple mobile robots formation control," in *IEEE/RSJ International Conference on Intelligent Robots and Systems*, 2005, pp. 887–892.

[16] J. Sanchez and R. Fierro, "Sliding mode control for robot formations," in *IEEE International Symposium on Intelligent Control*, 2003, pp. 438–443.

[17] L. Grüne and J. Pannek, *Nonlinear Model Predictive Control: Theory and Algorithms*, ser. Communications and Control Engineering. Springer London Dordrecht Heidelberg New York, 2011.

[18] R. Scattolini, "Architectures for distributed and hierarchical model predictive control – a review," *Journal of Process Control*, vol. 19, no. 5, pp. 723–731, 2009.

[19] W. B. Dunbar and R. M. Murray, "Distributed receding horizon control for multi-vehicle formation stabilization," *Automatica*, vol. 42, no. 4, pp. 549–558, 2006.

[20] M. W. Mehrez, G. K. I. Mann, and R. G. Gosine, "An optimization based approach for relative localization and relative tracking control in multi-robot systems," *Journal of Intelligent & Robotic Systems*, pp. 1–24, 2016.

[21] K. Hashimoto, S. Adachi, and D. V. Dimarogonas, "Distributed aperiodic model predictive control for multi-agent systems," *IET Control Theory Applications*, vol. 9, no. 1, pp. 10–20, 2015.

[22] A. Venkat, J. Rawlings, and S. Wright, "Stability and optimality of distributed model predictive control," *Proceedings of the 44th IEEE Conference on Decision and Control*, pp. 6680–6685, 2005.

[23] A. Richards and J. P. How, "Robust distributed model predictive control," *International Journal of Control*, vol. 80, no. 9, pp. 1517–1531, 2007.

[24] L. Grüne and K. Worthmann, *Distributed Decision Making and Control*. Springer Verlag, 2012, ch. A distributed NMPC scheme without stabilizing terminal constraints, pp. 261–287.

[25] K. Worthmann, M. W. Mehrez, M. Zanon, G. K. I. Mann, R. G. Gosine, and M. Diehl, "Model predictive control of nonholonomic mobile robots without stabilizing constraints and costs," *IEEE Transactions on Control Systems Technology*, vol. 24, no. 4, pp. 1394–1406, 2016.

[26] K. Worthmann, M. W. Mehrez, M. Zanon, G. K. Mann, R. G. Gosine, and M. Diehl, "Regulation of differential drive robots using continuous time mpc without stabilizing constraints or costs," in *Proceedings of the 5th IFAC Conference on Nonlinear Model Predictive Control (NPMC'15), Sevilla, Spain*, 2015, pp. 129–135.

[27] A. Richards and J. How, "A decentralized algorithm for robust constrained model predictive control," in *American Control Conference*, vol. 5, June 2004, pp. 4261–4266 vol.5.

[28] J. Pannek, "Parallelizing a State Exchange Strategy for Noncooperative Distributed NMPC," *System & Control Letters*, vol. 62, pp. 29–36, 2013.

[29] M. Fernández-Guasti, "Analytic geometry of some rectilinear figures," *Int. Jour. of Math. Ed. in Sci. & Tech.*, vol. 23, no. 6, pp. 895–901, 1992.

[30] A. Wächter and T. L. Biegler, "On the implementation of an interior-point filter line-search algorithm for large-scale nonlinear programming," *Mathematical Programming*, vol. 106, no. 1, pp. 25–57, 2006.

[31] J. Andersson, "A General-Purpose Software Framework for Dynamic Optimization," PhD thesis, Arenberg Doctoral School, KU Leuven, Department of Electrical Engineering (ESAT/SCD) and Optimization in Engineering Center, Kasteelpark Arenberg 10, 3001-Heverlee, Belgium, October 2013.

# Chapter 6

# An Optimization Based Approach for Relative Localization and Relative Tracking Control in Multi-Robot Systems

## 6.1 Abstract

In this chapter, an optimization based method is used for relative localization and relative trajectory tracking control in multi-robot systems (MRS's). In this framework, one or more robots are located and commanded to follow time varying trajectories with respect to another (possibly moving) robot reference frame. Such systems are suitable for a considerable number of applications, e.g. patrolling missions, searching operations, perimeter surveillance, and area coverage. Here, the nonlinear and constrained motion and measurement models in an MRS are incorporated to achieve an accurate state estimation algorithm based on nonlinear moving horizon estimation (MHE) and a tracking control method based on nonlinear model predictive control (MPC). In order to fulfill the real-time requirements, a fast and efficient algorithm based on a Real Time Iteration (RTI) scheme and automatic C-code generation, is adopted. Numerical simulations are conducted to: first, compare the performance of MHE against the traditional estimator used for relative localization, i.e. extended Kalman filter (EKF); second, evaluate the utilized relative localization and tracking control algorithm when applied to a team of multiple robots; finally, laboratory experiments were performed, for real-time performance evaluation. The conducted simulations validated the adopted algorithm and the experiments demonstrated its practical applicability.

## 6.2 Introduction

The increased number of applications of autonomous robotics in the recent decade led to a rapidly growing interest in multi-robot systems (MRS's). This is primarily because some given tasks may not be feasible by a single working robot, or multiple robots need a shorter time to execute the same task with a possibly higher performance. Potential application areas of MRS include: area coverage and patrolling missions [1], aerial/perimeter surveillance [2–4], search and rescue missions [5], searching operations [6], and situational awareness [7].

In the previously highlighted applications, the precise localization of the robotic members is a necessity in order to guarantee the success of such missions. Indeed, highly accurate localization can be achieved by equipping all the involved robots with the state-of-the-art sensory means. However, this increases the overall cost of such systems. Moreover, some robotic members, e.g. flying/hovering robots may not have enough payload/computational capacity to operate such means. Therefore, relative localization (RL) has been developed as a practical solution for effective and accurate execution of multi-robot collaborative tasks. The objective of relative localization is to detect and locate robots, with limited sensory capabilities (observed robots), with respect to another robots with accurate localization means (observing robots). This is achieved by using the relative observation between the two robots categories, see [8–12], for more details.

The most commonly used inter-robot relative measurements between a pair of robots in an MRS are the range and bearing measurements. Moreover, the extended Kalman filter (EKF) is the most commonly adopted nonlinear estimator for achieving the RL [9]. However, inappropriate initialization of EKF can generally lead to instability of the estimator as well as longer estimation settling time, i.e. the time required to reach acceptable estimation error levels. As a result, relative pose estimates become misleading. Indeed, this erroneous localization causes undesirable behaviour and possibly failure in collaborative missions. In order to avoid these issues, the majority of the past work assumed a known initial relative pose between two arbitrary robots, see, e.g. [13]. In addition to EKF, particle filter (PF) [14], unscented Kalman filter (UKF) [15], and pseudolinear Kalman filter [10, 16] have been also used to achieve RL in collaborative tasks.

The previously highlighted estimators, except the PF, use the Gaussian probability distribution to approximate a given state noise; the PF instead approximates the distribution via Monte Carlo sampling [17]. Moreover, in order to reduce the computational complexity, these filters employ Markov property, i.e the current state estimate is based only on the most recent measurement and the previous state estimate. In contrast to the previous methods, it is proposed in this chapter to employ a moving horizon estimation (MHE) to solve the RL problem. MHE considers the evolution of a constrained and possibly nonlinear model on a fixed time horizon, and minimizes the deviation of the model from a past window of measurements, see [18] for details. Therefore, MHE relaxes the Markov assumption mentioned above. Although MHE does not generally rely on any specific error distribution, tuning the estimator becomes easier when a probabilistic insight is considered, see [19, 20] for details, and [21] for a relative localization example using MHE in autonomous underwater vehicles.

Once the RL problem is accurately solved, the control objectives within a given MRS are classified under virtual structure, behaviour based, and leader follower approaches. In the virtual structure approach, the individual robots in the MRS formation are treated as points on a rigid body, where their individual motion are determined based on the formation overall-motion. In behaviour-based control, several desired behaviours are assigned to each robot in the MRS; the behaviours include formation keeping (stabilization), goal seeking, and obstacle avoidance. Finally, in the approach of the leader follower structure, a follower robot is assigned to follow the pose of a leader robot with an offset, see the overview articles [22–24] for more details.

The solution to the virtual structure approach is presented in [25], where authors designed an iterative algorithm that fits the virtual structure to the positions of the individual robots in the formation, then displaces the virtual structure in a prescribed direction and finally updates the robot positions. Solutions to the behaviour based approach include motor scheme control [26], null-space-based behaviour control [27], and social potential fields control [28]. The leader-follower control approach is achieved using feedback linearization [29], backstepping control [30], sliding mode control [31], and model predictive control [32].

In this chapter, based on relative localization achieved using MHE, it is proposed to use a centralized nonlinear model predictive control (MPC) to perform the relative trajectory tracking task, where one or more robots in a given MRS are commanded to follow time varying trajectories with respect to another robot reference frame. In addition, interrobot possible collisions must be avoided using the adopted controller. As this control problem can be hardly seen in the literature, it can be classified under the leader follower structure, with the exception that the relative references of the follower robots are time varying. In MPC, a cost function is formulated to minimize the tracking error for each of the follower robots. Then, an optimal control is obtained by solving a discrete nonlinear optimization problem over a pre-set prediction horizon [33].

Practical implementation of MHE and MPC requires online solutions of dynamic optimization problems, which can form a computational burden especially for systems with fast dynamics such as MRS's. In addition, classical implementation techniques of MHE and MPC to fast systems can cause a considerable time delay between taking measurements and applying new control commands. In order to overcome these problems, the real-time iteration (RTI) scheme originally proposed in [34] and implemented in an open-source software, i.e. ACADOtoolkit [35], has been used in this study. The used software exports tailored, efficient, and self contained C-code implementing fast MHE and MPC routines [20].

The majority of the work considering MRS in the literature treats the relative localization and control problems separately. In contrast, in this chapter, we combine the two objectives and adapt an optimization based solution for the relative localization and tracking control in multi-robot systems using fast MHE and MPC, respectively. Therefore, in this study, we give a particular attention to the practical implementation and evaluation of the adopted method to the considered system. Herein, the advantages of using MHE against EKF are presented first via numerical simulations. Then, the adopted estimation-control algorithm is evaluated, when applied to a given MRS. Finally, laboratory experiments were conducted to demonstrate a practical proof-of-concept of the used algorithm.

The chapter is organized as follows: in Section 6.3, an overview of multi-robot systems and the details of the relative localization and relative tracking control objectives are presented. A brief description of the employed estimation-control algorithm is provided in Section 6.4. Numerical simulation results are shown in Section 6.5 while experimental results are presented in Section 6.6. Finally, conclusions are drawn in Section 6.7.

## 6.3 Problem Setup

In this section, we provide an overview of the constructing elements in a typical multi-robot system (MRS), and introduce a detailed description of motion and measurement models involved in such systems. Moreover, an insightful explanation of relative localization and relative tracking control objectives is provided.

### 6.3.1 Multi-robot systems (MRS)

In a given MRS, robotic members, which include both aerial and ground robots, are classified into two classes, i.e. *observing robots* and *observed robots*. At a given time, a robot that takes inter-robot relative measurements of an arbitrary robot is defined as an observing robot while the robot that comes into the measurement range of the observing robot is referred to as an observed robot. For the considered MRS, the robotic members are assumed to possess the following characteristics [10–12, 36]:

i All robots in the MRS are equipped with dead-reckoning measurement means.

ii The observing robot is equipped with exteroceptive sensors in order to uniquely identify each of the observed robots in its field of view and measure their relative range and bearing.

iii The observing robot uses high precision sensors such that it can accurately perform its self-localization task while the observed robots have limited sensory/computational capabilities.

iv Robots can share data and information via communication devices and protocols.

v MRS' ground robots are assumed to navigate on flat surfaces. Moreover, the aerial robots perform hovering control such that the pitch and roll angles are kept stable via low velocity maneuvers, i.e. their 6-DOF are reduced to 4-DOF for simplicity.

Assume at a given time instant $k \in \mathbb{N}_0$, a team of robots contains $M \in \mathbb{N}$ observed robots and $L \in \mathbb{N}$ observing robots, where $M$ and $L$ are unknown to each robot, i.e. observing or observed, and also may change with time. In general, 6-DOF are required to describe a robot navigating a 3 dimensional space, see [37]; nonetheless, assumption (v) given above allows for the simplification of the 6-DOF into a 4-DOF [10]. Therefore, the relative pose of the $i^{th}$ observed robot in the $j^{th}$ observing robot's local Cartesian coordinate frame at a given time step $k$ is represented by the state vector $\mathbf{x}_i^j(k) \in \mathbb{R}^4$, $\mathbf{x}_i^j(k) = \left(x_i^j(k), y_i^j(k), z_i^j(k), \theta_i^j(k)\right)^\top$, where $\left(x_i^j(k), y_i^j(k), z_i^j(k)\right)^\top$ (m,m,m)$^\top$ is the relative position and $\theta_i^j(k)$ (rad) is the relative orientation, i.e. the relative yaw angle.

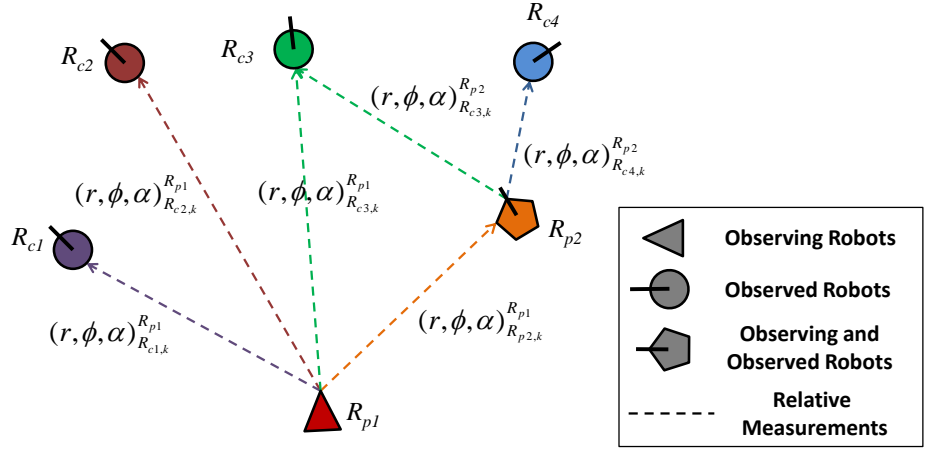

Figure 6.1: Configuration space of a multi-robot system at a given time instant.

Figure 6.1 shows a sample MRS configuration space at a given time instant $k$. Label $R_{pj}, j \in \{1, 2\}$ represents observing robots while label $R_{ci}, i \in \{1, 2, 3, 4\}$ represents observed robots. In this configuration, the relative measurements are indicated by the measurements vector $\mathbf{y}_i^j(k) \in \mathbb{R}^3$, $\mathbf{y}_i^j(k) = \left(r_i^j(k), \phi_i^j(k), \alpha_i^j(k)\right)^\top$ (m,rad,rad)$^\top$, which con-tains the relative range $r$, the relative azimuth angle $\phi$, and the relative elevation angle $\alpha$ of the $i^{th}$ robot measured by the $j^{th}$ robot[1]. As can be deduced from Figure 6.1, a robot can act simultaneously as an observing robot and an observed robot, e.g. $R_{p2}$.

### 6.3.2 Motion and measurement models

At a given time $t \geq 0$, the noise free continuous time motion-model of the relative-state $\mathbf{x}_i(t) = (x_i, y_i, z_i, \theta_i)^\top \in X \subset \mathbb{R}^4$ of the $i^{th}$ observed robot is given by [10, 12]

$$\dot{\mathbf{x}}_i(t) = \mathbf{f}_i\left(\mathbf{x}_i(t), \mathbf{u}_i(t), \mathbf{w}_i(t)\right)$$

$$\begin{pmatrix} \dot{x}(t) \\ \dot{y}(t) \\ \dot{z}(t) \\ \dot{\theta}(t) \end{pmatrix}_i = \begin{pmatrix} v_x(t)\cos(\theta(t)) - v_y(t)\sin(\theta(t)) - v_x^o(t) + y(t)\omega_z^o(t) \\ v_y(t)\cos(\theta(t)) + v_x(t)\sin(\theta(t)) - v_y^o(t) - x(t)\omega_z^o(t) \\ v_z(t) - v_z^o(t) \\ \omega_z(t) - \omega_z^o(t) \end{pmatrix}_i \tag{6.1}$$

where $\mathbf{f}_i(\cdot) : \mathbb{R}^4 \times \mathbb{R}^4 \times \mathbb{R}^4 \to \mathbb{R}^4$, is the analytic vector mapping function. $\mathbf{u}_i(t) = (v_{x,i}, v_{y,i}, v_{z,i}, \omega_{z,i})^\top \in U \subset \mathbb{R}^4$ is the control input velocities vector, of an observed robot, which contains the linear velocity components $(v_{x,i}, v_{y,i}, v_{z,i})^\top$ (m/s, m/s, m/s)$^\top$ and the angular velocity component $\omega_{z,i}$ (rad/s). $\mathbf{w}_i(t) = \left(v_{x,i}^o, v_{y,i}^o, v_{z,i}^o, \omega_{z,i}^o\right)^\top \in W \subset \mathbb{R}^4$ is the local velocities vector of an observing robot containing its linear and angular velocities, respectively. The sets $X, U$, and $W$ are specified by

$$X := [\underline{x}, \bar{x}] \times [\underline{y}, \bar{y}] \times [\underline{z}, \bar{z}] \times \mathbb{R},$$

$$U := [\underline{v}_x, \bar{v}_x] \times [\underline{v}_y, \bar{v}_y] \times [\underline{v}_z, \bar{v}_z] \times [\underline{\omega}_z, \bar{\omega}_z],$$

$$W := [\underline{v}_x^o, \bar{v}_x^o] \times [\underline{v}_y^o, \bar{v}_y^o] \times [\underline{v}_z^o, \bar{v}_z^o] \times [\underline{\omega}_z^o, \bar{\omega}_z^o].$$

---

[1]To keep the presentation simple later in the chapter, all terms denoted by $(\cdot)_i^j$ will be replaced by $(\cdot)_i$.

Since the rest of the chapter is based on the discrete time formulation, we define the *disturbed* Euler discrete model of (6.1), for a given sampling period $\delta > 0$ and a sampling instant $k \in \mathbb{N}_0$, as

$$\mathbf{x}_i(k+1) = \mathbf{f}_{\delta,i}(\mathbf{x}_i(k), \mathbf{u}_i(k), \mathbf{w}_i(k)) + \boldsymbol{\nu}_{\mathbf{x}} \tag{6.2}$$

where $\mathbf{f}_{\delta,i}(\cdot) : \mathbb{R}^4 \times \mathbb{R}^4 \times \mathbb{R}^4 \to \mathbb{R}^4$, is the discrete time dynamics given by

$$\mathbf{f}_{\delta,i}(\cdot) = \begin{pmatrix} x(k) \\ y(k) \\ z(k) \\ \theta(k) \end{pmatrix}_i + \delta \cdot \begin{pmatrix} v_x(k)\cos(\theta(k)) - v_y(k)\sin(\theta(k)) - v_x^o(k) + y(k)\omega_z^o(k) \\ v_y(k)\cos(\theta(k)) + v_x(k)\sin(\theta(k)) - v_y^o(k) - x(k)\omega_z^o(k) \\ v_z(k) - v_z^o(k) \\ \omega_z(k) - \omega_z^o(k) \end{pmatrix}_i . \tag{6.3}$$

Moreover, $\boldsymbol{\nu}_{\mathbf{x}}$ is a zero mean additive white Gaussian noise vector, which represents model uncertainties. To this end, ground observed robots are modeled by setting $z = v_z = 0$ in (6.3).

In terms of the relative Cartesian coordinates of an observed robot $(x_i, y_i, z_i)$, the general 3-D relative observation model is defined by $\mathbf{h}(\mathbf{x}_i(k))$ and given by (6.4).

$$\mathbf{y}_i(k) = \mathbf{h}(\mathbf{x}_i(k)) + \boldsymbol{\nu}_{\mathbf{y}}$$

$$\begin{pmatrix} r(k) \\ \phi(k) \\ \alpha(k) \end{pmatrix}_i = \begin{pmatrix} \sqrt{x(k)^2 + y(k)^2 + z(k)^2} \\ \arctan\left(\dfrac{y(k)}{x(k)}\right) \\ \arctan\left(\dfrac{z(k)}{\sqrt{x(k)^2 + y(k)^2}}\right) \end{pmatrix}_i + \begin{pmatrix} \nu_r \\ \nu_\phi \\ \nu_\alpha \end{pmatrix}, \tag{6.4}$$

where the measurement vector $\mathbf{y}_i = (r_i, \phi_i, \alpha_i)^\top \in Y \subset \mathbb{R}^3$ contains the relative range,

relative azimuth angle, and relative elevation angle, respectively. Parameters $\nu_r, \nu_\phi$ and $\nu_\alpha$ are zero mean additive Gaussian noises with standard deviations $\sigma_r, \sigma_\phi$, and $\sigma_\alpha$, respectively.

As can be noticed from the measurement model (6.4), the relative position components $(x, y, z)$ are indirectly measured by the model; however, the relative orientation, i.e. $\theta$, is neither directly nor indirectly involved in the measurement model. The full observability of the motion-measurement models pair, i.e. (6.2) and (6.4), is guaranteed provided that the observed robot linear velocity $v_x$ is non-zero, see [10, 37] for details.

### 6.3.3 Estimation-control objective

The main objectives for the MRS considered in this study can be summarized by the following two iterative tasks for all closed-loop instants $n \in \mathbb{N}_0$:

**a) Localize:** using the odometery measurements $\tilde{\mathbf{u}}_i(n) \in U$ and $\tilde{\mathbf{w}}_i(n) \in W$, as well as the sensory readings $\tilde{\mathbf{y}}_i(n) \in Y$, the localization task finds the best estimate $\hat{\mathbf{x}}_i(n) \in X$ of the relative-pose state vector $\mathbf{x}_i(n)$ for each observed robot, i.e. for all $i \in \{1, 2, \ldots, M\}$.

**b) Control:** based on the first task, the control task is to calculate a feedback control law $\boldsymbol{\mu}_i(\hat{\mathbf{x}}_i(n)) : \mathbb{R}^4 \to \mathbb{R}^4$, such that the $i^{th}$ observed robot is tracking a *reference* trajectory with respect to another observing robot with the reference relative state vector $\mathbf{x}_i^r(n) \in \mathbb{R}^4$ and reference speeds $\mathbf{u}_i^r(n) \in \mathbb{R}^4$. Moreover, while robots are following their relative reference trajectories, they must not be jeopardized by collision with other robots in the given MRS.

Traditional techniques solving the localization task are based on EKF. However, due to linearizations of models (6.2) and (6.4) involved in EKF implementation, loss of information and bias problems arise [16]. Moreover, proper initialization of the *estimated state* in EKF is crucial because any misleading initialization can lead to unstable estimation and/or longer estimation settling time, see, e.g., [9, 36], for more details. To overcome these problems, MHE is proposed as the main estimator to satisfy the first task, i.e. RL. To achieve the second task (control), MPC is proposed. The main reason for choosing MPC is its ability to handle constrained nonlinear systems, e.g. (6.2), in a straight forward manner. Thus, the collision avoidance task of the controller is implicitly included in its associated optimal control problem as well be seen in the following section, see also [23, 24, 38] for more details.

## 6.4 Estimator-Controller Synthesis

In this section, the mathematical formulation of MHE and MPC for the localization and tracking control objectives is presented first. Then, a brief description of the code generation toolkit used is given.

### 6.4.1 MHE formulation

For each observed robot in a given MRS, the MHE relative-state estimator is formulated as a least squares (LSQ) cost function $J_{N_E} : X \times U^{N_E} \times W^{N_E} \to \mathbb{R}_{\geq 0}$ shown in (6.5). Hence, using (6.5), the MHE repeatedly solves the constrained nonlinear dynamic optimization problem (6.6) over a fixed estimation horizon of length $N_E \in \mathbb{N}$ [20]. In (6.6), for $q \in \mathbb{N}_0$, the control sequences $\mathbf{u}_i$ and $\mathbf{w}_i$ are defined by

$$\mathbf{u}_i = (\mathbf{u}_i(q - N_E), \mathbf{u}_i(q - N_E + 1), \ldots, \mathbf{u}_i(q - 1)) \in U^{M \cdot N_E},$$
$$\mathbf{w}_i = (\mathbf{w}_i(q - N_E), \mathbf{w}_i(q - N_E + 1), \ldots, \mathbf{w}_i(q - 1)) \in U^{L \cdot N_E}.$$

$$J_{N_E}(\mathbf{x}_i(q - N_E), \mathbf{u}_i, \mathbf{w}_i) = \|\mathbf{x}_i(q - N_E) - \mathbf{x}^{\text{est}}_{i,q-N_E}\|^2_{\mathbf{A}_i}$$

$$+ \sum_{k=q-N_E}^{q} \|\mathbf{h}(\mathbf{x}_i(k)) - \tilde{\mathbf{y}}_i(k)\|^2_{\mathbf{B}_i} + \sum_{k=q-N_E}^{q-1} \left\| \begin{matrix} \mathbf{u}_i(k) - \tilde{\mathbf{u}}_i(k) \\ \mathbf{w}_i(k) - \tilde{\mathbf{w}}_i(k) \end{matrix} \right\|^2_{\mathbf{C}_i} \quad (6.5)$$

$$\min_{\mathbf{x}_i(q-N_E), \mathbf{u}_i \in \mathbb{R}^{4 \times N_E}, \mathbf{w}_i \in \mathbb{R}^{4 \times N_E}} J_{N_E}(\mathbf{x}_i(q - N_E), \mathbf{u}_i, \mathbf{w}_i) \quad (6.6)$$

subject to:

$$\mathbf{x}_i(k + 1) = \mathbf{f}_{\delta,i}(\mathbf{x}_i(k), \mathbf{u}_i(k), \mathbf{w}_i(k)) \quad \forall k \in \{q - N_E, \ldots, q - 1\},$$

$$\mathbf{x}_i(k) \in X \quad \forall k \in \{q - N_E, \ldots, q\},$$

$$\mathbf{u}_i(k) \in U \quad \forall k \in \{q - N_E, \ldots, q - 1\},$$

$$\mathbf{w}_i(k) \in W \quad \forall k \in \{q - N_E, \ldots, q - 1\}.$$

In (6.5), $\tilde{\mathbf{y}}_i$, $\tilde{\mathbf{u}}_i$, and $\tilde{\mathbf{w}}_i$ denote the actually measured system outputs and inputs, respectively. The first term in (6.5) is known as the arrival cost and it penalizes the deviation of the first state in the moving horizon window and its priori estimate $\mathbf{x}^{\text{est}}_{i,q-N_E}$ by the diagonal positive definite matrix $\mathbf{A}_i \in \mathbb{R}^{4\times 4}$. Normally, the estimate $\mathbf{x}^{\text{est}}_{i,q-N_E}$ is adopted from the previous MHE estimation step. Moreover, the weighting matrix $\mathbf{A}_i$ is chosen as a smoothed EKF-update based on sensitivity information gathered while solving the most recent MHE step. Therefore, $\mathbf{x}^{\text{est}}_{i,q-N_E}$ and $\mathbf{A}_i$ have to be only initialized, see [20, 39] for more details.

The second term in (6.5) penalizes the change in the system predicted outputs $\mathbf{h}(\mathbf{x}_i)$ from the actually measured outputs $\tilde{\mathbf{y}}_i$ by the diagonal positive-definite matrix $\mathbf{B}_i \in \mathbb{R}^{3\times 3}$. Similarly, the change in the applied control inputs ($\mathbf{u}_i$, and $\mathbf{w}_i$) from the measured inputs ($\tilde{\mathbf{u}}_i$, and $\tilde{\mathbf{w}}_i$) is penalized using the diagonal positive-definite matrix $\mathbf{C}_i \in \mathbb{R}^{8\times 8}$. The latter term is included in the cost function (6.5) to account for actuator noise and/or inaccuracy, see [19] for details. $\mathbf{B}_i$ and $\mathbf{C}_i$ are chosen to match the applied motion and measurement noise covariances. Furthermore, all estimated quantities in (6.6) are subject to bounds, which signify the system physical limitations. The solution of the optimization problem (6.6) leads mainly to the relative-state estimate sequence $\hat{\mathbf{x}}_i(k), k = (q - N_E, \cdots, q)$, where $\hat{\mathbf{x}}_i(0) := \hat{\mathbf{x}}_i(q)$ denotes the current estimate of the relative state vector of a given observed robot; $\hat{\mathbf{x}}_i(0)$ is later used as a measurement feedback to the MPC controller presented in the following subsection. Moreover, estimate sequences for the actually applied control inputs, i.e. $\hat{\mathbf{u}}_i$ and $\hat{\mathbf{w}}_i$, are byproducts of the MHE estimator.

### 6.4.2 MPC formulation

Motivated by the assumptions presented in Section 6.3.1, in particular that the observing robots have higher computational capabilities than the observed robots, a centralized MPC controller is proposed to achieve the relative trajectory tracking task. To this end, we define an MRS overall augmented nominal model, for a sampling instant $k \in \mathbb{N}_0$, as

$$\mathbf{x}_a(k+1) = \mathbf{f}_{\delta,a}(\mathbf{x}_a(k), \mathbf{u}_a(k), \mathbf{w}_a(k)), \tag{6.7}$$

where $\mathbf{x}_a = (\mathbf{x}_1^\top, \mathbf{x}_2^\top, \ldots, \mathbf{x}_M^\top)^\top \in X^M$ is the augmented relative state vector of $M$ observed robots; $\mathbf{u}_a = (\mathbf{u}_1^\top, \mathbf{u}_2^\top, \ldots, \mathbf{u}_M^\top)^\top \in U^M$ is the observed robots' augmented controls; $\mathbf{w}_a = (\mathbf{w}_1^\top, \mathbf{w}_2^\top, \ldots, \mathbf{w}_L^\top)^\top \in W^L$ is the augmented controls of the observing robots; and $\mathbf{f}_{\delta,a} : \mathbb{R}^{4 \cdot M} \times \mathbb{R}^{4 \cdot M} \times \mathbb{R}^{4 \cdot L} \to \mathbb{R}^{4 \cdot M}$ is the augmented nonlinear mapping given as

$$\mathbf{f}_{\delta,a} = \left(\mathbf{f}_{\delta,1}^\top, \mathbf{f}_{\delta,2}^\top, \ldots, \mathbf{f}_{\delta,M}^\top\right)^\top,$$

where $\mathbf{f}_{\delta,i}, i \in \{1, 2, \ldots, M\}$ are given by (6.3).

To this end, we define the augmented *reference* trajectories and the augmented *refer-*

*ence* speeds for the observed robots as

$$\mathbf{x}_a^r = \left(\mathbf{x}_1^{r\top}, \mathbf{x}_2^{r\top}, \ldots, \mathbf{x}_M^{r\top}\right)^\top \in X^M \text{ and } \mathbf{u}_a^r = \left(\mathbf{u}_1^{r\top}, \mathbf{u}_2^{r\top}, \ldots, \mathbf{u}_M^{r\top}\right)^\top \in U^M,$$

respectively.

Similar to the MHE formulation, the MPC controller is formulated using a least squares (LSQ) cost function $J_{N_C} : X^M \times U^{N_C \cdot M} \to \mathbb{R}_{\geq 0}$ penalizing the future deviation of the observed robots relative state vectors and control inputs from their time varying references as shown in (6.8).

$$\begin{aligned} J_{N_C}(\mathbf{x}_a(0), \mathbf{u}_a) &= \|\mathbf{x}_a(N_C) - \mathbf{x}_a^r(N_C)\|_{\mathbf{P}_a}^2 \\ &+ \sum_{k=0}^{N_c-1} \left(\|\mathbf{x}_a(k) - \mathbf{x}_a^r(k)\|_{\mathbf{Q}_a}^2 + \|\mathbf{u}_a(k) - \mathbf{u}_a^r(k)\|_{\mathbf{R}_a}^2\right) \end{aligned} \quad (6.8)$$

For each sampling instant ($n = 0, 1, \cdots$), MPC solves repeatedly the nonlinear dynamic optimization problem (6.9) over a future prediction horizon specified by $N_C \in \mathbb{N}$ [40]. In (6.9), the control sequence $\mathbf{u}_a$ is defined by

$$\mathbf{u}_a = (\mathbf{u}_a(0), \mathbf{u}_a(1), \cdots, \mathbf{u}_a(N_C - 1))_n \in U^{M \cdot N_C}.$$

The first term of the cost function (6.8) accounts for the terminal state penalty. This is achieved by using the diagonal positive-definite matrix $\mathbf{P}_a \in \mathbb{R}^{4 \cdot M \times 4 \cdot M}$, which penalizes the norm of the augmented terminal state error. This matrix is defined as $\mathbf{P}_a = \text{diag}(\mathbf{P}_1, \mathbf{P}_2, \ldots, \mathbf{P}_M)$, where $\mathbf{P}_i \in \mathbb{R}^{4 \times 4}$ is the terminal state penalty of the $i^{th}$ observed robot for all $i \in \{1, 2, \ldots, M\}$. The remaining terms in (6.8) account for the running costs of the future deviation between the robot states and controls, and their references.

$$\min_{\mathbf{u}_a \in \mathbb{R}^{2 \cdot M \times N_C}} J_{N_C}(\mathbf{x}_a(0), \mathbf{u}_a) \tag{6.9}$$

$$\text{subject to: } \mathbf{x}_a(0) = \hat{\mathbf{x}}_a(0),$$

$$\mathbf{x}_a(k+1) = \mathbf{f}_{\delta,a}(\mathbf{x}_a(k), \mathbf{u}_a(k), \mathbf{w}_a(k)), \quad \forall k \in \{0, 1, \ldots, N_C - 1\},$$

$$G(\mathbf{x}_a(k)) \leq 0 \quad \forall k \in \{1, 2, \ldots, N_C\},$$

$$\mathbf{x}_a(k) \in X^M \quad \forall k \in \{1, 2, \ldots, N_C\},$$

$$\mathbf{u}_a(k) \in U^M \quad \forall k \in \{0, 1, \ldots, N_C - 1\}.$$

Indeed, the reference control values, i.e. $\mathbf{u}^r_{\{\cdot\}}$, represent the desired speeds along the reference trajectories.

The norms of the running costs are weighted with the positive definite matrices $\mathbf{Q}_a \in \mathbb{R}^{4 \cdot M \times 4 \cdot M}$ and $\mathbf{R}_a \in \mathbb{R}^{4 \cdot M \times 4 \cdot M}$, where $\mathbf{Q}_a = \text{diag}(\mathbf{Q}_1, \mathbf{Q}_2, \ldots, \mathbf{Q}_M)$ and $\mathbf{R}_a = \text{diag}(\mathbf{R}_1, \mathbf{R}_2, \ldots, \mathbf{R}_M)$, respectively. $\mathbf{Q}_i \in \mathbb{R}^{4 \times 4}$ and $\mathbf{R}_i \in \mathbb{R}^{4 \times 4}, i \in \{1, 2, \ldots, M\}$ are the individual running costs weighing matrices defined for each observed robot. The weighting matrices $\mathbf{P}_a$, $\mathbf{Q}_a$, and $\mathbf{R}_a$ are chosen based on a simulation tuning procedure [19, 39].

As can be noticed from (6.9), the initial condition of the prediction of $\mathbf{x}_a$ is set to the augmented estimate vector $\hat{\mathbf{x}}_a(0) = (\hat{\mathbf{x}}_1^\top(0), \hat{\mathbf{x}}_2^\top(0), \ldots, \hat{\mathbf{x}}_M^\top(0))^\top \in \mathbb{R}^{4 \cdot M}$ resulting from solving problem (6.6) for all observed robots. Moreover, the collision avoidance constraints at instant $k$ of (6.9), i.e. $G(\mathbf{x}_a(k))$, is defined as

$$G(\mathbf{x}_a(k)) = (g_1(\mathbf{x}_a(k)), g_2(\mathbf{x}_a(k)))^\top.$$

$g_1(\mathbf{x}_a(k))$ accounts for the collision avoidance between the observed robots and, for all

$j \in \{1, 2, \ldots, M\}$ and $k \in \{1, 2, \ldots, N_C\}$, is defined as

$$g_1(\mathbf{x}_a(k)) = \Big(g_1(\mathbf{x}_a(k))|_{1,j}, \ldots, g_1(\mathbf{x}_a(k))|_{i-1,j}, g_1(\mathbf{x}_a(k))|_{i+1,j}, \ldots, g_1(\mathbf{x}_a(k))|_{M,j}\Big),$$
(6.10)

where

$$g_1(\mathbf{x}_a(k))|_{i,j} := 1 - \frac{\|(x_i(k), y_i(k), z_i(k)) - (x_j(k), y_j(k), z_j(k))\|^2}{4r_c^2}, \qquad i \neq j.$$

On the other hand, $g_2(\mathbf{x}_a(k))$ accounts for the collision avoidance between the observed robots and the observing robots and, for all $j \in \{1, 2, \ldots, L\}$ and $k \in \{1, 2, \ldots, N_C\}$, is defined as

$$g_2(\mathbf{x}_a(k)) = \Big(g_2(\mathbf{x}_a(k))|_{1,j}, \ldots, g_2(\mathbf{x}_a(k))|_{i,j}, \ldots, g_1(\mathbf{x}_a(k))|_{M,j}\Big).$$
(6.11)

where

$$g_2(\mathbf{x}_a(k))|_{i,j} = 1 - \frac{\|(x_i(k), y_i(k), z_i(k)) - (x_j(k), y_j(k), z_j(k))\|^2}{(r_c + r_p)^2}.$$

$r_c$ is the minimum Euclidean distance that all the robots in the given MRS must maintain from any observed robot while $r_p$ is the minimum Euclidean distance that all the robots in the given MRS must maintain from any observing robot.

The solution of problem (6.9), at a sampling instant $n \in \mathbb{N}_0$, leads to a minimizing control sequence

$$\mathbf{u}_a^\star|_n = (\mathbf{u}_a^\star(0), \mathbf{u}_a^\star(1), \cdots, \mathbf{u}_a^\star(N_C - 1))_n \in U^{M \cdot N_C}$$

where $\mathbf{u}_a^\star(0)|_n = \Big(\mathbf{u}_1^\star(0)^\top, \mathbf{u}_2^\star(0)^\top, \ldots, \mathbf{u}_M^\star(0)^\top\Big)_n^\top \in U^M$. Thus, the feedback control

law, for each observed robot in the given MRS, is defined by $\boldsymbol{\mu}_i(\hat{\mathbf{x}}_i(n)) = \mathbf{u}_i^\star(0)|_n$, $i \in \{1, 2, \ldots, M\}$. It has to be mentioned here that the observing robots speeds $\mathbf{w}_a$ are considered as constant parameters when solving (6.9). Algorithm 6.1 summarizes the overall estimation-control scheme, for a given MRS in which the term estimation-control refers to the main two iterative steps of the algorithm.

---

**Algorithm 6.1** Estimation-control scheme in MRS based on MHE and MPC

---
1: Set the number of the observed robots $M$
2: **for** every sampling instant $t_n, n = 0, 1, 2, \ldots$ **do**
3:     **for** each observed robot $R_{ci}, i = 1, 2, \ldots, M$ **do**
4:         Get the past $N_E$ odometery measurements, i.e. $\tilde{\mathbf{u}}_i$ and $\tilde{\mathbf{w}}_i$.
5:         Get the past $N_E + 1$ relative measurements, i.e. $\tilde{\mathbf{y}}_i$.
6:         Solve the optimization problem (6.6) over estimation horizon $N_E$, and find the best current estimate of the relative state $\mathbf{x}_i(0)|_n$, i.e. $\hat{\mathbf{x}}_i(0)|_n$.
7:     **end for**
8:     Set the vector $\hat{\mathbf{x}}_a(0)|_n = \left(\hat{\mathbf{x}}_1^\top(0), \hat{\mathbf{x}}_2^\top(0), \ldots, \hat{\mathbf{x}}_M^\top(0)\right)_n^\top$.
9:     Solve the optimization problem (6.9) over prediction horizon $N_C$, and find the minimizing control sequence $\mathbf{u}_a^\star|_n = (\mathbf{u}_a^\star(0), \mathbf{u}_a^\star(0), \cdots, \mathbf{u}_a^\star(N_C - 1))_n$.
10:     Extract $\mathbf{u}_a^\star(0)|_n = \left(\mathbf{u}_1^\star(0)^\top, \mathbf{u}_2^\star(0)^\top, \ldots, \mathbf{u}_M^\star(0)^\top\right)_n^\top$ from $\mathbf{u}_a^\star|_n$.
11:     Define the MPC-feedback control law $\boldsymbol{\mu}_i(\hat{\mathbf{x}}_i(n)) = \mathbf{u}_i^\star(0)|_n$, $i \in \{1, 2, \ldots, M\}$.
12:     **for** each observed robot $R_{ci}, i = 1, 2, \ldots, M$ **do**
13:         Apply $\boldsymbol{\mu}_i(\hat{\mathbf{x}}_i(n))$ to the observed robot $R_{ci}$.
14:     **end for**
15: **end for**

---

### 6.4.3 ACADO code generation

The open-source ACADO Code Generation toolkit [40, 41] is used in the implementation of the adopted algorithm. Based on a symbolic representation of MHE and MPC problems, ACADO generates an optimized and self-contained C code, which limits the associated computations to the most essential steps. This is done by examining the structure of the defined optimization problems, and avoiding the dynamic memory allocation; thus, eliminating any unnecessary and/or redundant computations [35, 41]. The exported C code is based on the Direct Multiple Shooting method, utilizing a condensing algorithm, and a real-time iteration (RTI) scheme. The RTI scheme reduces the time delay of the optimization problems (6.6) and (6.9) by *preparing* the associated computations in advance even before the measurements/feedback become available, see [35, 40, 41], for details.

Using ACADO code generation tool has been recently reported in several studies for different applications, e.g. autonomous driving [19], autonomous tractor–trailer systems [39, 42, 43], and wind energy systems [44, 45]. In this chapter, the use of fast MHE and MPC schemes is extended to MRS's relative localization and tracking control.

In the following two sections, the generated code has been tailored to solve the minimization problems (6.6) and (6.9) based on the Gauss Newton method, a real time iteration (RTI) scheme, a multiple-shooting discretization, and an Euler integrator for state prediction. The number of Gauss Newton iterations was limited to 10 while the Karush-Kuhn-Tucker (KKT) tolerance (a first order optimality condition) was set to $10^{-4}$. This particular choice provided a compromise between the execution time and the convergence of Algorithm 6.1 when implemented to the considered system, see [40, 41], for details. Indeed, an integrator with more than one integration step can be used, e.g. Runge-Kutta integrator, but the Euler integrator was chosen mainly because it is faster in state integration. Moreover, its performance was satisfactory for the considered problem. It has to be mentioned here that the same state integrator was used in the implementation of EKF in Section 6.5.1.

## 6.5 Simulation Results

In this section, two simulation studies are presented. First, a comparison simulation between MHE and EKF, which shows the main advantages of using MHE over EKF in RL. Second, a closed-loop simulation to evaluate Algorithm 6.1 when applied to a given MRS.

### 6.5.1 Comparison between MHE and EKF for RL

The performance of MHE against EKF is presented by a series of Monte Carlo numerical simulations. At a given sampling instant $n \in \mathbb{N}_0$, EKF finds the current estimate of an observed robot, i.e. $\hat{\mathbf{x}}_i(n), i \in \{1, 2, \ldots, M\}$ through two successive steps, i.e. a prediction step and a measurement update step, see [46] for details.

The conducted simulations here compromised a non-holonomic MRS with a ground observing robot and three aerial observed robots, i.e. $L = 1$ and $M = 3$. The observed robots are commanded to navigate along trajectories having different elevations. The state and control sets are selected as

$$X := [0, 6] \times [-3, 3] \times [0, 1] \times \mathbb{R}, \tag{6.12}$$
$$U := [-0.25, 0.6] \times [0, 0] \times [-0.25, 0.6] \times [-0.7, 0.7],$$
$$W := [-0.25, 0.6] \times [0, 0] \times [-0.25, 0.6] \times [-0.7, 0.7].$$

Note that, as we are considering only non-holonomic robots, all linear speeds in $y$-direction are set to 0. Moreover, the limits of the control sets are chosen based on the actuator limits of the robots used in real-time experiments.

In order to examine the effect of the sensory precision on the accuracy of the estimation, three noise configurations for $\boldsymbol{\nu}_\mathbf{y}$ in (6.4) and shown in Table 6.1, are used [9, 10]. Under all sensory noise configurations, the odometery measurements noise standard deviations specified by (6.13) are also used for all robots in the MRS.

$$\begin{aligned}
\sigma_{v_x} = \sigma_{v_y} = \sigma_{v_z} &= 0.01 \text{ m/s}, \\
\sigma_{v_x^o} = \sigma_{v_y^o} = \sigma_{v_z^o} &= 0.01 \text{ m/s}, \\
\sigma_{\omega_z} = \sigma_{\omega_z^o} &= 0.0873 \text{ rad/s}.
\end{aligned} \tag{6.13}$$

Table 6.1: Sensory noise configurations.

| Case | Relative measurements noise configuration |
|------|-------------------------------------------|
| 1 | $\sigma_r = 0.0068$ (m), $\sigma_\phi = \sigma_\alpha = 0.0036$ (rad) |
| 2 | $\sigma_r = 0.0167$ (m), $\sigma_\phi = \sigma_\alpha = 0.0175$ (rad) |
| 3 | $\sigma_r = 0.1466$ (m), $\sigma_\phi = \sigma_\alpha = 0.1000$ (rad) |

Following [20, 39], the estimation horizon length of MHE is chosen as ($N_E = 30$) using a trial and error method with the following tuning matrices

$$\mathbf{B} = \text{diag}(\sigma_r, \sigma_\theta, \sigma_\alpha)^{-1}, \text{ and}$$

$$\mathbf{C} = \text{diag}(\sigma_{v_x}, \sigma_{v_y}, \sigma_{v_z}, \sigma_{\omega_z}, \sigma_{v_x^o}, \sigma_{v_y^o}, \sigma_{v_z^o}, \sigma_{\omega_z^o})^{-1}.$$

The comparison between MHE and EKF is investigated through 15 Monte-Carlo simulations running at 10 Hz, i.e. $\delta = 0.1$ (seconds). The simulations presented here are conducted using the code generation tool of ACADO for MHE [20]. The generated C code has been compiled into a MEX function and all simulations were run using MATLAB.

The performance of the two estimators is evaluated for all observed robots by computing the root-mean-square-error (RMSE) of the relative position and relative orientation; thus, the conducted simulations are equivalent to 45 Monte-Carlo simulations. In order to meet the practical requirements of state estimation, the adopted estimators do not assume any knowledge of observed robots *initial* positions and/or orientations. Nonetheless, by inspecting the measurement model (6.4), the initial guess of the observed robots positions can be set using (6.14) for all $i \in \{1, 2, 3\}$.

$$\begin{aligned}
x_{i,o}^- &= \tilde{r}_{i,o} \cos(\tilde{\alpha}_{i,o}) \cos(\tilde{\phi}_{i,o}), \\
y_{i,o}^- &= \tilde{r}_{i,o} \cos(\tilde{\alpha}_{i,o}) \sin(\tilde{\phi}_{i,o}), \\
z_{i,o}^- &= \tilde{r}_{i,o} \sin(\tilde{\alpha}_{i,o}),
\end{aligned} \quad (6.14)$$

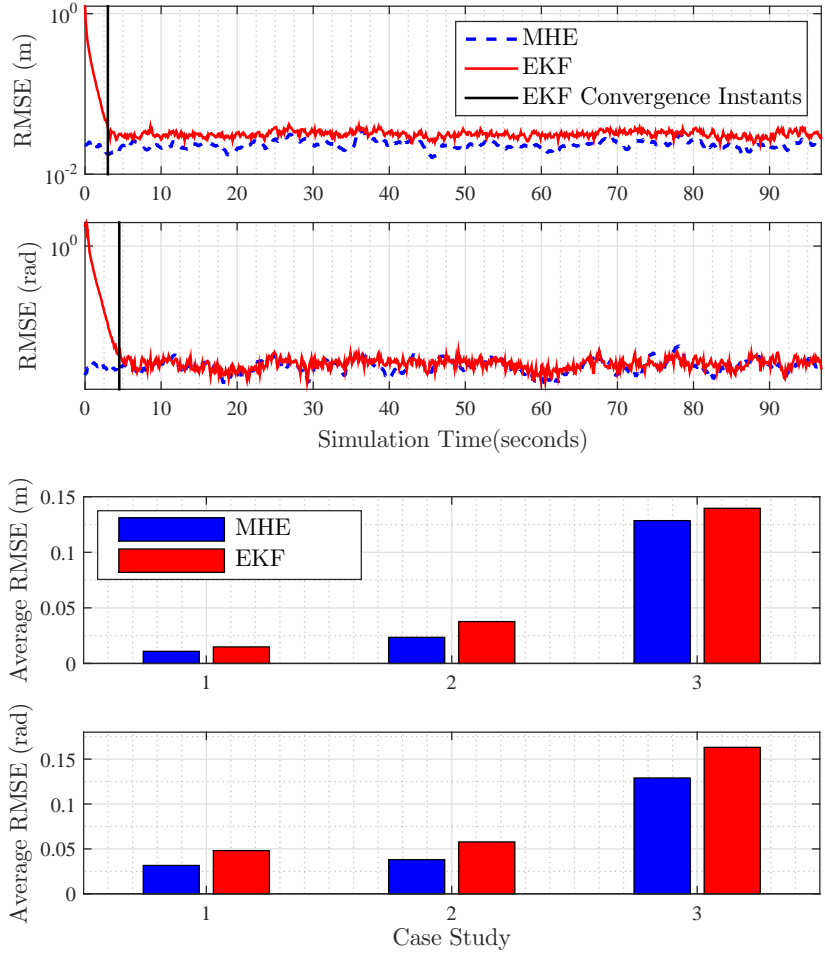

Figure 6.2: MHE and EKF performance for RL: Top: relative position and relative orientation estimation RMSE for case 2 shown in Table 6.1. Bottom: Relative position and relative orientation estimation *average* RMSE for all sensory noise configurations cases shown in Table 6.1.

where $x^-_{i,o}$, $y^-_{i,o}$ and $z^-_{i,o}$ are the initial guess of an observed robot position, and $\tilde{r}_{i,o}$, $\tilde{\phi}_{i,o}$ and $\tilde{\alpha}_{i,o}$ are the corresponding initial relative measurements. Additionally, the initial guess of an observed robot orientation is set to zero. The simulation results of the two estimators are summarized in Figure 6.2. Under the sensory noise configuration of Table 6.1 (case 2), the RMSE of the relative Euclidean position and relative orientation are presented in

Figure 6.2 (top); the RMSE of the relative position is calculated by

$$\left(\frac{1}{M \cdot n_{mc}} \sum_{i=1}^{M \cdot n_{mc}} \left(\hat{x}_e^2 + \hat{y}_e^2 + \hat{z}_e^2\right)\right)^{\frac{1}{2}},$$

where $\hat{x}_e, \hat{y}_e$ and $\hat{z}_e$ are the relative state position error Cartesian components for all observed robots, and $n_{mc} = 15$ is the number of Monte-Carlo simulations. As can be noticed from Fig 6.2 (top), MHE converges to acceptable estimation error limits immediately after the first estimation step is performed; however, EKF takes longer number of iterations until it reaches the same estimation error, i.e. 30 iterations (3 seconds) for relative position estimation error convergence and 45 iterations (4.5 seconds) for relative orientation estimation error convergence. Moreover, a more accurate relative positional estimation is observed under MHE while the relative orientation estimation accuracy of the two estimators are comparable. Indeed, it is crucial to achieve the state estimation with acceptable error limits as fast as possible in order to guarantee a satisfactory controller performance; thus, for RL, MHE outperforms EKF in this regard. Additionally, under the three cases of noise configurations shown in Table 6.1, the *average* estimation RMSE for relative position and orientation are presented in Figure 6.2 (bottom), which shows an increase in the estimation error as the uncertainty of the relative measurements increases. Nonetheless, in all cases, MHE showed more accurate estimation results when compared with EKF.

The conducted simulations shows MHE advantages over EKF in achieving a more accurate state estimation with fast convergence, despite the MHE relatively high computational intensity. Although the slow convergence of EKF can be treated by designing open-loop controls for the initial period of the control process to safely wait for the estimator convergence, EKF has additionally shown, in the literature, unsatisfactory performance in the existence of measurements' outliers when compared with MHE, see, e.g. [47, 48], for several numerical examples. Therefore, the high accuracy and the fast convergence of MHE, as well as the easy incorporation of state constraints in MHE formulation as shown in (6.6), suggests its use as a solution for relative state estimation in multi-robot systems.

### 6.5.2 Estimation and control using MHE and MPC

Algorithm 6.1 is evaluated numerically when applied to a team of non-holonomic robots with a ground observing robot and three aerial observed robots, i.e. $L = 1, M = 3$. Here, the observed robots are required to track reference trajectories in the observing robot frame of reference, which can be moving; such a mission can be seen in a number applications, e.g. area coverage, patrolling missions, searching operations, and perimeter surveillance, see, e.g. [1, 3, 6, 49], for more details. To this end, two case studies are considered depending on the observing robot speeds:

(a) The observing robot is moving in a straight line, i.e. $v_x^o = 0.1$ (m/s) and $\omega_z^o = 0$ (rad/s),

(b) The observing robot is moving in a circle, i.e. $v_x^o = 0.1$ (m/s) and $\omega_z^o = -0.05$ (rad/s).

See Table 6.2 for the detailed *relative* reference trajectories, of the observed robots, under these two cases. The simulations running time is set to 75 (seconds) with a sampling time $\delta = 0.2$ (seconds). Similar to the simulations presented in Section 6.5.1, the observed robots have odometery measurements noise standard deviations given by (6.13). In addition, the observing robot is adopting a relative measurement means with sensory noise configuration specified by Table 6.1 (case 2).

State and control limits used here are as given by (6.12). Indeed, the selected relative reference trajectories (Table 6.2) are chosen such that the corresponding reference speeds are not violating the limits (6.12). Moreover, to inspect the controller ability in fulfilling the collision avoidance task, the selected reference trajectories are designed such that there are locations along the mission, where at least two observed robots share the same reference position at the same time. The initial postures of all observed robots are also set such that a large initial tracking error is imposed on the controller.

Table 6.2: Simulation scenarios reference trajectories

|  | Case 1: $v_x^o = 0.1$ (m/s), $\omega_z^o = 0.0$ (rad/s). | Case 2: $v_x^o = 0.1$ (m/s), $\omega_z^o = -0.05$ (rad/s). |
|---|---|---|
| $R_{ci}$ | | |
| $R_{c1}$ | $x_1^r = 4 + 1.5\sin(0.22t)$ $y_1^r = 1.5\cos(0.22t)$ $z_1^r = 0.5$ | $x_1^r = 4 + 1.5\sin(0.2t)$ $y_1^r = 1.5\cos(0.2t)$ $z_1^r = 0.5$ |
| $R_{c2}$ | $x_2^r = 2 + 1.5\sin(0.13t)$ $y_2^r = -1 + 1.5\cos(0.13t)$ $z_2^r = 0.5$ | $x_2^r = 2 + 1.5\sin(0.1t)$ $y_2^r = -1 + 1.5\cos(0.1t)$ $z_2^r = 0.5$ |
| $R_{c3}$ | $x_3^r = 2 + 1.5\sin(0.14t)$ $y_3^r = 0.5 + 1.5\cos(0.14t)$ $z_3^r = 0.5$ | $x_3^r = 2 + 1.5\sin(0.12t)$ $y_3^r = 0.5 + 1.5\cos(0.12t)$ $z_3^r = 0.5$ |

For MHE, the estimation horizon is selected as ($N_E = 20$); estimator initialization and the selection of the associated tuning matrices are chosen as presented in Section 6.5.1. Additionally, following [39], MPC prediction horizon is set to ($N_C = 20$) using a trial and error method with the associated LSQ weighting matrices selected, for all $i \in \{1, 2, 3\}$, as

$$\begin{aligned}
\mathbf{Q}_i &= \text{diag}(10, 10, 10, 0.1), \\
\mathbf{R}_i &= \text{diag}(25, 25, 25), \text{ and} \\
\mathbf{P}_i &= \text{diag}(50, 50, 50, 0.5).
\end{aligned} \quad (6.15)$$

As can be noticed from (6.15), the weights that correspond to the robots position in $\mathbf{Q}$ are chosen larger than that for the orientation, to give more freedom in steering the robots to reach their desired positions. Moreover, due to the experimental requirements, the elements of the weighting matrix $\mathbf{R}$ are chosen larger than that for $\mathbf{Q}$ to obtain a well damped closed-loop response and compensate for the unmodeled dynamics of the robots used in the experiments. Since the terminal value of the tracking error shown in (6.8) is the most necessary future error value that needs to be minimized, elements of $\mathbf{P}$ are chosen 5 times larger than that for $\mathbf{Q}$, see, e.g. [39], for more details on MPC tuning guidelines. The radii used in the collision avoidance constraints (6.10) and (6.11) are chosen as $r_c = 0.225$ (m) and $r_p = 0.3$ (m). All the simulations presented here are conducted using the code generation tool of ACADO for MHE and MPC [20, 41], which resembles the necessary real-time implementation requirements. The generated C code has been compiled into MEX functions and all simulations are conducted using MATLAB.

Figure 6.3 and Figure 6.4 summarize the performance of Algorithm 6.1 when applied to the considered MRS under the two considered scenarios shown in Table 6.2. These figure show the 3D exhibited trajectories of the observed robots with respect to the observing robot frame of reference[2].

Figure 6.3 and Figure 6.4 (left) illustrate the performance of the estimator by depicting the estimated trajectories against their ground truth data. As can be noticed, despite the exhibited trajectories by the observed robots, MHE provided an accurate state estimation along all trajectories. Figure 6.5 shows the estimation error components and the $3\sigma$ error boundaries (where $\sigma$ is the estimation error standard deviation) for the aerial observed robot $R_{c1}$ under case (1) of Table 6.2. It can be seen from Figure 6.5 that the estimation errors did not go beyond the $3\sigma$ error boundaries. This indicates that the estimation process is conservative, i.e. no overconfident estimation has occurred. This observation was identical to the other observed robots. Moreover, the MHE estimator achieved the relative localization task with $\sim 10$ (cm) positional accuracy and $\sim 0.18$ (rad) orientation accuracy. Figure 6.3 and Figure 6.4 (right) show the performance of the controller

---

[2]Click here for a visualization video of the overall motion of the MRS considered here.

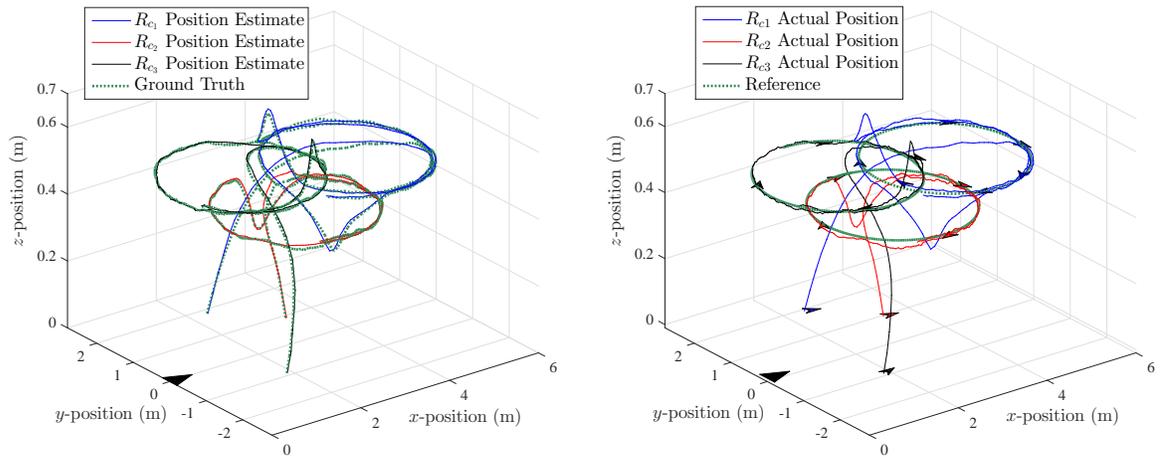

Figure 6.3: Algorithm 6.1 performance for case 1 given by Table 6.2: estimator performance (left); controller performance (right). Observing robot is highlighted by a black color filled triangle at the origin.

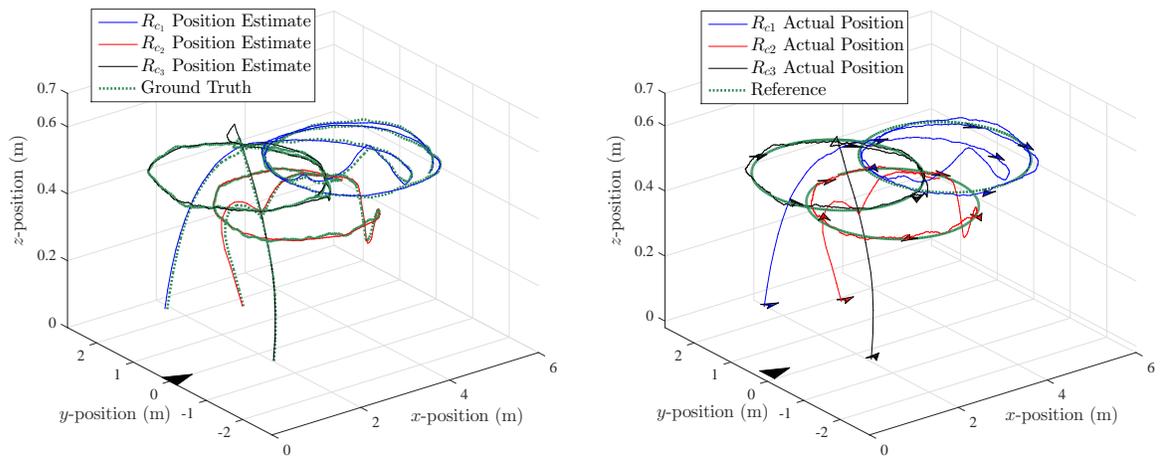

Figure 6.4: Algorithm 6.1 performance for case 2 given by Table 6.2: estimator performance (left); controller performance (right). Observing robot is highlighted by a black color filled triangle at the origin.

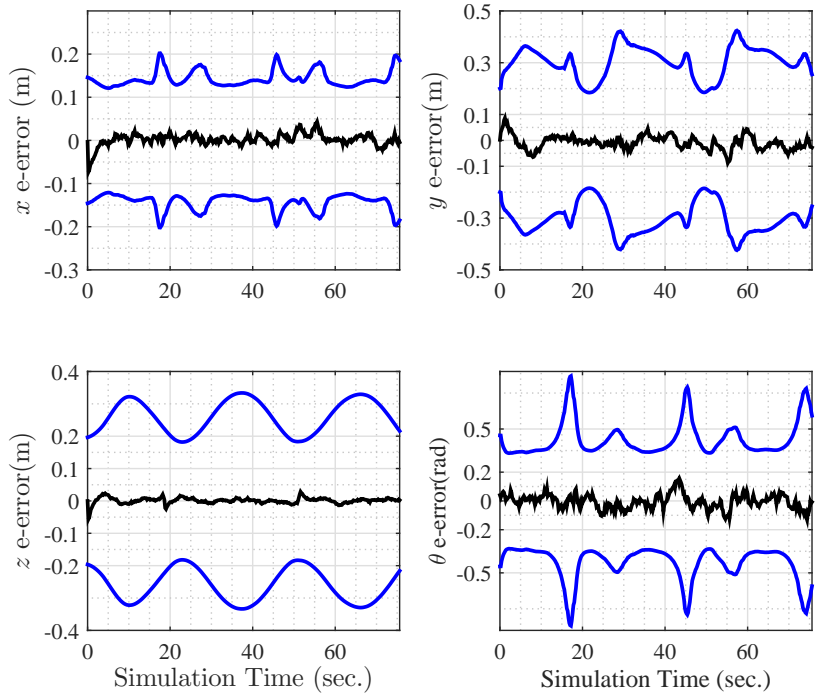

Figure 6.5: Relative state estimation error for ($R_{c1}$) in Case 1: estimation error (solid black), $3\sigma$ boundaries (solid blue).

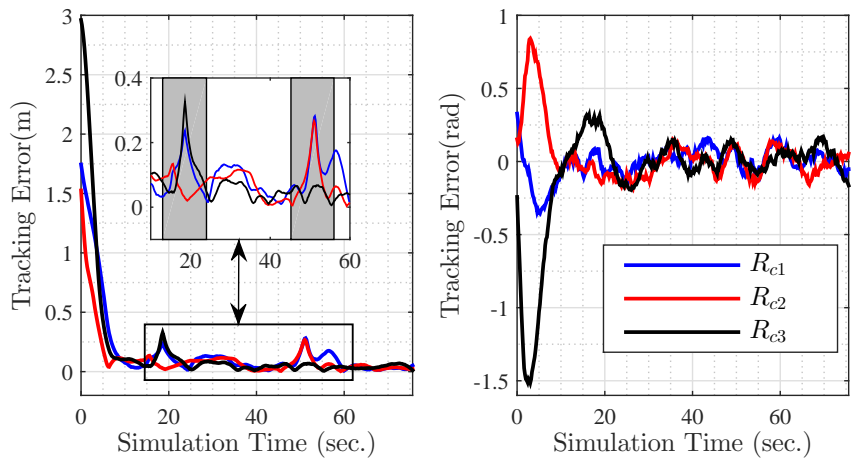

Figure 6.6: Relative tracking error for all observed robots (case 1 given by Table 6.2): relative position error (left), relative orientation error (right). Collision avoidance periods are highlighted in gray.

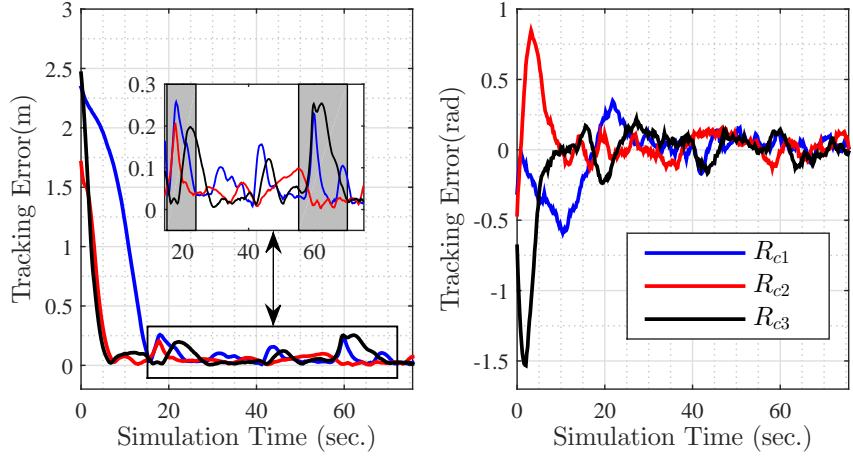

Figure 6.7: Relative tracking error for all observed robots (case 2 given by Table 6.2): relative position error (left), relative orientation error (right). Collision avoidance periods are highlighted in gray.

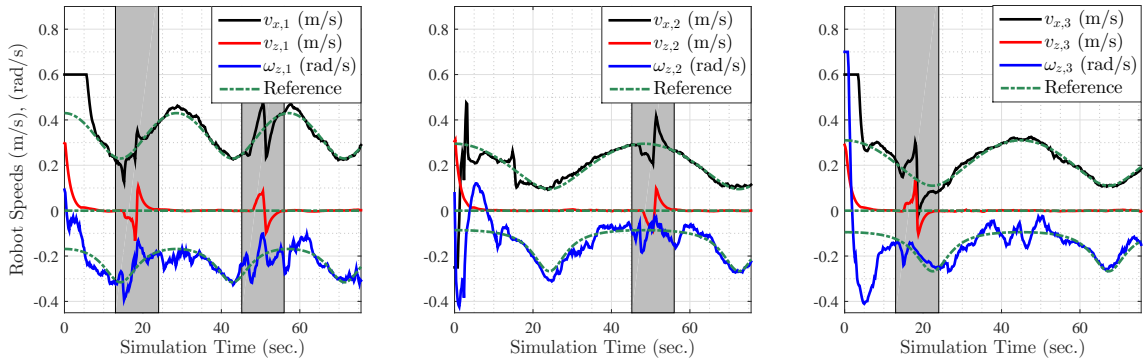

Figure 6.8: Control actions for all observed robots case 1. Collision avoidance periods are highlighted in gray.

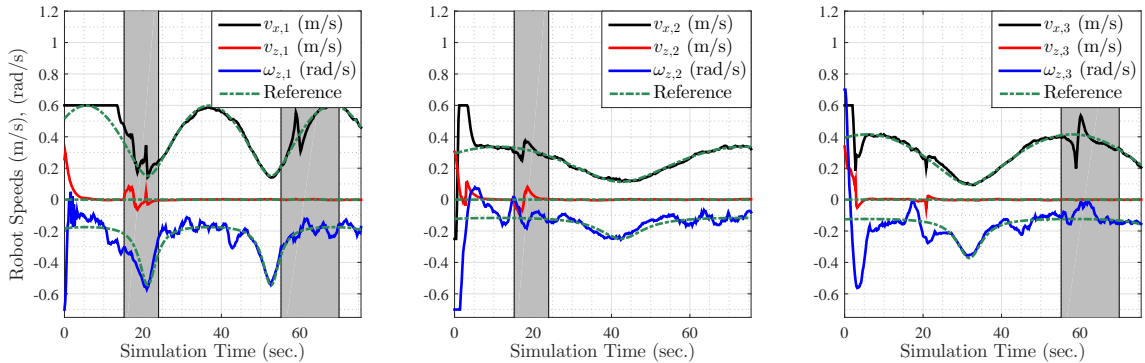

Figure 6.9: Control actions for all observed robots case 2. Collision avoidance periods are highlighted in gray.

Table 6.3: Algorithm 6.1 computational times (case 1).

| Computational Time | | $R_{c1}$ | $R_{c2}$ | $R_{c3}$ |
|---|---|---|---|---|
| MHE | maximum (ms) | 3.4 | 3.3 | 3.3 |
| | average (ms) | 3.1 | 3.1 | 3.1 |
| MPC | maximum (ms) | | 20 | |
| | average (ms) | | 15 | |
| Total | maximum (ms) | | 30 | |
| | average (ms) | | 19.3 | |

by depicting the actual trajectories of the observed robots and their references[3]. The position and orientation of the robots are indicated by color-filled triangles. These figures show that the MPC controller provided an accurate trajectory tracking regardless the initially imposed error. Moreover, at the moments, where the observed robots are about to collide, the controller provides collision free maneuvers in which the observed robots are commanded to deviate from their references to fulfill this task. Figure 6.6 and Figure 6.7 show the relative tracking errors of all observed robots in terms of the Euclidean relative position, and the relative orientation; the steady state relative tracking errors are observed to be within $\sim 12$ (cm) for the relative position and $\sim 0.16$ (rad) for the relative orientation. However, during two time intervals (highlighted in gray in Figure 6.6 and Figure 6.7) positional tracking errors deviated largely from the steady state values for collision avoidance purposes. This can be also seen in Figure 6.8 and Figure 6.9, where the control commands of all observed robots deviate also from their references to achieve collision free maneuvers. Moreover, these figures show that the control actions never went beyond their physical limits (6.12).

Since the adopted estimator-controller algorithm will be implemented experimentally, the associated computational requirements must be satisfied. The conducted simulations here have been run on a standard PC with a quad-core 2.5 GHz processor and 8 GB of memory. MHE and MPC computational times, for case 1 given by Table 6.2, are presented in Table 6.3, where it can be noticed that the adopted auto-generated C-code provided a fast implementation of the utilized algorithm, i.e. the feedback control actions $\boldsymbol{\mu}_i(\hat{\mathbf{x}}_i), i \in \{1,2,3\}$ were always ready to be sent to the observed robots within a maximum delay of $\sim 30$ (milliseconds) from the time the measurements $\tilde{\mathbf{y}}_i, i \in \{1,2,3\}$ are available. Similar computational requirements are observed for case 2 in Table 6.2.

---

[3]The actual trajectories of the observed robots in Figure 6.3 and Figure 6.4 (right) are the ground truth data in the same figures (left).

## 6.6 Experimental Results

Algorithm 6.1 was validated experimentally using a team of three ground non-holonomic robots, i.e. an Irobot Create®, which served as an observing robot, and two Pioneer 3-AT® research platforms, which served as observed robots; thus, $L = 1$ and $M = 2$, see, e.g. [50], for the technical details of Pioneer 3-AT mobile robots. To achieve high accuracy generation of the ground truth data, all robots were localized using the commercial motion capture system *OptiTrack*® with a setup of 4 cameras. Motivated by [10, 51], the ground truth data has been used to synthesis the relative measurements $\tilde{\mathbf{y}}_i, i \in \{1,2\}$ by computing the relative range $r_i$ and relative bearing $\phi_i$ using (6.4). Although synthesizing the relative measurements relaxes assumption (ii) presented in Section 6.3.1 in the conducted experiments, it provides a mean to simulate the behaviour of a relative measurement sensor and the freedom to choose its level of accuracy from Table 6.1, see, e.g. [52] for a possible relative measurement sensor that can be utilized here. The conducted experiments adopted a remote host[4] running a TCP/IP Matlab client, which communicates wirelessly with the robotic team to acquire the odometery measurements while remotely communicating with *OptiTrack*'s optical motion capture software (*Motive*®) in order to stream the ground truth data. Moreover, the same remote host sends the feedback control commands to servers running on the observed robots' platforms through the same TCP/IP connection. It has to be mentioned here that, in order to meet the real-time communication requirements, the client-server module of Pioneer-3AT mobile robots is based on an unacknowledged communication packet interface with simple means for dealing with ignored/lost packets. Moreover, the communication network used in the experiment was tested, for sending and receiving data, and the transmission delays were observed to be within the order of magnitude of ones of milliseconds and, thus, suitable for the sampling time used in the experiment.

---
[4]The used remote host in the experiments is the same PC used in simulations of Section 6.5.2.

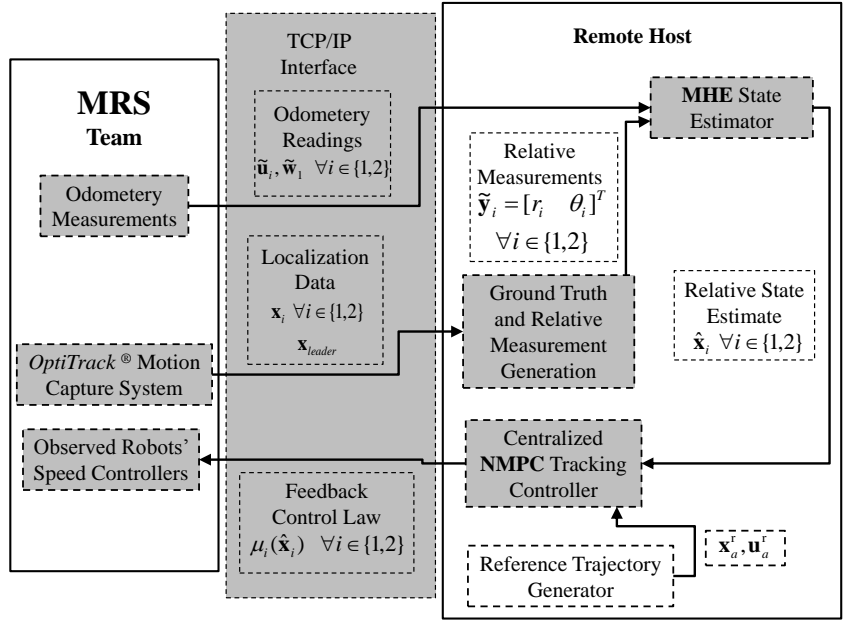

Figure 6.10: Block diagram of the experimental setup used to validate Algorithm 6.1.

The client performs the estimation-control associated computations by using the same MEX functions generated in Section 6.5 for MHE and MPC. The block diagram of the architecture used to validate the estimator-controller design is shown in Figure 6.10.

In the experiment, the observed robots were required to track the relative reference trajectories shown in Table 6.4 while the observing robot is stationary. Similar to Section 6.5.2, these trajectories were designed such that there is at least an instant at which the two references share the same location. The total experiment time was set to 65 (seconds) with a sampling rate of 10 Hz, i.e. $\delta = 0.1$ (seconds). Furthermore, the measurement noise configuration given in Table 6.1 (case 2) was used. MHE and MPC tuning matrices were chosen following the method presented in Section 6.5.2 with $N_E = N_C = 30$. Additionally, the radii shown in the inequality constraints (6.10) and (6.11) were chosen as $r_c = 0.5$ (m) and $r_p = 0.3$ (m).

Table 6.4: Reference trajectories in experiments.

| $R_{ci}$ | Reference |
|---|---|
| $R_{c1}$ | $x_1^r = 2 + 1.2\sin(0.15t)$<br>$y_1^r = 1.2\cos(0.15t)$ |
| $R_{c1}$ | $x_1^r = 2 + 1.2\sin(0.2t)$<br>$y_1^r = -0.9 + 1.2\cos(0.2t)$ |

Figure 6.11 summarizes the performance of MHE and MPC; to increase the readability of the figure, the performance has been visualized over two consequent time intervals[5]. As can be concluded from Figure 6.11, the utilized estimation-control design provided an accurate relative localization as well as relative tracking control. Figure 6.12 depicts the relative state estimation error of the observed robots poses, where it can be noticed that MHE provided an accurate state estimation from the first iteration, i.e. the relative localization accuracy has been observed to be within $\sim$ 15 (cm) for the relative position estimation and within $\sim 0.23$ (rad) for the relative orientation estimation.

Figure 6.13 shows the relative tracking errors of the observed robots in terms of the Euclidean relative position, and the relative orientation; the steady state relative tracking errors were observed to be within $\sim$ 11 (cm) for the relative position and $\sim 0.12$ (rad) for the relative orientation. In addition, it can be seen that during the obstacle avoidance period (highlighted in gray), the two robots deviated largely from their references.

---

[5]Click here for a video of the experiment (the experiment is at the second part of the video).

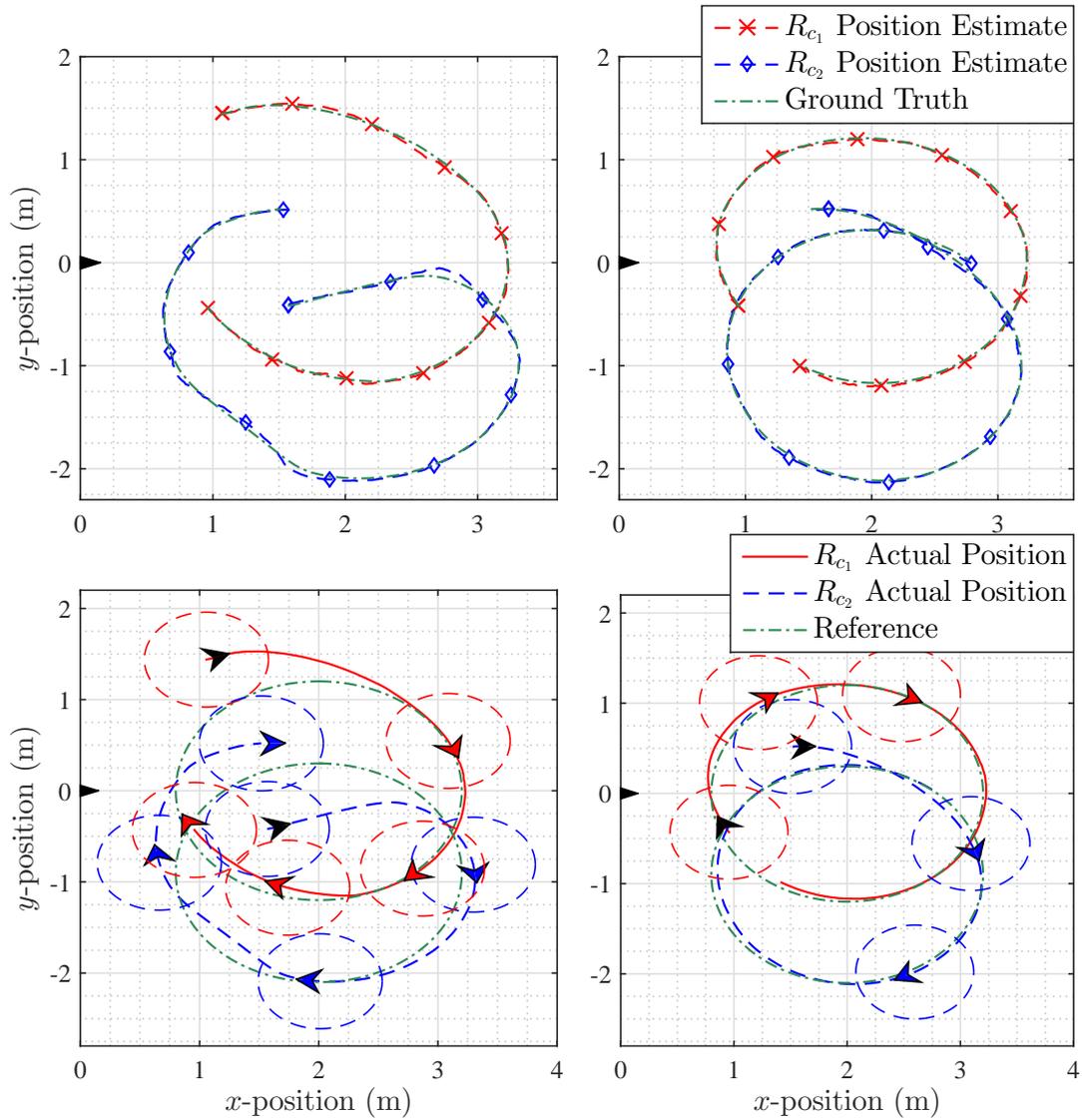

Figure 6.11: Experimental performance of Algorithm 6.1: relative localization results (top), relative tracking results (bottom). Left column subfigures show the performance during the time interval 0 to 30 (seconds) while right column subfigures are for the time interval 30.1 to 65 (seconds). The Observing robot is highlighted by a black color-filled triangle at the origin. Subfigures shown at the bottom highlight snapshots taken for the observed robots, at the same time instants, and illustrated by arrows with dotted circles indicating their minimum radius $r_c$; the first snapshot for each subfigure is highlighted by a black color-filled arrow.

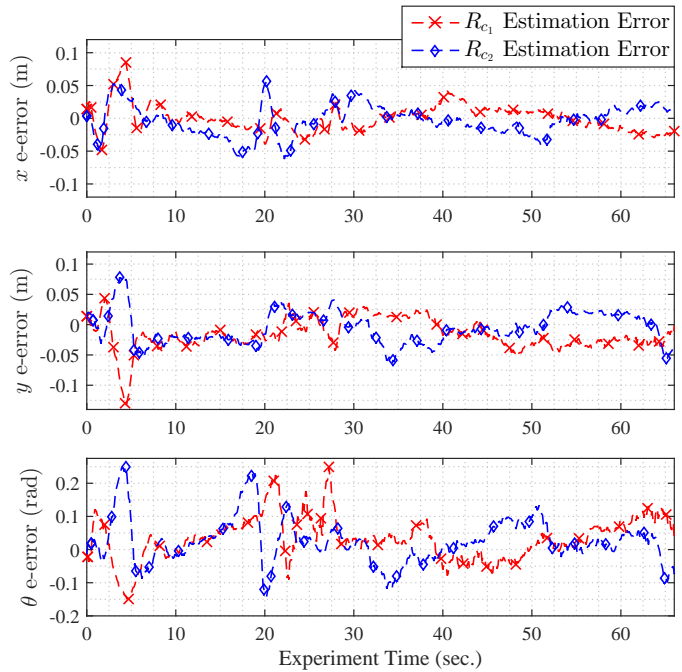

Figure 6.12: Relative state estimation error components.

This period is again highlighted in Figure 6.14, which shows the observed robots' speeds deviation from their references. Indeed, these deviations were necessary to satisfy the inequality constraints (6.10) and (6.11).

Algorithm 6.1 computational times needed in the experiments are summarized in Table 6.5. As can be noticed, the feedback control actions, i.e. $\boldsymbol{\mu}_i(\hat{\mathbf{x}}_i)$, $i \in \{1, 2\}$, were always ready to be sent to the observed robots speed controllers within a maximum time of $\sim 10.4$ (milliseconds) from the time the relative measurement $\tilde{\mathbf{y}}_i, i \in \{1, 2\}$ were available. In summary, Algorithm 6.1 was demonstrated to be practically successful at least at the scale of the considered experiment.

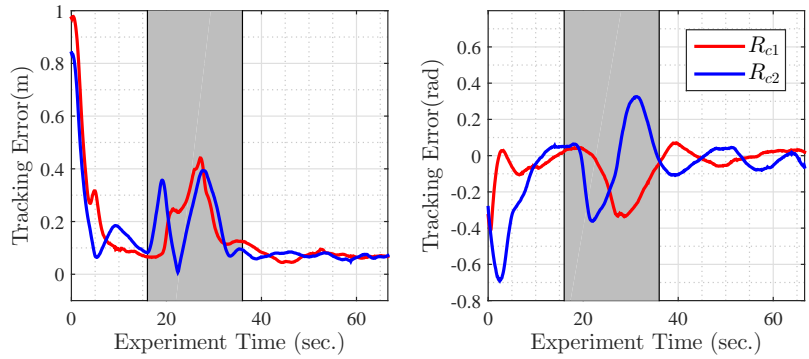

Figure 6.13: Relative trajectory tracking errors of the observed robots. Collision avoidance period is highlighted in gray.

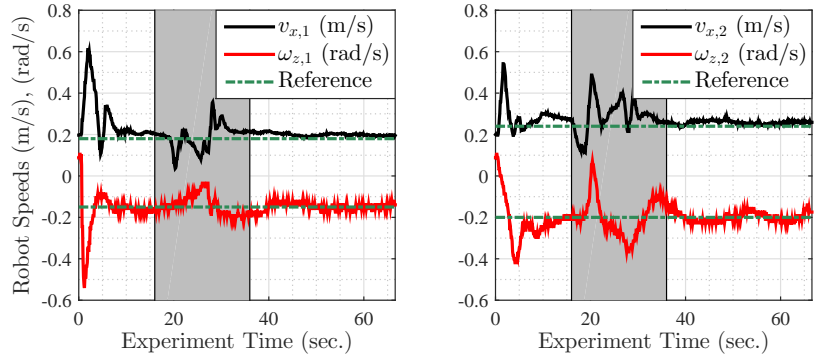

Figure 6.14: Speeds of the observed robots. Collision avoidance period is highlighted in gray.

Table 6.5: Algorithm 6.1 computational times in experiments.

|       | Computational Time | $R_{c1}$ | $R_{c2}$ |
|-------|--------------------|----------|----------|
| MHE   | maximum (ms)       | 1.8      | 2.2      |
|       | average (ms)       | 1.6      | 1.6      |
| MPC   | maximum (ms)       | 6.4      |          |
|       | average (ms)       | 4.4      |          |
| Total | maximum (ms)       | 10.4     |          |
|       | average (ms)       | 7.6      |          |

## 6.7 Conclusions

In this chapter, we utilized an optimization based estimation and control algorithm using MHE and MPC for relative localization and relative trajectory tracking control in multi-robots systems. Although MHE and MPC are known in the literature to be computationally intense, we adopted an auto-generated C-code implementing fast and efficient algorithms based on a real-time iteration (RTI) scheme. The used algorithm is evaluated through a series of simulations and real-time experiments. Firstly, a comparison between MHE and EKF has been conducted. This comparison showed the main advantages of the adopted relative localization scheme, which lie in the high accuracy and fast convergence of MHE. Secondly, the algorithm has been evaluated on an MRS mission and showed acceptable performance in fulfilling the mission requirements, as well as the real-time implementation requirements. Finally, an experimental evaluation of the presented algorithm has been conducted in which the practical requirements for estimation, control, and computations have been fulfilled.

## References

[1] J. Acevedo, B. Arrue, I. Maza, and A. Ollero, "Distributed approach for coverage and patrolling missions with a team of heterogeneous aerial robots under communication constraints," *International Journal of Advanced Robotic Systems*, vol. 10, no. 28, pp. 1–13, 2013.

[2] R. Beard, T. McLain, D. Nelson, D. Kingston, and D. Johanson, "Decentralized cooperative aerial surveillance using fixed-wing miniature uavs," *Proceedings of the IEEE*, vol. 94, no. 7, pp. 1306–1324, July 2006.

[3] D. Kingston, R. Beard, and R. Holt, "Decentralized perimeter surveillance using a team of uavs," *Robotics, IEEE Transactions on*, vol. 24, no. 6, pp. 1394–1404, Dec 2008.

[4] J. Acevedo, B. Arrue, I. Maza, and A. Ollero, "Cooperative large area surveillance with a team of aerial mobile robots for long endurance missions," *Journal of Intelligent and Robotic Systems*, vol. 70, no. 1-4, pp. 329–345, 2013.

[5] M. Bernard, K. Kondak, I. Maza, and A. Ollero, "Autonomous transportation and deployment with aerial robots for search and rescue missions," *Journal of Field Robotics*, vol. 28, no. 6, pp. 914–931, 2011.

[6] I. Maza and A. Ollero, "Multiple uav cooperative searching operation using polygon area decomposition and efficient coverage algorithms," in *Distributed Autonomous Robotic Systems 6*, R. Alami, R. Chatila, and H. Asama, Eds. Birkhäuser Basel, 2007, pp. 221–230.

[7] M. A. Hsieh, A. Cowley, J. F. Keller, L. Chaimowicz, B. Grocholsky, V. Kumar, C. J. Taylor, Y. Endo, R. C. Arkin, B. Jung, D. F. Wolf, G. S. Sukhatme, and D. C. MacKenzie, "Adaptive teams of autonomous aerial and ground robots for situational awareness," *Journal of Field Robotics*, vol. 24, no. 11-12, pp. 991–1014, 2007.

[8] F. Rivard, J. Bisson, F. Michaud, and D. Letourneau, "Ultrasonic relative positioning for multi-robot systems," in *Robotics and Automation, 2008. ICRA 2008. IEEE International Conference on*, 2008, pp. 323–328.

[9] O. De Silva, G. Mann, and R. Gosine, "Development of a relative localization scheme for ground-aerial multi-robot systems," in *Intelligent Robots and Systems (IROS), 2012 IEEE/RSJ International Conference on*, 2012, pp. 870–875.

[10] T. Wanasinghe, G. I. Mann, and R. Gosine, "Relative localization approach for combined aerial and ground robotic system," *Journal of Intelligent and Robotic Systems*, vol. 77, no. 1, pp. 113–133, 2015.

[11] T. Wanasinghe, G. Mann, and R. Gosine, "Distributed leader-assistive localization method for a heterogeneous multirobotic system," *Automation Science and Engineering, IEEE Transactions on*, vol. 12, no. 3, pp. 795–809, July 2015.

[12] M. W. Mehrez, G. K. Mann, and R. G. Gosine, "Nonlinear moving horizon state estimation for multi-robot relative localization," in *Electrical and Computer Engineering (CCECE), 2014 IEEE 27th Canadian Conference on*, May 2014, pp. 1–5.

[13] J. Fenwick, P. Newman, and J. Leonard, "Cooperative concurrent mapping and localization," in *Robotics and Automation, 2002. Proceedings. ICRA '02. IEEE International Conference on*, vol. 2, 2002, pp. 1810–1817.

[14] L. Carlone, M. Kaouk Ng, J. Du, B. Bona, and M. Indri, "Simultaneous localization and mapping using rao-blackwellized particle filters in multi robot systems," *Journal of Intelligent and Robotic Systems*, vol. 63, no. 2, pp. 283–307, 2011.

[15] S. Xingxi, W. Tiesheng, H. Bo, and Z. Chunxia, "Cooperative multi-robot localization based on distributed ukf," in *Computer Science and Information Technology (ICCSIT), 2010 3rd IEEE International Conference on*, vol. 6, 2010, pp. 590–593.

[16] C. Pathiranage, K. Watanabe, B. Jayasekara, and K. Izumi, "Simultaneous localization and mapping: A pseudolinear kalman filter (plkf) approach," in *Information and Automation for Sustainability, 2008. ICIAFS 2008. 4th International Conference on*, 2008, pp. 61–66.

[17] J. B. Rawlings and B. R. Bakshi, "Particle filtering and moving horizon estimation," *Computers and Chemical Engineering*, vol. 30, no. 10–12, pp. 1529 – 1541, 2006.

[18] C. Rao and J. Rawlings, "Nonlinear moving horizon state estimation," in *Nonlinear Model Predictive Control*, ser. Progress in Systems and Control Theory, F. Allgöwer and A. Zheng, Eds. Birkhäuser Basel, 2000, vol. 26, pp. 45–69.

[19] M. Zanon, J. Frasch, and M. Diehl, "Nonlinear moving horizon estimation for combined state and friction coefficient estimation in autonomous driving," in *Control Conference (ECC), 2013 European*, 2013, pp. 4130–4135.

[20] H. Ferreau, T. Kraus, M. Vukov, W. Saeys, and M. Diehl, "High-speed moving horizon estimation based on automatic code generation," in *Decision and Control (CDC), 2012 IEEE 51st Annual Conference on*, 2012, pp. 687–692.

[21] S. Wang, L. Chen, D. Gu, and H. Hu, "An optimization based moving horizon estimation with application to localization of autonomous underwater vehicles," *Robotics and Autonomous Systems*, vol. 62, no. 10, pp. 1581 – 1596, 2014.

[22] M. Saffarian and F. Fahimi, "Non-iterative nonlinear model predictive approach applied to the control of helicopters group formation," *Robotics and Autonomous Systems*, vol. 57, no. 6-7, pp. 749 – 757, 2009.

[23] M. W. Mehrez, G. K. Mann, and R. G. Gosine, "Formation stabilization of nonholonomic robots using nonlinear model predictive control," in *Electrical and Computer Engineering (CCECE), 2014 IEEE 27th Canadian Conference on*, May 2014, pp. 1–6.

[24] W. Dunbar and R. Murray, "Model predictive control of coordinated multi-vehicle formations," in *Decision and Control, 2002, Proceedings of the 41st IEEE Conference on*, vol. 4, 2002, pp. 4631–4636 vol.4.

[25] M. Lewis and K.-H. Tan, "High precision formation control of mobile robots using virtual structures," *Autonomous Robots*, vol. 4, no. 4, pp. 387–403, 1997.

[26] T. Balch and R. Arkin, "Behavior-based formation control for multirobot teams," *Robotics and Automation, IEEE Transactions on*, vol. 14, no. 6, pp. 926–939, 1998.

[27] G. Antonelli, F. Arrichiello, and S. Chiaverini, "Flocking for multi-robot systems via the null-space-based behavioral control," in *Intelligent Robots and Systems, 2008. IROS 2008. IEEE/RSJ International Conference on*, 2008, pp. 1409–1414.

[28] T. Balch and M. Hybinette, "Social potentials for scalable multi-robot formations," in *Robotics and Automation, 2000. Proceedings. ICRA '00. IEEE International Conference on*, vol. 1, 2000, pp. 73–80 vol.1.

[29] G. Mariottini, F. Morbidi, D. Prattichizzo, G. Pappas, and K. Daniilidis, "Leader-follower formations: Uncalibrated vision-based localization and control," in *Robotics and Automation, 2007 IEEE International Conference on*, 2007, pp. 2403–2408.

[30] X. Li, J. Xiao, and Z. Cai, "Backstepping based multiple mobile robots formation control," in *Intelligent Robots and Systems, 2005. (IROS 2005). 2005 IEEE/RSJ International Conference on*, 2005, pp. 887–892.

[31] J. Sanchez and R. Fierro, "Sliding mode control for robot formations," in *Intelligent Control. 2003 IEEE International Symposium on*, 2003, pp. 438–443.

[32] D. Gu and H. Hu, "A model predictive controller for robots to follow a virtual leader," *Robotica*, vol. 27, pp. 905–913, 10 2009.

[33] L. Grüne and J. Pannek, *Nonlinear Model Predictive Control: Theory and Algorithms*, ser. Communications and Control Engineering. Springer London Dordrecht Heidelberg New York, 2011.

[34] M. Diehl, H. Bock, J. P. Schlöder, R. Findeisen, Z. Nagy, and F. Allgöwer, "Real-time optimization and nonlinear model predictive control of processes governed by

differential-algebraic equations," *Journal of Process Control*, vol. 12, no. 4, pp. 577 – 585, 2002.

[35] B. Houska, H. Ferreau, and M. Diehl, "ACADO Toolkit – An Open Source Framework for Automatic Control and Dynamic Optimization," *Optimal Control Applications and Methods*, vol. 32, no. 3, pp. 298–312, 2011.

[36] T. Wanasinghe, G. Mann, and R. Gosine, "Pseudo-linear measurement approach for heterogeneous multi-robot relative localization," in *Advanced Robotics (ICAR), 2013 16th International Conference on*, Nov 2013, pp. 1–6.

[37] N. Trawny, X. Zhou, K. Zhou, and S. Roumeliotis, "Interrobot transformations in 3-d," *Robotics, IEEE Transactions on*, vol. 26, no. 2, pp. 226–243, April 2010.

[38] T. P. Nascimento, A. P. Moreira, and A. G. S. Conceição, "Multi-robot nonlinear model predictive formation control: Moving target and target absence," *Robotics and Autonomous Systems*, vol. 61, no. 12, pp. 1502 – 1515, 2013.

[39] E. Kayacan, E. Kayacan, H. Ramon, and W. Saeys, "Learning in centralized nonlinear model predictive control: Application to an autonomous tractor-trailer system," *Control Systems Technology, IEEE Transactions on*, vol. 23, no. 1, pp. 197–205, Jan 2015.

[40] M. Vukov, W. Van Loock, B. Houska, H. Ferreau, J. Swevers, and M. Diehl, "Experimental validation of nonlinear mpc on an overhead crane using automatic code generation," in *American Control Conference (ACC), 2012*, June 2012, pp. 6264–6269.

[41] B. Houska, H. J. Ferreau, and M. Diehl, "An auto-generated real-time iteration algorithm for nonlinear {MPC} in the microsecond range," *Automatica*, vol. 47, no. 10, pp. 2279 – 2285, 2011.

[42] E. Kayacan, E. Kayacan, H. Ramon, and W. Saeys, "Distributed nonlinear model predictive control of an autonomous tractor–trailer system," *Mechatronics*, vol. 24, no. 8, pp. 926 – 933, 2014.

[43] ——, "Robust tube-based decentralized nonlinear model predictive control of an autonomous tractor-trailer system," *Mechatronics, IEEE/ASME Transactions on*, vol. 20, no. 1, pp. 447–456, Feb 2015.

[44] M. Zanon, S. Gros, and M. Diehl, "Rotational start-up of tethered airplanes based on nonlinear mpc and mhe," in *Control Conference (ECC), 2013 European*, July 2013, pp. 1023–1028.

[45] M. Zanon, G. Horn, S. Gros, and M. Diehl, "Control of dual-airfoil airborne wind energy systems based on nonlinear mpc and mhe," in *Control Conference (ECC), 2014 European*, June 2014, pp. 1801–1806.

[46] S. Thrun, W. Burgard, and D. Fox, *Probabilistic Robotics (Intelligent Robotics and Autonomous Agents)*.   The MIT Press, 2005.

[47] A. Alessandri and M. Awawdeh, "Moving-horizon estimation for discrete-time linear systems with measurements subject to outliers," in *Decision and Control (CDC), 2014 IEEE 53rd Annual Conference on*, Dec 2014, pp. 2591–2596.

[48] M. J. Awawdeh, "Moving-horizon estimation for outliers detection and data mining applications," Ph.D. dissertation, Universita degli Studi di Genova, 2015.

[49] N. Hazon, F. Mieli, and G. Kaminka, "Towards robust on-line multi-robot coverage," in *Robotics and Automation, 2006. ICRA 2006. Proceedings 2006 IEEE International Conference on*, May 2006, pp. 1710–1715.

[50] M. W. Mehrez, G. K. Mann, and R. G. Gosine, "Stabilizing nmpc of wheeled mobile robots using open-source real-time software," in *Advanced Robotics (ICAR), 2013 16th International Conference on*, Nov 2013, pp. 1–6.

[51] M. Cognetti, P. Stegagno, A. Franchi, G. Oriolo, and H. Bulthoff, "3-d mutual localization with anonymous bearing measurements," in *Robotics and Automation (ICRA), 2012 IEEE International Conference on*, May 2012, pp. 791–798.

[52] O. D. Silva, G. K. I. Mann, and R. G. Gosine, "An ultrasonic and vision-based relative positioning sensor for multirobot localization," *IEEE Sensors Journal*, vol. 15, no. 3, pp. 1716–1726, March 2015.

# Chapter 7

# Summary and Outlook

In this chapter, we summarize the main contributions introduced in the thesis and discuss a number of possible extensions to the presented work.

The primary focus of this research study was to adopt optimization based solutions for control and state-estimation in non-holonomic mobile robots. Single robot control problems; multiple robots control problems; and relative-localization in ground-aerial multi-robot systems (MRS's) were the central problems discussed in the thesis. Optimization based solutions used are model predictive control (MPC) and moving horizon estimation (MHE). The conducted research formed four main objectives:

I. Designing an asymptotically stable MPC control scheme, without stabilizing constraints or cost, for the point stabilization control problem of a non-holonomic mobile robot.

II. Designing an asymptotically stable MPC control scheme, without stabilizing constraints or cost, for the path following control problem of a non-holonomic mobile robot.

III. Designing a reduced communication distributed model predictive control (DMPC)

scheme for the regulation control problem of a multi-robots system (MRS).

IV. Designing a relative localization scheme based on an MHE observer for a ground-aerial (MRS).

## 7.1 Research Summary based on Objective I

Within the first objective, the asymptotic stability of non-holonomic mobile robots' point stabilization under MPC was considered. After an extensive literature review, it has been observed that the available MPC controllers for the point stabilization problem can be shown to be asymptotically stable by the assistance of stabilizing constraints and/or costs. Adding these stabilizing enforcing ingredients introduces computational complexities and implementation challenges as highlighted in Sections 2.2 and 3.2. Therefore, in this study, we verified a novel approach for guaranteeing MPC asymptotic stability for the point stabilization problem of a non-holonomic mobile robot. This approach is based on deriving upper bounds on the value function and computing a prediction horizon length such that asymptotic stability of the MPC closed-loop is rigorously proven. To this end, open loop trajectories were constructed in order to derive the value function bounds. Additionally, we showed that these upper bounds can be attained by suitably choosing the running costs of the MPC controller, i.e. a non-quadratic running costs was adopted. The proposed running costs penalizes the direction orthogonal to the desired orientation more than other directions. The success of the theoretical findings in this part of the study was verified by a series of numerical simulations. First, the simulation results showed local mere convergence to the desired stabilization point using the proposed running costs in comparison with the quadratic running costs. Moreover, the simulation results showed the satisfaction of the relaxed-Lyapunov inequality–asymptotic stability condition–for several initial conditions.

## 7.2 Research Summary based on Objective II

Secondly, we considered the path following task of a non-holonomic mobile robot using MPC. The control objective here is to follow a state-space geometric reference. Here, the speed to move along the reference is not given a priori, but rather is obtained via a timing law, which adds an extra degree of freedom to the controller. Similar to the point stabilization case, it has been shown in the literature that asymptotic stability of the model predictive path following control (MPFC) can be guaranteed by adding stabilizing constraints and/or costs. As a result, the asymptotic stability of MPFC without stabilizing constraints or costs was considered in this research. Although we limited our analysis to a particular sub-set of the state-set, we showed that asymptotic stability can be guaranteed without stabilizing constraints or costs by appropriately choosing the prediction horizon length. To this end, we first formulated the MPFC control problem as a point stabilization problem of an augmented dynamics system. Then, we designed open-loop control maneuvers such that the robot is steered from any given initial condition to the end-point on a considered path. This led to bounds on the MPFC value function, which were later used to derive prediction horizon lengths such that closed-loop stability is verified. The theoretical development in this part of the thesis was verified by a set of numerical experiments. Initially, we showed the effect of the MPFC tuning parameters, e.g. running costs weights, on the stabilizing prediction horizon length. Moreover, our numerics verified the asymptotic stability of MPFC by conducting closed-loop simulations for a number of feasible initial condition.

## 7.3 Research Summary based on Objective III

Thirdly, the regulation control of an MRS to a prescribed pattern was considered. Here, we adopted a DMPC scheme to achieve this control task. In DMPC, each subsystem in a given MRS solves its own (local) optimal control problem (OCP) and communicate its own prediction trajectory to other robots in the MRS. This communicated information is later used to formulate coupling constraints, which accounts for inter-robot collision avoidance. As shown in the literature, the communication of the full prediction trajectories is the most common technique to formulate the coupling constraints. Therefore, in this thesis we focused on reducing the communication load in DMPC for non-holonomic MRS's. To this end, we introduced a quantization technique to achieve this purpose. The proposed method is based on partitioning the operating region into a grid and then projecting the predicted trajectories onto this grid to generate an occupancy grid prediction. The occupancy grid is later communicated via a proposed differential-communication scheme. Additionally, we introduced a coupling constraints formulation in DMPC, which is based on the introduced communication method and a squircle approximation. The success of the proposed approach was verified by numerical simulations, which showed a significant decrease in the data-exchange load. The simulation results showed also the convergence of an MRS to a prescribed reference pattern.

## 7.4 Research Summary based on Objective IV

Finally, we considered the state estimation problem in ground-aerial MRS's. Here, relative (range and bearing) measurements are used to achieve relative localization (RL) of robots with limited sensory capabilities with respect to robots with accurate localization means. The first class of robots is denoted by (observed robots) while the second class is referred to as (observing robots). It has been observed in the literature that this state estimation problem was not treated by an optimization based method, e.g. moving horizon estimation (MHE). Therefore, in the last phase of this thesis, we gave a particular attention to solving the RL problem in MRS's using MHE. To this end, we formulated the RL problem in the framework of MHE and showed the advantages of using MHE by comparing its performance against extended Kalman filter (EKF). As a complementary part of this localization scheme, we proposed a centralized MPC controller to achieve relative trajectory in MRS's based on the RL results obtained via MHE. More precisely, we introduced an algorithm in which relative localization is achieved via MHE, while the relative trajectory tracking is fulfilled by a centralized MPC. The proposed algorithm was verified by numerical simulations as well as laboratory experiments.

## 7.5 Summary of Contributions

In summary, the following are the key contributions of the thesis objectives in using optimization based solutions for control and state estimation in non-holonomic mobile robots.

1. **Contributions from Objective I:**

    - A novel implementation of MPC control, without stabilizing constraints or costs, for the point stabilization of a non-holonomic mobile robot.

    - A novel proof of the controllability assumptions for point stabilization based on a newly introduced running costs.

    - An algorithm for determining a stabilizing prediction horizon length.

2. **Contributions from Objective II:**

    - A novel implementation of MPC control, without stabilizing constraints or costs, for the path following control of a non-holonomic mobile robot.

    - A novel proof of the controllability assumptions for path following.

3. **Contributions from Objective III:**

- A novel design of a DMPC based on occupancy grid and differential communication.

- A novel approach for communication in DMPC for non-holonomic mobile robots pattern regulation.

- A Novel formulation of coupling constraints based on occupancy grid and squircle approximation.

4. **Contributions from Objective IV:**

    - A novel implementation of MHE observers in relative localization in ground-aerial MRS's.

    - Development of an optimization based algorithm in which relative localization as well as relative trajectory tracking are achieved by means of MHE and MPC, respectively.

    - Experimental validation of the developed localization-control algorithm.

### 7.5.1 List of Publications

The following list of scientific articles is a result of the research conducted in this thesis. The core chapters of the thesis, i.e. Chapters 3, 4, 5, and 6, are the publications 3, 4, 6, and 2 of the list, respectively.

**Journal Articles:**

1. Karl Worthmann, Mohamed W. Mehrez, George K.I. Mann, Raymond G. Gosine, and Jürgen Pannek, "Interaction of Open and Closed Loop Control in MPC", Accepted for publication in Automatica, 2017 [1].

2. Mohamed W. Mehrez, George. K.I. Mann, and Raymond G. Gosine, "An Optimization Based Approach for Relative Localization and Relative Tracking Control in Multi-Robot Systems", in Journal of Intelligent and Robotic Systems, 2017, [2].

3. Karl Worthmann, Mohamed W. Mehrez, Mario Zanon, George K.I. Mann, Raymond G. Gosine, and Moritz Diehl, "Model Predictive Control of Nonholonomic Mobile Robots Without Stabilizing Constraints and Costs", in IEEE Transactions on Control Systems Technology, 2016 [3].

**Conference Articles:**

4. Mohamed W. Mehrez, Karl Worthmann, George K.I. Mann, Raymond G. Gosine, and Timm Faulwasser, "Predictive Path Following of Mobile Robots without Terminal Stabilizing Constraints", in Proceedings of the IFAC 2017 World Congress, Toulouse, France, 2017, [4].

5. Mohamed W. Mehrez, Karl Worthmann, George K.I. Mann, Raymond G. Gosine, and Jürgen Pannek, "Experimental Speedup and Stability Validation for Multi-Step MPC", in Proceedings of the IFAC 2017 World Congress, Toulouse, France, 2017, [5].

6. Mohamed W. Mehrez, Tobias Sprodowski, Karl Worthmann, George K.I. Mann, Raymond G. Gosine, Juliana K. Sagawa, and Jürgen Pannek, "Occupancy Grid based Distributed MPC of Mobile Robots", [6].

7. Karl Worthmann, Mohamed W. Mehrez, Mario Zanon, George K. I. Mann, Raymond G. Gosine, and Moritz Diehl, "Regulation of Differential Drive Robots Using Continuous Time MPC Without Stabilizing Constraints or Costs", in Proceedings of the 5th IFAC Conference on Nonlinear Model Predictive Control (NPMC'15), 2015, [7].

8. Mohamed W. Mehrez, George K. I. Mann, and Raymond G. Gosine, "Comparison of Stabilizing NMPC Designs for Wheeled Mobile Robots: an Experimental Study", in Proceedings of the Moratuwa Engineering Research Conference (MERCON'15), 2015, [8].

9. Mohamed W. Mehrez, George K. I. Mann, and Raymond G. Gosine, "Formation stabilization of nonholonomic robots using nonlinear model predictive control", in Proceedings of the IEEE 27th Canadian Conference on Electrical and Computer Engineering (CCECE'14), 2014, [9].

10. Mohamed W. Mehrez, George K. I. Mann, and Raymond G. Gosine, "Nonlinear moving horizon state estimation for multi-robot relative localization", in Proceedings of the IEEE 27th Canadian Conference on Electrical and Computer Engineering (CCECE'14), 2014, [10].

11. Mohamed W. Mehrez, George K. I. Mann, and Raymond G. Gosine, "Stabilizing NMPC of wheeled mobile robots using open-source real-time software", in Proceedings of the 16th International Conference on Advanced Robotics (ICAR'13), 2013, [11].

12. Mohamed W. Mehrez, George K. I. Mann, and Raymond G. Gosine, "Control Profile Parameterized Nonlinear Model Predictive Control of wheeled Mobile Robots", in Proceedings of the Newfoundland Electrical and Computer Engineering Conference (NECEC'13), 2013, [12].

## 7.6 Future Research Directions

The research work presented in the thesis has a number of possible future extensions. These future developments aim at the generalization of the results and improving its practicality.

**Single robot control under MPC:** Here, we propose the following extensions

- More complex dynamic models can be used and tested under MPC without stabilizing constraints or costs. These models may consider uneven (rough) terrains, e.g. the models presented in [13]. The models may also be learned online as shown in the studies [14, 15].

- For both point stabilization and path following control problems, domains with obstacles can be considered when studying asymptotic stability under MPC.

- Less conservative open loop maneuvers can be also developed in order to find shorter prediction horizons, which guarantee closed loop asymptotic stability.

- For the path following control problem, stability can be studied for more generic state set rather than only a sub-set of it. Moreover, stability can be studied for periodic paths, i.e. paths that repeat, see, e.g. [16–18].

- Results presented here can be also extended to practical asymptotic stability of MPC closed-loop instead of the more conservative settings, i.e. asymptotic stability. In practical stability, convergence only the neighborhood of the reference is required rather than the mere convergence to the reference itself, see, e.g. [19], for details.

**Occupancy grid based DMPC:** Here, we propose the following extensions

- Verifying the proposed controller via laboratory experiments.

- Investigation of the proposed method closed-loop (practical) asymptotic stability by following the results presented in [20].

- Investigating dynamic-order priority rules of the distributed optimization instead of only a fixed-order priority rule.

**Relative localization in MRS:** Here, we propose the following extensions

- Implementing the proposed localization and control algorithm to actual systems with aerial vehicles. Moreover, considering the experimental implementation with actual relative range and bearing sensors similar to the one proposed in [21].

- Modifying the proposed algorithm to achieve relative tracking control via distributed MPC instead of centralized MPC.

- Considering a more general 6-DOF model of aerial vehicles dynamics instead of the 4-DOF model used.

# References

[1] K. Worthmann, M. W. Mehrez, G. K. Mann, R. G. Gosine, and J. Pannek, "Interaction of open and closed loop control in mpc," *Automatica*, vol. 82, pp. 243 – 250, 2017.

[2] M. W. Mehrez, G. K. I. Mann, and R. G. Gosine, "An optimization based approach for relative localization and relative tracking control in multi-robot systems," *Journal of Intelligent & Robotic Systems*, vol. 85, no. 2, p. 385–408, 2017.

[3] K. Worthmann, M. W. Mehrez, M. Zanon, G. K. I. Mann, R. G. Gosine, and M. Diehl, "Model predictive control of nonholonomic mobile robots without stabilizing constraints and costs," *IEEE Transactions on Control Systems Technology*, vol. 24, no. 4, pp. 1394–1406, 2016.

[4] M. W. Mehrez, K. Worthmann, G. K. I. Mann, R. Gosine, and T. Faulwasser, "Predictive path following of mobile robots without terminal stabilizing constraints," in *Proc. 20th IFAC World Congr., Toulouse, France*, 2017, pp. 10 268–10 273.

[5] M. W. Mehrez, K. Worthmann, G. K. Mann, R. G. Gosine, and J "Experimental speedup and stability validation for multi-step mpc," *IFAC-PapersOnLine*, vol. 50, no. 1, pp. 8698 – 8703, 2017, 20th IFAC World Congress.

[6] M. W. Mehrez, T. Sprodowski, K. Worthmann, G. K. I. Mann, R. G. Gosine, J. K. Sagawa, and J. Pannek, "Occupancy Grid based Distributed Model Predictive Control of Mobile Robots," in *IEEE/RSJ International Conference on Intelligent Robots and Systems (IROS)*, 2017, pp. 4842 – 4847.

[7] K. Worthmann, M. W. Mehrez, M. Zanon, G. K. Mann, R. G. Gosine, and M. Diehl, "Regulation of differential drive robots using continuous time mpc without stabilizing constraints or costs," in *Proceedings of the 5th IFAC Conference on Nonlinear Model Predictive Control (NPMC'15), Sevilla, Spain*, 2015, pp. 129–135.

[8] M. W. Mehrez, G. K. I. Mann, and R. G. Gosine, "Comparison of stabilizing nmpc designs for wheeled mobile robots: An experimental study," in *2015 Moratuwa Engineering Research Conference (MERCon)*, April 2015, pp. 130–135.

[9] M. W. Mehrez, G. K. Mann, and R. G. Gosine, "Formation stabilization of nonholonomic robots using nonlinear model predictive control," in *IEEE 27th Canadian Conference on Electrical and Computer Engineering (CCECE)*, May 2014, pp. 1–6.

[10] ——, "Nonlinear moving horizon state estimation for multi-robot relative localization," in *Electrical and Computer Engineering (CCECE), 2014 IEEE 27th Canadian Conference on*, May 2014, pp. 1–5.

[11] ——, "Stabilizing nmpc of wheeled mobile robots using open-source real-time software," in *Advanced Robotics (ICAR), 2013 16th International Conference on*, Nov 2013, pp. 1–6.

[12] M. W. Mehrez, G. Mann, and R. Gosine, "Control profile parameterized nonlinear model predictive control of wheeled mobile robots," in *The 22nd Annual Newfoundland Electrical and Computer Engineering Conference*.   IEEE Newfoundland and Labrador Section, 2013.

[13] A. Mandow, J. L. Martinez, J. Morales, J. L. Blanco, A. Garcia-Cerezo, and J. Gonzalez, "Experimental kinematics for wheeled skid-steer mobile robots," in *Proc. Int. Conf. Intelligent Robots and Syst.*, Oct 2007, pp. 1222–1227.

[14] C. J. Ostafew, A. P. Schoellig, and T. D. Barfoot, "Conservative to confident: Treating uncertainty robustly within learning-based control," in *2015 IEEE International Conference on Robotics and Automation (ICRA)*, May 2015, pp. 421–427.

[15] ——, "Learning-based nonlinear model predictive control to improve vision-based mobile robot path-tracking in challenging outdoor environments," in *2014 IEEE International Conference on Robotics and Automation (ICRA)*, May 2014, pp. 4029–4036.

[16] A. Alessandretti, A. P. Aguiar, and C. Jones, "Trajectory-tracking and Path-following Controllers for Constrained Underactuated Vehicles using Model Predictive Control," in *Proc. European Control Conf.*, Zürich, Switzerland, 2013, pp. 1371–1376.

[17] T. Faulwasser, *Optimization-based Solutions to Constrained Trajectory-tracking and Path-following Problems*.   Shaker, Aachen, Germany, 2013.

[18] T. Faulwasser and R. Findeisen, "Nonlinear model predictive control for constrained output path following," *IEEE Trans. Automat. Control*, vol. 61, no. 4, pp. 1026–1039, April 2016.

[19] L. Grüne and J. Pannek, "Practical NMPC suboptimality estimates along trajectories," *Syst. and Control Lett.*, vol. 58, no. 3, pp. 161 – 168, 2009.

[20] L. Grüne and K. Worthmann, *Distributed Decision Making and Control*. Springer Verlag, 2012, ch. A distributed NMPC scheme without stabilizing terminal constraints, pp. 261–287.

[21] O. De Silva, G. Mann, and R. Gosine, "Development of a relative localization scheme for ground-aerial multi-robot systems," in *Intelligent Robots and Systems (IROS), 2012 IEEE/RSJ International Conference on*, 2012, pp. 870–875.